\documentclass[11pt, a4paper]{book}
\pdfoutput=1

\usepackage{graphicx}
\usepackage{amsfonts}
\usepackage{amssymb}
\usepackage{amsmath}
\usepackage{amsthm}
\usepackage[all,2cell]{xy}
\usepackage{enumerate}
\usepackage[T1]{fontenc}
\usepackage{textcomp}
\usepackage{stmaryrd}  
\usepackage{eucal}
\usepackage{verbatim}
\usepackage{setspace}
\usepackage[safe]{tipa}
\usepackage{simplewick}
\usepackage{textpos}
\usepackage{hyperref}
\UseAllTwocells
\input xy

\TPGrid{1}{1} 

\theoremstyle{plain}
\newtheorem{thm}{Theorem}[chapter]
\newtheorem{cor}[thm]{Corollary}
\newtheorem{lem}[thm]{Lemma}
\newtheorem{prop}[thm]{Proposition}

\newtheorem*{conjecture}{Conjecture}
\newtheorem*{thmt}{Theorem}

\theoremstyle{definition}
\newtheorem{defn}[thm]{Definition}

\newcommand\mi{{-1}}
\newcommand\mmi{{\text{\rm -} 1}}
\newcommand{\la}[1]{\stackrel{\,\,\, #1}{\longleftarrow}}
\newcommand{\ra}[1]{\stackrel{\,\,\, #1}{\longrightarrow}}
\newcommand{\sra}[1]{\stackrel{ #1}{\rightarrow}}
\newcommand{\sla}[1]{\stackrel{ \, \,\, #1 }{\leftarrow}}

\newcommand\ba{\begin{aligned}}
\newcommand\ea{\end{aligned}}
\newcommand\ig{\includegraphics}

\newcommand\XXX[1]{}
\newcommand\Iso{\stackrel{\sim}{\Rightarrow}}

\newcommand\uIso{\operatorname{\text{\rm uIso}}}

\newcommand\id{\operatorname{\text{\rm id}}}
\newcommand\pre{\operatorname{\text{\rm pre}}}
\newcommand\ch{\operatorname{\text{\rm ch}}}
\newcommand\post{\operatorname{\text{\rm \rm post}}} 

\newcommand\ad{\operatorname{\text{\rm ad}}}

\newcommand\codim{\operatorname{\text{\rm codim}}}
\newcommand\proj{\operatorname{\text{\rm proj}}}
\newcommand\CYau{\operatorname{\mathcal{C}\mathcal{Y}\text{\rm au}}}

\newcommand\ccdot{\smash{\cdot}}
\newcommand\be{\begin{equation}}
\newcommand\Td{\operatorname{\text{\rm Td}}}

\newcommand\vol{\operatorname{\text{\rm vol}}}
\newcommand\ee{\end{equation}}
\newcommand{\ttr}[1] {{\mathbb{T}\text{\rm r}(#1)}}

\newcommand{\threearrow}{\; \, \smash{\ba \ig{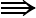} \ea} \; \; \,}

\newcommand\Sections{\operatorname{\text{\rm Sections}}}
\newcommand\sign{\operatorname{\text{\rm sign}}}
\newcommand\EvenHanded{\operatorname{\text{\rm Even-Handed}}}
\newcommand\piv{{\operatorname{\text{\rm piv}}}}
\newcommand\EG{\mathcal{E}G}
\newcommand\A{\mathcal{A}}
\newcommand\B{\mathcal{B}}
\newcommand\BG{\mathbf{B}G}
\newcommand\Paths{\operatorname{\text{\rm Paths}}}

\newcommand\fusion{{\operatorname{\text{\rm fusion}}}}

\newcommand{\Adj}{\operatorname{\text{\rm Adj}}}

\newcommand{\Tr}{\operatorname{\text{\rm Tr}}}

\newcommand{\Ind}{\operatorname{\text{\rm Ind}}}

\newcommand{\LBun}{\operatorname{\text{\rm LBun}}}
\newcommand{\LCat}{\operatorname{\mathcal{L}\text{\rm Cat}}}

\newcommand{\SLCat}{\operatorname{\mathcal{SL}\text{\rm Cat}}}
\newcommand{\Set}{\operatorname{\text{\rm Set}}}

\newcommand{\TRep}{\operatorname{\text{\rm 2}\mathcal{R}\text{\rm ep}}}
\newcommand{\Class}{\operatorname{\text{\rm Class}}}
\newcommand{\Amb}{\operatorname{\text{\rm Amb}}}
\newcommand{\Rep}{\operatorname{\text{\rm Rep}}}

\newcommand{\Hom}{\operatorname{\text{\rm Hom}}}
\newcommand{\Hilb}{\operatorname{\text{\rm Hilb}}}
\newcommand{\Pic}{\operatorname{\text{\rm Pic}}}
\newcommand{\Stab}{\operatorname{\text{\rm Stab}}}
\newcommand{\Aut}{\operatorname{\text{\rm Aut}}}

\newcommand{\End}{\operatorname{\text{\rm End}}}
\newcommand{\FinSpaces}{\operatorname{\mathcal{F}\text{\rm inSpaces}}}
\newcommand{\Vect}{\operatorname{\text{\rm Vect}}}
\newcommand{\THilb}{\operatorname{\text{\rm 2}\mathcal{H}\text{\rm ilb}}}

\newcommand{\Cat}{\operatorname{\mathcal{C}\text{\rm at}}}

\newcommand{\X}{\mathcal{X}}
\newcommand{\Y}{\mathcal{Y}}
\newcommand{\Z}{\mathcal{Z}}
\newcommand{\E}{\mathcal{E}}
\newcommand{\T}{\mathsf{T}}
\newcommand{\ttt}{\mathsf{t}}
\newcommand{\sss}{\mathsf{s}}
\newcommand{\Ss}{\mathsf{S}}
\newcommand{\Ob}{\operatorname{\text{\rm Ob}}}

\newcommand{\Nat}{\operatorname{\text{\rm Nat}}}

\newcommand{\V}{\mathcal{V}}
\newcommand{\Dim}{\operatorname{\text{\rm Dim}}}

\newcommand{\Gerbes}{\operatorname{\mathcal{G}\text{\rm \rm{erbes}}}}
\newcommand{\FFix}{\operatorname{\text{\rm Fix}}}
\newcommand{\UOneTor}{U(1)\text{\rm -Tor}}
\newcommand{\comments}[1]{}

 \DeclareMathOperator{\Fun}{Fun}

\def\clap#1{\hbox to 0pt{\hss#1\hss}}

\def\mathclap{\mathpalette\mathclapinternal}

\def\mathclapinternal#1#2{%
\clap{$\mathsurround=0pt#1{#2}$}}

\newcommand{\lag}[2]{ \xy {\ar@{->}^(0.45){\scriptscriptstyle #2}_(0.4){\scriptscriptstyle
#1} (6,0)*{}; (0,0)*{}};
\endxy}

\newcommand{\biglag}[2]{ \xy {\ar@{->}^(0.45){\scriptscriptstyle #2}_(0.4){\scriptscriptstyle
#1} (15,0)*{}; (0,0)*{}};
\endxy}

\newcommand{\ratone}{\sqrt{\frac{\phantom{i}\mathclap{\smash{\ba \ig{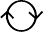} \ea}}\phantom{T}}{\phantom{i}\mathclap{\smash{\ba \ig{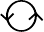} \ea}}\phantom{T}}}}

\newcommand{\rattwo}{\sqrt{\frac{\phantom{i}\mathclap{\smash{\ba \ig{e227.pdf} \ea}}\phantom{T}}{\phantom{i}\mathclap{\smash{\ba \ig{e226.pdf} \ea}}\phantom{T}}}}

\newcommand{\nlag}[3]{{ \begin{xy}   (4,-0.7)*{\scriptscriptstyle #1}="d"; (0,-0.7)*{\scriptscriptstyle #2}="a"; (1.2,-0.7)*{\scriptscriptstyle <}  {\ar@{-}_{\scriptscriptstyle \; \; #3} "d"; "a"} \end{xy}}}
\newcommand{\mlag}[3]{{\scriptscriptstyle #2 \sla{\smash{\, #3}} \scriptscriptstyle #1}}
\newcommand{\mmlag}[3]{{\scriptscriptstyle #2 \xleftarrow{\;\;\; \smash{\mathclap{\scriptscriptstyle #3} \; }} \scriptscriptstyle #1}}
\newcommand{\blag}[3]{{#2 \xleftarrow{\, #3} #1}}

 \reversemarginpar
\SelectTips{cm}{}

\newcommand{\Fix}[2] {
    \xy
        (0,0)*{#1};
        (0,2)*{\ba \ig{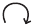} \ea};
        (2.8,3)*{\scriptstyle #2};
    \endxy}

\newcommand{\FixP}[2] {
    \xy
        (0,0)*{#1};
        (0.3,2)*{\ba \ig{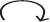} \ea};
        (3.6,3)*{\scriptstyle #2};
    \endxy}

\newcommand{\Fixx}[4] {
    \xy
        (0,0)*{#1};
        (-0.2,2.5)*{\ba \ig{loop.pdf}\ea};
        (#3,#4)*{\scriptstyle #2};
    \endxy}

\newcommand{\sFix}[2] {
    \xy
        (0,-1.7)*{\scriptstyle #1};
        (0.2,-0.3)*{\ba \ig{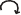}\ea};
        (2,0)*{\scriptscriptstyle #2};
    \endxy}

 \reversemarginpar
\DeclareFontFamily{U}{rsf}{}
\DeclareFontShape{U}{rsf}{m}{n}{
  <5> <6> rsfs5 <7> <8> <9> rsfs7 <10-> rsfs10}{}
\DeclareMathAlphabet{\mathscr}{U}{rsf}{m}{n}

\begin{document}

\title{On unitary 2-representations of finite groups and topological quantum field theory}

\author{Bruce Bartlett\\
 \\
 \\
 \\
 A thesis submitted in fulfilment of the requirements\\
  for the degree of Doctor of Philosophy\\
  to the\\
  University of Sheffield\\
  Department of Pure Mathematics \\
  Oct 2008
  \date{}}

\maketitle
\chapter*{Abstract}
\addcontentsline{toc}{chapter}{Abstract}

This thesis contains various results on unitary 2-representations of finite groups and their 2-characters, as well as on pivotal structures for fusion categories.  The motivation is extended topological quantum field theory (TQFT), where the 2-category of unitary 2-representations of a finite group is thought of as the `2-category assigned to the point' in the untwisted finite group model. 

The first result is that the braided monoidal category of transformations of the identity on the 2-category of unitary 2-representations of a finite group computes as the category of conjugation equivariant vector bundles over the group equipped with the fusion tensor product. This result is consistent with the extended TQFT hypotheses of Baez and Dolan, since it establishes that the category assigned to the circle can be obtained as the `higher trace of the identity' of the 2-category assigned to the point.

The second result is about 2-characters of 2-representations, a concept which has been introduced independently by Ganter and Kapranov. It is shown that the 2-character of a unitary 2-representation can be made functorial with respect to morphisms of 2-representations, and that in fact the 2-character is a unitarily fully faithful functor from the complexified Grothendieck category of unitary 2-representations to the category of unitary equivariant vector bundles over the group. 

The final result is about pivotal structures on fusion categories, with a view towards a conjecture made by Etingof, Nikshych and Ostrik.  It is shown that a pivotal structure on a fusion category cannot exist unless certain involutions on the hom-sets are plus or minus the identity map, in which case a pivotal structure is the same thing as a twisted monoidal natural transformation of the identity functor on the category. Moreover the pivotal structure can be made spherical if and only if these signs can be removed. 

\cleardoublepage
\vspace{20cm}
\section*{\begin{center}Dedication\end{center}}
\addcontentsline{toc}{chapter}{Dedication}
\thispagestyle{empty}
\begin{center}
{\em To mom and dad, granny and gramps, nanny and pa!}
\end{center}
\chapter*{Acknowledgements}
\addcontentsline{toc}{chapter}{Acknowledgements} Firstly my heartfelt thanks goes to my supervisor Simon Willerton whom I had the great fortune of learning so many ideas from these last three years; he has been my mentor and my friend. I wish to thank all my fellow Houses of Maths inmates for the good times that we had, and all the postgrads, present and past, the young 'uns and the old guns.

As far as work is concerned, I wish to especially thank the `locals' David Gepner, James Cranch, Paul Buckingham, Almar Kaid, Dave Stern, Eugenia Cheng, Tom Bridgeland, Vic Snaith and Richard Hepworth for stimulating conversations and for helping me at various points during this thesis.

I would like to acknowledge Michael Shulman , Nick Gurski and Tom Leinster who have helped me with various categorical notions. I am also grateful to Scott Carter for patiently explaining to me the movie moves needed in Appendix \ref{movieapp}.

I am sincerely grateful to John Baez whose column `This Week's Finds' has had such a profound effect on me, and to Urs Schreiber, whom I count myself very lucky to know.

I would also like to thank Jamie Vicary, Mathieu Anel, Frank Neumann, Hellen Colman and Steve Lack.

Finally, I would like to acknowledge the support of the Excellence Exchange Scheme, as well as the Glasgow Mathematical Journal Trust. 
\tableofcontents
\chapter{Introduction}

In this thesis we prove three main results. The first has to do with verifying the `crossing with the circle' equation in the extended finite group TQFT, the second has to do with characterizing pivotal structures on fusion categories, and the third is a proof that the 2-character functor is unitarily fully faithful.

\section{Background}
Before we explain these results, we first recall some of the broad themes of extended topological quantum field theory, so as to establish a context.

\subsection*{Chern-Simons theory and the birth of TQFT}
In 1989, Witten gave a beautiful geometric description of the then recently discovered knot invariant of Jones in terms of a three-dimensional quantum field theory known as {\em Chern-Simons theory} \cite{ref:witten, ref:jones}. From the point of view of ordinary physics this theory is rather odd in that the action is {\em purely topological} --- it does not depend on the metric on spacetime. To obtain a knot invariant this feature is critical for otherwise the path integral would not be invariant under ambient isotopy of the knot. Quantum field theories of this type are known as {\em topological quantum field theories} (TQFT's), and they are particularly amenable to mathematical analysis because many of the infinities which plague full-blown theories do not arise; for instance the Hilbert spaces of the theory are {\em finite-dimensional}. In the years since Witten published his paper, Chern-Simons theory has become the poster-child of TQFT, and a focus point for the recent renaissance in the interaction between geometry and physics. For an introduction to these ideas, we recommend \cite{ref:atiyah, ref:kohno} as well as the original paper of Witten \cite{ref:witten}. For more on the background and history of this interaction, see \cite[Chapter 11]{ref:woit}.

\subsection*{Extended topological quantum field theory}
From a mathematical perspective, the name of the game has been to try and understand the formal aspects of topological quantum field theory. Just what {\em is} a TQFT really, and how does one construct one? What are the geometric ingredients which come into play? 

An answer to the first question was apparently provided right in the beginning by Atiyah and Segal \cite{ref:atiyah2, ref:segal}. From their work it emerged that an $n$-dimensional TQFT can be thought of a {\em functor}
 \[
  Z \colon \text{nCob} \rightarrow \text{Vect}
 \]
from the category nCob of $(n-1)$-dimensional closed manifolds and cobordisms between them to the category of vector spaces. Moreover, this functor is required to preserve the monoidal structures on these categories, so that a disjoint union of manifolds must get sent to the tensor product of their corresponding vector spaces, and similarly the cobordism which switches two manifolds around should get sent to the linear map which swaps the factors in the tensor product of their corresponding vector spaces. Now a closed $n$-manifold $M$ can be interpreted as a cobordism from the empty set to itself, so the fact that the functor $Z$ preserves the monoidal structure implies that $Z$ will send $M$ to a map from the complex numbers to itself, which is the same thing as a {\em number}. So in this picture, a TQFT assigns a {\em number} to a closed $n$-manifold and a {\em vector space} to a closed $(n-1)$-manifold.

Over time though, it was realized that such a formalization failed to capture some important aspects of the main examples of TQFT's. For instance, Reshetikhin, Turaev and others found a setting where the ideas of Witten could be made mathematically rigorous, by interpreting everything in terms of the category of representations of a {\em quantum group} \cite{ref:resh}, which turned out to be a braided monoidal category. It was found that this description and the Atiyah-Segal monoidal functor description were somehow different sides of the same coin, and the problem was to find a formalism where both descriptions had a common home.

Another problem was that the formalization of a TQFT as a monoidal functor only captured a small subset of the gluing laws which actually held in practice. The action in a quantum field theory is usually of a {\em local} nature, which suggests that the theory should be natural with respect to all possible gluing laws of all codimensions, not just gluing two $(n-1)$-manifolds along an $n$-dimensional cobordism.

These considerations gave rise to the notion of an {\em extended} TQFT as one which behaves well with respect to these extra gluing laws \cite{ref:baez_dolan_hda0, ref:freed, ref:freed3, ref:freed_quinn, ref:lawrence, ref:quinn, ref:yetter}. It was realized that just as an ordinary $n$-dimensional TQFT assigns a number to a closed $n$-manifold and a vector space to a closed $(n-1)$-manifold, an {\em extended} TQFT goes further and assigns a {\em linear category} to a closed $(n-2)$-manifold, and a {\em linear 2-category} to a closed $(n-3)$-manifold, and so on. In Chern-Simons theory for instance, where $n=3$, the braided monoidal category of representations of a quantum group --- the approach to the subject adopted by Reshitkin, Turaev and others --- suddenly found a home as the {\em category assigned to the circle}. In other words, the language of {\em higher categories} began to offer a unified approach to the subject, although admittedly this formalism was not to everybody's taste.

\subsection*{The Baez-Dolan hypotheses}
The problem remained to find a concise description of the idea of an extended TQFT in terms of this new language, in a way which captured all the gluing laws. This was the point at which Baez and Dolan released their seminal paper \cite{ref:baez_dolan_hda0}, the first of a long series on `Higher dimensional algebra and topological quantum field theory'. In this paper they made two hypotheses about the nature of extended TQFT's, or at least about the nature of {\em unitary} TQFT's, which are the kind of theories which generally arise in physics. We will state each individual conjecture slightly differently here, but the `total hypothesis' amounts to the same thing. The first hypothesis is as follows.

 \begin{center}
\framebox{\parbox[b]{12cm}{ \textbf{Extended TQFT Hypothesis (Baez-Dolan \cite{ref:baez_dolan_hda0}).} An $n$-dimensional unitary extended TQFT is a weak $n$-functor, preserving all levels of duality, from the $n$-category $n\mathcal{C}$ob of cobordisms to $n\mathcal{H}$ilb, the $n$-category of $n$-Hilbert spaces.}
}
 \end{center}
There are three important ingredients here. The first is the implicit assertion that the mathematical structure whose objects are closed $0$-manifolds (i.e. collections of points), whose `1-morphisms' are 1-manifolds with boundary, whose `2-morphisms' are manifolds with corners, and so on up to $n$, actually forms an {\em $n$-category.} In other words, the implicit assertion is that {\em higher categories are precisely the right language to describe the cutting and pasting behaviour of manifolds with corners}. The theory of higher categories is still young and undeveloped, but this certainly represents a major motivation for the subject; see \cite{ref:cheng_gurski2} for a recent approach in this regard.

The second and most crucial ingredient is the phrase `preserving all levels of duality', which is a hypothesis for what the word {\em unitary} means in an extended TQFT context. The idea is that the $n$-category $n\mathcal{C}$ob is the archetypal `n-category with duals', by which is meant that for every $k$-morphism $f \colon a \rightarrow b$ there is a dual morphism $f^* \colon b \rightarrow a$, together with $(k+1)$-morphisms $\eta \colon \id_a \Rightarrow f^* f$ and $\epsilon \colon ff^* \Rightarrow \id_b$, which themselves have dual morphisms and coherence data of their own, and so on, up to the top level $n$, where we have run out of higher morphisms. The basic idea is that $n$-categories with duals describe the geometry of `how things happen' in a manner akin to the `{\em catastrophes}' of singularity theory \cite{ref:arnold, ref:thom}. That is, performing $f$ and then $f^*$ introduces a `kink', and $\eta$ represents the process of removing this `kink', and so on. This is a profound idea which lies at the very heart of the Extended TQFT Hypothesis, as we are about to see in the next hypothesis, and the reader is urged to read \cite{ref:baez_dolan_hda0} in order to understand it better.

The third ingredient is the target $n$-category, namely the {\em $n$-category of $n$-Hilbert spaces}. An `$n$-Hilbert space' is meant to be an $(n-1)$-category which has the structure and properties of a Hilbert space, appropriately categorified. For instance, a {\em 2-Hilbert space} was defined by Baez \cite{ref:baez_2_hilbert_spaces} essentially as a $\mathbb{C}$-linear category where the {\em hom-sets} are finite-dimensional Hilbert spaces (this `categorifies' the idea that the inner product of two vectors in a Hilbert space is a complex number) which are equipped with  compatible  involutions $* \colon \Hom(x,y) \rightarrow \overline{\Hom(y,x)}$  (this categorifies the equation $(x,y) = \overline{(y,x)}$ in a Hilbert space). In fact, this is the least important ingredient of the Extended TQFT Hypothesis. All that really matters is that the $n$-category which the extended TQFT $Z$ takes values in is some kind of linear category with duals at all levels; the examples will ultimately decide what is the best notion. For instance, although Chern-Simons theory can be formulated as taking values in these higher Hilbert spaces, most other examples of extended TQFT's, such as Rozansky-Witten theory (see \cite{ref:roberts_willerton} and the references therein) and the theories arising in the `Geometric Langlands programme' \cite{ref:benzvi} rather take values in some kind of {\em derived} higher category.

The second hypothesis of Baez and Dolan about the nature of extended TQFT's is as follows.

 \begin{center}
\framebox{\parbox[b]{12cm}{ \textbf{Cobordism Hypothesis (Baez-Dolan \cite{ref:baez_dolan_hda0}.)} The $n$-category $n\mathcal{C}$ob of cobordisms is the free stable $n$-category with duals on one object.}
}
 \end{center}
This is the remarkable assertion that the concept `free $n$-category with duals' is rich enough to capture, in one foul swoop, and in a concise algebraic fashion, all the cutting and pasting behaviour of manifolds! Said differently, it is the statement that {\em $n$-categories with duals are the grammar of space}. An $n$-category is called {\em stable} if it can be regarded as the top $n$-dimensional part of an $(k+n)$-category with $k \gg n$; this condition is there simply to ensure that we are talking about {\em abstract} cobordisms and not `unstable' embedded objects such as {\em tangles}.

It is not our task here to talk about the evidence for these hypotheses; for that the reader is referred to \cite{ref:baez_dolan_hda0} and the rest of the `Higher Dimensional Algebra' series. What we wish to stress is the {\em power} of these hypotheses. For instance, by combining them together we see that they assert the following:
 \begin{center}
\framebox{\parbox[b]{12cm}{\textbf{The primacy of the point$\xy (0,-0.8)*{\bullet}; \endxy$  }  An $n$-dimensional unitary extended TQFT is completely described by the $n$-Hilbert space it assigns to a point.} }
 \end{center}
This is a striking assertion; the Ultimate Statement with respect to the locality of the theory: it says that {\em everything}, including all the quantum operations etc., is determined by `abstract higher-categorical nonsense' from what the theory assigns to a point! Considering that the quantum invariants involved in theories such as  Chern-Simons theory, Rozansky-Witten theory, Seiberg-Witten theory and Geometric Langlands represent some of the deepest connections between physics and geometry, this represents a powerful entry point for higher-categorical ideas into mainstream mathematics.

\subsection*{Verifying the primacy of the point}
In Chapter \ref{ppoint} we will actually {\em `verify'} the `primacy of the point' in the case $n=2$. It has long been known that a unitary two-dimensional TQFT is characterized by $n$ positive real numbers $k_1, \ldots, k_n$, which are the eigenvalues of the hermitian handle creation operator \cite{ref:durhuus_jonsson}. We will see that this is precisely the data which characterizes a 2-Hilbert space up to strong unitary isomorphism. (We have put quotation marks around the word `verify' here because most of the terms of the Extended TQFT Hypothesis have still not been made entirely precise, even for $n=2$; in particular it is not completely clear what it means to `preserve duals at all levels', though a precise formulation of the other aspects of the theory for $n=2$ can be found in the work of Morton \cite{ref:morton})

In fact, Costello has actually proved a far more powerful result along these lines \cite{ref:costello}. He has shown that specifying an extended two-dimensional open-closed topological conformal field theory is the same thing as specifying a {\em Calabi-Yau $A_\infty$-category}. A topological conformal field theory is a very rich structure that lies between a TQFT and a conformal field theory, and the main examples are the $A$ and $B$ models introduced by Witten \cite{ref:witten2}. An $A_\infty$-category can be thought of as a collection of objects with graded vector spaces between them which forms a `category up to coherent homotopy', and a {\em Calabi-Yau} $A_\infty$-category is, roughly speaking, an $A_\infty$-category equipped with a non-degenerate invariant pairing on the spaces of morphisms. One can think of a Calabi-Yau $A_\infty$-category as a `derived version' of a 2-Hilbert space, and thus one can interpret Costello's result as affirming the primacy of the point for $n=2$ in the world of topological conformal field theories.


\subsection*{The `crossing with the circle' equation}
If the extended TQFT hypothesis is indeed the correct framework for capturing the higher gluing laws of an extended TQFT, another spin-off would be that one would expect that the quantum invariants $Z(M)$ assigned to closed manifolds $M$ would satisfy
 \begin{center}
\framebox{\parbox[b]{7cm}{ \[
  Z(M \times S^1) \cong \mbox{Dim} \, Z(M)
 \]}}
 \end{center}
where `Dim' refers to the {\em higher-categorical dimension} of $Z(M)$, a sort of higher-categorical trace operation, defined as the collection of transformations of the identity functor on $Z(M)$. We call this the `crossing with the circle' equation and it makes some interesting predictions. Firstly, let us be clear: this equation is claiming that the {\em entire plethora of topological operations} on the higher-categorical structure $Z(M \times S^1)$ --- such as the product arising from the higher pair of pants, the braiding, the Frobenius pairing etc. --- can be computed purely from  `higher-categorical abstract nonsense' on $\Dim Z(M)$. This gets to the heart of the entire motivation for studying extended TQFTs: {\em topological operations translate into higher-categorical operations}.

Consider the implications that this equation would hold for Chern-Simons theory. When $M = S^1$ it says that the Hilbert space assigned to the torus --- known as the {\em Verlinde algebra} --- can be calculated as the algebra of natural transformations of the identity functor on the category assigned to the circle. Now the category assigned to the circle is the category of positive energy representations of the loop group at level $k$. Since this category is semisimple, the algebra of natural transformations has a natural basis consisting of the natural transformations supported at the irreducible representations. We thus conclude that `the characters diagonalize the fusion rules'. This latter statement is in fact a celebrated result known as the {\em Verlinde conjecture} --- and we have just seen how it would arise in principle as an elementary byproduct of formulating the theory in higher categorical terms.

We also remark that some recent computations in this regard, in the context of geometric representation theory, can be found in the work of Ben-Zvi, Francis and Nadler \cite{ref:benzvi_francis_nadler}.

\subsection*{The Freed picture of extended field theory}
So far we have been answering the question ``What {\em is} an extended TQFT?'', but we have had nothing to say about the geometric ingredients one needs to actually {\em construct} such a theory. Freed has established a beautiful picture in this regard, and everything we do in this thesis has been heavily influenced by his ideas \cite{ref:freed, ref:freed2, ref:freed3}. In particular we warmly recommend the recent \cite{ref:freed_remarks}, which gives an elegant overview of the entire framework of extended Chern-Simons theory along these lines.

In the framework of Freed, the geometric data needed to construct an $n$-dimensional extended TQFT consists of the following (it should be understood that this is a `working hypothesis'; much still needs to be done in order to make it precise):
\begin{itemize}
 \item For every manifold $M$ with $\dim M \leq n$, a {\em space of fields} $\mathcal{P}_M$ on $M$. The main requirement is that these fields must be {\em local}, so that they behave well under cutting and pasting. For instance in Chern-Simons theory based on a gauge group $G$ the space of fields $\mathcal{P}_M$ is the space of flat $G$-bundles with connection on $M$.
 \item For every such manifold $M$, a {\em higher-line bundle} $\mathcal{L} \rightarrow \mathcal{P}_M$ over the space of fields $\mathcal{P}_M$. The fibers of the line bundle $\mathcal{L}$ are one-dimensional $(n- \dim M)$-Hilbert spaces, so that when $\dim M = n$ it is a `bundle of numbers of unit norm', i.e. a {\em function} $e^{i S[\cdot]} \colon \mathcal{P}_M \rightarrow U(1)$, while if $\dim M = n-1$ it is a `bundle of hermitian lines' i.e. a conventional {\em hermitian line bundle}, while if $\dim M = n-2$ it is a `2-line bundle', i.e. an assignment of a one-dimensional 2-Hilbert space to every field $P \in \mathcal{P}_M$, and so on down to the point, where $\dim M = 0$, which has a one-dimensional $n$-Hilbert space sitting above it.
 \end{itemize}
The beauty of this picture is the startling new interpretation given to the age-old idea of the `action' of a field theory: it is now reinterpreted as {\em a machine which constructs higher line bundles over the spaces of fields}. The {\em quantum} side of the story is just as attractive: at the level of closed manifolds the quantum theory simply assigns to a closed manifold $M$ the space of sections of these higher line bundles:
 \[
  Z(M) = \Gamma(\mathcal{L}).
 \]
This picture suggests that many of our conceptual hang-ups about the path integral in quantum field theory might have arisen from the fact that we gave too much credence to the idea of the action as a {\em number} --- this is just the top-dimensional `degenerate case' of what the action really is, and things look a lot more palatable at higher codimension! For instance, most mathematicians feel queasy when instructed to `add up the value of the action over all the fields', but have far less psychological qualms when instructed to `take the space of sections of a line bundle'.

The question now arises as to how one actually comes up with this geometric data of the higher line bundles over the spaces of fields. In most cases of interest these bundles arise from a mechanism of mapping into some fixed `target space' $X$ which has a `master higher-line bundle' $\mathcal{L} \rightarrow X$ sitting above it, and then transferring this master line bundle to the various spaces of fields $\mathcal{P}_M$ by a process known as {\em transgression}. In physics this mechanism is known as a {\em $\sigma$-model}. Schreiber and Waldorf have recently made much progress towards a higher-categorical understanding of this mechanism \cite{ref:schreiber_waldorf1, ref:schreiber_waldorf2, ref:schreiber}. In particular, they showed how the process of transgression can be elegantly interpreted higher-categorically simply as postcomposition with the higher transport functor. In addition, Sati, Schreiber, Sk\u{o}da and Stevenson have recently outlined a theory of twisted nonabelian cohomology in a smooth $\omega$-groupoid setting \cite{ref:sati_et_al}, providing a concrete toolbox for understanding the geometry of these `master higher-line bundles'.

At this point, we turn our attention to a specific toy model where everything can be made quite explicit: the {\em finite group model}.

\subsection*{The finite group model}
For each dimension $n$, Dijkgraaf, Witten and Freed have shown how there is an extended
TQFT called the {\em finite group
model} \cite{ref:dijkgraaf_witten, ref:freed}. The input for this theory is a finite group $G$ and an $n$-cocycle $\omega \in
Z^n(G, U(1))$. When $n=3$, it is a finite group analogue
of Chern-Simons theory, and one thinks of the 3-cocycle $\omega$ as a low-level description of a `3-line bundle' over $\BG$, or alternatively a `2-gerbe' over $\BG$. Here and elsewhere in this thesis, we use the notation that $\BG$ refers to the group $G$ thought of as a one-object category (the use of boldface avoids confusion with  `$BG$', the geometric realization of the nerve of $\BG$).

Willerton has described an elegant way to set up the geometric data of the finite group model via gerbes and finite groupoids \cite{ref:simon}, and we recall this description here. In the finite group model, the space of fields $\mathcal{P}_M$ over a
manifold $M$ (with $\dim M \leq n$) is the groupoid whose objects are the principal
$G$-bundles over $M$ and whose morphisms are equivariant isomorphisms of $G$-bundles. A finite model for $\mathcal{P}_M$ is obtained by choosing a finite number of points on $M$ (at least one for each component), and
denoting by $\Pi_1(M)$ the restriction of the fundamental groupoid of $M$ to
these points. It is convenient to allow more than one point on each
component, as this makes things more natural if we are considering
spaces with boundaries. Nevertheless, since a finite principal
$G$-bundle is determined by its holonomy, we have
 \[
  \mathcal{P}_M = \Fun (\Pi_1 (M) , \BG),
 \]
the category of functors and natural transformations from $\Pi_1(M)$ into $\BG$.

The action is defined by transgression, and works as follows. One
uses the diagram
 \[
  \mathcal{P}_M \la{\pi} \mathcal{P}_M \times \Pi_1(M)
  \ra{\text{ev}} \BG,
 \]
to pull back and push down the cocycle $\omega$ onto $\mathcal{P}_M$
to give a $(d-\dim M)$-cocycle $\tau_M (\omega) \in
Z^{(d-n)}(\mathcal{P}_M, U(1))$.

The quantum theory assigns to each closed manifold $M$ the space of sections of the higher line bundle over $\mathcal{P}_M$ represented by the transgressed cocycle $\tau_M(\omega)$.

In Table 1 we tabulate the resulting quantum invariants in the three-dimensional theory, which is the context of this thesis. We mention again the disclaimer here that the notion of a `3-dimensional extended TQFT' is not yet a precise well-defined concept --- working towards this point is the entire purpose of the Baez-Dolan extended TQFT hypotheses.

\vspace{0.5cm}
\begin{table}[t]
\begin{tabular}{cccl}
 $\dim$ & space & $Z_3$ & explanation \\
 \hline \vspace{0.2cm}
 3 & $M$ & $\int_{\mathcal{P}_M} \tau_M (\omega)$ & no. of $G$-bundles $P$ on $M$, each weighted by $\frac{\tau_M(\omega)}{|\Aut(P)|}$ \vspace{0.2cm}  \\
 2 & $\Sigma$ & $\Gamma_{\mathcal{P}_\Sigma} (\tau_\Sigma(\omega))$ & flat
 sections of $\tau_\Sigma(\omega)$-twisted line bundle on $\mathcal{P}_\Sigma$ \vspace{0.2cm} \\
 1 & $S^1$ & $\Hilb_G^{\tau(\omega)} (G)$ & \parbox{7cm}{category of $\tau(\omega)$-twisted equivariant vector bundles over $G$} \vspace{0.2cm} \\
 0 & pt & $\TRep^\omega(G$) & 2-category of $\omega$-twisted 2-representations of $G$
\end{tabular}

 \caption{The quantum invariants assigned to closed manifolds in the twisted three-dimensional finite group model.}
 \end{table}

For more details about the invariants assigned to 1, 2 and 3-dimensional closed manifolds and the notation we are using in Table 1, we refer the reader to the work of Willerton \cite{ref:simon}, though we will review this technology in Chapter \ref{GerbesChap}. The entry we wish to explain in the above table is the quantum invariant $Z(\text{pt})$ assigned to a {\em point}. The space of fields $\mathcal{P}_\text{pt}$ over a point is the groupoid of $G$-bundles over a point, which can be thought of simply as $\BG$. The `higher 3-line bundle' $\mathcal{L} \rightarrow \mathcal{P}_\text{pt}$ over the space of fields is obtained by using the 3-cocycle $\omega$ to twist the {\em associator} on $\THilb$, the 2-category of 2-Hilbert spaces, so as to form a `bundle of one-dimensional 3-Hilbert spaces over $\mathcal{P}_\text{pt} = \BG$':
 \[
  \THilb \times_\omega \BG \rightarrow \BG.
 \]
The quantum invariant assigned to the point is then the space of flat sections of this bundle. This procedure is equivalent to thinking of the 3-cocycle $\omega$ on $\B G$ as a recipe for constructing a {\em 2-group} $(G, \omega)$ (see \cite{ref:baez_lauda_2-groups, ref:martins_porter}), and then setting $Z(\text{pt})$ to be the 2-category of `2-representations' of this 2-group, as in \cite{ref:elgueta}. In fact, because we want to ultimately construct a {\em unitary} extended TQFT, we must require that these 2-representations be `unitary'. Let us emphasize the final result:

 \begin{center}
\parbox[b]{11cm}{ In the twisted 3d finite group theory, the quantum invariant assigned to the point is expected to be the 2-category of unitary 2-representations of the 2-group $(G,\omega)$.}
 \end{center}

Of course, in the rest of the thesis we will make precise what we mean by `unitary 2-representations' of a group. We will essentially be restricting ourselves to the {\em untwisted} model, where $\omega = 1$, not because of a lack of technology to deal with the twisted case but simply because the untwisted case has proved a project enough in itself.

\section{Verifying $Z(S^1) \simeq  \text{Dim }  Z(\text{pt})$ in the finite group model}
Having a concrete toy model of extended TQFT like the finite group model allows us to test some of the predictions of the Extended TQFT Hypothesis. For instance, we can use the `cocycles on finite groupoids' technology of Willerton \cite{ref:simon} to check that the `crossing with the circle' equation indeed holds --- in the full twisted model --- in all dimensions where $Z(M)$ is a category or a vector space. We have collected some of our computations in this regard in Appendix \ref{AppendixTQFTFacts}.

The next case to check is when $Z(M)$ is a {\em 2-category}. For instance, in the three-dimensional model this requires that we take $M$ to be a point. As we have seen, in the three-dimensional untwisted model $Z(\text{pt})$ is expected to be the 2-category of unitary 2-representations of the group $G$, while $Z(S^1)$ is the braided monoidal category $\Hilb^\text{fusion}_G(G)$ of conjugation equivariant vector bundles over the group under the fusion tensor product (see Appendix \ref{AppFreed}). In Chapter \ref{Dim2RepGchap} we verify that the `crossing with the circle' equation indeed holds in this case, and this represents our first main theorem in this thesis --- though it will be proved last, in Chapter \ref{Dim2RepGchap}.

\begin{thmt} The higher-categorical dimension of the 2-category of unitary 2-representations of a finite group $G$ is equivalent, as a braided monoidal category, to the category of conjugation-equivariant hermitian vector bundles over $G$ equipped with the fusion tensor product. In symbols,
 \[
  \Dim \TRep(G) \simeq \Hilb_G^\fusion (G).
 \]
\end{thmt}
As we have intimated, this result is significant because the braided monoidal structure on $\Hilb_G^\fusion (G)$ has hitherto been obtained from topological considerations on the pair of pants cobordism, as explicated by Freed \cite{ref:freed}. Our result shows that we can also compute this braided monoidal structure directly from `abstract nonsense' on $\TRep(G)$.

\section{The 2-character functor is unitarily fully faithful}
In the previous sections we have outlined at great length our motivation for being interested in the 2-category of unitary 2-representations of a finite group --- it is the {\em 2-category assigned to the point in the three-dimensional untwisted finite group TQFT}. Except for Chapters \ref{EvenHandedchap} and \ref{FusionChap} which are about even-handed structures and fusion categories and can be read separately from the rest, this 2-category $\TRep(G)$ is the entire focus of this thesis: {\em defining} it, explicating the {\em geometry} behind how it works, and showing how the main results about {\em ordinary} representations of groups have categorified analogues for {\em 2-representations}.

Our main result in this latter regard is the following. We define what it means to take the {\em 2-character} of a unitary 2-representation (this definition has also been given independently by Ganter and Kapranov \cite{ref:ganter_kapranov_rep_char_theory}), so as to obtain a {\em unitary conjugation equivariant vector bundle over the group $G$}. This is to be thought of as the categorification of a class function on $G$.

Then we show how one can also take the 2-character of a {\em morphism} of unitary 2-representations, so as to obtain a {\em morphism} between their corresponding vector bundles over $G$ (this latter notion was not considered in \cite{ref:ganter_kapranov_rep_char_theory}). Just as the ordinary character of a representation does not depend on the isomorphism class of the representation, the 2-character of a {\em morphism} of 2-representations does not depend on its isomorphism class, hence the 2-character descends to a functor from the {\em Grothendieck category} of unitary 2-representations to the category of conjugation equivariant vector bundles over the group. Our main result is that after one tensors the hom-sets in $[\TRep(G)]$ with $\mathbb{C}$, the 2-character functor becomes {\em unitarily fully faithful}:

\begin{thmt} The 2-character functor
 \[
  \chi  \colon [\TRep(G)]_\mathbb{C} \rightarrow \Hilb_G(G)
  \]
is a unitarily fully faithful functor from the complexified Grothendieck category of unitary 2-representations of $G$ to the category of unitary conjugation equivariant vector bundles over $G$.
\end{thmt}
This is to be thought of as a categorification of the fact that the ordinary character of ordinary unitary representations
 \[
  \chi \colon [\Rep(G)]_\mathbb{C} \rightarrow \text{Class}(G)
 \]
is a unitary isomorphism from the complexified Grothendieck group of the category of unitary representations of $G$ to the space of class functions on $G$.

To prove this result, we establish a geometric correspondence between unitary 2-representations and {\em finite equivariant gerbes equipped with metrics}, and we show that under this correspondence the `2-character' of a 2-representation corresponds to the {\em `geometric character'} of its associated equivariant gerbe. It turns out that the behaviour of the geometric character on the hom-sets is a rearrangement of the {\em twisted character map} of Willerton \cite{ref:simon}, and thus it follows from the technology of \cite{ref:simon} that the geometric character map is unitarily fully faithful, and hence so also is the 2-character.

Of course, we have swept a great deal under the carpet in this terse description, and we now briefly expand in a bit more detail on two points which we have glossed over in the above argument: the geometric correspondence for {\em ordinary} representations of groups and their characters, and the notion of an {\em even-handed structure} on a 2-category.

\subsection*{Geometry of ordinary representations of groups and their characters}
The basic idea of geometric quantization in the equivariant context is that every representation of a Lie group $G$ arises as the `quantization' of a classical geometric system having symmetry group $G$. Normally this is expressed in the language of symplectic geometry and polarizations \cite{ref:woodhouse}, and it can become quite intricate. However, we found it conceptually useful to set up a simpler but more formal version of this correspondence, in order that our later results on the correspondence between unitary 2-representations and `finite equivariant gerbes equipped with metrics' could be viewed in the right context. We showed that the category of unitary representations of a group is equivalent to a certain category of {\em equivariant hermitian holomorphic line bundles over compact hermitian manifolds}, and moreover that under this correspondence, the {\em character} of a representation corresponds to the {\em geometric character} of its associated equivariant line bundle, defined as an integral of a certain kernel over the base manifold. When the line bundle is sufficiently positive, this integral can be re-expressed using the Atiyah-Segal equivariant index theorem as an integral over the {\em fixed points} of the group action. These are precisely the ideas which get categorified for equivariant gerbes in Chapter \ref{GerbesChap}.

We admit this is mostly a formal result, simply revolving around the idea  every finite dimensional Hilbert space $V$ identifies antilinearly as the space of sections of a holomorphic line bundle --- the {\em hyperplane line bundle} over the projective space of $V$. But it is useful to have this correspondence at hand, in order to compare it with its `categorified version', the correspondence between 2-representations and equivariant gerbes. We also add the disclaimer here that whereas the geometric correspondence for ordinary unitary representations of Lie groups (finite, compact or noncompact) involves `genuine' geometric objects such as complex manifolds and holomorphic line bundles, the geometric correspondence we set up for unitary 2-representations of finite groups is only a {\em toy model} since the `equivariant gerbes' we talk about are {\em finite}. This has to do with the fact that a 2-Hilbert space is semisimple and is thus too discrete to carry continuous geometry; it is not the final word on what a `categorified finite-dimensional inner product space' should be, but only a first approximation. We will discuss this point further in Chapter \ref{2HilbChap}.

\subsection*{Even-handed structures on 2-categories}
In order to define the 2-character of a morphism of unitary 2-representations, one needs to have good control over adjoint functors in the 2-category of 2-Hilbert spaces, or geometrically speaking, one needs to ensure that one can choose a flat section of the `ambidextrous adjunction bundle'. We call this an {\em even-handed structure}, and the behaviour of the 2-character on morphisms will in general {\em depend} on the choice of this structure. However, we show that the 2-category of 2-Hilbert spaces has a canonical even-handed structure, which uses the inner products and duality on the hom-sets in an essential way --- this is analogous to the way that the adjoint of a linear map between vector spaces requires an inner product. An algebraic geometer will recognize this as the assertion that a 2-Hilbert space comes equipped with a canonical `trivial Serre functor'. This is a  feature which makes the theory of unitary 2-representations richer than the theory of unadorned 2-representations on the 2-vector spaces of Kapranov and Voedodsky \cite{ref:kapranov_voevodsky} --- the inner products enter in an essential way.

The idea of an even-handed structure on a 2-category is a general one though, and can be applied to any 2-category where the 1-morphisms have `ambidextrous duals', in the sense that for every morphism $F \colon A \rightarrow B$ there exists a morphism $F^* \colon B \rightarrow A$ that is simultaneously left and right adjoint to $F$. When the 2-category has only one object, so that it can be regarded as a monoidal category, then an even-handed structure is essentially the same thing as a {\em pivotal structure} on the monoidal category in the sense of Joyal and Street \cite{ref:joyal_street_bmc2} and Freyd and Yetter \cite{ref:freyd_yetter}.

However, it is our contention that the notion of an even-handed structure has a number of conceptual advantages over the notion of a pivotal structure, and we explain these advantages in Chapter \ref{EvenHandedchap}. Chiefly, it actually fits the examples we have in mind, such as the 2-category of 2-Hilbert spaces, and also it meshes well with the {\em string diagram notation} for working with 2-categories which we use throughout this thesis. This latter reason might seem purely aesthetic, but it is important because by translating pivotal structures into the language of even-handed structures we have been able to precisely characterize the set of pivotal structures on a {\em fusion category}, making progress on a conjecture that was made by Etingof, Nikshych and Ostrik \cite{ref:eno}. This is our third main result in this thesis, to which we now turn.

\section{Characterizing pivotal structures on fusion categories}
A {\em fusion category} is a semisimple linear monoidal category where every object has a dual. They have been much studied in the field of quantum algebra (see the lecture notes of M\"{u}ger \cite{ref:müger_lectures} for a recent overview). A seminal paper in this regard has been that of Etingof, Nikshych and Ostrik \cite{ref:eno}, and one of the conjectures made there was that {\em every fusion category admits a pivotal structure}. Since an `even-handed structure' on a monoidal category is roughly the same thing as a pivotal structure, except that it meshes better with the string diagram notation, we wanted to see if using the framework of even-handed structures would help make progress on this conjecture.

Indeed, using a string diagram observation by Hagge and Hong \cite{ref:hagge_hong}, we were able to translate many of the results of \cite{ref:eno} into a purely string diagrammatic framework. This enabled us to identify a pivotal structure on a fusion category as a {\em twisted monoidal natural transformation of the identity functor on the category}. The twisting is governed by a collection of signs $\{\epsilon^i_{jk}\}$, one for each triple of simple objects in the category, which arise from the nature of the duality operation in the category. Namely, after choosing signs for the square roots of the {\em paired dimensions} of M\"{u}ger \cite{ref:müger_from_subfactors_to_categories_and_topology_I}, one can canonically define certain {\em involution operators} $T^i_{jk}$ on the hom-sets of the form $\Hom(X_i, X_j \otimes X_k)$ where $X_i, X_j$ and $X_k$ are simple objects. We have shown that one of the requirements for the existence of a pivotal structure on the category is that $T^i_{jk} = \pm \id$ for all $i,j,k$; if this is the case we record these signs in the symbols $\epsilon^i_{jk} = \sign(T^i_{jk})$. We call these signs the {\em pivotal symbols} since they are somewhat analogous to the {\em 6j symbols} \cite{ref:carter_flath_saito} which govern the associator information in a fusion category.

Making a different choice of signs for the square roots of the paired dimensions would have resulted in different signs for the involution operators $T^i_{jk}$, and we formalize this by saying that the equivalence class of the signs $\{ \epsilon^i_{jk} \}$ of the involution operators $T^i_{jk}$ in a fusion category gives rise to a class $[\epsilon] \in H_\text{piv}([C], \mathbb{Z}/2)$ in the {\em `pivotal cohomology'} of the category. We show that $C$ can carry a {\em spherical} pivotal structure if and only if this class is trivial (a `spherical' pivotal structure is a notion due to Barrett and Westbury \cite{ref:barrett_westbury} and refers to a pivotal structure where the `left and right traces' of endomorphisms coincide).

Given a collection of signs $\{\epsilon^i_{jk}\}$, we define an {\em $\epsilon$}-twisted monoidal natural transformation of the identity on a fusion category $C$ as a collection of numbers $\{t_i\}_{i \in I}$ where $I$ indexes the simple objects satisfying
 \[
 t_j t_k = \epsilon^i_{jk} t_i \quad \text{whenever $X_i$ appears in $X_j \otimes X_k$.}
 \]
We write the collection of solutions to these equations as $\Aut_\otimes^\epsilon(\id_C)$. If $\epsilon^i_{jk} \equiv 1$ for all $i,j,k$ then this is the same thing as a monoidal natural transformation of the identity on $C$ due to the semisimplicity of $C$.

Let us record our result; it is phrased in the formalism of {\em even-handed structures} but these are equivalent to pivotal structures, as we show in Chapter \ref{pivssec}.
\begin{thmt} Let $C$ be a fusion category over $\mathbb{C}$ with representative simple objects $X_i$. Suppose that a choice of roots $d_i^2 = d_{\{i,i^*\}}$ of the paired dimensions has been made, with resulting involution operators $T^i_{jk} \colon \Hom(X_i, X_j \otimes X_k) \rightarrow \Hom(X_i, X_j \otimes X_k)$.  Then:
 \begin{enumerate}
  \item Unless $T^i_{jk}=\pm \id$ for all $i,j$ and $k$, the fusion category $C$ cannot carry an even-handed structure.
  \item Suppose that $T^i_{jk} = \epsilon^i_{jk} \id$ for all $i,j$ and $k$, where $\epsilon^i_{jk} = \pm 1$. Then an even-handed structure on $C$ is the same thing as an $\epsilon$-twisted monoidal natural transformation of the identity on $C$. That is, there is a canonical bijection of sets
       \[
        \EvenHanded(C) \cong \Aut_\otimes^\epsilon (\id_C).
       \]
   \item Furthermore, the even-handed structure can be made spherical if and only if $[\epsilon] = 0$ in $H_\text{piv}([C], \mathbb{Z}/2)$.
  \end{enumerate}
\end{thmt}
We emphasize that this result does not yet settle the conjecture of Etingof, Nikshych and Ostrik; it might turn out for instance that $T^i_{jk} = \pm \id$ {\em always} holds in a fusion category, and moreover that the class $[\epsilon]$ is {\em always} trivial. But we hope that it does at least clarify some of the issues involved.

\section{Comparison with previous work} There have already been a number of works on 2-representations of groups on semisimple linear categories, such as that of Barrett and Mackaay \cite{ref:barrett_mackaay}, Crane and Yetter \cite{ref:crane_yetter} and Ostrik \cite{ref:ostrik},  with the most relevant for this thesis being that of Elgueta \cite{ref:elgueta} and that of Ganter and
Kapranov \cite{ref:ganter_kapranov_rep_char_theory}. Elgueta performed a thorough
and careful investigation of the 2-category of 2-representations of a 2-group (a 2-group is a monoidal category with structure and properties analogous to that of a group \cite{ref:baez_lauda_2-groups}) acting on Kapranov and Voevodsky's 2-vector spaces, and his motivation was therefore to work with co-ordinatized versions of 2-vector spaces amenable for direct computation, and to classify the various structures which appear.

Ganter and Kapranov were motivated by equivariant homotopy theory, namely to try and find a categorical construction which would produce the sort of generalized group characters which crop up in Morava $E$-theory; they (independently) introduced the {\em categorical character} (which we call the 2-character) of a 2-representation and showed that it indeed achieves this purpose. Since they had no reason not to, they also worked with co-ordinatized 2-vector spaces (this time of the form $\Vect^n$); also they did not investigate in any depth morphisms and 2-morphisms of 2-representations.

As we have explained, we have been motivated by extended TQFT, where the 2-category of unitary 2-representations of a group appears as the 2-category assigned to the point in the untwisted three-dimensional finite group model. The finite group model itself is just a stepping stone, since our real goal is to try and understand more deeply the extended nature of {\em Chern-Simons theory} \cite{ref:freed_remarks, ref:freed_teleman_hopkins}. The language of Chern-Simons theory is the geometric language of moduli-stacks, line-bundles, equivariant structures, flat sections and such like, and this has therefore motivated our approach to 2-representations and is what distinguishes our approach from previous approaches, though we remark that related ideas do appear in \cite{ref:ganter_kapranov_rep_char_theory}.   For instance, as far as possible we try to work directly with the underlying 2-Hilbert spaces of the 2-representations themselves as opposed to some `co-ordinitization' of them, a strategy which might be important in a more intricate geometric setting.  We hope that some of the ideas we have developed in this paper will also translate into the more advanced geometric contexts of the paper of Ganter and Kapranov \cite{ref:ganter_kapranov_rep_char_theory}, as well as the work of Freed, Teleman and Hopkins \cite{ref:freed_teleman_hopkins}.

\section{Overview of thesis}
In Chapter \ref{GeometryOrdinaryRepChap} we establish the geometric correspondence between unitary representations of groups and equivariant holomorphic hermitian line bundles over compact hermitian manifolds, and we show that the character of a representation corresponds to the geometric character of its associated equivariant line bundle.

In Chapter \ref{2HilbChap} we review the notion of a 2-Hilbert space due to Baez \cite{ref:baez_2_hilbert_spaces}. We will be framing some of the basic definitions of the theory a bit differently, attempting to steer clear of using the concept of an `abelian category', since it is likely that more refined notions of `2-Hilbert spaces' in the future will have more of a {\em derived} flavour. We also show that the 2-category of 2-Hilbert spaces is equivalent to a certain 2-category whose objects are finite sets equipped with positive real scale factors; this is the basis for the geometric correspondence between unitary 2-representations and finite equivariant gerbes equipped with `metrics'.

In Chapter \ref{stchap} we review the string diagram notation for working with 2-categories, and we prove that this notation is perfectly well-defined even in the fully weak setting; it is not only a notation for {\em strict} 2-categories, which is a general misconception.

In Chapter \ref{EvenHandedchap} we define even-handed structures on 2-categories, and we compare this notion to other similar notions of duality. When the 2-category consists of linear categories, linear functors and natural transformations, we show how an even-handed structure arises from the data of a {\em trace} on the hom-sets in the category. In this way we establish that the 2-category $\THilb$ of 2-Hilbert spaces, and the 2-category $\CYau$ --- whose objects are the graded derived categories $D(X)$ of Calabi-Yau manifolds --- are both equipped with canonical even-handed structures.

In Chapter \ref{FusionChap} we apply our notion of an even-handed structure to fusion categories, and we obtain our aforementioned result, that a pivotal structure on a fusion category corresponds to a twisted monoidal natural transformation of the identity on the category, and that moreover the twist must vanish if the category is to admit a spherical pivotal structure.

In Chapter \ref{2RepsChap} we define unitary 2-representations of finite groups on 2-Hilbert spaces, and we give some examples. Using string diagrams, we define how to take the 2-character of a unitary 2-representation, and using the even-handed structure on $\THilb$ we also define how to take the 2-character of a {\em morphism} of unitary 2-representations, so that the 2-character descends to functor from the Grothendieck category of unitary 2-representations to the category of conjugation equivariant unitary vector bundles over the group $G$.

In Chapter \ref{GerbesChap} we define finite equivariant gerbes equipped with metrics, and the 2-category they constitute. We define the geometric character of an equivariant gerbe in terms of the transgressed line bundle, and we show how the geometric character also becomes a functor from the Grothendieck category of equivariant gerbes to the category of equivariant vector bundles over the group $G$. Finally we show that after tensoring the hom-sets with $\mathbb{C}$ this functor is unitarily fully faithful.

In Chapter \ref{GerbesCharChap} we establish the correspondence between unitary 2-representations and finite equivariant gerbes equipped with metrics, and we show that the 2-character of a 2-representation corresponds to the geometric character of its associated equivariant gerbe. This enables us to show that the 2-character is also unitarily fully faithful.

Finally in Chapter \ref{Dim2RepGchap} we compute the `higher-categorical dimension' of the 2-category of unitary 2-representations of $G$, and we show that it is equivalent as a braided monoidal category to the category of conjugation equivariant vector bundles over the group under the fusion tensor product. We do this by working in the framework of equivariant gerbes, where one can write down canonical formulas without having to make choices. This implies the corresponding result for 2-representations because of the equivalence between these two categories, at least under a certain technical assumption which is believed indeed to hold.

In Appendix \ref{AppendixTQFTFacts} we verify the `crossing with the circle' equation in the twisted finite group model for low codimension. In Appendix \ref{AppendixFusion} we prove a result about duals in fusion categories that we use in Chapter \ref{FusionChap}. In Appendix \ref{movieapp} we show that if there is a coherent adjoint equivalence between 2-categories then there is a well-defined braided monoidal equivalence between their higher-categorical dimensions. In Appendix \ref{AppFreed} we recall the fusion product on the category of conjugation equivariant vector bundles over a finite group, as in Freed \cite{ref:freed}.

\chapter[Geometry of ordinary representations]{The geometry of ordinary representations of groups\label{GeometryOrdinaryRepChap}}
We begin this thesis by laying down a context for our results on unitary 2-representations and their 2-characters in later chapters. Namely, in this chapter we show how the category of unitary representations of a Lie group $G$ can be regarded as being equivalent to a certain category of {\em equivariant line bundles}, and that under this equivalence the {\em character} of a representation corresponds to the {\em geometric character} of the line bundle. The entire purpose of this chapter is to show that the geometric correspondence we develop in later chapters between unitary 2-representations and equivariant gerbes, and between 2-characters of unitary 2-representations and geometric characters of equivariant gerbes, also has an analogue at the level of {\em ordinary} representations.

For the most part, what we do in this chapter is quite formal, and simply revolves around the idea that every vector space can be thought of as the space of sections of a line bundle. Elementary as it is, this idea can nevertheless be quite profound, for instance to a category theorist it is the {\em Yoneda lemma}. Nevertheless, we do not claim that the results in this chapter are particularly important; we only wish to set the scene for our later results on the geometry of {\em categorified} representations of groups.

In the first two sections we will review known material. In Section \ref{B1} we recall the definitions of some basic geometric notions, such as holomorphic line bundles and kernels, and in Section \ref{B2} we review the {\em Bergman kernel}, the canonical kernel on a holomorphic hermitian line bundle over a compact hermitian manifold. We present three ways of thinking about it: as a sum over the fundamental modes of propagation, as the large time limit of the heat kernel, and as a limit of a path integral. We also relate the Bergman kernel to the {\em coherent states} framework. After having reviewed this material, the rest of the chapter is our own contribution, except where we say otherwise.

In Section \ref{B4} we define the category of line bundles and kernels. In Section \ref{B5} we define the {\em geometric line bundle} over projective space, which is an alternative notation for the {\em hyperplane line bundle} available in the presence of inner products, and we show how the vector space can be recovered antilinearly as the space of sections of the geometric line bundle. In Section \ref{B6} we establish that the category of finite dimensional Hilbert spaces and the category of line bundles and kernels are equivalent. In Section \ref{B7} we show that this correspondence continues to hold in the presence of a $G$-action; this establishes the equivalence between the category of unitary representations of $G$ and the category of equivariant line bundles.

Finally in Section \ref{B8} we define the {\em geometric character} of an equivariant line bundle, and we show that that the character of a representation computes as the geometric character of its associated equivariant line bundle. It is precisely this theorem which we will `categorify' in Theorem \ref{transthm} in Chapter \ref{GerbesCharChap}. We close the chapter by giving a concrete formula for the geometric character using {\em index theory}, valid in the case when the curvature of the line bundle is sufficiently positive. This formula computes the geometric character as an integral over the {\em fixed points} of the group action. In equation \ref{needlabel} of Chapter \ref{GerbesChap} we will give the analogue of this formula for the geometric character of an {\em equivariant gerbe}.

\section{Holomorphic hermitian line bundles and kernels\label{B1}}
In this section we recall some basic geometric notions. A basic reference is \cite{ref:wells}.

A {\em hermitian metric} on a complex manifold is a Riemannian metric which preserves the complex structure. A {\em hermitian manifold} is a complex manifold equipped with a hermitian metric. Note that the inner products on the tangent spaces are only required to vary smoothly with the basepoint and not holomorphically; in particular every complex manifold can be equipped with a hermitian metric.

A {\em holomorphic line bundle} $L$ over a complex manifold $X$ is a line bundle over $X$ whose total space $L$ is a complex manifold and whose projection map $\pi \colon L \rightarrow X$ is holomorphic. A {\em holomorphic section} of $L$ is a holomorphic map $s \colon X \rightarrow L$ such that $\pi \circ s = \id$. We will write the vector space of holomorphic sections of a line bundle $L$ as $\Gamma(L)$. Note that $\Gamma(L)$ must be finite dimensional if $X$ is compact.

A holomorphic line bundle $L$ is called a {\em hermitian holomorphic line bundle} if each fiber is equipped with a hermitian inner product and if these inner products vary smoothly with the basepoint. Note that every holomorphic line bundle can be equipped with a hermitian inner product.

All of the line bundles we deal with in this chapter will be holomorphic. We will often simply say `line bundle' or `hermitian line bundle' instead of `holomorphic line bundle' and `holomorphic hermitian line bundle'.

If $L$ is a hermitian line bundle over a compact hermitian manifold $X$, then the space of holomorphic sections $\Gamma(L)$ of $L$ can be given an inner product by integrating the fibrewise inner product with respect to the volume form defined by the hermitian metric on $X$:
 \[
  \langle s, s' \rangle = \int_X (s(x), s'(x)) \, \vol_x.
 \]

\begin{figure}[t]
\centering
\ig{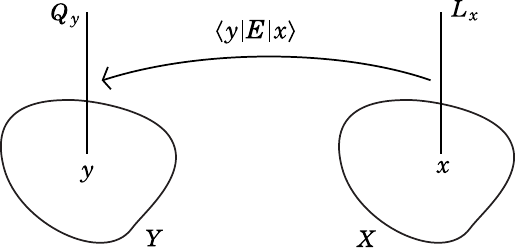}
\caption{\label{kernel1}  A holomorphic kernel from a line bundle $L \rightarrow X$ to another line bundle $Q \rightarrow Y$.}
\end{figure}

Suppose that $L$ and $Q$ are hermitian holomorphic line bundles over compact hermitian manifolds $X$ and $Y$. The {\em homomorphism line bundle} $\hom (L,Q)$ is the line bundle over $Y \times \overline{X}$ whose fiber at $(y, x)$ is the collection of linear maps from $L_x$ to $Q_y$. (We must use the complex conjugate manifold $\overline{X}$ here because $\overline{L}$ is a holomorphic line bundle over $\overline{X}$, not $X$.) Its structure as a holomorphic line bundle is obtained from using the inner products in the fibers to identify it with the tensor product line bundle $Q \boxtimes \overline{L} \rightarrow Y \times \overline{X}$. A {\em holomorphic kernel} from $L$ to $Q$ is a holomorphic section
 \[
 \langle E \rangle \in \Gamma(\hom(L, Q))
\]
of the homomorphism line bundle. In other words a holomorphic kernel consists of a linear map
 \[
  \langle y | E | x \rangle \colon L_x \rightarrow Q_y
 \]
from every fiber of $L$ to every fiber of $Q$ which varies antiholomorphically with respect to $x$ and holomorphically with respect to $y$ (see Figure \ref{kernel1}).
Note that the collection of holomorphic kernels from $L$ to $Q$ is canonically isomorphic to the space of linear maps from $\Gamma(L)$ to $\Gamma(Q)$ via the following sequence of canonical isomorphisms:
 \begin{align} \label{seqisos}
 \Gamma(\hom(L, Q)) &\cong \Gamma(\overline{L} \boxtimes Q) \nonumber \\
 &\cong \Gamma(\overline{L}) \otimes \Gamma(Q) \quad \text{(standard fact about tensor product bundles \cite{ref:wells})} \nonumber \\
 &\cong \overline{\Gamma(L)} \otimes \Gamma(Q) \quad \text{(by definition)} \\
 &\cong \Gamma(L)^\vee \otimes \Gamma(Q)\quad \text{(Hilbert space identifies antilineary with its dual)} \nonumber\\
 &\cong \Hom(\Gamma(L), \Gamma(Q)) \nonumber.
 \end{align}
We remind the reader that in this thesis we write $V^\vee := \Hom(V, \mathbb{C})$ for the linear dual of a vector space $V$.

\section{The Bergman kernel\label{B2}}
Every hermitian line bundle $L$ over a compact complex hermitian manifold $X$ comes equipped with a canonical kernel to itself called the {\em Bergman} (or {\em reproducing}) {\em kernel}, which we write as
  \[
   \langle \;\; \rangle \in \Gamma (\hom(L, L)).
  \]
For references on the Bergman kernel, we refer the reader to \cite{ref:kirwin, ref:berndtsson, ref:berman_et_al, ref:berman}. The Bergman kernel is a fundamental geometric object and in this section we describe it in three different ways. We also describe how it relates to the {\em coherent states} formalism.

\subsection{The Bergman kernel via modes of propagation}
\begin{figure}[t]
\centering
\ig{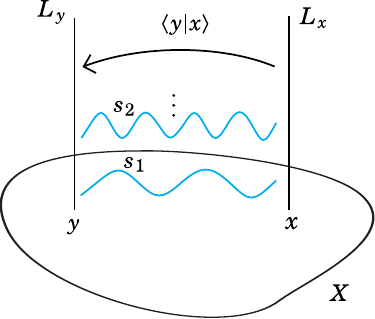}
\caption{\label{kernel2}  The Bergman kernel $\langle y | x \rangle$ sums over the fundamental modes of propagation from $L_x$ to $L_y$.}
\end{figure}

The most elementary description of the Bergman kernel is obtained by applying the sequence of isomorphisms \eqref{seqisos} in the reverse direction to the identity operator on $\Gamma(L)$. The result can be described as follows. Choose an orthonormal basis $\{s_i\}$ for $\Gamma(L)$. The Bergman kernel can then be described fibrewise as the map
 \begin{align*}
  \langle y | x \rangle \colon L_x &\rightarrow L_y \\
   v & \mapsto \sum_i ( s_i(x), v ) \, s_i (y).
 \end{align*}
In other words, {\em the Bergman kernel transports $v$ from the fiber at $x$ to the fiber at $y$ by summing over the various modes of propagation from $x$ to $y$} (see Figure \ref{kernel2}). We use the term `mode of propagation' to refer to these orthonormal sections $s_i$ because they are the zero eigenmodes of the Kodaira-Laplace operator. In this way we think of the Bergman kernel as providing the information of how the various lines in the line bundle $L$ are {\em correlated} with respect to each other. Note that this description is independent of the choice of orthonormal basis of sections $\{s_i\}$, since if $\{s'_i\}$ is another orthonormal basis, then $s'_i = \sum_j U_{ij} s_j$ for a unitary matrix $U$, and hence
 \[
  \sum_i (s'_i (x), v) \, s_i' (y) = \sum_{i,j,k} \overline{U}_{ij} U_{ik} (s_j(x), v) \, s_k (y) = \sum_k (s_k (x), v) \, s_k (y).
 \]

As an example, the Bergman kernel on the torus is constructed from the theta-functions, and on higher genus surfaces from automorphic forms (see \cite{ref:kirwin}). Also, the Bergman kernel $\langle L' | L \rangle $ on the hyperplane line bundle on the projective space of a Hilbert space is simply the orthogonal projection operator from the line $L$ onto the line $L'$, as we will see in Section \ref{B5}.

The importance of the Bergman kernel lies in its {\em reproducing property}.
\begin{lem}\label{repprop} For every holomorphic section $s \in \Gamma(L)$, we have
   \[
    s(y) = \int_X \langle y | x \rangle (s(x)) \, \vol_x.
   \]
\end{lem}
\begin{proof}
Choose an orthonormal basis $s_i$ for $\Gamma(L)$. We have:
  \begin{align*}
   \int_X \langle y | x \rangle (s(x)) \vol_x &= \int_X \sum_i (s_i(x), s(x)) s_i(y) \vol_x \\
    &= \sum_i \langle s_i, s \rangle s_i(y) \\
    &= s(y).
   \end{align*}
\end{proof}
Note that the Bergman kernel $\langle y | x \rangle $ is a {\em nonlocal object} in the sense that it is sensitive to the geometry of $X$ and $L$ not only near the points $x$ and $y$ but also to the geometry at other points $z$ very far away. This is because locally deforming the geometry around $z$ while leaving the geometry around $x$ and $y$ unchanged will alter the basis of sections $s_i$, and this information will thus `propagate' to the Bergman kernel $\langle y | x \rangle$. On the other hand, if the curvature of the line bundle $L$ is positive everywhere, then there do exist asymptotic expansions of the diagonal Bergman kernel $\langle x | x \rangle$ in terms of the local geometry at $x$ (see for example \cite{ref:berman, ref:berman_et_al, ref:ma_marinescu, ref:kirwin}) --- but these formulas only apply in the limit of high powers $L^p$ of $L$ (the so-called {\em semiclassical limit}). In this limit the Bergman kernel becomes more and more localized along the diagonal, so that there are no longer any correlations between distant points $x$ and $y$.

\subsection{The Bergman kernel as the large time limit of the heat kernel}
A more intrinsic description of the Bergman kernel is that it is the {\em large time limit of the heat kernel} on the space of sections of $L$. Let us review the notion of the heat kernel; see \cite{ref:ma_marinescu, ref:roe, ref:berline_getzler_vergne, ref:kirwin} for further details.

Since $X$ is a complex manifold, the complexification of the tangent bundle splits into holomorphic and antiholomorphic pieces
 \[
  TX \otimes_\mathbb{R} \mathbb{C} = T^{1,0}X \oplus T^{0,1}X.
 \]
We write the vector bundle of antiholomorphic differential forms with values in $L$ as
 \[
 S =  (T^{0,1}X)^\vee \otimes L.
 \]
We define the differential operator $D \colon C^\infty (S) \rightarrow C^\infty(S)$ as
 \[
 D = \sqrt{2} ( \overline{\partial} + \overline{\partial}^* ).
 \]
Here $\partial$ is the twisted Dolbeault operator --- the differential operator given locally by the holomorphic charts on the base manifold $X$ and the canonical connection on the hermitian line bundle $L$ ---  and $\overline{\partial}^*$ is its adjoint with respect to the inner products induced from the hermitian inner product on $X$ and the fibrewise inner products on $L$.

The {\em heat equation} for $D$ is the partial differential equation
 \[
  \frac{\partial s}{\partial t} = -D^2 s.
 \]
One can express the solutions to this equation in terms of the {\em heat kernel} --- a smooth time-dependent collection of linear maps
 \[
  \langle y | e^{-tD^2} | x \rangle \colon S_x \rightarrow S_y
 \]
having the property that
 \[
  s_t(y) = (e^{-tD^2}s_0) (y) = \int_X \langle y | e^{-tD^2} | x  \rangle s_0 (x)\, \vol_x.
 \]
Indeed, the heat kernel is the unique time-dependent section of $S \boxtimes S^\vee$ over $X \times X$ which is continuous in $t$, continuously differentiable in $x$ and $y$, satisfies the heat-equation in the variable $y$, and `tends to a $\delta$-function' as $t \rightarrow 0$.

Observe that as $t \rightarrow \infty$, the operator $e^{-tD^2}$ becomes the projector onto the space of {\em harmonic forms} $\Gamma(S)$ --- the forms satisfying $Ds = 0$:
 \[
  \lim_{t \rightarrow \infty} e^{-tD^2} = \pi \colon C^\infty (S) \rightarrow \Gamma (S).
 \]
If we restrict ourselves to $0$-forms (i.e., smooth sections of the line bundle $L$), then this means that the Bergman kernel is the large-time limit of the heat kernel, restricted to the space of holomorphic sections:
 \[
  \langle y | x \rangle = \lim_{t \rightarrow \infty} \langle y | e^{-tD^2} | x \rangle  |_{\Gamma(L)}.
 \]
We interpret this intuitively as the statement that `at large times, only the holomorphic modes of propagation from $x$ to $y$ survive'.

\subsection{The Bergman kernel as a path integral}
A good geometric way to think of the heat kernel
 \[
 \langle y | e^{-tD^2} | x \rangle \colon S_x \rightarrow S_y
 \]
is as a `weighted average of over all paths from $x$ to $y$ of the parallel transport from $S_x$ to $S_y$'. This is the {\em path integral} approach, which has long been employed by physicists but which has actually been made mathematically rigorous, notably in the recent work of B\"{a}r and Pf\"{a}ffle \cite{ref:bar_pfaffle}. In the limit $t \rightarrow \infty$, this provides a path integral interpretation for the Bergman kernel.

Let us review the result of B\"{a}r and Pf\"{a}ffle, specialized to our setting of a hermitian line bundle $L$ over a compact hermitian manifold $X$. We need to make some definitions. A {\em partition of length $t$} is a sequence of positive numbers $T = (t_1, \ldots, t_r)$ which add up to $t$. We write $|T|$ for the largest difference $t_{i+1} - t_i$ in $T$ and $\sigma_j(T) := t_1 + \cdots + t_j$. Let $x$ and $y$ be points in $X$, and let $\gamma \colon [0,t] \rightarrow X$ be a continuous curve such that $\gamma(0) = x$ and $\gamma(t) = y$, and put $x_j := \gamma(\sigma_j(T)$. We say that $\gamma$ is a {\em geodesic polygon} from $x$ to $y$ with respect to $T$ if for any two subsequent points $x_j$ and $x_{j+1}$ the curve $\gamma|_{[\sigma_j(T), \sigma_{j+1}(T)]}$ is the unique shortest geodesic joining them (this implies that they are not cut points of each other). Note in particular that when restricted to one of the subintervals $[\sigma_j(T), \sigma_{j+1}(T)]$ the curve is smooth and parameterized proportionally to arc-length. We write
 \[
  S(\gamma, t) :=  \int_0^t \left[\frac{1}{4} |\dot{\gamma}(s)|^2 - \frac{1}{3} \text{scal}_{\gamma(s)}\right] ds
  \]
for the {\em action} of a geodesic polygon of length $t$, where $\text{scal}_x$ refers to the scalar curvature of $X$ at a point $x$.

We denote the set of all geodesic polygons from $x$ to $y$ with respect to $T$ as $\Paths_T(x,y)$. Since a geodesic is uniquely determined by its end points (unless they are cut points of each other), we can identify $\Paths_T(x,y)$ with an open set of $X^{r+1}$ via the correspondence
 \[
  \gamma \leftrightarrow (x_0, \ldots, x_r).
 \]
The Riemannian metric on $X^{r+1}$ induces a measure on $\Paths_x^y (T) $ which we denote by $\mathcal{D}\gamma$.

Since $L$ is a hermitian line bundle, there exists a unique holomorphic covariant derivative $\nabla$ acting on sections of $L$ which preserves the metric. We write parallel transport along a curve $\gamma$ with respect to $\nabla$ as $\tau(\gamma) \colon L_\gamma(0) \rightarrow L_{\gamma(t)}$. Finally, we set
 \[
  Z(T) := \prod_{j=1}^r (4 \pi t_j)^{\frac{\dim X}{2}}.
 \]
We are now ready to state their result.

\begin{thm}[{refer B\"{a}r and Pf\"{a}ffle \cite[Thm 6.1]{ref:bar_pfaffle}}] Suppose $L$ is a hermitian line bundle over a compact hermitian manifold $X$. Then one can define a sequence of partitions $T_n$ of length $t$ with $|T_n| \rightarrow 0$, such that
 \[
  \langle y |  e^{-tD^2} | x \rangle = \lim_{n \rightarrow \infty} \frac{1}{Z(T_n)} \int_{\Paths_{T_n}(x,y)} e^{-S(\gamma, t)} \tau(\gamma)_0^t \mathcal{D}\gamma
 \]
where the $n \rightarrow \infty$ limit refers to convergence in $C^0(X \times \overline{X}, L \boxtimes \overline{L})$.
\end{thm}
\noindent We will write
 \[
  \frac{1}{Z_t} \int_{\Paths_t(x,y)} e^{-S(\gamma, t)} \tau(\gamma)_0^t \mathcal{D}\gamma
 \]
as shorthand for the limit above.
\begin{cor} The Bergman kernel can be expressed as the large time limit of a path integral:
 \[
   \langle y | x \rangle = \lim_{t \rightarrow \infty} \frac{1}{Z_t} \int_{\Paths_t(x,y)} e^{-S(\gamma, t)} \tau(\gamma)_0^t \mathcal{D}\gamma.
 \]
\end{cor}
This description of the Bergman kernel $\langle y | x \rangle$ makes it very clear why it is a nonlocal object, since it sums over all paths from $x$ to $y$ and hence is sensitive to the geometry of the entire manifold.

\subsection{The Bergman kernel and coherent states}
The Bergman kernel often manifests itself via the {\em coherent states} formalism; we recommend the paper of Kirwin \cite{ref:kirwin} as an introduction to this approach. For each $x \in X$ and $v \in L_x$, define the {\em coherent state at $x$ with basevector $v$} as the holomorphic section $|x\rangle_v$ of $L$ obtained by applying the Bergman kernel to $v$:
 \[
  |x\rangle_v (y) = \langle y | x \rangle (v).
 \]
For instance, when $L$ is the hyperplane line bundle over the projective space of a finite-dimensional Hilbert space $V$, the coherent state arising from a vector $v \in V$ is simply the section of the hyperplane line bundle which orthogonally projects $v$ onto every line $L$, as we will see in Section \ref{B5}. In general, coherent states have many marvelous properties from the viewpoint of geometric quantization: they form an overcomplete basis for $\Gamma(L)$, they are `maximally peaked' around $x$, and they `evolve classically'. We refer the reader to \cite{ref:kirwin} for a definition of these terms; we will not delve into these aspects of the theory.

Suppose that $A \colon \Gamma(L) \rightarrow \Gamma(L)$ is a linear map. Using the coherent states, we can think of $A$ as a holomorphic kernel via the prescription
 \begin{align*}
  \langle y | A | x \rangle \colon L_x & \rightarrow L_y \\
   v & \mapsto (A |x\rangle_v)(y).
  \end{align*}
The trace of $A$ can be computed as an integral over $X$ as follows. We use the following convention: on the diagonal the map $\langle x | A | x \rangle \colon L_x \rightarrow L_x$ is a map from a one-dimensional vector space to itself, and we identify it with its {\em trace}, that is with the complex number it represents.
 \begin{lem} The trace of a linear operator $A \colon \Gamma(L) \rightarrow \Gamma(L)$ can be computed as
  \[
   \Tr A = \int_X \langle x | A | x \rangle \, \vol_x.
  \]
 \end{lem}
 \begin{proof} Choose an orthonormal basis of sections $s_i$ for $\Gamma(L)$, and expand the coherent state $|x\rangle_v$ as
  \[
   |x \rangle_v = \sum_i \langle s_i, |x \rangle_v \rangle s_i.
  \]
We calculate these expansion coefficients as:
 \begin{align*}
  \langle s_i, |x \rangle_v \rangle &= \int_X (s_i(y), \langle y | x \rangle v ) \, \vol_y \\
   &= \sum_k \int_X (s_i(y), (s_k(x), v)s_k(y)) \, \vol_y \\
   &= (s_i(x), v).
  \end{align*}
Therefore, choosing each $v \in L_x$ to have unit norm, we have:
 \begin{align*}
  \int_X \langle x | A | x \rangle \, \vol_x &= \sum_{i,j} \langle s_i, |x\rangle_v \rangle \langle s_j, As_i \rangle (v, s_j (x)) \\
   &= \sum_{i,j} \langle s_j, A s_i \rangle \int_X (s_i(x), s_j(x)) \\
   &= \sum_i \langle s_i, A s_i \rangle \\
   &= \Tr A.
  \end{align*}
\end{proof}

\section{The category of line bundles and kernels\label{B4}}
In this section we define the category $\LBun$ whose objects are holomorphic line bundles and whose morphisms are kernels between them. First we give the definition, then we expand out what it means, and finally we prove that the resulting structure indeed forms a category.

 \begin{defn} The category $\LBun$ is defined as follows. An object is a holomorphic hermitian line bundle $L \rightarrow X$ over a compact hermitian manifold $X$. A morphism is a holomorphic kernel. The identity morphisms are the Bergman kernels, and composition is defined by integration with respect to the volume form defined by the metric.
 \end{defn}
Let us write out explicitly the definition of composition. Suppose that $L \rightarrow X$, $Q \rightarrow Y$ and $R \rightarrow Z$ are hermitian line bundles over $X$, $Y$ and $Z$ respectively, and that $\langle E \rangle \colon L \rightarrow Q$ and $\langle F \rangle \colon Q \rightarrow Z$ are holomorphic kernels, that is
 \[
  \langle E \rangle \in \Gamma(\hom(L, Q)) \quad \text{and} \quad \langle F \rangle \in \Gamma(\hom(Q, R)).
 \]
Then the composite kernel
 \[
  \langle F \rangle \circ \langle E \rangle \in \Gamma(\Hom(L, R))
 \]
is defined by integrating over $y \in Y$:
 \be \label{formkernel}
  \langle z | \langle F\rangle \circ \langle E\rangle | x \rangle v = \int_Y \langle z | F | y \rangle \langle y | E | x \rangle v \, \vol_y.
 \ee
This new kernel is indeed holomorphic in $x$ and antiholomorphic in $z$ (think of the approximations to the integral by finite sums). We now confirm that the resulting structure indeed forms a category.

\begin{lem} The structure $\LBun$ indeed forms a category.
\end{lem}
\begin{proof} Firstly, composition is clearly associative, by Fubini's theorem:
 \begin{align*}
 \langle t | (\langle G \rangle \circ \langle F \rangle) \circ \langle E \rangle | x \rangle v &= \int_Z \langle t | G | z \rangle \circ \langle z | \langle F \rangle \circ \langle E \rangle | x \rangle v \, \vol_z \\
 &= \int_Z \int_Y \langle t | G | z \rangle \circ \langle z | F | y \rangle \circ \langle y | F | x \rangle v \, \vol_y \vol_z \\
 &= \int_Y \int_Z \langle t | G | z \rangle \circ \langle z | F | y \rangle \circ \langle y | F | x \rangle v \, \vol_z \vol_y \\
 &=  \langle t | \langle G \rangle \circ (\langle F \rangle \circ \langle E \rangle ) | x \rangle v.
 \end{align*}
We must check that the Bergman kernel $\langle \; \; \rangle$ is the left and right unit for composition. Suppose that $\langle E \rangle$ is a holomorphic kernel from a line bundle $L \rightarrow X$ to another line bundle $Q \rightarrow Y$. Then
 \begin{align}
  \langle y | \langle \; \; \rangle \circ \langle E \rangle | x \rangle v &= \int_X \langle y | x' \rangle \circ \langle x' | E | x \rangle v \, \vol_{x'} \nonumber \\
  &= \int_X \langle y | x' \rangle s (x') \label{bobuncle}
 \end{align}
where $s$ is the holomorphic section of $L$ defined by
 \[
  s(x') = \langle x' | E | x \rangle v.
 \]
Therefore the expression \eqref{bobuncle} must equal $\langle y | E | x \rangle v$ by the reproducing property (Lemma \ref{repprop}) of the Bergman kernel. This establishes that the Bergman kernel is the left unit for composition.

To see that it is the right unit for composition, firstly observe that if $w \in Q_y$ and $x' \in X$ then the section $t$ of $L$ defined by
 \[
  t(x) = \langle y | E | x \rangle^* w
 \]
is holomorphic.\XXX{I'm not 100 sure why it is holomorphic}
 Therefore we have:
  \begin{align*}
   \bigl( w, \langle y | \langle E \rangle \circ \langle \rangle | x \rangle v \bigr) &= \int_X \bigl( w, \langle y | E | x' \rangle \circ \langle x' | x \rangle v \bigr) \, \vol_{x'} \\
   &= \int_X ( t(x'), \langle x' | x \rangle v) \, \vol_{x'} \\
   &= \sum_i \int_X (t(x'), s_i(x')) \, \vol_{x'} \, (s_i(x), v) \\
   &= \sum_i \bigl( \langle s_i, t \rangle s_i (x), v \bigr) \\
   &= (t(x), v) \\
   &= (\langle y | E | x \rangle^* w, v) \\
   &= (w, \langle y | E | x \rangle v).
  \end{align*}
Thus $\langle y | \langle E \rangle \circ \langle \rangle | x \rangle v =  \langle y | E | x \rangle v$.
\end{proof}

\section{The geometric line bundle over projective space\label{B5}}
In this section we introduce the {\em geometric line bundle} $\tau$ over the projective space of a finite dimensional Hilbert space, and we show how the Hilbert space can be recovered antilinearly as the space of holomorphic sections of $\tau$.

Suppose that $V$ is an $(n+1)$-dimensional vector space. Recall that the {\em projective space} $\mathbb{P}(V)$ of $V$ is the compact complex manifold consisting of all one-dimensional subspaces $l \subseteq V$. The {\em tautological line bundle} $\mathcal{O}(-1)$ over $\mathbb{P}(V)$ is the holomorphic line bundle with total space
 \[
 \mathcal{O}(-1) = \big\{ (l, v) \colon l \text{ is a line in $V$, $v \in l$}\big\},
 \]
so that the fiber at a line $l \in \mathbb{P}(V)$ is the line $l$ itself:
 \[
  \mathcal{O}(-1)_l = l.
 \]
The {\em hyperplane line bundle} $\mathcal{O}(1)$ is defined as the dual of the tautological line bundle, so that its fibers are given by
 \[
  \mathcal{O}(1)_l = l^\vee
 \]
where $l^\vee := \Hom(l, \mathbb{C})$ is the linear dual of $l$. If $V$ is equipped with an inner product then one can antilinearly identify each line $l$ with its linear dual $l^\vee$. We define the {\em geometric line bundle} $\tau$ over $\mathbb{P}(V)$ as the line bundle whose fibers are the complex conjugate of the tautological line bundle fibers,
 \[
  \tau_l = \overline{l},
 \]
but whose whose holomorphic structure is derived from the hyperplane line bundle $\mathcal{O}(1)$ by making the fibrewise identifications $\overline{l} \cong l^\vee$. In other words, the geometric line bundle $\tau$ is just another notation for the hyperplane line bundle $\mathcal{O}(1)$, available in the presence of an inner product, and we use it because oftentimes (but not always) it makes for shorter and more intuitive formulas.

To be completely explicit here, we are using the notation that $\overline{l}$ refers to the vector space whose underlying abelian group is $l$, the only difference being that scalar multiplication acts antilinearly. Convention dictates that when a vector $u \in l$ is thought of as a vector in $\overline{l}$ one should really write it as $\overline{u}$. However, this convention tends to clutter up the formulas, and so we will rarely use it unless the distinction between $\overline{u} \in \overline{l}$ and $u \in l$ is important in a specific calculation. In other words, if a formula requires as input a certain vector $u \in \overline{l}$, but as written only receives a vector $u \in l$, then it is understood that one should first write a bar above the $u$ in order to actually apply the formula, and the same convention goes for the output side of the formula.

Now, the fibers of the geometric line bundle carry a natural inner product inherited from $V$,
 \[
  (\overline{v}, \overline{w})_{\overline{l}} = (w, v)_V.
 \]
(This is one place where we need to be careful about writing $\overline{v}$ instead of $v$.) This means that the geometric line bundle is a {\em hermitian line bundle}.

Moreover, the projective space $\mathbb{P}(V)$ carries a natural hermitian metric called the {\em Fubini-study metric}, which is characterized up to a scalar multiple as the unique metric on $\mathbb{P}(V)$ which is invariant under the unitary group $U(V)$. We will normalize this metric by demanding that if $v \in V$ is a unit vector, then
 \[
  \int_{\mathbb{P}(V)} |(u, v)|^2 \vol_l = 1 \quad \text{(where $u \in l, (u,u)=1$)}
 \]
where $\vol_l$ is the corresponding volume form of the metric. More informally, we are normalizing the Fubini-Study metric in such a way that `the probabilities for finding a vector $v$ in the various states $u$ add up to unity'. Note that this normalization does not depend on the choice of the unit vector $v \in V$ since the Fubini-Study metric is invariant under unitary transformations.

This gives the vector space $\Gamma(\tau)$ of holomorphic sections of the geometric line bundle an inner product via
 \[
  \langle s, s' \rangle = \int_{\mathbb{P}(V)} \left(s(l), s'(l) \right) \, \vol_l
 \]
Note that the space of sections $\Gamma(\tau)$ is finite-dimensional since $\mathbb{P}(V)$ is compact. We now show that $\Gamma(\tau)$ identifies with $\overline{V}$.

\begin{figure}
\centering
\ig{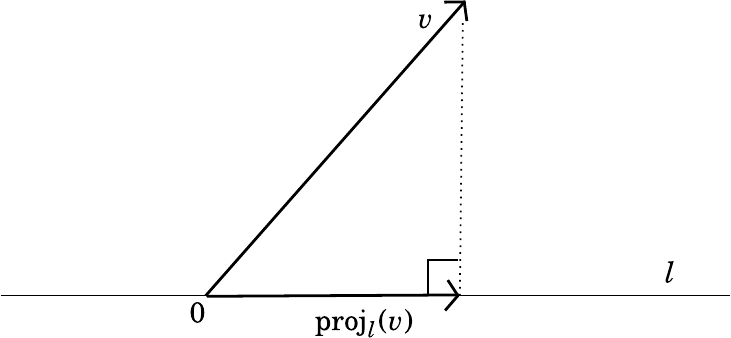}
\caption{\label{projfig}The section $s_v$ of the geometric line bundle over projective space associated to a vector $v \in V$.}
\end{figure}

Each vector $v \in V$ gives rise to a holomorphic section $s_v \in \Gamma(\tau)$, defined by orthogonally projecting $v$ onto every line $l$ as in  Figure \ref{projfig}, and then thinking of the result as an element of $\overline{l}$:
 \[
  s_v (l) = (l, \proj_l (v) \in \overline{l}).
 \]
This section $s_v$ is indeed holomorphic, as we show in the following lemma.

\begin{lem} For each $v \in V$, the section $s_v$ of the geometric line bundle $\tau$ described above is holomorphic.
\end{lem}
\begin{proof} Recall that the geometric line bundle $\tau$ is just another notation for the hyperplane line bundle $\mathcal{O}(1)$, in that the charts for the former are defined in terms of the charts for the latter by making the fiberwise identifications $\overline{l} \cong l^\vee$. Thus it suffices to check that $s_v$ is holomorphic when thought of as a section of the hyperplane line bundle.

In terms of the hyperplane line bundle, the section $s_v$ translates into the section
 \begin{align*}
  s'_v \colon \mathbb{P}(V) & \rightarrow  \mathcal{O}(1) \\
   l & \mapsto \left(l, (v, \cdot)|_l \right)
  \end{align*}
since for all $u \in l$ we have
 \[
  \left( \proj_l(v), u\right) = (v,u).
 \]
by the definition of what it means to orthogonally project a vector.

To see that $s'_v$ is holomorphic, let us recall the standard holomorphic charts for projective space and the hyperplane line bundle. (The author thanks Neil Strickland for this section, which gives a co-ordinate free presentation).

For any
nonzero linear map $\alpha \colon V \rightarrow \mathbb{C}$, put $U_\alpha = \alpha^{-1}\{1\}$ and
$P_\alpha=\{l\in \mathbb{P}(V) : \alpha(l) = \mathbb{C}\}$.
Define the chart maps
$\phi_\alpha \colon U_\alpha \rightarrow P_\alpha$ by $\phi_\alpha(v)=\mathbb{C} v$.  (Of course, $\phi_\alpha$ is only really a `chart' once we fix a holomorphic identification between the affine space $U_\alpha$ and $\mathbb{C}^{\dim(V) -1}$, but this easily can be done by fixing an inner product on $V$ in such a way that $\alpha = (a, \cdot)$ for some $a \in V$ with $||a||=1$, thereby obtaining $U_i = a + a^\perp$, and then choosing a linear identification $a^\perp \cong \mathbb{C}^{\dim(V)-1}$. ) The transition map
$\phi_\beta^{-1}\phi_\alpha$ is given by $v\mapsto v/\beta(v)$ and so is
holomorphic where it is defined.

The maps $\phi_\alpha$ are the
charts for the standard holomorphic atlas on $\mathbb{P}(V)$, after one makes the identifications $U_\alpha \cong \mathbb{C}^{\dim(V)-1}$ in the manner described in the parenthesis above.

The standard charts for the hyperplane line bundle $\mathcal{O}(1)$ are given by
 \begin{align*}
  \xi_\alpha \colon U_\alpha \times \mathbb{C} & \rightarrow \mathcal{O}(1)|_{P_\alpha} \\
   (v, z) & \mapsto (\mathbb{C}v, z \alpha|_{\mathbb{C}v}) \\
   \left(\text{unique } v \in l \text{ s.t. } \alpha(v)=1, \lambda(v)\right)  & \mapsfrom (l, \lambda)
 \end{align*}
Now, note that for any linear map $\gamma \colon V \rightarrow \mathbb{C}$ we have a section $t_\gamma \colon \mathbb{P}(V) \rightarrow \mathcal{O}(1)$ given by $t_\gamma(l)=(l,\gamma|_l)$.  This satisfies $\xi_\alpha^{-1}(t_\gamma(\phi_\alpha(v)))=\gamma(v)$, which proves that $t_\gamma$ is holomorphic. This finally shows that $s'_v$ is holomorphic, since $s'_v = t_{(v, \cdot)}$.

For a co-ordinate description, fix an orthonormal basis $\{e_i\}$ of $V$ and write $v = \sum_{j=0}^n v^i e_i$ and $l = [z^0 : \cdots : z^n]$. Then the charts for projective space described above translate into
 \begin{align*}
  \phi_i \colon  \mathbb{C}^n & \rightarrow \{ [z^0 : \cdots : z^n] : z^i \neq 0\} \\
 (w^1, \ldots, w^n) & \mapsto [w^1 : \cdots : w^i : 1 : w^{i+1} : \cdots : w^n] \\
( \frac{z^0}{z^i}, \ldots, \frac{z^{i-1}}{z^i}, \frac{z^{i+1}}{z^i}, \ldots, \frac{z^n}{z^i}) & \mapsfrom [z^0 : \cdots : z^n].
 \end{align*}
Thus the section $s_v$ (equivalently, $s'_v$) computes locally as
 \[
  (w^1, \ldots, w^n) \mapsto \overline{v^0}w^1 + \cdots + \overline{v^{i-1}}w^i + \overline{v^i} + \overline{v^{i+1}}w^{i+1} + \cdots + \overline{v^n}w^n
 \]
which is clearly holomorphic in $w$.
\end{proof}
In this way one can recover $V$ antilinearly as the space of sections of the geometric line bundle. This is a well-known fact in algebraic geometry in the guise of the hyperplane line bundle, but here we give a proof which emphasizes the inner products involved. A category theorist should think of this as the `decategorified Yoneda lemma' since it says that knowing a vector $v \in V$ is the same thing as knowing the inner products of $v$ with all the other vectors.
 \begin{lem}\label{myoneda} The map
 \begin{align*}
   \overline{V} & \rightarrow \Gamma(\tau) \\
   v & \mapsto s_v
 \end{align*}
is a unitary isomorphism of finite dimensional Hilbert spaces.
 \end{lem}
 \begin{proof}
To see that the map $v \mapsto s_v$ is antiunitary, one needs to remember that complex scalars are defined to act on the fibers of the geometric line bundle antilinearly. Hence for $v, v' \in V$ we have
 \begin{align*}
  \langle s_v, s_{v'} \rangle  &= \int_{\mathbb{P(V)}} \left( (u,v) \cdot u, \, (u,v') \cdot u ) \right) \; \vol_l \quad ( u \in l, ||u||=1) \\
  &= \int_{\mathbb{P}(V)} (v',u)(u,v) \; \vol_l.
 \end{align*}
Note that this formula defines an alternative inner product on $V$ via
 \[
  (v, v')_2 := \langle s_{v'}, s_v \rangle = \int_{\mathbb{P}(V)} (v,u)(u,v') \; \vol_l.
 \]
By linear algebra, we must therefore have $(v, v')_2 = (v, Av')$ for some positive self-adjoint operator $A \colon V \rightarrow V$. Now $A$ commutes with $U(V)$ because $(v, v')_2$ is $U(V)$ invariant. Hence $A = r \id$ for some positive constant $r > 0$, so that $(v, v')_2 = r (v,v')$ for all $v, v' \in V$ (the author thanks Neil Strickland for this argument). But $r$ must be equal to unity, because for $v \in V$ with $(v,v) = 1$, we have
 \begin{align*}
  (v, v)_2 &= \int_{\mathbb{P}(V)} |(u, v)|^2 \vol_l \quad ( u \in l, ||u||=1) \\
           &= 1
 \end{align*}
by our normalization convention for the Fubini-Study metric. Thus we have $\langle s_v, s_{v'} \rangle = (v', v)$ so that the map $v \mapsto s_v$ is antiunitary.

In addition, {\em all} holomorphic sections of $\tau$ are of the form $s_v$ for some $v \in V$, by the following argument. Suppose that $s \in \Gamma(\tau)$. Observe that the maps
 \begin{align*}
  f_i \colon \mathbb{P}(V) & \rightarrow \mathbb{C} \\
  l & \mapsto (e_i, s(l))
 \end{align*}
are holomorphic functions \XXX{check this} on the compact manifold $\mathbb{P}(V)$, and hence constant. We claim that
 \[
  s = s_{\sum f_i e_i}.
 \]
This is one place where we need to be careful about complex conjugation signs. We should really be writing
 \[
  f_i(l) = (e_i, \overline{s(l)})
 \]
since $s(l) \in \overline{l}$. Now, we calculate:
\begin{align*}
  s_{\sum f_i e_i} (l) &= \sum_i \overline{f_i} s_{e_i} (l) \\
  &= \sum_i (\overline{s(l)}, e_i) \overline{(u, e_i)u} \quad \text{(where $u \in l, (u,u)=1$)}\\
   &= \overline{\sum_i (e_i, \overline{s(l)})(u, e_i) u} \\
   &= \overline{(u, \overline{s(l)}) u} \\
   &= \overline{\overline{s(l)}} \\
   &= s(l).
  \end{align*}
Hence $s = s_{\sum f_i e_i}$ as we claimed, and the map $\overline{V} \rightarrow \Gamma(\tau)$ is surjective.
\end{proof}
This result is the basis of this entire chapter: if we can think of every vector space as the space of sections of a line bundle, then it is not hard to imagine that the category of representations of a group is equivalent to a category of equivariant line bundles, which is the result we are working towards. For now, let us record a result we promised in the earlier sections.
\begin{cor} The Bergman kernel of the geometric line bundle over the projective space of a finite-dimensional Hilbert space $V$ is given fibrewise by orthogonally projecting one line onto another, ie.
 \[
  \langle l' | l \rangle v = \proj_{l'} (v)
 \]
where $l$ and $l'$ are lines in $V$.
\end{cor}
\begin{proof} Choose an orthonormal basis $\{e_i\}$ for $V$. By the previous lemma, the corresponding sections $s_{e_i}$ of the geometric line bundle over $\mathbb{P}(V)$ form an orthonormal basis for $\Gamma(L)$. Suppose $l$ and $l'$ are lines in $\mathbb{P}(V)$, and $v \in l$. By definition the Bergman kernel computes as
 \begin{align*}
  \langle l' | l \rangle v &= \sum_i \left( s_{e_i}(l), \, v\right) s_{e_i}(l') \\
   &= \sum_i \frac{1}{(v,v)} \left( (v, e_i) v, v \right) \left(w, e_i\right) w \quad \text{(where $w \in l', (w,w)=1$)} \\
   &= \sum_i (w,e_i)(e_i, v) w \\
   &= (w,v)w \\
   &= \proj_{l'} (v).
  \end{align*}
 \end{proof}

\section{An equivalence of categories\label{B6}}
In this section we show that the category $\Hilb$ of finite-dimensional Hilbert spaces and linear maps is equivalent to the category $\LBun$ of holomorphic line bundles over compact hermitian manifolds and kernels between them. We do this by defining functors $\Hilb \rightarrow \LBun$ and $\LBun \rightarrow \Hilb$ and showing that they establish an equivalence of categories.

\begin{figure}[t]
\centering
\ig{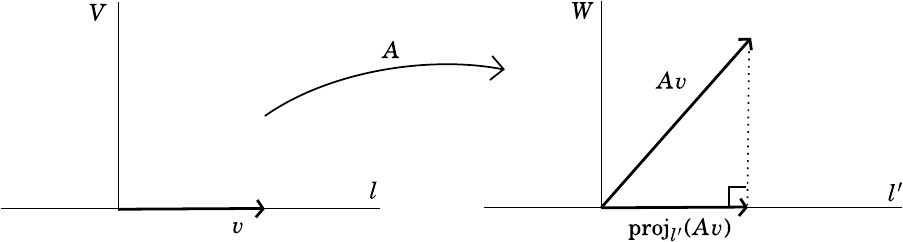}
\caption{\label{projfig2} The kernel $\langle l' | A | l \rangle$ is defined by applying $A$ and then projecting onto the line $l'$.}
\end{figure}

\begin{prop} The map $\tau \colon \Hilb \rightarrow \LBun$ defined as follows is a functor.
 \begin{itemize}
  \item A finite-dimensional Hilbert space $V$ gets sent to the geometric line bundle $\tau_V \rightarrow \mathbb{P}(V)$.
  \item A linear map $A \colon V \rightarrow W$ gets sent to the holomorphic kernel $\langle A \rangle \colon \tau_V \rightarrow \tau_W$ defined fibrewise by applying $A$ and then projecting (see Figure \ref{projfig2}),
 \begin{align*}
  \langle l' | A | l \rangle \colon \overline{l} & \rightarrow \overline{l'} \\
   v & \mapsto \proj_{l'}(Av) \in \overline{l'}.
  \end{align*}
 \end{itemize}
\end{prop}

\begin{proof} If $A \colon U \rightarrow V$ and $B \colon V \rightarrow W$ are linear maps, and $u \in U$, then:
 \begin{align*}
  \langle l'' | \langle B \rangle \circ \langle A \rangle | l \rangle u & = \int_{\mathbb{P}(V)} \langle l'' | B | l' \rangle \circ \langle l' | A | l \rangle u \, \vol_{l'} \\
   &=\int_{\mathbb{P}(V)}  \proj_{l''}(B \proj_{l'} (Au)) \, \vol_{l'} \\
   &= \int_{\mathbb{P}(V)} (w, Bv)(v, Au) w \, \vol_{l'} \quad \left( \ba \begin{array}{c}\text{where $v \in l', (v,v)=1$}\\ \text{and $w \in l'', (w,w)=1$} \end{array} \ea \right) \\
   &= \langle s_{B^* w}, s_{Au} \rangle w   =  (B^*w, Au) w \\
   &= (w, B Au) w  =   \proj_{l''} (BAu).
 \end{align*}
Moreover the identity linear map $\id \colon V \rightarrow V$ gets sent to the Bergman kernel. To see this, suppose $l$ and $l'$ are lines in $V$ and $v \in l$. Choose an orthonormal basis $\{e_i\}$ for $V$, and let $\{s_{e_i} = \proj_{\cdot} (e_i)\}$ be the corresponding orthonormal basis of sections of the geometric line bundle. Then:
 \begin{align*}
  \langle l' | \id | l \rangle v &= (w,v) w \quad \text{(where $w \in l', (w,w)=1$)} \\
                                 &= \sum_i (e_i, v)(w, e_i) w = \sum_i (s_{e_i}(l), v) s_{e_i} (l') \\
                                 &= \langle l' | l \rangle v.
  \end{align*}
\end{proof}
We now write down the functor which goes in the other direction.

\begin{prop} The map $\overline{\Gamma} \colon \LBun \rightarrow \Hilb$ defined as follows is a functor.
 \begin{itemize}
  \item A line bundle $L \rightarrow X$ gets sent to the Hilbert space $\overline{\Gamma(L)}$, with the inner product defined by integrating the fibrewise inner products over $X$.
  \item A holomorphic kernel $\langle E \rangle$ from $L \rightarrow X$ to $Q \rightarrow Y$ gets sent to the linear map
 \[
  E \colon \overline{\Gamma(L)} \rightarrow \overline{\Gamma(Q)}
 \]
defined by integrating $\langle E \rangle$ over $X$,
 \[
  (Es)(y) = \int_X \langle y | E | x \rangle s(x) \vol_x.
 \]
 \end{itemize}
\end{prop}
\begin{proof} This construction is functorial, by Fubini's theorem:
 \begin{align*}
  (F (Es))(z) &= \int_Y \langle z | F | z \rangle (Es)(y) \, \vol_y \\
                &= \int_Y \int_X \langle z | F | y \rangle \circ \langle y | E | x \rangle s(x) \, \vol_x \vol_y \\
                &= \int_X \int_Y \langle z | F | y \rangle \circ \langle y | E | x \rangle s(x) \, \vol_y \vol_x \\
                &= \int_X \langle z | \langle F \rangle \circ \langle E \rangle | x \rangle s(x) \, \vol_x.
  \end{align*}
Moreover, since by Lemma \ref{repprop}
 \[
 s(y) = \int_X \langle y | x \rangle s(x) \, \vol_x
 \]
we see that the functor sends Bergman kernels (the identity morphisms in $\LBun$) to identity linear maps, as it should.
\end{proof}
Finally we show that these constructions establish an equivalence of categories.

\begin{prop} The functors $\tau\colon \Hilb \rightarrow \LBun$ and $\overline{\Gamma} \colon \Hilb \rightarrow \LBun$ defined above establish an equivalence of categories between the category of finite-dimensional Hilbert spaces and the category of hermitian line bundles.
\end{prop}
\begin{proof} To see that these categories simply {\em are} indeed equivalent, observe that the functor $\overline{\Gamma} \colon \LBun \rightarrow \Hilb$ is essentially surjective by Lemma \ref{myoneda}, and fully faithful due to the sequence of isomorphisms given in \eqref{seqisos}. However, we need to show that the functors defined above are mutual inverses of each other up to isomorphism. To do this, we define unit and counit natural transformations $\eta \colon \id \Rightarrow \overline{\Gamma} \circ \tau$ and $\text{ev} \colon \overline{\Gamma} \circ \tau \Rightarrow \id$ and prove that they are isomorphisms.

Indeed, for every line bundle $L \rightarrow X$ there is a canonical kernel
 \[
  \langle \text{ev} \rangle_L \colon \tau_{\overline{\Gamma(L)}} \rightarrow L
 \]
from the geometric line bundle over the space of sections of $L$ to $L$, defined as follows. Let $[s]$ be the one-dimensional line in $\overline{\Gamma(L)}$ through the section $s$, and $s'$ is a member of this line. Then the kernel $\langle \text{ev} \rangle_L$ is defined by simply evaluating a section at a point $x \in X$:
 \begin{align*}
  \langle x | \text{ev} | [s] \rangle_L \colon [s] & \rightarrow L_x \\
   s' & \mapsto s'(x).
  \end{align*}
To show that this collection of kernels defines a natural transformation, we must show that for every kernel $\langle E \rangle \colon L \rightarrow Q$ the following diagram commutes:
 \[
  \xymatrix{ \tau_{\overline{\Gamma(L)}} \ar[r]^-{\langle \tau(E) \rangle} \ar[d]_{\langle \text{ev} \rangle_L}  & \tau_{\overline{\Gamma(Q)}} \ar[d]^{\langle \text{ev}\rangle_Q} \\
  L \ar[r]_{\langle E \rangle} & Q}
 \]
Indeed, by expanding out the definitions we compute:
 \begin{align*}
  \langle y | \langle \text{ev} \rangle_Q \circ \langle \tau(E) \rangle | [s] \rangle (s') &= \int_{\overline{\Gamma(Q)}} \langle y | \text{ev} | [t]] \rangle \langle [t] | \tau(E) | [s] (s') \; \vol_{[t]} \\
  &= \int_{\overline{\Gamma(Q)}} \int_Y \int_X \left(t'(y'), \langle y' | E | x \rangle s'(x)\right) t'(y) \; \vol_x \vol_y \vol_{[t]} \\
 &= \int_X \langle y | E | x \rangle s'(x) \; \vol_x \\
 &= \langle y | \langle E \rangle \circ \langle \text{ev} \rangle_L | [s] \rangle (s').
\end{align*} \XXX{I cheated slightly above. Work this out before the vivo!}
We claim that this kernel has an inverse
 \[
  \langle \text{ev}^\mi \rangle \colon \tau_{\overline{\Gamma(L)}} \rightarrow L
 \]
given by orthogonally projecting the coherent state supported at $x$ onto the line $[s]$:
 \begin{align*}
  \langle  [s] | \text{ev}^\mi | x \rangle \colon L_x & \rightarrow [s]  \\
   v & \mapsto \proj_{[s]} |x\rangle_v.
 \end{align*}
Indeed, this is the `kernel' version of the well-known {\em Kodaira embedding map} (see \cite[pg 8]{ref:kirwin}). We have
 \begin{align*}
  \langle y | \langle \text{ev} \rangle \circ \langle \text{ev}^\mi \rangle | x \rangle v &= \int_{\mathbb{P}(\overline{\Gamma(L)})}  \langle y | \text{ev}| [s] \rangle \circ \langle [s] | \text{ev}^\mi | x \rangle v \, \vol_{[s]} \\
  &= \int_{\mathbb{P}(\overline{\Gamma(L)})} (\sum_i \langle s_i, s \rangle s_i (x), v) s(y) \, \vol_{[s]} \\
  &= \sum_{i,j} (s_i(x), v) s_j(y) \, \int_{\mathbb{P}(\overline{\Gamma(L)})} \langle s_j, s \rangle \langle s, s_i \rangle \, \vol_{[s]} \\
  &= \langle y | x \rangle v,
  \end{align*}
because if we write $V = \overline{\Gamma(L)}$ and $e_i \in V$ for $s_i \in \overline{\Gamma(L)}$ and also $v$ for $s$, then the integral term on the right of the second-last equation above computes as
 \begin{align*}
  \int_{\mathbb{P}(\overline{\Gamma(L)})} \langle s_j, s \rangle \langle s, s_i \rangle \, \vol_{[s]} &= \int_{\mathbb{P}(V)} ( e_j, v )(v, e_i) \, \vol_l \\
  &= \int_{\mathbb{P}(V)} (s_{e_j}(l), s_{e_i}(l)) \vol_l \\
  &= \langle s_j, s_i \rangle \\
  &= \delta_{ij}.
 \end{align*}
A similar calculation shows that $\langle \text{ev}^\mi \rangle \circ \langle \text{ev} \rangle$ also equals the Bergman kernel $\langle \rangle$. This establishes that the counit map $\langle \text{ev} \rangle$ is an isomorphism. We have already proved In Lemma \ref{myoneda} that every Hilbert space $V$ is antilinearly isomorphic to the space of sections of its geometric line bundle; that is we established an isomorphism
 \[
  \eta_V \colon V \rightarrow \overline{\Gamma(\tau_V)}.
  \]
These isomorphisms are natural with respect to $V$, because if $A \colon V \rightarrow W$ is a linear map, with $l \in \mathbb{P}(V)$ and $l' \in \mathbb{P}(W)$, and $v \in V$, then
 \begin{align*}
  \Big[\overline{\Gamma} (\langle f \rangle) \circ \eta_V (v) \Big] (l') &= \int_{\mathbb{P}(V)} \langle l' | A | l \rangle s_v (l) \; \vol_l \\
   &= \langle l' | \langle A \rangle \circ \langle \; \rangle | l \rangle v \\
   &= \proj_{l'} (Av) = \Big[\eta_W \circ A(v) \Big] (l').
  \end{align*}
Hence the functors described above are mutual inverses of each other up to natural isomorphism.
\end{proof}
One caveat to this result is that the functors $\Hilb \rightarrow \LBun$ and $\LBun \rightarrow \Hilb$ defined above are {\em antilinear} at the level of hom-sets. It is not possible to avoid a `wrinkle' like this. This is because regarding a vector $v$ in a Hilbert space $V$ as a section $s_v$ of the geometric line bundle
 \begin{align*}
  V & \rightarrow \overline{\Gamma(L)} \\
  v & \mapsto s_v
  \end{align*}
is analogous to the Yoneda embedding of a category $C$ into its category of presheaves,
\begin{align*}
  C &\hookrightarrow [C^\text{op}, \Set] \\
  v &\mapsto \Hom(-, v),
\end{align*}
where the `op' is unavoidable. However, we could have formulated our result differently, for instance we could have defined the functor $\Hilb \rightarrow \LBun$ to be {\em contra}variant, in which case it would indeed be a linear functor. But that would undermine the entire purpose of this chapter, which was to show how the geometric picture for `categorified vector spaces' we develop in Chapter \ref{2HilbChap} also has a counterpart at the level of {\em ordinary} vector spaces. In Chapter \ref{2HilbChap}, all of our functors are covariant, which has caused us to make this choice here.

\section{Unitary representations and equivariant line bundles\label{B7}}
In this section we show that the relationship between finite dimensional Hilbert spaces and holomorphic line bundles continues to hold in the presence of a group action.

Let $G$ be a Lie group --- it could be finite, compact or noncompact --- and $X$ a compact hermitian manifold on which $G$ acts by orientation preserving isometries. A {\em unitary equivariant holomorphic hermitian line bundle over $X$} is a holomorphic hermitian line bundle $L \rightarrow X$ equipped with a lift of the action of $G$ on $X$ --- that is, for each $g \in G$ and $x \in X$, there is a unitary map
 \[
  g \colon L_x \rightarrow L_{g \cdot x},
 \]
and these maps vary smoothly with $g$ and holomorphically in $x$. If $L \rightarrow X$ and $Q \rightarrow Y$ are two such equivariant line bundles, then an {\em equivariant holomorphic kernel} from $L$ to $Q$ is a holomorphic kernel $E \in \Gamma(Q \boxtimes \overline{L})$ such that the following diagram commutes for all $x \in X$, $y \in Y$ and $g \in G$:
 \be \label{eqdm}
  \ba \xymatrix{ L_x \ar[rr]^{\langle y | E | x \rangle} \ar[d]_g && Q_y \ar[d]^g \\ L_{g \cdot x} \ar[rr]_{\langle g \cdot y | E | g \cdot x \rangle} && Q_{g \cdot y}} \ea
 \ee
In this way, $G$-equivariant line bundles assemble into a category $\LBun(G)$.

\begin{defn} Given a Lie group $G$, the category $\LBun(G)$ of $G$-equivariant line bundles is defined as follows. An object is a unitary equivariant holomorphic hermitian line bundle $L \rightarrow X$ over a compact hermitian manifold $X$ on which $G$ acts by orientation preserving isometries, and a morphism is an equivariant holomorphic kernel.
\end{defn}

We write $\Rep(G)$ for the category of unitary representations of $G$ on finite-dimensional Hilbert spaces and intertwiners between them. Just as in the previous section, we now write down functors $\Rep(G) \rightarrow \LBun(G)$ and $\LBun(G) \rightarrow \Rep(G)$ and prove they establish an equivalence of categories.

If $V$ is a unitary representation of $G$ then $G$ acts by orientation preserving isometries on $\mathbb{P}(V)$ and hence the geometric line bundle $\tau_V$ becomes a unitary equivariant holomorphic hermitian line bundle over $\mathbb{P}(V)$. Moreover, if $A \colon V \rightarrow W$ is an intertwiner of representations, then the corresponding kernel $\langle A \rangle \colon \tau_V \rightarrow \tau_W$ is equivariant --- because if $l$ and $l'$ are lines in $V$ and $W$ respectively, with $v \in l$ and $w \in l'$ of unit norm, then
 \begin{align*}
  g \cdot \proj_{l} (Av) &= g \cdot (w, Av) w \\
   &= (w, Av) g \cdot w \\
   &= (g \cdot w, g \cdot Av) (g \cdot w) \\
   &= \proj_{g \cdot l}(A (g \cdot v)).
 \end{align*}
This establishes a functor $\Rep(G) \rightarrow \LBun(G)$.

In the reverse direction, suppose that $L$ is a unitary $G$-equivariant holomorphic hermitian line bundle over a compact hermitian manifold $X$. Then the space of sections $\Gamma(L)$ (and hence also its complex conjugate $\overline{\Gamma(L)}$) inherits a representation of $G$ via the prescription
 \[
  (g \cdot s)(x) := g \cdot s(g^\mi \cdot x).
 \]
This representation is unitary, because
 \begin{align*}
  \langle g \cdot s, g \cdot s' \rangle &= \int_X  (( g \cdot s)(x), (g \cdot s' )(x)) \, \vol_x \\
   &= \int_X  (g \cdot s(g^\mi \cdot x), g \cdot s'(g^\mi \cdot x)) \, \vol_x \\
    &= \int_X ( s(g^\mi \cdot x), s'(g^\mi \cdot x)) \, \vol_x \\
    &= \int_X (s(x), s'(x)) \, \vol_x \\
    &= \langle s, s' \rangle
 \end{align*}
where the second-last step uses the fact that $G$ acts on $X$ by orientation preserving isometries. A similar argument shows that the linear map $A \colon \overline{\Gamma(L)} \rightarrow \overline{\Gamma(Q)}$ induced by an equivariant holomorphic kernel $\langle A \rangle \colon L \rightarrow Q$ is an intertwiner. This establishes a functor $\LBun(G) \rightarrow \Rep(G)$.

\begin{thm} The above constructions establish an equivalence $\Rep(G) \simeq \LBun(G)$ between the category of unitary representations of a Lie group $G$ and the category of $G$-equivariant holomorphic hermitian line bundles over compact hermitian manifolds.
\end{thm}
\begin{proof} We show that the functor $\LBun(G) \rightarrow \Rep(G)$ is essentially surjective and fully faithful. To show that it is essentially surjective, all we need to do is show that if $V$ is a unitary representation of $G$, then the isomorphism
 \begin{align*}
   V &\rightarrow \overline{\Gamma(\tau_V)} \\
   s &\mapsto s_v = \proj_{\cdot} (v)
 \end{align*}
from Lemma \ref{myoneda} is equivariant. Indeed, we have
 \begin{align*}
  s_{g \cdot v} (l) &= \proj_l (g \cdot v) \\
                    &= (u, g \cdot v) u \quad \text{(where $u \in l, (u,u)=1$)} \\
                    &= (g^\mi \cdot u, v) u \quad \text{(since $g$ is a unitary map)} \\
                    &= g \cdot (g^\mi \cdot u, v) (g^\mi \cdot u) \\
                    &= g \cdot \proj_{g^\mi \cdot l} (v) \\
                    &= (g \cdot s_v) (l).
  \end{align*}
To see that the functor is fully faithful, suppose that $L$ and $Q$ are equivariant line bundles. By the sequence of canonical isomorphisms presented in \eqref{seqisos}, we can identify the space of intertwiners from $\Gamma(L)$ to $\Gamma(Q)$ as the $G$-invariant sections of the homomorphism line bundle:
 \[
  \Hom_{\Rep(G)} (\Gamma(L), \Gamma(Q)) \cong \Gamma(\Hom(L, Q))^G.
  \]
But a $G$-invariant section of the homomorphism line bundle is precisely a kernel satisfying the commutative diagram \eqref{eqdm}. This completes the proof.
\end{proof}

\section{Geometric characters of equivariant line bundles\label{B8}}
In this section we define the {\em geometric character} of an equivariant line bundle $L$ over a space $X$ as a certain integral over $X$, and we show that it computes as the character of the associated representation of $G$ on $\Gamma(L)$. This is a purely formal result, but by the Atiyah-Singer equivariant index theorem this integral can actually be expressed as a certain integral over the {\em fixed points} of the action of $G$ on $X$. This is significant because it is precisely this idea that gets `categorified' in Chapters \ref{GerbesChap} and \ref{GerbesCharChap}.

If $L$ is an equivariant line bundle over a compact hermitian manifold $X$, we define the {\em associated kernel} of a group element $g \in G$ as the collection of maps
 \[
 \langle y | g | x \rangle \colon L_x \rightarrow L_y
\]
obtained by postcomposing with the Bergman kernel:
\[
  L_x \ra{g} L_{g \cdot x} \ra{\langle y | g \ccdot x \rangle} L_y.
 \]
We define the {\em geometric character} $\ch_L$ of $L$ as the function on $G$ defined by integrating the associated kernel over the diagonal:
 \[
  \ch_L (g) = \int_X \langle x | g | x \rangle \, \vol_x
 \]
We can compute the geometric character as follows.
\begin{prop} The geometric character of an equivariant line bundle $L$ computes as the character of the associated representation of $G$ on $\Gamma(L)$:
 \[
 \ch_L (g) = \Tr(g \colon \Gamma(L) \rightarrow \Gamma(L)).
 \]
\end{prop}
\begin{proof}
 We have
  \begin{align*}
  \ch_L (g) &= \int_X \langle x | g | x \rangle \\
   &= \int_X \left(v, \sum_i (s_i(g \ccdot x), g \ccdot v) s_i(x) \right) \; \vol_x \quad \text{(where $v \in L_x, (v,v) =1$)} \\
   &= \sum_i \int_X \left(s_i (g \ccdot x), g \ccdot v\right) \; \left(v, s_i (x)\right) \; \vol_x \\
   &= \sum_i \int_X \left(s_i (g \ccdot x), g \ccdot v\right) \; \left(g \ccdot v, g \ccdot s_i (x)\right) \; \vol_x \\
   &= \sum_i \int_X  \left(s_i(g \ccdot x), g \ccdot s_i (x)\right) \; \vol_x \\
   &= \sum_i \int_X  \left(s_i(x), g \ccdot s_i (g^\mi \ccdot  x)\right) \; \vol_x \\
   &= \sum_i \langle s_i, g \ccdot s_i \rangle \\
   &= \Tr(g \colon \Gamma(L) \rightarrow \Gamma(L)).
    \end{align*}
  \end{proof}
We now write this result in such a way as to make the analogy with our corresponding result in Chapter \ref{GerbesCharChap} more explicit. To do this, let us assume that that the Lie group $G$ is compact (we are allowing here that it be finite). In that case we can sensibly speak about the Hilbert space $L^2(\text{Class}(G))$ of square integrable class functions on $G$, because we can use the bi-invariant Haar measure on $G$. We also have the Grothendieck ring $[\Rep(G)]_\mathbb{C}$, the vector space generated by the isomorphism classes of unitary representations of $G$, which is equipped with the inner product
 \[
  ([\rho], [\sigma]) = \dim \Hom(\rho, \sigma).
 \]
It is well known \cite{ref:brocker_dieck} that the character map provides a unitary isomorphism from the Hilbert space completion of the Grothendieck ring to the Hilbert space of class functions:
 \[
  \chi \colon \widehat{[\Rep(G)]_\mathbb{C}} \rightarrow \Class (G).
 \]
All these facts might have analogues for noncompact groups, but this is beyond the expertise of the present author.
\begin{cor}[{Compare Chapter \ref{GerbesCharChap}, Thm \ref{transthm}}] If the Lie group $G$ is compact, we have a commutative diagram of unitary isomorphisms of Hilbert spaces:
 \[
  \ba \xymatrix @1 @C=0.1in { \widehat{[\Rep(G)]_\mathbb{C}} \ar[rr]^\cong \ar[dr]_{\chi}  && \widehat{[\LBun(G)]_\mathbb{C}} \ar[dl]^{\ch} \\ & \Class(G)} \ea.
 \]
\end{cor}
We close this chapter by providing a concrete formula for the geometric character using index theory. The idea is that if the curvature of the line bundle $L$ is sufficiently positive, then all higher cohomology groups vanish, so that we can compute things about $\Gamma(L)$ using formulas involving the entire cohomology group $H^* (X, L)$. We will also assume in addition that the base manifold $X$ is {\em K\"{a}hler}. It is not actually necessary to make this assumption in order to write the geometric character as an integral over the fixed points, but without it we do not have a direct relationship between the Bergman kernel and the large time limit of the kernel of the Dirac operator \cite[pg 21]{ref:duistermaat}.

Let us explain the terms in the proposition below; our main reference is \cite[chap. 6]{ref:berline_getzler_vergne}, with input from
\cite[chap. 12]{ref:duistermaat} and \cite{ref:atiyah_singer}. The geometric data is a $G$-equivariant hermitian holomorphic line bundle $L$ over a compact $n$-dimensional K\"{a}hler manifold $X$ and an element $g \in G$. The canonical line bundle of $X$ is the holomorphic line bundle $K = \Gamma^n (T^{1,0} X)^\vee$. The fixed point locus $X^g$ forms a submanifold of $X$, and the normal bundle along $X^g$ is denoted $\mathcal{N}$. The Riemannian curvature of $X$ restricted to the fixed point submanifold $X^g$ splits up as
 \[
  R|_{X^g} = R^0 \oplus R^1
 \]
with $R^0 \in \mathcal{A}^2(X^g, \mathfrak{so}(X^g))$ and $R^1 \in \mathcal{A}^2 (X^g, \mathfrak{so}(\mathcal{N}))$, where the notation $\mathcal{A}^k(M, \mathcal{E}) \equiv \Gamma(M, \Lambda T^k M \otimes \mathcal{E})$ refers to the space of differential $k$-forms on a manifold $M$ with values in a vector bundle $\mathcal{E}$. The Todd genus form of $X^g$ is
 \[
  \Td(X^g) = \text{det}\left(\frac{R_0}{e^{R_0} - 1}\right) \in \mathcal{A}^{2 \bullet} (X^g, \mathbb{C}).
 \]
The curvature of the line bundle $L$ will be denoted as $c(L) \in \mathcal{A}(X, \mathbb{C})$.

Since $X^g$ is fixed by $g$, each  tangent space $TX_x$ for $x \in X^g$ is a representation of $G$; moreover with respect to the orthogonal splitting
 \[
  TX_x = TX^g_x \oplus \mathcal{N}_x
 \]
we see that $g$ acts as the identity on $TX^g_x$ so that the only interesting representation theoretic data is the representation of $G$ on the normal tangent space $\mathcal{N}_x$. Finally, the fixed point submanifold $X^g$ splits up into a finite number of connected components which we write as $X^g_m$, and for $x \in X^g_m$ the action of the group element $g$ on the complex line $L_x$ is by definition unitary and hence can be written as multiplication by $e^{i \theta_m} \in U(1)$.

We can now state the formula.

\begin{prop}[{Compare Chapter \ref{GerbesChap}, Eqn \ref{needlabel}}] \label{Lastprop} Suppose $L \in \LBun(G)$ is an equivariant line bundle over a compact K\"{a}hler manifold $X$. If the line bundle $K^\vee \otimes L$ is positive, then the geometric character computes as an integral over the fixed points:
 \[
  \ch_L (g) = i^{-\dim (X)/2} \sum_m \frac{e^{i \theta_m}}{(2 \pi)^{\dim(X^g_m)/2}}\int_{X^g_m} \frac{ \Td(X^g) e^{-c(L)}}{\det \left(1 - ge^{-R^1} \right)}.
 \]
\end{prop}
\begin{proof} Firstly, the Kodaira vanishing theorem asserts that the condition that $K^\vee \otimes L$ is positive is precisely the condition that the higher cohomology groups $H^i(X, L)$ vanish (see \cite{ref:poon} and the references therein) so the only nonvanishing cohomology group is $\Gamma(L) \equiv H^0(X, L)$. Now, recall that the Bergman kernel can be expressed as the large time limit of the heat kernel. Unravelling our definitions, from first principles we thus have
 \begin{align*}
  \ch_L (g) = \int_X \langle x | g | x \rangle \, \vol_x &= \lim_{t \rightarrow \infty} \int_X \Tr \; \langle x | g e^{-tD^2} | x \rangle|_{\Gamma(L)} \, \vol_x \\
   &= \lim_{t \rightarrow \infty} \int_X \text{STr} \; \langle x | g e^{-tD^2} | x \rangle_{H^*(X, L)} \, \vol_x.
  \end{align*}
By the McKean-Singer formula the integral in the last line does not actually depend on $t$, thus we may evaluate it in the limit $t \rightarrow 0$. In this limit, the heat kernel $\langle y | e^{-tD^2} | x \rangle$ tends towards a $\delta$-function, so that the integral localizes over the fixed points $x \in X^g$. The local version \cite[Thm 6.11]{ref:berline_getzler_vergne} of the equivariant Atiyah-Segal-Singer index theorem \cite{ref:atiyah_singer} precisely computes the resulting top-degree differential form. At a fixed point $x \in X^g_m$, this formula gives:
 \[
   \lim_{t \rightarrow 0} \; \text{STr} \; \langle x | g e^{-tD^2} | x \rangle = \frac{i^{-\dim (X)/2} e^{i \theta_m}}{(2 \pi)^{\dim(X^g_m)/2}} \frac{ \Td(X^g) e^{- c(L)}}{\det \left(1 - ge^{-R^1} \right)} .
 \]
This is the integrand we used in the statement of the proof.
\end{proof}

\chapter{2-Hilbert spaces\label{2HilbChap}}

In this chapter we review the notion of a {\em 2-Hilbert space} due to Baez \cite{ref:baez_2_hilbert_spaces}. A 2-Hilbert space can be thought of as the `categorification' of the notion of a finite-dimensional Hilbert space, and we need them because we want to study {\em unitary} 2-representations of groups. Although much of the material in this chapter is taken from \cite{ref:baez_2_hilbert_spaces}, we will be framing some of the basic definitions and the theory a bit differently. In particular, we will define a 2-Hilbert space as a {\em Cauchy complete} $H^*$-category instead of as an {\em abelian} $H^*$-category; these notions are equivalent but we would like to steer clear of using `abelian categories' because we would like to align the theory more closely with its `derived' cousins, such as the Calabi-Yau $A_\infty$-categories of Costello \cite{ref:costello}.

In addition we wish to emphasize the {\em geometric} picture of 2-Hilbert spaces analogous to that presented in Chapter \ref{GeometryOrdinaryRepChap}: just as a Hilbert space may be regarded as the space of sections of a hermitian line bundle over a compact complex manifold equipped with a metric, a 2-Hilbert space can be regarded as the space of sections of a `categorified hermitian line bundle' over a finite space equipped with a `metric'.

In Section 1 we define $H^*$-categories. In Section 2 we define 2-Hilbert spaces as Cauchy complete $H^*$-categories, and we establish the equivalence of this definition to that of Baez. We also compare 2-Hilbert spaces to the {\em 2-vector spaces} of Kapranov and Voevodsky. In Section 3 we define the 2-category $\THilb$ of 2-Hilbert spaces, characterize the 1-morphisms and 2-morphisms in $\THilb$, and define what it means for a pair of 2-Hilbert spaces to be {\em strongly unitarily equivalent}. We also show that {\em every unitary commutative Frobenius algebra is isomorphic to the algebra of natural transformations of the identity on a 2-Hilbert space}, which follows up on our discussion in this regard from the Introduction. Finally in Section 4 we define a 2-category $\FinSpaces$ of finite sets `equipped with metrics', and we prove the equivalence $\THilb \simeq \FinSpaces$.

\section{$H^*$-categories\label{hcats}}
Consider the category $\Rep(G)$ of finite-dimensional unitary representations of a compact Lie group $G$ (more generally, one could consider {\em twisted} unitary representations of a compact Lie groupoid; for the finite case see \cite{ref:simon}). Since its objects are {\em unitary} representations, the category $\Rep(G)$ has two salient features which the ordinary category of `mere' representations on vector spaces does not have. Firstly, every morphism $f \colon V \rightarrow W$ has an adjoint $f^* \colon W \rightarrow V$, and secondly there is an inner product on the hom-sets\footnote{The factor $\frac{1}{|G|}$ arises due to the fact that we should really think of the inner product as an integral over the groupoid $\BG$, and the correct measure to use on a groupoid is that obtained by dividing by the size of the automorphism groups as in \cite{ref:simon, ref:baez_dolan_1}.},
\be \label{definnp}
 (f,g) = \frac{1}{|G|}\Tr (f^* g).
\ee
We abstract these two features into the definition of an {\em $H^*$-category}. Write $\Hilb$ for the category of finite dimensional Hilbert spaces and arbitrary linear maps; a Hilb-category is a category enriched over Hilb.

\begin{defn}[{Baez \cite[Prop 3]{ref:baez_2_hilbert_spaces}}]An {\em $H^*$-category} is a $\Hilb$-category $H$ equipped with antilinear maps $* \colon \Hom(x,y) \rightarrow \Hom(y,x)$ for all $x,y \in H$, such that
  \begin{itemize}
   \item $f^{* *} = f$,
   \item $(f g)^* = g^* f^*$,
   \item $(f g, h) = (g, f^* h)$,
   \item $(f g, h) = (f, h g^*)$
   \end{itemize}
  whenever both sides of the equation are defined.
\end{defn}
The reason for the terminology `$H^*$-category' is because there is an already established notion of an {\em $H^*$-algebra} (due to Ambrose \cite{ref:ambrose}), and an $H^*$-category can be thought of as a `many-object $H^*$-algebra'. In the finite-dimensional setting, Baez has recast the definition of Ambrose as follows.
\begin{defn}[{Baez \cite{ref:baez_2_hilbert_spaces}}, {Ambrose \cite{ref:ambrose}}] An {\em $H^*$-algebra} is $A$ is a finite dimensional Hilbert space that is also an associative algebra with unit, equipped with an antilinear involution $* \colon A \rightarrow A$ satisfying
 \begin{align*}
  (ab, c) &= (b, a^*c) \\
  (ab, c) &= (a, cb^*)
 \end{align*}
for all $a,b,c \in A$. An isomorphism of $H^*$-algebras is a unitary map that is also an involution-preserving algebra isomorphism.
\end{defn}
The basic example of an $H^*$-algebra is the space of linear operators on a finite-dimensional Hilbert space $H$, where the product is the composite of operators, the involution is the usual adjoint operation, and the inner product is given by
 \[
  (a, b) = k \Tr (a^*b)
 \]
for some real number $k > 0$. It follows from the work of Ambrose that {\em all} finite-dimensional $H^*$-algebras are isomorphic to an orthogonal direct sum of $H^*$-algebras of this form. The scale factor $k$ is important and highlights the difference between $H^*$-algebras and (finite-dimensional) $C^*$-algebras, where such a factor does not enter.

The link between $H^*$-algebras and $H^*$-categories is that an $H^*$-algebra is clearly the same thing as a one-object $H^*$-category, while from a finite set $S$ of objects in an $H^*$-category one may obtain an $H^*$-algebra by forming the `arrow algebra'
 \[
  A = \bigoplus_{x,y \in S} \Hom(x,y).
 \]

\section{2-Hilbert spaces}
The notion of an $H^*$-category should be thought of as the `categorification' of the notion of an inner product space. Now, recall that a {\em Hilbert space} is defined to be a Cauchy complete inner-product space. There is an elegant way to formulate the notion of `Cauchy completeness' for arbitrary enriched categories due to Lawvere \cite{ref:lawvere}. If one thinks of a metric space as a category enriched in $[0, \infty]$, one recovers the ordinary notion of Cauchy completeness.

Lawvere's notion of Cauchy completeness in an enriched setting works as follows. Let $\mathcal{V}$ be the monoidal category over which the categories in question are enriched (in our setting of $H^*$-categories, $\mathcal{V} = \Hilb$). We will use the convention that a {\em bimodule} $E \colon A \rightsquigarrow B$ between $\mathcal{V}$-enriched categories is a $\mathcal{V}$-functor
 \[
  E \colon B \otimes A^{\text{op}} \rightarrow \mathcal{V}
 \]
where $B \otimes A^{\text{op}}$ is the category whose object set is $\text{Ob} \,B \times \text{Ob} \, A$ and whose hom-sets are the tensor products of the hom-sets between the factors. Note that every $\V$-functor $F \colon A \rightarrow B$ gives rise canonically to a pair of bimodules
 \[
  F_* \colon A \rightsquigarrow B \quad \text{and} \quad F^* \colon B \rightsquigarrow A
 \]
defined by setting
 \[
  F_* (b, a) = \Hom(F(a), b) \quad \text{and} \quad F^*(a, b) = \Hom(b, F(a)).
 \]
The essential characteristic of such a pair of bimodules coming from a $\mathcal{V}$-functor $F$ is that $F^*$ is the {\em right adjoint} of $F_*$. We can now give the definition of Cauchy completeness.
\begin{defn}[{Lawvere \cite{ref:lawvere}}] A $\mathcal{V}$-category $B$ is {\em Cauchy complete} if every bimodule $E \colon A \rightsquigarrow B$ having a right adjoint is isomorphic to a bimodule of the form $F_*$ for a $\mathcal{V}$-functor $F \colon A \rightarrow B$.
\end{defn}
We encourage the reader to read the paper of Lawvere to obtain a richer understanding of this notion. Recalling that $H^*$-categories are Hilb-enriched categories with some extra structure, this allows us to make the following definition.

\begin{defn} A {\em 2-Hilbert space} is a Cauchy complete $H^*$-category.
\end{defn}
Intuitively speaking, this says that a 2-Hilbert space is an $H^*$-category which does not have any `holes' in it. This definition is on the face of it different to that of Baez in \cite{ref:baez_2_hilbert_spaces}, where a 2-Hilbert space is defined as an {\em abelian} $H^*$-category. We prefer to follow M\"{u}ger \cite[page 6]{ref:müger_from_subfactors_to_categories_and_topology_I} and avoid using the language of abelian categories, because it has been our experience that the notion of kernels and cokernels plays no conceptual role when performing calculations in 2-Hilbert spaces, being replaced instead by the more important notions of {\em subobject} and {\em direct sum}. Moreover, we will see in Chapter \ref{yausec} that a reasonable notion of a `nondiscrete' 2-Hilbert space is the graded derived category of coherent sheaves $\mathbf{D}(X)$ on a compact Calabi-Yau manifold $X$. These categories are definitively not abelian, so we have steered clear of using the term `abelian category' in our definition of a 2-Hilbert space. We hope that our definition also provides a cleaner way to see a 2-Hilbert space as the `categorification' of a Hilbert space.

Nevertheless, the following proposition establishes the equivalence of these definitions, and shows that being Cauchy complete is equivalent to requiring that the category has all direct sums and subobjects.

\begin{prop} The following conditions on an $H^*$-category are equivalent:
 \begin{enumerate}
  \item It is Cauchy complete (that is, it is a 2-Hilbert space).
  \item Every pair of objects has a direct sum, and all idempotents split.
  \item It is abelian.
 \end{enumerate}
 \end{prop}
\begin{proof} (i) $\Rightarrow$ (ii). Street has shown that {\em an enriched category is Cauchy complete if and only if it admits all absolute colimits} \cite{ref:street_colimits}. An absolute colimit is one which is preserved by all $\mathcal{V}$-functors, that is it is one which is {\em defined by equations}. Since direct sums and idempotent splittings are defined by equations, this proves (i) $\Rightarrow$ (ii).

(ii) $\Rightarrow$ (iii) In Proposition 10 of \cite{ref:baez_2_hilbert_spaces}, Baez used the classification of $H^*$-algebras to show that if an $H^*$ category has a zero object, every pair of objects has a direct sum, and every morphism has a kernel, then it is an abelian category which is moreover semisimple (see below for the definition of `semisimple'). We obtain a zero object by splitting the zero morphism regarded as an idempotent. We only need to show that the existence of kernels is guaranteed by the fact that all idempotents split, and this follows from the classification of $H^*$-algebras. For suppose we are given a morphism $f \colon x \rightarrow y$. Then the $H^*$-algebra
 \[
  A = \Hom(x,x) \oplus \Hom(x,y) \oplus \Hom(y,x) \oplus \Hom(y,y)
 \]
can be written as
 \[
  A = \bigoplus_{i=1}^n L^2(H_i^x \oplus H_i^y, k_i)
 \]
for finite-dimensional Hilbert spaces $H_i^x$ and $H_i^y$ (one should think of $H^i_x$ as $\Hom(e_i, x)$ where the $e_i$ are a representative set of simple objects). Thus we can write $f$ as
\[
 f = \bigoplus_i f_i
\]
where $f_i \colon H^i_x \rightarrow H^i_y$. Let $p_i$ be the orthogonal projection onto $\text{ker}(f_i) \subset H^i_x$. Then $p = \sum_i p_i$ is an idempotent on $x$ and a splitting of $p$ provides a kernel for $f$.

(iii) $\Rightarrow$ (i) Baez has shown that if the category is abelian then it must be semisimple. And if it is semisimple then it is equivalent to $\Hilb^n$ and hence must be Cauchy complete since the completion of $\mathbb{C}$ (thought of as a one object Hilb-enriched category) is equivalent to $\Hilb$; this is a general result about categories enriched in vector spaces \cite{ref:borceaux_dejean}.
\end{proof}
We would like to emphasize that the notion of Cauchy completeness is the natural friend of mathematics in a `derived' setting, due to the following observation by Street \cite[pg 3]{ref:street_colimits} which we used above:
 \begin{quotation}
 Hence a {\em category is Cauchy complete if and only if it admits all absolute colimits.} As an example, since mapping cones and suspensions can be defined equationally for dg-categories (categories enriched in complexes of abelian groups), Cauchy complete dg-categories admit suspensions and mapping cones (as well as the expected direct sums and splittings for idempotents).
 \end{quotation}

We follow the formulation of M\"{u}ger \cite{ref:müger_from_subfactors_to_categories_and_topology_I} and say that a linear category is {\em semisimple} if it has direct sums and subobjects and if there exist objects $e_i$ labelled by a set $I$ such that $\Hom(e_i, e_j) \cong \delta_{ij}\mathbb{C}$ (such objects are called {\em simple}) and such that for any two objects $x$ and $y$ the composition map
  \[
   \bigoplus_{i \in I} \Hom(x, e_i) \otimes \Hom(e_i, y) \rightarrow \Hom(x,y)
  \]
is an isomorphism. In other words, a semisimple linear category is one which is equivalent to $\Vect^I$ for some set $I$. It is called {\em finite dimensional} if the index set $I$ is finite --- that is, if there are only finitely many isomorphism classes of simple objects. The following fact, which we have encountered already, is very important.
\begin{prop}[{Baez \cite[Prop 10]{ref:baez_2_hilbert_spaces}}] Every 2-Hilbert space is semisimple.
\end{prop}
This proof uses the classification of $H^*$-algebras, which shows that the algebra generated by the arrows between any fixed set of objects in a 2-Hilbert space must essentially break up as a sum of matrix algebras, from which the result follows. However, it is important to understand that a 2-Hilbert space $H$ has more information in it than merely being a semisimple linear category, because every simple object $e_i$ comes equipped with a positive real {\em scale factor}
 \[
  k_i := (\id_{e_i}, \id_{e_i}).
 \]
We think of these numbers as equipping the index set $I$ with a `Riemannian metric'. The set $I$ of simple objects together with their scale factors  characterizes the 2-Hilbert space $H$ up to {\em strong unitary equivalence}, a concept we will define in the next section.

Given an object $x \in H$, we say that a basis
 \[
   \{ a_p \colon e_i \rightarrow x \}_{p=1}^{\dim \Hom(e_i, x)}
 \]
for $\Hom(e_i, x)$ is a {\em $*$-basis} if
 \[
 a_p^* a_q = \delta_{pq} \id_{e_i}.
 \]
Let us agree to use the following notation to strip off scalars from scalar multiples of identity morphisms on simple objects:
 \[
  \langle \lambda \id_{e_i} \rangle := \lambda \in \mathbb{C}.
 \]
Then since
 \[
 (a_p, a_q) = (\id_e, a_p^* a_q) = (\id_e, \id_e) \langle a_p^* a_q \rangle = k_i \langle a_p^* a_q \rangle,
 \]
a $*$-basis is the same thing as an orthogonal basis where each basis element has norm $\sqrt{k_i}$; in particular a $*$-basis always exists. The fact that 2-Hilbert spaces are semisimple means that we can expand any morphism $f \colon x \rightarrow y$ in terms of its `matrix elements' as
\[
  f = \sum_{i,q,p} \langle b_{i,q}^* f a_{i,p} \rangle \,  b^*_{i,q} a_{i,p}
  \]
where $\{a_{i,p} \colon e_i \rightarrow x\}$ and $\{b_{i,p} \colon e_i \rightarrow y\}$ are $*$-bases for the hom-sets $\Hom(e_i, x)$ and $\Hom(e_i, y)$ as $e_i$ runs over the representative simple objects of $H$.

\subsection*{Decategorifying a 2-Hilbert space}
If $H$ is a 2-Hilbert space, we will write $[H]_\mathbb{C}$ for the {\em complexified Grothendieck group} of $H$ --- the tensor product of $\mathbb{C}$ with the abelian semigroup generated by the isomorphism classes of objects $[v]$ in $H$ under the relations $[v \oplus w] = [v] + [w]$. A basis for $[H]_\mathbb{C}$ is given by the isomorphism classes $[e_i]$ of simple objects of $H$. We regard $[H]_\mathbb{C}$ as a Hilbert space with inner product defined on the generating elements by
 \[
 ([v], [w]) = \dim \Hom(v, w).
\]
Note that the Grothendieck group $[H]_\mathbb{C}$ only uses the structure of $H$ as a semisimple category; it does not `see' the scale factors on the simple objects.

\subsection*{Comparison with 2-vector spaces}
The usual notion of `categorified vector space' that has been used in a 2-representation theory context has been that of a a {\em 2-vector space} due to Kapranov and Voevodsky \cite{ref:kapranov_voevodsky}, as in the works of Elgueta \cite{ref:elgueta}, Barrett and Mackaay \cite{ref:barrett_mackaay} and Ganter and Kapranov \cite{ref:ganter_kapranov_rep_char_theory}; see also the work of Morton \cite{ref:morton} and also the additional paper by Elgueta \cite{ref:elgueta3}. A 2-vector space is defined essentially as an `abelian $\Vect$-module category equivalent to $\Vect^n$ for some $n$'. One can think of a finite-dimensional 2-Hilbert space as a 2-vector space `equipped with an inner product'. However it is important to note that their definitions have entirely different flavours: a 2-Hilbert space is defined {\em intrinsically} via the duality on the hom-sets and the inner products, with semisimplicity --- the fact that the resulting category is equivalent to $\Vect^n$ --- as a consequence, whereas for a 2-vector space semisimplicity is added in by hand.

This underscores the considerable constraint that a $*$-structure places on a linear category. For instance, the category $\Rep(A)$ of finite-dimensional representations of an algebra $A$ is always an abelian $\Vect$-module category, but it fails to be a 2-vector space in general precisely because of the lack of duality on the hom-sets in $\Rep(A)$.

However, one advantage of a 2-vector space is that it explicitly includes a prescription for categorified scalar multiplication --- to say that it is a $\Vect$-module category is to say that one can `tensor' an object of a 2-vector space with a vector space to form another object. This construct is missing from the definition of a 2-Hilbert space. On the other hand, since every finite-dimensional 2-Hilbert space is linearly equivalent to $\Vect^n$ for some $n$, it can always be added in artificially, but this is an unnatural way to proceed and it may prove advantageous to add this ability explicitly to the definition of a 2-Hilbert space.

In any event, we wish to underscore once more our main reason for working with 2-Hilbert spaces: our motivation is the Extended TQFT Hypothesis of Baez and Dolan \cite{ref:baez_dolan_hda0}, where we are required to work with {\em unitary} structures. A 2-Hilbert space is a richer structure than a 2-vector space --- a 2-vector space is classified up to equivalence simply by the number of nonisomorphic simple objects, while we will see that a 2-Hilbert space is classified up to strong unitary equivalence by the simple objects {\em and their scale factors}. This is crucial in order to make the link with extended TQFT, because we will shortly see it is precisely these scale factors which encode the Frobenius algebra information on the vector space assigned to the circle in a 2d extended TQFT, and hence confirm the `primacy of the point' as we explained in the introduction to this thesis.

\section{The 2-category of 2-Hilbert spaces\label{BigTCat}} In this section we recall the definition of the 2-category $\THilb$ of 2-Hilbert spaces, as in \cite{ref:baez_2_hilbert_spaces}. Then we characterize the morphisms and 2-morphisms in $\THilb$; in the latter case we also define the {\em inner product} between 2-morphisms. We end the section by classifying 2-Hilbert spaces up to {\em strong unitary equivalence}.

A functor $F \colon H \rightarrow H'$ between 2-Hilbert spaces is called {\em linear} if it is linear on the level of hom-sets. It is called a {\em $*$-functor} if $F(f^*) = F(f)^*$ for all morphisms $f$ in $H$.

\begin{defn} The 2-category $\THilb$ of 2-Hilbert spaces has finite-dimensional 2-Hilbert spaces for objects,
linear $*$-functors for morphisms, and natural transformations for 2-morphisms.
\end{defn}
This definition is slightly different to that of Baez \cite{ref:baez_2_hilbert_spaces}, where the functors between 2-Hilbert spaces are additionally required to preserve direct sums (that is, they are required to be exact functors between the underlying abelian categories, which is equivalent). But since direct sums are defined by {\em equations} and not by a universal property, any linear functor automatically preserves them, so this requirement is actually superfluous.

In the remainder of this section we characterize the morphisms and 2-morphisms in $\THilb$.

\subsection{Characterizing linear $*$-functors}
Our first result is that a linear $*$-functor between 2-Hilbert spaces is roughly the same thing as a `matrix of vector spaces'. This is well-known for linear functors between semisimple categories (see eg. \cite{ref:bakalov_kirillov}), but we need to show that the requirement $F(f^*) = F(f)^*$ does not introduce significant additional constraints.

We will find it convenient here and elsewhere to work with {\em marked} 2-Hilbert spaces, by which we mean a 2-Hilbert space where each isomorphism class of simple objects has been endowed with a distinguished representative $e_i$. In practice many 2-Hilbert spaces actually arrive in this way, for instance it is common to have certain preferred choices for the irreducible representations of a group from the outset. A morphism of marked 2-Hilbert spaces is just an ordinary morphism of 2-Hilbert spaces which has no regard for the distinguished simple objects, and similarly for the 2-morphisms; we write the 2-category of marked 2-Hilbert spaces as $\THilb_\text{m}$. Clearly the 2-functor $\THilb_\text{m} \rightarrow \THilb$ which forgets the marked simple objects is an equivalence of 2-categories.

\begin{lem} \label{charstlem} \begin{enumerate}
 \item A linear $*$-functor $F \colon H \rightarrow H'$ between marked 2-Hilbert spaces is determined up to natural isomorphism by the nonnegative integers \\ $\dim \Hom(e_\mu, Fe_i)$ where $e_i$ and $e_\mu$ run over the marked simple objects in $H$ and $H'$ respectively.
 \item Conversely, given a matrix of nonnegative integers $F_{\mu, i}$ there exists a unique (up to natural isomorphism) linear $*$-functor $F \colon H \rightarrow H'$ with $\dim \Hom(e_\mu, Fe_i) = F_{\mu, i}$.
 \end{enumerate}
\end{lem}
\begin{proof}
(i) Suppose $G \colon H \rightarrow H'$ is another linear $*$-functor with $\dim \Hom(e_\mu, Ge_i)=\dim \Hom(e_\mu, Fe_i)$. The fact that 2-Hilbert spaces are semisimple means that we have a canonical isomorphism
 \begin{align*}
  \Hom(Fe_i, Ge_i) & \stackrel{\cong}{\rightarrow} \bigoplus_\mu \Hom_{\Hilb} (\Hom(e_\mu, Fe_i), \Hom(e_\mu, Ge_i)) \\
  f & \mapsto \post(f)
 \end{align*}
where $\post(f)$ refers to the operation of postcomposing with $f$.  Hence we see that for each $i$ there must exist an isomorphism $\gamma_{e_i} \colon Fe_i \rightarrow Ge_i$; by Lemma \ref{unique} this extends (uniquely) to a natural isomorphism $\gamma \colon F \Rightarrow G$.

\noindent (ii) We must define the behaviour of $F$ on objects and on morphisms. For each marked simple object $e_i$, choose $F(e_i) \in H'$ arbitrarily but satisfying
 \[
  F(e_i) = \bigoplus_\mu F_{\mu,i} e_\mu
  \]
where as always the $\bigoplus$ symbol is interpreted in the sense of `{\em a} direct sum' and not `{\em the} direct sum' which wouldn't make any sense. Then for every other object $x \in H$, choose $F(x) \in H'$ arbitrarily but satisfying
\[
  F(x) = \bigoplus_i \dim \Hom(e_i, x) F(e_i).
  \]
The task is now to define $F$ on morphisms and to ensure that the result is a linear $*$-functor. Begin by choosing, for each $x \in H$, a $*$-basis
 \[
  \{ v^{(x)}_{i, p} \colon e_i \rightarrow x \}_{p = 1}^{\dim \Hom(e_i, x)}.
 \]
The value of $F$ on morphisms is then determined by what it does to these basis vectors
 \[
  F \colon \Hom(e_i, x) \rightarrow \Hom(Fe_i, Fx)
 \]
because for $f \colon x \rightarrow y$ we must then set
 \be \label{formF}
  F(f) = \sum_{i,p,q} \langle v^{(y)*}_{i,p} f v^{(x)}_{i,q} \rangle \; F(v^{(y)}_{i,p}) F(v^{(x)}_{i,q})^*.
 \ee
To proceed, make a choice of pairwise orthogonal isometric embeddings of Hilbert spaces
 \[
  A_{i, \mu, p} \colon \Hom(e_\mu, Fe_i) \hookrightarrow \Hom(e_\mu, Fx), \quad p = 1 \ldots \dim \Hom(e_i, x)
 \]
which decompose the Hilbert space $\Hom(e_\mu, Fx)$ into $\dim \Hom(e_i, x)$ orthogonal subspaces of dimension $\dim \Hom(e_\mu, Fe_i)$. Each isometry $A_{i, \mu, p}$ is represented by a morphism $a_{i, \mu, p} \colon Fe_i \rightarrow Fx$ defined uniquely by the requirement that $A_{i, \mu, p} = \post(a_{i, \mu, p})$. Note that $a_{i, \mu, p}^* a_{i, \mu, q} = \delta_{pq} \id_{Fe_i}$ because
 \begin{align*}
  \post(a_{i, \mu, p}^* a_{i, \mu, q}) &= \post(a_{i, \mu, p}^* ) \post( a_{i, \mu, q}) \\
   &= \post(a_{i, \mu, p})^* \post( a_{i, \mu, q})\\
   &= A_{i, \mu, p}^* A_{i, \mu, q} \\
   &= \delta_{pq} \id.
  \end{align*}
Now define the value of $F$ on the basis vectors as
 \[
  F(v^{(x)}_{i,p}) = \sum_\mu a_{i, \mu, p}.
 \]
Equation \eqref{formF} then defines the value of $F$ on all morphisms. Note that we have $F(f^*) = F(f)^*$ by construction, and $F$ is clearly linear. The reader will also easily verify that $F(g \circ f) = F(g) \circ F(f)$ and $F(\id) = \id$, so that we have indeed created a linear $*$-functor $F$ with $\dim \Hom(e_\mu, Fe_i) = F_{\mu, i}$. The fact that $F$ is independent of the choices we made up to natural isomorphism follows from part (i).
\end{proof}
\subsection{The Hilbert space of natural transformations}
In this subsection we show that the vector space $\Nat(F, G)$ of natural transformations between linear $*$-functors $F, G \colon H \rightarrow H'$ is naturally equipped with an inner product. We will need this inner product in Chapters \ref{2RepsChap} and \ref{GerbesCharChap} when we investigate 2-characters of 2-representations. Firstly, we record the basic fact that a natural transformation is freely determined by its components on the simple objects.
\begin{lem} \label{unique} Suppose $F, G : H \rightarrow H'$ are morphisms between 2-Hilbert spaces and that $H$ is marked, with
distinguished simple objects $e_i$. Then given any arbitrary collection of morphisms
 \[
  \theta_i : F(e_i) \rightarrow G(e_i)
 \]
in $H'$, there exists a unique natural transformation $\theta : F \Rightarrow G$ such that $\theta_{e_i} =
\theta_i$.
\end{lem}
\begin{proof} We show uniqueness first. Suppose that $\theta \colon F \Rightarrow G$ is a natural transformation. Given $x \in H$, choose a $*$-basis $a_{i,p} \in \Hom(e_i, x)$ for each simple object $e_i$. We have the following diagram:
  \[
   \xymatrix{ F(e_i) \ar[r]^{F(a_{i,p})} \ar[d]_{\theta_i} & F(x) \ar[d]^{\theta_x} \\ G(e_i) \ar[r]_{G(a_{i,p})} &
   G(x)}
  \]
Thus \be \label{natform}
\begin{aligned}
 \theta_x F(a_{i,p}) = G(a_{i,p}) \theta_i \quad \Rightarrow \quad \theta_x F(a_{i,p}) F(a_{i,p}^*) & = G(a_{i,p}) \theta_i
 F(a_{i,p}^*) \\
  \quad \Rightarrow \quad \sum_{i,p} \theta_x F(a_{i,p}) F(a_{i,p}^*) & = \sum_{i,p} G(a_{i,p}) \theta_i
 F(a_{i,p}^*) \\
  \quad \Rightarrow \quad \theta_x & = \sum_{i,p} G(a_{i,p}) \theta_i
 F(a_{i,p}^*).
\end{aligned}
\ee
 This establishes uniqueness. To establish existence, we must check naturality. Let $f : x\rightarrow y$ be a
morphism in $H$, and choose a $*$-basis $b_{i,p} \in \Hom(e_i, y)$ for each $i$. Then since we can expand $f$
in terms of $a_{i,p}$ and $b_{i,q}$, to check naturality we need only check it for the basis morphisms $b_{i,q}
a_{i,p}^* : x \rightarrow y$. And indeed,
 \[
  \begin{aligned}
   G(b_{i,q} a_{i,p}^*) \theta_x &= \theta_y F(b_{i,q} a_{i,p}^*) \\
   \Longleftrightarrow  \contraction{\sum_{j,r} G(b_{i,q}}{a_{i,p}^*}{)G(}{a_{j,r}}
  \sum_{j,r} G(b_{i,q} a_{i,p}^*) G(a_{j,r}) \theta_i F(a_{j,r}^*)
   &=  \contraction{\sum_{j,r} G(b_{j,r}) \theta_i F(}{b_{j,r}^*}{)F(}{b_{i,q}}
\sum_{j,r} G(b_{j,r}) \theta_i F(b_{j,r}^*) F(b_{i,q} a_{i,p}^*) \\
   \Longleftrightarrow G(b_{i,q}) \theta_i F(a_{i,p}^*) &= G(b_{i,q}) \theta_i F(a_{i,p}^*).
  \end{aligned}
\]
\noindent Moreover, it is clear from the explicit formula \eqref{natform} that $\theta_{e_i} = \theta_i$.
   \end{proof}
With this lemma in  mind, a natural inner product on the space $\Nat(F, G)$ of natural transformations is
\be \label{theinnpro}
      \langle \theta,  \theta' \rangle = \sum_i k_i (\theta_{e_i}, \theta'_{e_i})
 \ee
where $e_i$ runs over representative simple objects for $H$, and $k_i$ are their scale factors. This formula should be compared with \cite[Eqn 3.16]{ref:freed2}.

Note that the Hilbert space $\Nat(\id_H, \id_H)$ of natural transformations of the identity functor on a 2-Hilbert space $H$ is isomorphic to the complexified Grothendieck group $[H]_\mathbb{C}$, since it is spanned by those transformations $\id_i$ whose components are supported exclusively on the class of the simple object $e_i$ where they equal the identity map. However the map which sends $\id_i \mapsto [e_i]$ {\em is not a unitary isomorphism} since
 \[
  ([e_i], [e_i]) = 1 \text{ in } [H]_\mathbb{C} \quad \text{but} \quad (\id_i, \id_i) = k_i^2 \text{ in } \Nat(\id_H, \id_H).
 \]
This is an important caveat to keep in mind.

\subsection{Unitary equivalence for 2-Hilbert spaces}
A morphism $f \colon x \rightarrow y$ in a 2-Hilbert space is called {\em unitary} if $f^*f = \id_x$ and $ff^* = \id_y$, and a natural transformation $\theta \colon F \Rightarrow G$ in $\THilb$ is called {\em unitary} if all its components are unitary. We will follow the terminology of Baez \cite{ref:baez_2_hilbert_spaces} and say that a pair of 2-Hilbert spaces $H, H'$ are {\em equivalent} if they are equivalent in the 2-category $\THilb$, that is if there exist linear $*$-functors $F \colon H \rightarrow H'$ and $G \colon H' \rightarrow H$ together with natural isomorphisms $\eta \colon \id \Rightarrow F^*F$ and $\epsilon \colon FF^*$. The 2-Hilbert spaces $H$ and $H'$ are called {\em unitarily equivalent} if $\eta$ and $\epsilon$ can be chosen to be unitary.

One might ask for a stronger condition, not mentioned in \cite{ref:baez_2_hilbert_spaces}. Let us say that a linear $*$-functor $F \colon H \rightarrow H'$ is {\em unitary} if is a unitary linear map at the level of the hom-sets
 \[
  F \colon \Hom(x,y) \rightarrow \Hom(Fx, Fy).
 \]
We shall say that the 2-Hilbert spaces $H$ and $H'$ above are {\em strongly unitarily equivalent} if $F$ and $G$ are unitary functors. The following proposition classifies 2-Hilbert spaces up to these various notions of equivalence.

\begin{prop} \begin{enumerate} \item {\rm (Baez \cite[Prop 27, Cor 28]{ref:baez_2_hilbert_spaces})}. A pair of 2-Hilbert spaces $H$ and $H'$ are equivalent if and only if they are unitarily equivalent if and only if they have the same dimension.
\item A pair of 2-Hilbert spaces $H$ and $H'$ are strongly unitarily equivalent if and only if there is a bijective map between their isomorphism classes of simple objects which preserves the scale factors.
\end{enumerate}
\end{prop}
\begin{proof} (i) Since 2-Hilbert spaces are semisimple, $H$ and $H'$ are equivalent if and only if they have the same number of isomorphism classes of simple objects --- that is, if they have the same dimension. Moreover, the unit and counit maps establishing the equivalence can always be made unitary, because the categories are semisimple, so a natural isomorphism can be seen as a collection of linear isomorphisms between Hilbert spaces, and these can always be made unitary.

(ii) A strong unitary equivalence $F \colon H \rightarrow H'$ must send simple objects in $H$ to simple objects in $H'$, and $F$ will be unitary if and only if it is unitary at the level of the simple objects. But $(\id_{F(e_i)}, \id_{F(e_i)}) = (F(\id_{e_i}), F(\id_{e_i})$ so the statement follows.
\end{proof}

\subsection{Frobenius algebras from 2-Hilbert spaces\label{ppoint}}
In the introduction to this thesis, we pointed out that the data needed to define a unitary two-dimensional TQFT \cite{ref:durhuus_jonsson} --- namely, a list of positive real numbers --- is precisely the same data which characterizes a 2-Hilbert space, which serves to confirm the `primacy of the point' from an extended TQFT point of view. The positive real numbers represent the eigenvalues of the hermitian handle creation operator from the viewpoint of TQFT, or the scale factors on the simple objects from the viewpoint of a 2-Hilbert space.

On the other hand, it is well known that a two-dimensional TQFT corresponds to a commutative Frobenius algebra. A {\em unitary} two-dimensional TQFT corresponds to a `unitary commutative $\dagger$-Frobenius algebra' in the sense of Vicary \cite{ref:vicary}. We can reformulate this notion for our purposes here as follows.

\begin{defn} A {\em unitary} commutative Frobenius algebra is a commutative $H^*$-algebra. Regarded as a Frobenius algebra, the Frobenius form is given by $\epsilon(a) = (1, a)$.
\end{defn}
By the classification of $H^*$-algebras from Section \ref{hcats}, we see that a unitary commutative Frobenius algebra consists of a bunch of orthogonal idempotents $1_i$ having positive real scale factors $\epsilon(1_i) = \lambda_i$.

We now show how one can obtain every unitary commutative Frobenius algebra by the operation of taking the `higher-categorical dimension' of some 2-Hilbert space $H$, defined as the algebra of natural transformations of the identity functor:
 \[
  \Dim H = \Nat(\id_H, \id_H).
 \]
This algebra has a natural inner product which we defined in \eqref{theinnpro}. We saw that it has a natural basis $\id_i$ consisting of transformations whose components are supported exclusively on the class of the simple object $e_i$ where they equal the identity map, and by definition the Frobenius form computes as
 \[
  \epsilon(\id_i) = k_i^2.
 \]
In other words, we have established the following result, which shows that in principle every `ordinary' two-dimensional unitary TQFT $Z \colon \text{2Cob} \rightarrow \Hilb$ can be extended to an extended TQFT $Z \colon \text{2}\mathcal{C}\text{ob} \rightarrow \THilb$.

 \begin{prop} Every unitary commutative Frobenius algebra is isomorphic to $\Dim H$ for some 2-Hilbert space $H$.
 \end{prop}
For instance, in the 2-Hilbert space $\Rep(G)$ of unitary representations of a finite group, we find that $\Dim \Rep(G)$ is the commutative Frobenius algebra spanned by the idempotents $\id_i$ where $i$ runs over the irreducible representations, and whose Frobenius form computes as
 \[
  \epsilon(\id_i) = \frac{|\dim V_i|^2}{|G|^2}.
 \]
This is precisely the Frobenius algebra arising in the two-dimensional finite group model, as in the lecture notes of Segal \cite{ref:segal2}.

\section{The geometry of 2-Hilbert spaces\label{AGeomSec}} In Chapter \ref{GeometryOrdinaryRepChap}, we showed how the category Hilb of finite-dimensional Hilbert spaces and linear maps is equivalent to a geometric category $\LBun$ whose objects are hermitian line bundles over compact hermitian manifolds and whose morphisms are kernels. In this section we show how this picture `categorifies' for 2-Hilbert spaces.

\subsection*{The 2-category of finite spaces}

\begin{figure}[t]
\centering
\ig{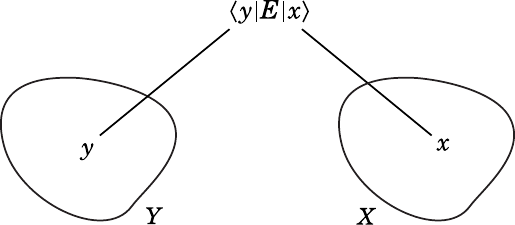}
\caption{\label{buggerfig} A morphism between finite spaces is a Hilbert bundle over their product. }
\end{figure}

\begin{defn} The 2-category $\FinSpaces$ is defined as follows. An object is a finite set $X$ equipped with a positive real scale factor $k_x$ for each $x \in X$. A morphism $E \colon X \rightarrow Y$ is a bundle of Hilbert spaces over $Y \times X$, with composition given by convolution. A 2-morphism is a map of vector bundles.
\end{defn}
Let us explain this definition. One should think of a finite set $X$ equipped with positive real scale factors $k_x$ for each $x \in X$ as a compact discrete smooth manifold equipped with a `metric'.  A morphism $E \colon X \rightarrow Y$ is to be thought of as the categorification of the idea of a {\em kernel} from Chapter \ref{GeometryOrdinaryRepChap} --- instead of the amplitude for going from $x \in X$ to $y \in Y$ being a mere {\em number} (more precisely, an element in a complex line), the amplitude is now a {\em Hilbert space}, which we will write as
  \[
  \langle y | E | x \rangle.
 \]
See Figure \ref{buggerfig}.
We compose morphisms $E \colon X \rightarrow Y$ and $F \colon Y \rightarrow Z$ by setting $F \circ E$ to be the vector bundle over $Z \times X$ whose space of ways of going from $x \in X$ to $z \in Z$ is given by the weighted direct sum over all $y \in Y$ of the tensor product amplitudes:
 \be \label{mycompform}
  \langle z | F \circ E | x \rangle := \hat{\displaystyle \bigoplus_{y \in Y}} \, \frac{1}{k_y} \; \langle z | F | y \rangle \otimes \langle y | E | x \rangle.
 \ee
By this notation, we mean that the inner product is given on the homogenous components by
 \[
  (v_1 \otimes w_1, v_2 \otimes w_2) = \frac{1}{k_y} (v_1, v_2)(w_1, w_2)
 \]
where $k_y$ is the scale factor on $y \in Y$. The formula \eqref{mycompform} for composition of morphisms in $\FinSpaces$ should be compared with the formula \eqref{formkernel} for composition of kernels in $\LBun$ from Chapter \ref{GeometryOrdinaryRepChap}:
  \[
    \langle z | F \circ E | x \rangle  = \int_Y \vol_y \; \langle z | F | y \rangle \circ \langle y | E | x \rangle.
  \]
The hat on the direct sum is there to indicate that $\langle z | F \circ E | x \rangle$ is not defined as the {\em formal} direct sum over $y \in Y$ of $\langle z | F | y \rangle \otimes \langle y | E | x \rangle$, which would require an ordering on $Y$, but is rather defined geometrically as the `space of sections of the fibrewise tensor product bundle $F \otimes E$ over $Y \subset Z \times Y \times Y \times X$'. That is to say, we are defining composition via a geometric pull-push formula of vector bundles identical to the formula for composition of kernels in the derived category context (see \cite{ref:ganter_kapranov_rep_char_theory, ref:caldararu_willerton}).

The identity morphisms $\id \colon X \rightarrow X$ are the vector bundles over $X \times X$ which have weighted copies of $\mathbb{C}$ along the diagonal:
 \[
  \langle y | \id | x \rangle := \delta_{x,y} k_x \mathbb{C}.
 \]
The inner product is chosen to ensure that the natural choice of coherence isomorphisms
 \[
  E \circ \id \rightarrow E,
  \]
namely the ones which are given fibrewise by the map
 \begin{align*}
  \langle y | E | x \rangle \otimes \langle x | \id | x \rangle & \rightarrow \langle y | E | x \rangle \\
  v \otimes 1 & \mapsto v
 \end{align*}
are unitary. We have now described composition of 1-morphisms; horizontal and vertical composition of 2-morphisms works in the obvious way. The resultant structure $\FinSpaces$ forms a 2-category, though not a strict one. To see this, observe that if $K \colon Z \rightarrow A$ is another morphism of finite spaces, then although the vector bundle $(K \circ F) \circ E$ is not truly {\em equal} to $K \circ (F \circ E)$, there is a canonical isomorphism between them, which is fibrewise just the ordinary associator of vector spaces:
 \[
  (\langle a | K | z \rangle \otimes \langle z | F | y \rangle) \otimes \langle y | E | x \rangle \cong \langle a | K | z \rangle \otimes (\langle z | F | y \rangle \otimes \langle y | E | x \rangle).
 \]
This is one advantage of insisting on the geometric definition of convolution of vector bundles and not using {\em formal} direct sums --- the orderings necessitated by the latter would clutter up this argument. In any event, since the ordinary associator on vector spaces satisfies the appropriate coherence diagrams, we see that $\FinSpaces$ is indeed a 2-category.

\subsection*{The inner product on the 2-morphisms}
The 2-morphisms in $\FinSpaces$ can be equipped with a natural inner product. Suppose that $E, F \colon X \rightarrow Y$ are 1-morphisms in $\FinSpaces$, and that $f, g \colon E \Rightarrow F$ are 2-morphisms between them. In other words, $f$ and $g$ are a collection of fibrewise linear maps of Hilbert spaces
 \[
  f_{y, x} \colon \langle y | E | x \rangle \rightarrow \langle y | F | x \rangle.
 \]
The we define their inner product by
 \be \label{myip}
  \langle f, g \rangle := \sum_{x \in X, y \in Y} k_x k_y \Tr(f_{x,y}^* g).
 \ee

\subsection*{An equivalence of 2-categories}
We now show that the 2-category of 2-Hilbert spaces is equivalent to the 2-category of finite spaces.

We define a 2-functor
 \[
 \langle \cdot \rangle \colon \THilb_\text{m} \rightarrow \FinSpaces
 \]
as follows. To a marked 2-Hilbert space $H$ it assigns the set of marked simple objects, which we will write as $\langle H \rangle$. A linear $*$-functor $F \colon H \rightarrow H'$ gives rise to a bundle $\langle F \rangle$ of Hilbert spaces over $\langle H' \rangle \times \langle H \rangle$, defined fibrewise by
 \[
  \langle \mu | F | i \rangle := \Hom(e_\mu, Fe_i).
 \]
Finally, a natural transformation $\theta \colon F \Rightarrow G$ gives rise to a map of vector bundles $\langle \theta \rangle \colon \langle F \rangle \rightarrow \langle G \rangle$, defined by postcomposing with the components of $\theta$ on the marked simple objects $e_i \in H$:
 \[
  \langle \theta \rangle_{\mu, i} := \post(\theta_{e_i}) \colon \langle \mu | F | i \rangle \rightarrow \langle \mu | G | i \rangle.
  \]
\begin{prop}\label{THilbprop} The map $\langle \cdot \rangle \colon \THilb_\text{m} \rightarrow \FinSpaces$ is a weak 2-functor, with unitary coherence 2-cells, and is an equivalence of 2-categories. Moreover it is a unitary linear map at the level of the 2-morphisms.
\end{prop}
\begin{proof} We need to establish the coherence data. Suppose that $F \colon H \rightarrow H'$ and $G \colon H' \rightarrow H''$ are morphisms of marked 2-Hilbert spaces. We define the compositor
 \[
  \langle G \rangle \circ \langle F \rangle \rightarrow \langle G \circ F \rangle
  \]
fibrewise as
 \begin{align} \label{mycompositor}
    \hat{\mathop{\bigoplus}_{\mu \in \langle H' \rangle}} \frac{1}{k_\mu} \langle \lambda | G | \mu \rangle
    \otimes \langle \mu |F | i \rangle & \rightarrow \langle \lambda | G F | i \rangle  \\
  f \otimes g & \rightarrow G(f) \circ g
  \end{align}
where $f \in \langle \mu | F | i\rangle$ and $g \in \langle \lambda | G | \mu \rangle$. This map is an isomorphism, since if we choose a $*$-basis $\{ a^{(\mu, i)}_p \colon e_\mu \rightarrow \sigma(e_i)\}$ for each space $\langle \mu | \sigma | i \rangle$ then we can describe its inverse as the map which sends
 \be \label{expinverse}
  \langle \lambda | G F | i \rangle \ni h \mapsto \sum_{\mu, p}
  G(a_p^{(\mu, i)})^* \circ h \,
  \otimes \, a^{(\mu, i)}_p.
  \ee
Moreover the compositor is unitary, because if we choose a $*$-basis $\{b^{(\lambda, \mu)}_r \colon e_\lambda \rightarrow G(e_\mu)\}$ for $\langle \lambda | G | \mu \rangle$ then
 \begin{align*}
  ( b^{(\lambda, \mu)}_r \otimes b^{(\lambda, \mu)}_s, \, a^{(\mu, i)}_p \otimes a^{(\mu, i)}_q) &= \frac{1}{k_\mu} (b^{(\lambda, \mu)}_r, b^{(\lambda, \mu)}_s) (a^{(\mu, i)}_p, a^{(\mu, i)}_q) \\
  &= k_\lambda \delta_{rs} \delta_{pq} \\
  &= (b^{(\lambda, \mu)}_r \otimes b^{(\lambda, \mu) *}_s, \, G(a^{(\mu, i)}_p)^* G(a^{(\mu, i)}_q) ) \\
   &= ( G(a^{(\mu, i)}_p) \circ b^{(\lambda, \mu)}_r, \, G(a^{(\mu, i)}_q) \circ b^{(\lambda, \mu)}_s).
  \end{align*}
Similarly we define the unit coherence isomorphism
 \begin{align*}
  \id_{\langle H \rangle} & \rightarrow \langle \id_H \rangle \\
  \intertext{fibrewise by}
  \mathbb{C} & \rightarrow \langle i | \id | i \rangle \\
    1 & \mapsto \id_{e_i}.
  \end{align*}
This map is unitary, by the definition of the inner product on the fibers of the vector bundle $\id_{\langle H \rangle}$.

Taken together, this coherence data exhibits $\langle \cdot \rangle \colon \THilb_m \rightarrow \FinSpaces$ as a weak 2-functor. This 2-functor is clearly essentially surjective, because from any finite set of points $X$ one can define a corresponding 2-Hilbert space $\Hilb^|X|$. We need only check that it is a local equivalence, that is that at the level of hom-categories the functor
 \[
 \Hom_{\THilb_\text{m}} (H, H') \rightarrow \Hom_{\FinSpaces} (\langle H \rangle, \langle H' \rangle)
  \]
is an equivalence of categories. Lemma \ref{unique} establishes that this functor is fully faithful, while Lemma \ref{charstlem} establishes that it is essentially surjective.

Finally, we need to show that the 2-functor is a unitary linear map at the level of 2-morphisms. Suppose $F, G \colon H \rightarrow H'$ are linear $*$-functors, and $\theta, \phi \colon F \Rightarrow G$ are natural transformations between them. Then:
 \begin{align*}
   (\langle \theta \rangle, \, \langle \phi \rangle) & := \sum_{i \in \langle H \rangle, \mu \in \langle H' \rangle} k_i k_\mu \Tr (\langle \theta \rangle_{\mu, i}^* \circ \langle \phi \rangle_{\mu, i}) \\
   &= \sum_{i, \mu} k_i k_\mu \Tr (\post(\theta_{e_i})^* \post(\phi_{e_i})) \\
   &= \sum_{i, \mu} k_i k_\mu \Tr (\post(\theta_{e_i}^* \phi_{e_i})) \\
   &= \sum_i k_i (\id_{F(e_i)}, \, \theta_{e_i}^* \phi_{e_i})\\
   &= \sum_i k_i (\theta_{e_i}, \, \phi_{e_i}) \\
   &= (\theta, \phi).
 \end{align*}
\end{proof}
In Chapter \ref{GerbesCharChap} we will prove an equivariant version of this result, namely that the 2-category of unitary 2-representations of a group $G$ on marked 2-Hilbert spaces is equivalent to the 2-category of {\em equivariant gerbes}.

\chapter{String diagrams\label{stchap}} \XXX{rewrite this chapter!!}
In this chapter we recall for the reader's benefit the string diagram notation for 2-categories. This notation is particularly suited to describe structures such as adjunctions and monads. The use of string diagrams as a convenient notational device has its roots in the physics literature where they appear as {\em Feynman diagrams} or the {\em diagrammatic tensor calculus} of Penrose (see \cite{ref:baez_stay_rosetta} and the references therein). In time mathematicians realized that this notation works in the general context of any monoidal category or 2-category, and soon applications were found in knot theory and quantum algebra. We highlight here especially the work of Joyal and Street \cite{ref:joyal_street, ref:street_string_diagrams}, as well as Freyd and Yetter \cite{ref:freyd_yetter_2} and Turaev \cite{ref:turaev2}. Because we will use it over and over again in this thesis, and because there are some misconceptions about the notation we wish to clear up, we have given this section a chapter of its own. \XXX{Don't like this last sentence}

\begin{figure}
 \centering
 \ig{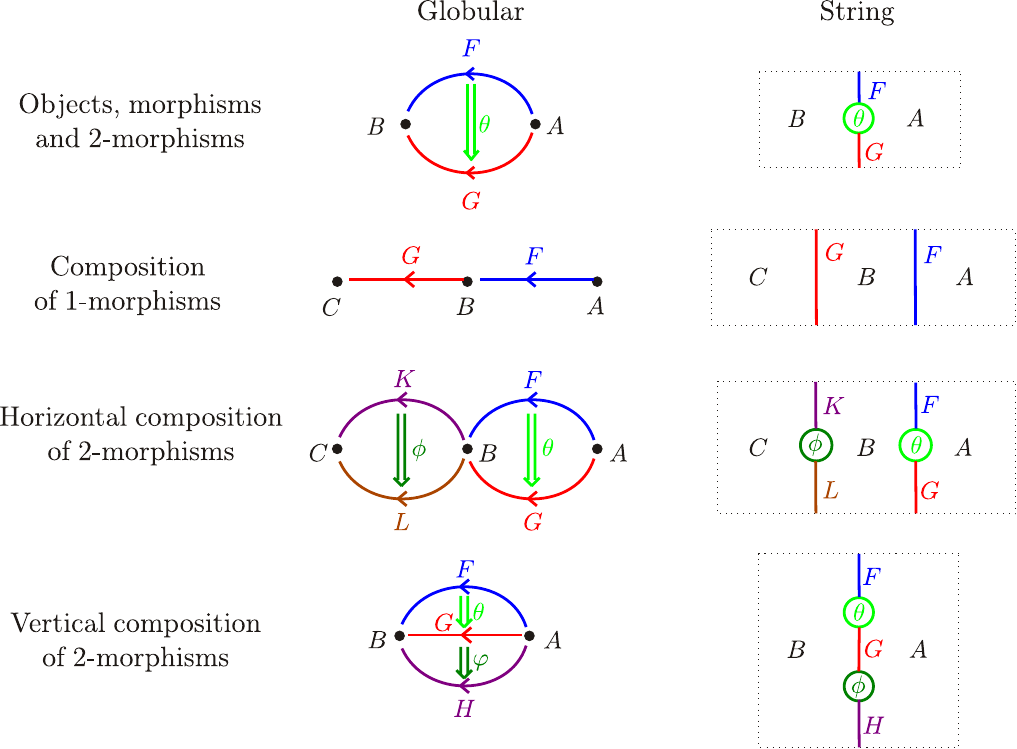}
 \caption{\label{fig:sd} Globular notation for 2-categories versus string diagram notation. }
 \end{figure}

String diagrams are a two-dimensional graphical notation for working with 2-categories, and may be regarded as the `Poincar\'{e} duals' of the ordinary globular notation. The basic idea is summarized in Figure \ref{fig:sd}. The reader who does not understand this figure and who is still confused by these diagrams is referred to Section 2.2 of \cite{ref:lauda} or Section 1.1 of \cite{ref:caldararu_willerton} for more details. In our diagrams, composition of 1-morphisms runs from right to left (so a $G$ to the left of an $F$ means $G$ after $F$), and composition of 2-morphisms runs from top to bottom. We stress that string diagrams are not merely a mnemonic but are a perfectly rigorous notation.

We remind the reader that by a `2-category' we mean a {\em not-necessarily-strict} 2-category (also called a {\em bicategory}). We wish to strongly emphasize the following point, which is often misunderstood in the literature.

\begin{quote}
 String diagrams are not just a notation for {\em strict} 2-categories. When interpreted correctly, they are a perfectly rigorous notation for {\em weak} 2-categories.
\end{quote}

\begin{prop} After specifying a parenthesis scheme for the source and target 1-morphisms of a string diagram, the resultant 2-morphism represented by the diagram is independent of the choice of parentheses, associators and unit 2-isomorphisms used to interpret the interior of the diagram.
\end{prop}
\begin{proof}Consider a typical string diagram:
 \[
   \ig{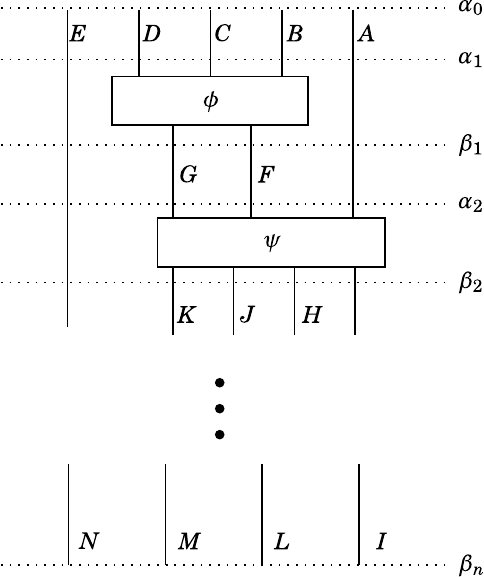}
 \]
Divide the diagram up using a timeline as shown in such a way that each 2-morphism occurs in its own time segment with each such segment having an `initial' state, which occurs just before the step of vertically composing with that 2-morphism, and a `final' state, which occurs just after. The source and target 1-morphisms of the diagram are written as $\alpha_0$ and $\beta_n$ respectively. Their parenthesis scheme is required to be given; suppose for instance that they are
 \[
 \alpha_0 = E \circ [(D \circ \id) \circ [(C \circ B) \circ A]] \quad \text{and} \quad \beta_n = [(N \circ M) \circ L] \circ [I \circ \id].
 \]
Another aspect of the problem which is understood to be given beforehand is the nature of the 2-morphisms which make up the diagram --- their source and target 1-morphisms, as well as the parentheses schemes on these 1-morphisms. For instance:
 \[
 \phi \colon (D \circ C)\circ B \Rightarrow G \circ F \quad \text{and} \quad \psi \colon (F \circ \id) \circ A \Rightarrow J \circ (H \circ K).
 \]
Now, suppose we have two different interpretations of this string diagram:
\[
 \ba \xymatrix{E \circ [(D \circ \id) \circ [(C \circ B) \circ A]] \ar@{=>}[d]^{a_1} \\ [E \circ ((D \circ C)\circ B)] \circ A \ar@{=>}[d]^{[\id_E \ast \phi] \ast \id_A} \\ [E \circ (G \circ F)] \circ A \ar@{=>}[d]^{a_2} \\ (E \circ G) \circ ((F \circ \id) \circ A) \ar@{=>}[d]^{\id_{E \circ G} \ast \psi} \\ \vdots}\ea \;\; \text{and} \;\;
 \ba \xymatrix{E \circ [(D \circ \id) \circ [(C \circ B) \circ A]] \ar@{=>}[d]^{a'_1} \\ E \circ [((D \circ C)\circ B) \circ A] \ar@{=>}[d]^{\id_E \ast [\phi \ast \id_A]} \\ E \circ [(G \circ F) \circ A] \ar@{=>}[d]^{a'_2} \\ E \circ [G \circ ((F \circ \id) \circ A)] \ar@{=>}[d]^{\id_E \ast [\id_G \ast \psi]} \\ \vdots} \ea
  \]
Here $a_1$ refers to the chosen sequence of associators and unit isomorphisms which transform the initial string at $\alpha_0$ into the string at $\alpha_1$; the same goes for $a'_1$, $a_2$, $a'_2$ etc. Using coherence for 2-categories (in the form which says `all diagrams of constrains commute' as in the thesis of Gurski \cite{ref:gurski}) and naturality of the associator, we can transform the initial segment of the second interpretation of the diagram into the initial segment of the first interpretation as follows:
 \[
\!\!\!\!\!\!\!\!\!\! \ba \xymatrix{E \circ [(D \circ \id) \circ [(C \circ B) \circ A]] \ar@{=>}[d]^{a'_1} \\ E \circ [((D \circ C)\circ B) \circ A] \ar@{=>}[d]^{\id_E \ast [\phi \ast \id_A]} \\ E \circ [(G \circ F) \circ A] \ar@{=>}[d]^{\cdots} \\ {}} \ea  =   \ba \xymatrix{E \circ [(D \circ \id) \circ [(C \circ B) \circ A]] \ar@{=>}[d]^{a_1} \\ [E \circ ((D \circ C)\circ B)] \circ A \ar@{=>}[d]^{a_1' a_1^\mi} \\ E \circ [((D \circ C)\circ B) \circ A] \ar@{=>}[d]^{\id_E \ast [\phi \ast \id_A]} \\ E \circ [(G \circ F) \circ A] \ar@{=>}[d]^{\cdots} \\ {}} \ea  =
 \ba \xymatrix{E \circ [(D \circ \id) \circ [(C \circ B) \circ A]] \ar@{=>}[d]^{a_1} \\ [E \circ ((D \circ C)\circ B)] \circ A \ar@{=>}[d]^{[\id_E \ast \phi] \ast \id_A} \\ [E \circ (G \circ F)] \circ A \ar@{=>}[d]^{a_1''} \\  E \circ [(G \circ F) \circ A] \ar@{=>}[d]^{\cdots} \\ {}} \ea
\]
We can repeat this procedure all the way down the diagram; at the end we will have transformed the second interpretation into the first one.
\end{proof}
Now, whenever a string diagrams occurs in this paper it will always be clear from the context precisely what the parenthesis scheme on the input and output 1-morphisms is intended to be. For instance:
\[
\text{``The 2-morphism $\eta\colon \id \Rightarrow (G \circ F) \circ (F^* \circ G^*)$ is defined as } \ba \ig{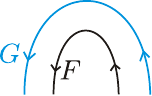} \ea \text{ .''}
   \]
In other words, every string diagram in this paper has a precise and rigorous meaning, even in a fully weak setting. 
\chapter{Even-handed structures}\label{EvenHandedchap}

The main purpose of this chapter is to introduce an alternative paradigm of duality for 1-morphisms in a 2-category which we call an {\em even-handed structure}, and to investigate this notion for a variety of 2-categories of interest. We will need the notion of an even-handed structure in Chapter \ref{2RepsChap} in order to describe the `push-forward' map for 2-characters of unitary 2-representations of finite groups.

The main issue when considering duals for 1-morphisms in a 2-category is that a 1-morphism can have both a {\em right} and a {\em left} dual (we use the words `dual' and `adjoint' interchagably), and these two notions are apriori independent of each other. But in many applications, such as the push-forward map for 2-characters in Chapter \ref{2RepsChap}, and in calculating the `quantum dimensions' of objects in fusion categories as in Chapter \ref{FusionChap}, one needs to go beyond the `snake-style' diagrams and be able to draw {\em closed circles} in the string diagram notation, which requires one to have a rule which partners right dual structure maps with left dual structure maps.

There have been various approaches to addressing this need for coupling right and left duals together; we are proposing yet another. The distinguishing feature of our approach is that we do not require that a permanent fixed choice of dual $F^*$ be made for every morphism $F$ from the outset. We only require that for {\em each} choice of right dual $F^*$, there is an assigned way to portray $F^*$ also as a {\em left} dual (whence the name even-handed structure). The two main advantages of this slightly relaxed approach are the following: it meshes well with the string diagram notation (which cannot be said for a straightforward implementation of a `fixed dual' approach) and also it {\em actually fits the examples we have in mind}. In most 2-categories made up of categories, functors and natural transformations, adjoints are not canonically given but are only specified up to isomorphism. However there is often a canonical {\em function} which {\em turns} right adjoints into left adjoints and vice versa, and getting the theory to match this kind of situation has been our main motivation.

In Section \ref{defesec} we define the notion of an even-handed structure. In Section~\ref{geominterp} we give a geometric interpretation of an even-handed structure as an `even-handed trivialization' of a bundle-like structure over the morphisms in the 2-category which we call the {\em ambijunction gerbe}. In Section~\ref{srel} we compare the notion of an even-handed structure to other conceptions of `duals for morphisms' in monoidal categories and 2-categories, such as the {\em pivotal categories} approach of Joyal and Street \cite{ref:joyal_street_bmc2} and Freyd and Yetter \cite{ref:freyd_yetter}, as well as the {\em 2-categories with duals} approach of Baez and Langford \cite{ref:baez_langford}. In Section \ref{catsection} we specialize the notion of an even-handed structure to the case of a 2-category consisting of categories, functors and natural transformations, where it is more convenient to adopt the hom-set isomorphisms approach to adjunctions as opposed to the abstract approach in terms of unit and counit maps. We use this reformulation in Section \ref{popq} to show how one can obtain an even-handed structure on a collection of linear categories from the data of a {\em trace} on each category, and that in fact when the categories involved are semisimple an even-handed structure is {\em the same thing}, up to a global scale factor, as providing a trace on each category. Finally in Sections \ref{EvenTHilb} and \ref{yausec} we give our two main examples of even-handed structures arising from traces: the even-handed structure on the 2-category $\THilb$ of {\em 2-Hilbert spaces}, where the trace arises from the inner products on the hom-sets, and the even-handed structure on the 2-category $\CYau$ (whose objects are the graded derived categories $D(X)$ of Calabi-Yau manifolds $X$, and whose morphisms are functors naturally isomorphic to integral kernels), where the trace arises from the Serre duality on the hom-sets.

\section{Definition of an even-handed structure\label{defesec}}
Recall that a morphism $F^* \colon B \rightarrow A$ in a 2-category is said to be a {\em right adjoint of} another morphism $F \colon A \rightarrow B$ (equivalently $F$ is a {\em left adjoint} of $F^*$) if there exist unit and counit 2-morphisms
 \begin{align*}
  \eta \colon \id_A \Rightarrow F^* F, \quad &\text{drawn as} \quad \ba \ig{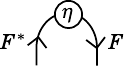} \ea \\
  \epsilon \colon FF^* \Rightarrow \id_B, \quad & \text{ drawn as } \quad\ba \ig{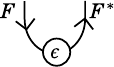} \ea
 \end{align*}
which satisfy the adjunction equations (`snake diagrams'):
 \[
 \ba \ig{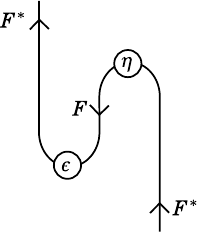} \ea = \ba \ig{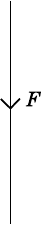} \ea \quad \text{and} \quad \ba \ig{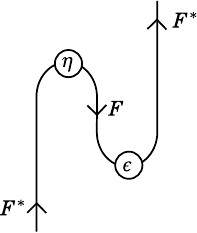} \ea = \ba \ig{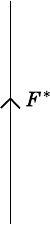} \ea \; .
 \]
The arrows on these diagrams are not really necessary but they are there to provide a visual aid. We will write
 \[
 \Adj(F \dashv F^*)
 \]
for the set\footnote{We are tacitly assuming here that $\Aut(F)$ is a {\em set} and not a proper class. This is the case in most reasonable situations and certainly holds in all our examples.} of all pairs of unit and counit 2-morphisms $(\eta, \epsilon)$  which exhibit $F^*$ as a right adjoint of $F$.
It is well known that the group $\Aut(F)$ of automorphisms of $F$ acts freely and transitively on $\Adj(F \dashv F^*)$ by `twisting' the unit and counit maps:
 \[
  \delta \cdot   \left( \ba \ig{e308.pdf} \ea, \;\; \ba \ig{e309.pdf} \ea \right) = \left(\ba \ig{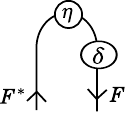} \ea , \;\; \ba \ig{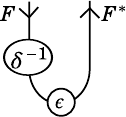} \ea \right).
 \]
We say that a 2-category is {\em has ambidextrous adjoints} if every morphism has an adjoint and if all adjoints are two-sided --- that is, for every morphism $F$ there exists a morphism $F^*$ which is simultaneously a right and left adjoint of $F$. Note that having ambidextrous adjoints is a {\em property} of the 2-category and not a structure, because we are not specifying the choice of these adjoints in any way.

\begin{defn}\label{evsdefn} An {\em even-handed structure} on a 2-category with ambidextrous adjoints consists of, for each pair $(F, F^*)$ of adjoint morphisms, a map
 \begin{align*}
  \Psi \colon \Adj(F \dashv F^*) &\rightarrow \Adj(F^* \dashv F) \\
   \left( \ba \ig{e308.pdf} \ea, \;\; \ba \ig{e309.pdf} \ea \right) & \mapsto \left(\ba \ig{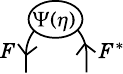} \ea, \;\; \ba \ig{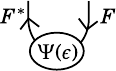} \ea \right)
 \end{align*}
satisfying:
 \begin{enumerate}
  \item $\Psi^2 = \id$.
  \item For all composable morphisms $F$ and $G$,
   \[
    \Psi\left( \ba \ig{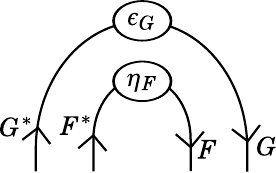} \ea, \ba \ig{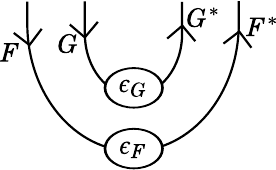} \ea \right) = \left( \ba \ig{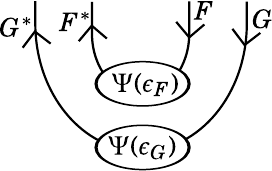} \ea, \ba \ig{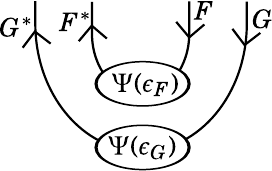} \ea \right).
   \]
   Moreover for all objects $A$, $\Psi(\mathbf{1}_A) = \mathbf{1}_A$, where $\mathbf{1}_A = (\lambda^\mi, \lambda)$ is the trivial adjunction on $\id_A$, with $\lambda \colon \id_A \circ \id_A \Rightarrow \id_A$ the left/right (they are equal in this setting) identity coherence 2-morphism. With this in mind, this would be drawn in string diagrams as $\Psi\left( \ba \ig{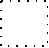} \ea, \ba \ig{e262.pdf} \ea \right) = \left( \ba \ig{e262.pdf} \ea, \ba \ig{e262.pdf} \ea \right)$.
  \item For all 2-morphisms $\theta \colon F \Rightarrow G$ and all choices of right adjoints $(F^*, \eta_F, \epsilon_G)$ and $(G^*, \eta_G, \epsilon_G)$ for $F$ and $G$ respectively, the following 2-morphisms $G^* \Rightarrow F^*$ agree:
   \be \label{adjoints1}
     \theta^\dagger := \ba \ig{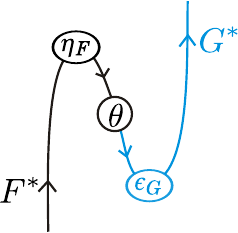} \ea = \ba \ig{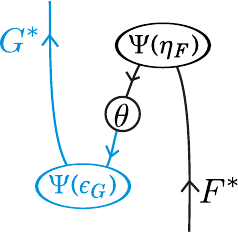} \ea =: {}^\dagger \! \theta_\Psi.
   \ee
  \end{enumerate}
 \end{defn}
We call the 2-morphism $\theta^\dagger$ on the left-hand side of (iii) above the {\em right dagger} of $\theta$, because it is computed using the structure of $F^*$ as a {\em right} adjoint of $F$, and we call the 2-morphism on the right hand-side the {\em left dagger of $\theta$ computed according to $\Psi$}, because it is computed using the structure of $F^*$ as a {\em left} adjoint of $F$. The requirement that $\theta^\dagger = {}^\dagger \! \theta_\Psi$, which we call the {\em even-handed equation}, therefore asks that $\Psi$ turns right adjoints into left adjoints in such a way that the right and left daggers of 2-morphisms always agree, and should be regarded as a naturality condition on $\Psi$.

We remark that we have chosen the notation $\theta^\dagger \colon G^* \Rightarrow F^*$ as opposed to $\theta^* \colon G^* \Rightarrow F^*$ because, with a view to TQFT, we would like to follow Baez and Langford's approach \cite{ref:baez_langford} and reserve the $*$-symbol exclusively for the notion of dual in a higher category which means `a morphism which goes in the opposite direction'. In other words, the notation $\theta^*$ is reserved for a 2-morphism which goes from $G$ to $F$ and not from $G^*$ to $F^*$.

Another remark is that in all of our examples, the requirement $\Psi^2 = \id$ already follows from the other axioms, but we have been unable to show that this holds in general (see Section \ref{pivssec} for further discussion).

In this thesis we give three classes of examples of 2-categories equipped with even-handed structures:
  \begin{itemize}
  \item Braided monoidal categories, for instance representation categories of ribbon Hopf algebras. See the next two subsections.
  \item Monoidal categories equipped with a pivotal structure. See Section \ref{pivssec} and Chapter \ref{FusionChap}.
  \item Collections of linear categories equipped with traces, thought of as a 2-category where the morphisms are linear functors and the 2-morphisms are natural transformations. For instance, the 2-category of 2-Hilbert spaces (Section \ref{EvenTHilb}) and the 2-category of derived categories of Calabi-Yau manifolds (Section \ref{yausec}).
  \end{itemize}
The reader is invited to jump to these sections to get the flavour of the theory.

We end this section by recording the following useful lemma about how an even-handed structure
 \[
  \Psi \colon \Adj(F \dashv F^*) \rightarrow \Adj(F^* \dashv F)
 \]
is natural with respect to changing $F$ and $F^*$ by an invertible 2-morphism.
\begin{lem}\label{narlem} \begin{enumerate}
 \item If $\delta \colon F \Rightarrow F'$ is an invertible 2-morphism, then
  \[
  \Psi \left( \ba \ig{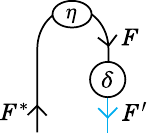} \ea, \ba \ig{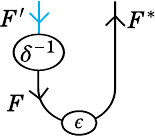} \ea \right) = \left(\ba \ig{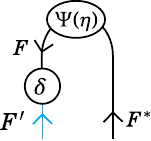} \ea, \ba \ig{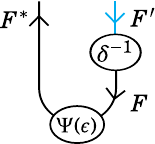} \ea \right).
  \]
 \item If $\gamma \colon F^* \Rightarrow (F^*)'$ is an invertible 2-morphism, then
   \[
    \Psi \left( \ba \ig{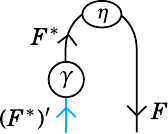} \ea, \ba \ig{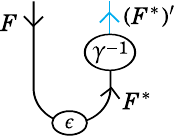} \ea \right) = \left(\ba \ig{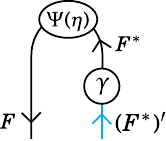} \ea, \ba \ig{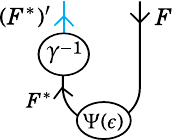} \ea \right).
    \]
    \end{enumerate}
\end{lem}
\begin{proof} We prove (i), the proof of (ii) is similar. Let us write
 \[
  (\eta', \epsilon') \equiv \left(\ba \ig{e384.pdf} \ea, \ba \ig{e385.pdf} \ea \right).
 \]
To prove (i) it is sufficient to check that
 \[
  \ba \ig{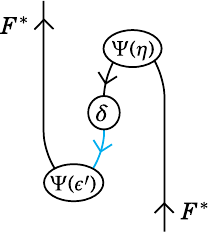} \ea = \ba \ig{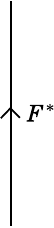} \ea.
 \]
But this follows from the even-handed equation, because the left-hand side is ${}^\dagger \! \delta_\Psi$ while the right hand side is $\delta^\dagger$.
\end{proof}

\subsection*{Even-handed structures on monoidal categories}
A monoidal category can be thought of as a one-object 2-category by reinterpreting the objects as 1-morphisms from the object to itself (with the composition being tensor product) and the morphisms as 2-morphisms; in this picture duals for morphisms gets reinterpreted as duals for {\em objects}. So it makes sense to talk about an `even-handed structure on a monoidal category with ambidextrous duals'. The point about having an even-handed structure is that it allows one to take `quantum traces' of endomorphisms $f \colon V \rightarrow V$, an operation which takes values in $\Hom(1,1)$. In fact one can form the `right trace' or the `left trace' of $f$:
  \[
   \Tr_r (f) =  \ba \ig{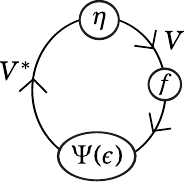} \ea \quad \quad \Tr_l (f) = \ba \ig{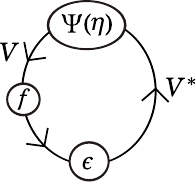} \ea.
  \]
Here $(V^*, \eta, \epsilon)$ is {\em any} choice of right dual for $V$. The trace does not depend on this choice, because if $((V^*)', \eta', \epsilon')$ is another choice of right dual of $V$, then it is well-known that there is a canonical isomorphism $\delta \colon V^* \rightarrow (V^*)'$ transforming $(\eta, \epsilon)$ into $(\eta', \epsilon')$, so that for the right trace we have
 \[
  \ba \ig{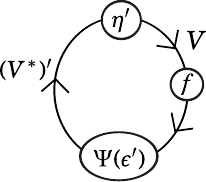} \ea = \ba \ig{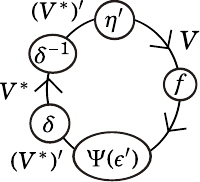} \ea = \ba \ig{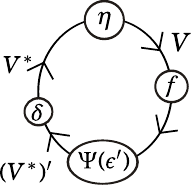} \ea = \ba \ig{e390.pdf} \ea
 \]
where the final equality uses Lemma \ref{narlem} together with the fact that $\Psi^2 = \id$. A similar argument holds for the left trace. We say that the even-handed structure $\Psi$ is {\em spherical} if the right and left traces always coincide. This notion is due to Barrett and Westbury \cite{ref:barrett_westbury}, although more precisely they worked in the framework of {\em pivotal structures} and not even-handed structures; we will compare these two frameworks in Section \ref{pivssec}. The word `spherical' refers to the fact that one may imagine the closed loop involved in the trace map as living on the surface of a sphere, so that one is requiring the trace to be invariant under the isotopy of the sphere which transforms the right trace into the left one.

\subsection*{Even-handed structures on braided monoidal categories}
Our first class of examples of even-handed structures comes from monoidal categories equipped with a braiding. Firstly recall a {\em rigid} monoidal category is one with the property that every object has an (unspecified) right dual, and that a {\em twist} on a braided monoidal category $C$ is a natural isomorphism $\theta$ of the identity which satisfies $\theta_{V \otimes W} = \sigma_{W,V} \circ (\theta_W \otimes \theta_V) \circ \sigma_{V, W}$ for all objects $V,W$. If $C$ is rigid, then we will say that the twist is {\em compatible with duals in $C$} if $\theta_{V^*} = \theta_V^\dagger$ for some choice (and hence for {\em all} choices) of right dual $(V^*, \eta, \epsilon)$ of $V$.

The following fact is well known in quantum algebra in the framework of pivotal structures (see for example Cor 2.2.3 of \cite{ref:bakalov_kirillov}); here we show how it works in the framework of even-handed structures.
\begin{lem}\label{natbugger} Let $C$ be a rigid braided monoidal category. Then there is a canonical bijection
 \[
 \left\{ \ba \begin{array}{c} \text{even-handed structures on $C$}  \end{array} \ea \right\} \cong \left\{ \ba \begin{array}{c} \text{twists on $C$} \end{array} \ea \right\}.
 \]
Moreover, if the twist is compatible with duals then the corresponding even-handed structure is spherical.
\end{lem}
\begin{proof}
Given a twist $\theta$, we obtain an even-handed structure by setting
 \[
 \Psi \left( \ba \ig{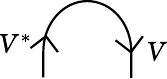} \ea, \; \ba \ig{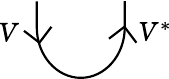} \ea \right) =  \left( \ba \ig{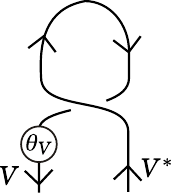} \ea, \; \ba \ig{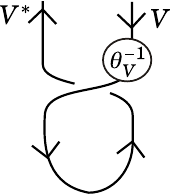} \ea \right).
\]
Similarly given an even-handed structure $\Psi$, we obtain a twist by setting
 \[
  \theta_V = \ba \ig{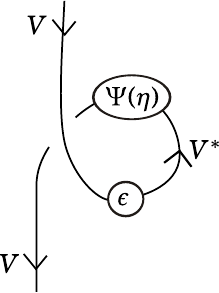} \ea
 \]
where $(V^*, \eta, \epsilon)$ is {\em any} choice of right dual for $V$ (the construction is independent of this choice). A quick calculation shows that these constructions are inverse to each other, for instance:
  \[
  \theta_V \mapsto \left( \ba \ig{e327.pdf} \ea \mapsto \ba \ig{e329.pdf} \ea \right) \mapsto \ba \ig{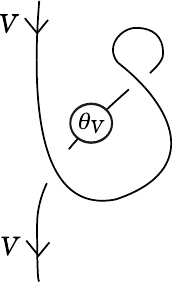} \ea = \theta_V.
\]
Moreover if $\theta_{V^*} = \theta_V^\dagger$ then $\Psi$ is spherical:
 \[
  \Tr_l (f) = \ba \ig{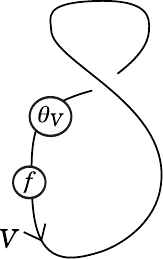} \ea = \ba \ig{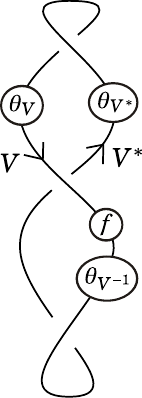} \ea = \ba \ig{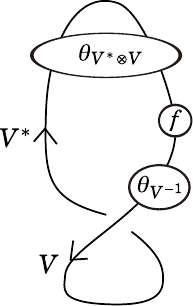} \ea = \ba \ig{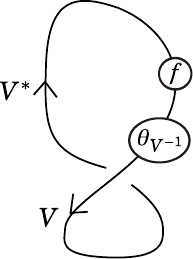} \ea = \Tr_r(f).
 \]
We remark that if the category is semisimple, the converse also holds (see Chapter \ref{FusionChap}).
\end{proof}

\section{Geometric interpretation\label{geominterp}}
\begin{figure}[t]
 \centering
 \[
 \begin{array}{c} \ba \ig{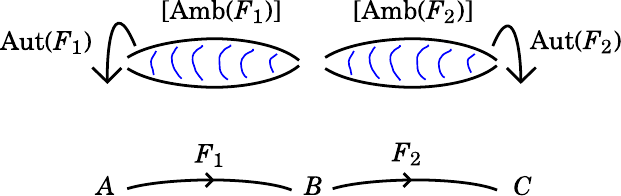} \ea \quad \quad \ba \ig{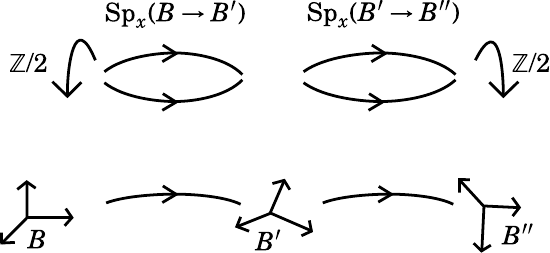} \ea \\ \text{(a)} \quad \quad \qquad \qquad\qquad \qquad \qquad  \qquad \quad \text{(b)} \end{array}
 \]
 \caption{\label{ambgerbe} (a) The `ambijunction gerbe' of a 2-category with ambidextrous adjoints. (b) The spin gerbe of a Riemannian manifold. }
 \end{figure}

One can think of an even-handed structure geometrically as a `flat section' of a certain `ambijunction bundle' over the 1-morphisms in the 2-category; this is in fact the approach we adopted in some previous work \cite{ref:bartlett}. In this viewpoint, instead of thinking of an even-handed structure as a rule for partnering right adjoints with left adjoints, one thinks of it as a coherent choice $F \mapsto F^{[*]}$ of {\em isomorphism class of ambidextrous adjoint} for every morphism $F$ in the 2-category. To understand this viewpoint, let us be precise about what we mean by an `ambidextrous adjoint'.
\begin{defn}
 An {\em ambidextrous adjoint} of a morphism $F \colon A \rightarrow B$ in a 2-category is a quintuple $\langle F^* \rangle \equiv (F^*, \eta, \epsilon, n, e)$ where  $F^* \colon B \rightarrow A$ is a morphism, $\eta \colon \id_A \Rightarrow F^*F$ and $ \epsilon \colon FF^* \Rightarrow \id_B$ are unit and counit maps exhibiting $F^*$ as a right adjoint of $F$, and $n \colon \id_B \Rightarrow F F^*$ and $e \colon F^* F \Rightarrow \id_A$ are unit and counit maps which exhibit $F^*$ as a left adjoint of $F$.
 \end{defn}
We write the data of a particular ambidextrous adjoint of $F$ in string diagrams as
 \[
   \langle F^* \rangle  = \left( \ba \ig{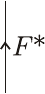} \ea, \ba \ig{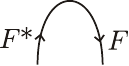} \ea, \ba
   \ig{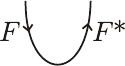} \ea, \ba \ig{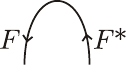} \ea, \ba \ig{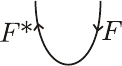} \ea \right) \, .
 \]
The choices of ambidextrous adjoints for $F$ organize themselves into a groupoid $\Amb(F)$ which we call the {\em ambijunction groupoid} of $F$. An object is a choice of ambidextrous adjoint $\langle F^* \rangle$ of $F$. A morphism $\gamma \colon \langle F^* \rangle \rightarrow  \langle (F^*)' \rangle$ of ambidextrous adjoints of $F$ is an invertible 2-morphism $\gamma \colon  F^* \Rightarrow (F^*)'$, drawn as
 \[
  \ba \ig{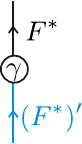} \ea \,
 \]
such that `twisting' the unit and counits of $\langle F^* \rangle$ by $\gamma$ results in $\langle (F^*)' \rangle$, that is,
 \[
  \left( \ba \ig{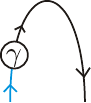} \ea, \ba \ig{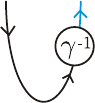} \ea, \ba
  \ig{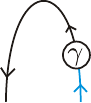} \ea, \ba \ig{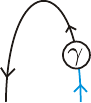} \ea \right) = \left( \ba
  \ig{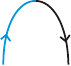} \ea, \ba \ig{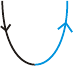} \ea, \ba \ig{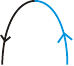} \ea, \ba
  \ig{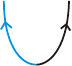} \ea \right).
 \]
We write $[\Amb(F)]$ for the set of isomorphism classes in the ambijunction groupoid of $F$, and we write the class of a particular ambidextrous adjoint $\langle F^* \rangle$ as $[F^*]$. The properties of the ambijunction groupoid are summarized in the following lemma, whose proof involves elementary string diagram considerations.
\begin{lem} \label{amblem1} Suppose $F$ is a morphism in a 2-category and that the groupoid $\Amb(F)$ is nonempty. Then:
 \begin{enumerate}
  \item There is at most one arrow between any two ambidextrous adjunctions in $\Amb(F)$.
  \item The group $\Aut(F)$ of automorphisms of $F$ acts freely and transitively on $[\Amb(F)]$ by twisting the unit and counit maps which display $F^*$ as a left adjoint of $F$,
      \[
  [F^*] \stackrel{\alpha}{\mapsto} \left[ \ba \ig{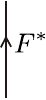} \ea, \ba \ig{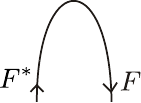} \ea, \ba \ig{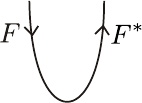} \ea, \ba
  \ig{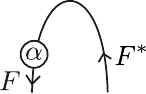} \ea, \ba \ig{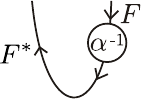} \ea \right].
 \]
 \end{enumerate}
\end{lem}
Geometrically, this lemma says that the sets $[\Amb(F)]$ of isomorphism classes of ambidextrous adjoints for each 1-morphism $F$ in a 2-category $\mathcal{C}$ can be thought of as forming a principal bundle-like structure
\be \label{gerbe1}
 [\Amb(\mathcal{C})] \rightarrow \text{1-Mor}(\mathcal{C})
\ee
over the 1-morphisms in $\mathcal{C}$, which we call the {\em ambijunction gerbe} (see Figure \ref{ambgerbe}). The idea is that the fiber of a 1-morphism $F$ in the ambijunction gerbe is the $\Aut(F)$-torsor $[\Amb(F)]$ of isomorphism classes of ambidextrous adjoints for $F$. An even-handed structure $\Psi$ can then be thought of as `even-handed trivialization' of this gerbe, defined as the section which sends
 \[
  F \mapsto \left[\ba \ig{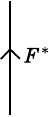} \ea, \; \ba \ig{e308.pdf} \ea, \; \ba \ig{e309.pdf} \ea, \; \ba \ig{e310.pdf} \ea, ;\ba \ig{e311.pdf} \ea \right] \in [\Amb(F)]
 \]
where $(F^*, \eta, \epsilon)$ is any choice of right adjoint for $F$; the fact that this is well-defined follows from the naturality properties of $\Psi$. Conversely, one could {\em define} an even-handed structure as `an even-handed trivialization' $[*]$ of the ambijunction gerbe,
 \[
  F \mapsto F^{[*]} \in [\Amb(F)],
 \]
from which one could recover $\Psi$ as
 \begin{multline*}
 \Psi\left(\ba \ig{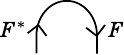} \ea, \; \ba \ig{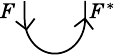} \ea\right) := \text{unique unit and counit maps } \left(\ba \ig{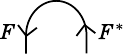} \ea, \; \ba \ig{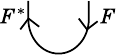} \ea \right) \\
 \text{ such that } F^{[*]} = \left[\ba \ig{e312.pdf} \ea, \; \ba \ig{e313.pdf} \ea, \; \ba \ig{e314.pdf} \ea, \; \ba \ig{e315.pdf} \ea, \;\ba \ig{e316.pdf} \ea \right].
 \end{multline*}
Further details on this point of view are contained in \cite{ref:bartlett}, but we will not need them here.

The idea of an even-handed structure as a trivialization of the `ambijunction gerbe' is analagous to Murray and Singer's reformulation \cite{ref:murray_singer} of the notion of a spin structure on a Riemannian manifold $M$ as a trivialization of the {\em spin gerbe}, which we briefly recall for the reader's benefit. At each point $x \in M$, we may regard the set $\text{Fr}_x(M)$ of orthonormal frames $B$ of the tangent space at $x$ as a codiscrete groupoid (in other words, we imagine that there is precisely one arrow $B \rightarrow B'$ between any two frames). For each arrow $B \rightarrow B'$ in $\text{Fr}_x(M)$ we have the set $\text{Sp}_x (B \rightarrow B')$ of homotopy classes of paths of orthogonal rotations which transform $B$ into $B'$ (it is a $\mathbb{Z}/2$-torsor). Performing this construction for all $x\in M$ gives rise to the {\em spin gerbe}
 \be \label{gerbe2}
  \text{Sp}(M) \rightarrow \text{Fr}(M)
 \ee
and Murray and Singer showed that a trivialization of this gerbe is precisely a spin structure on $M$ in the classical sense. The idea is that \eqref{gerbe2} is completely analagous to \eqref{gerbe1}, as we have indicated in Figure \ref{ambgerbe}.

In Section \ref{tracelinsec} we will offer two further facts which tentatively link even-handed structures with spin structures: firstly, that the set of even-handed structures on a collection of semisimple categories determines a weighting (or `metric') on each category, but only up to a global scale factor, which is analagous to the fact that the dimension of the space of harmonic spinors on a spin manifold does not depend upon rescaling the metric \cite{ref:hitchin}, and secondly that the formula $\Psi(\phi) = * \,\, \phi^* *$ for the even-handed structure on the 2-category of 2-Hilbert spaces is roughly analagous to the formula for the adjoint of the derivative operator, $d^* = * d *$.

\section{Relationship with other approaches\label{srel}}
In this section we compare the notion of an even-handed structure to the more established approaches to duality in 2-categories: {\em pivotal structures} and {\em monoidal 2-categories with duals}.

\subsection{Pivotal structures on monoidal categories\label{pivssec}}
To the author's knowledge, the notion of a {\em pivotal structure} on a monoidal category (or at least, a {\em pivotal category}, which is a monoidal category equipped with a pivotal structure) goes back to Freyd and Yetter \cite{ref:freyd_yetter}, who actually credit this notion to unpublished work of Joyal and Street \cite{ref:joyal_street_bmc2}. For a more recent account of the notion of a pivotal structure, we recommend the lecture notes of Boyarchenko \cite{ref:boyarchenko}.

Suppose that every object $V$ in a monoidal category $C$ comes equipped with a {\em chosen} right dual $(V^\star, \eta_V, \epsilon_V)$; we will call this a {\em system of right duals} on $C$. We use the notation $V^\star$ instead of our more usual $V^*$ in order to make the conceptual distinction between duals picked `freely at the point of calculation' written as $V^*$, and distinguished `fixed-for-all-time' duals written as $V^\star$. At any event, it is not hard to see that choosing a system of right duals gives rise to a contravariant monoidal functor $\star \colon C \rightarrow C$.
\begin{defn} A {\em pivotal structure} on monoidal category with a chosen system of right duals is a monoidal natural transformation $\gamma \colon \id \Rightarrow \star \star$ satisfying $\gamma_{V^\star} = (\gamma_V^{-1})^\dagger$.
\end{defn}
We now make a few remarks about this definition. Firstly it is deceptively simple: if we expand out what it means for $\gamma$ to be a {\em monoidal} natural isomorphism, we find that it is the requirement that
 \be \label{monism}
 \ba \ig{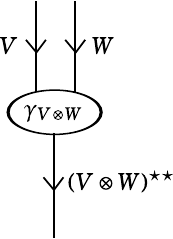} \ea \; = \; \ba \ig{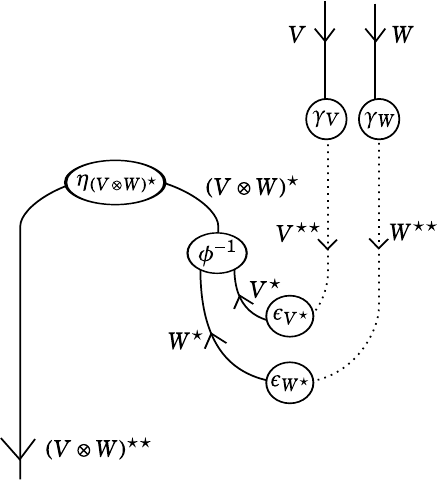} \ea
 \ee
where $\phi \colon W^* \otimes V^* \rightarrow V \otimes W$ is the coherence isomorphism presenting $\star$ as a monoidal functor,
\[
 \phi = \ba \ig{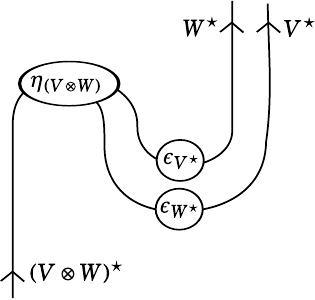} \ea.
\]
One should compare this tricky diagram to the corresponding diagram in the definition of an even-handed structure:
 \[
  \Psi\left(\ba \ig{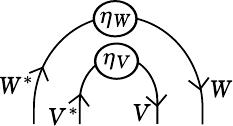} \ea\right) = \ba \ig{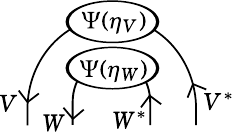} \ea.
 \]
Secondly, the last requirement on the definition --- that $\gamma_{V^\star} = (\gamma_V^{-1})^\dagger$ --- was present in the early definition of Freyd and Yetter \cite{ref:freyd_yetter} as well as in the treatment of Bakalov and Kirillov \cite{ref:bakalov_kirillov} but in more recent literature appears to have been dropped (see eg. \cite{ref:calaque_etingof, ref:müger_lectures, ref:boyarchenko}). This is probably because in many situations, such as for monoidal linear categories which are semisimple, this condition is in fact a consequence of $\gamma$ being monoidal and hence unnecessary. Under the corrrespondence between even-handed structures and pivotal structures which we spell out below, this condition corresponds to the condition $\Psi^2 = \id$ which we noted in Section \ref{defesec} that we have been unable to show follows from the other axioms.

\begin{prop}\label{pivev} Let $C$ be a monoidal category with ambidextrous duals. Then for each choice of system of right duals $\star$, there is a canonical bijection
 \[
  \left\{ \text{even-handed structures on $C$} \right\} \cong \left\{ \ba \begin{array}{c} \text{pivotal structures }\\ \text{on $C$ with respect to $\star$} \end{array} \ea \right\}.
 \]
Moreover, this bijection is natural with respect to altering the system of right duals $\star$.
\end{prop}
\begin{proof}
Given an even-handed structure $\Psi$ we obtain a pivotal structure $\gamma \colon \id \Rightarrow \star \star$ by setting
\[
 \gamma_V := \ba \ig{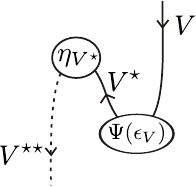} \ea, \text{ which we'll often draw simply as } \ba \ig{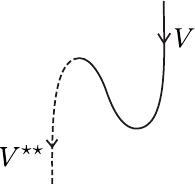} \ea.
\]
These isomorphisms are natural, because if $\theta \colon V \rightarrow W$ is a morphism, then the even-handed equation $\theta^\dagger = ^\dagger \! \theta_\Psi$ implies that
 \[
  \theta^{\ast \ast} \circ \gamma_V = \ba \ig{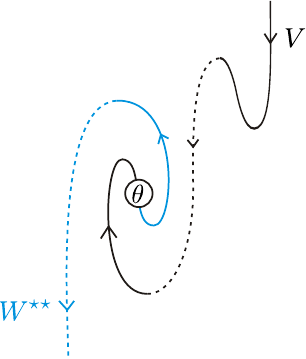} \ea = \ba \ig{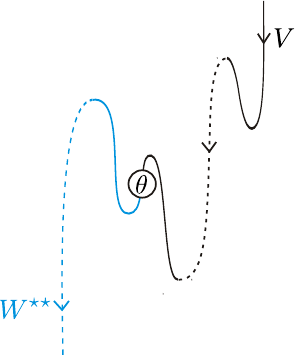} \ea = \ba \ig{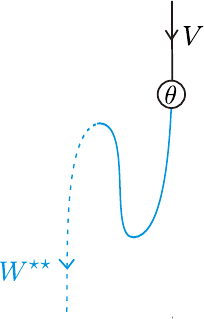} \ea = \gamma_W \circ \theta.
 \]
We also need to show that $\gamma$ is a {\em monoidal} natural isomorphism, that is that equation \eqref{monism} is satisfied. If we expand this out, we see that it is the requirement that
 \[
 \ba \ig{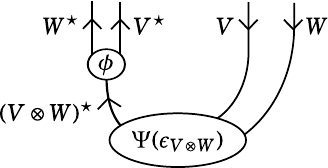} \ea = \ba \ig{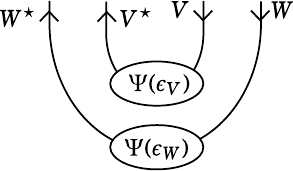} \ea.
\]
But this follows from the naturality of $\Psi$ as expressed in Lemma \ref{natbugger} and the fact that it preserves tensor products:
 \begin{align*}
 \ba \ig{e398.pdf} \ea &= \Psi\left(\ba \ig{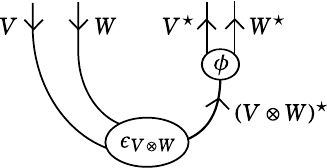} \ea\right)
 \\ &= \Psi \left(\ba \ig{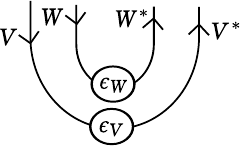} \ea \right) \\
 &= \ba \ig{e399.pdf} \ea .
\end{align*}
Finally, the requirement $\gamma_{V^\star} = (\gamma_V^\mi)^\dagger$ boils down to the requirement that
 \[
  \ba \ig{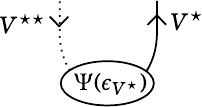} \ea = \ba \ig{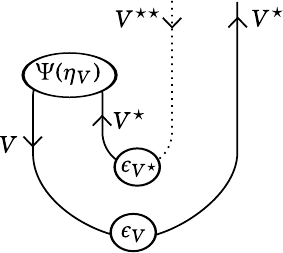} \ea.
 \]
And this follows from applying $\Psi$ to both sides and using $\Psi^2 = \id$:
 \[
  \Psi \left(\ba \ig{e403.pdf} \ea\right) = \ba \ig{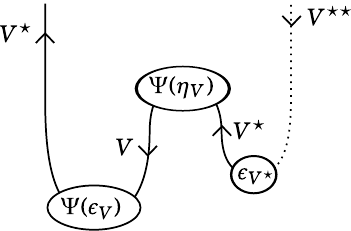} \ea = \ba \ig{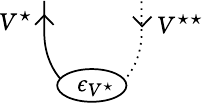} \ea.
 \]

In the reverse direction, from a chosen system of right duals $(V^\star, \eta_V, \epsilon_V)$ and a pivotal structure $\gamma \colon \id \Rightarrow \star \star$ we obtain an even-handed structure by defining its action on the chosen right duals $(V^\star, \eta_V, \epsilon_V)$ as
 \[
 \left( \ba \ig{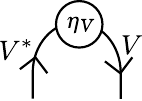} \ea, \; \ba \ig{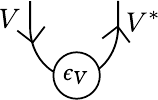} \ea \right) \stackrel{\Psi}{\mapsto} \left( \ba \ig{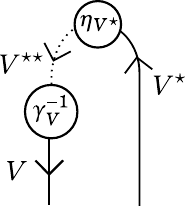} \ea, \; \ba \ig{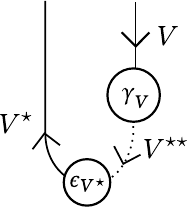} \ea \right)
 \]
and then extending this by naturality of $\Psi$ to {\em any} right dual $(V^*, \eta, \epsilon)$:
 \[
\left( \ba \ig{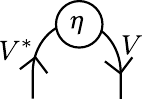} \ea, \; \ba \ig{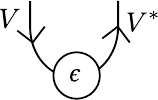} \ea\right) \stackrel{\Psi}{\mapsto} \left( \ba \ig{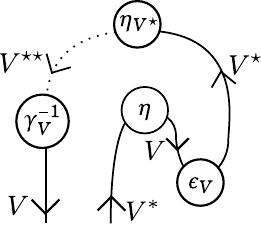} \ea, \; \ba \ig{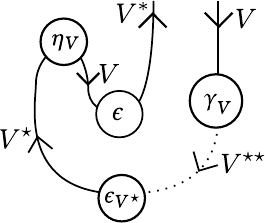} \ea \right)
 \]
(This illustrates one of the difficulties of the pivotal structure framework and of having {\em fixed} right duals $V^\star$: it makes for awkward expressions when one wishes to use some other choice of right dual $V^*$ in a calculation.) We have $\Psi^2 = \id$ because we have seen that this is equivalent to $\gamma_{V^\star} = (\gamma_V^{-1})^\star$. Moreover, the fact that $\Psi$ respects the monoidal structure follows from the fact that $\gamma$ is monoidal, by running our argument above backwards. Finally the even-handed equation corresponds precisely to the naturality of $\gamma$:
 \[
  ^\dagger \! \theta_\Psi = \ba \ig{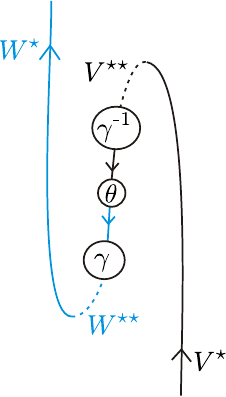} \ea = \ba \ig{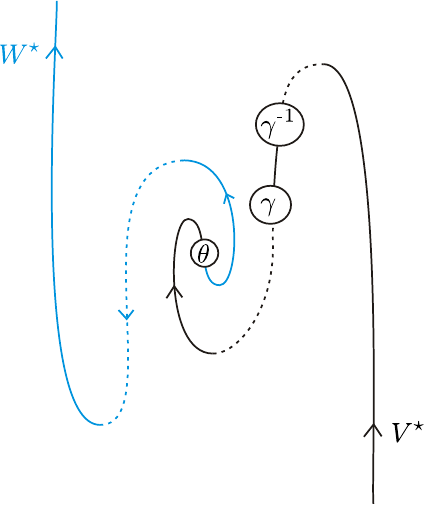} \ea = \ba \ig{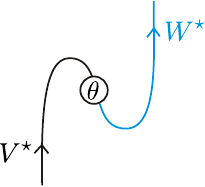} \ea = \theta^\dagger.
 \]
These constructions are easily seen to be mutual inverses of each other, so the first statement of the proposition follows.

Finally, suppose that $\star'$ is another system of choices of right dual $(V^{\star'}, \eta'_V, \epsilon'_V)$ for each object $V$, drawn using different colours as \XXX{Sort out colours!!}
 \[
  \left( \ba \ig{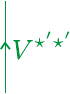} \ea, \ba \ig{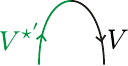} \ea, \ba \ig{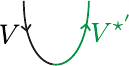} \ea \right),
 \]
and suppose that $\gamma$ and $\gamma'$ are the pivotal structures with respect to $\star$ and $\star'$ respectively, arising from the even-handed structure on $C$ via the above prescription. If we write $\xi \colon \star \star \Rightarrow \star' \star'$ for the canonical natural isomorphism between the two systems of right duals, then for each object $V$ we have
 \begin{align*}
  \xi_V \circ \gamma_V = \ba \ig{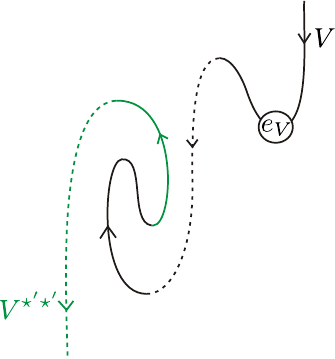} \ea &= \ba \ig{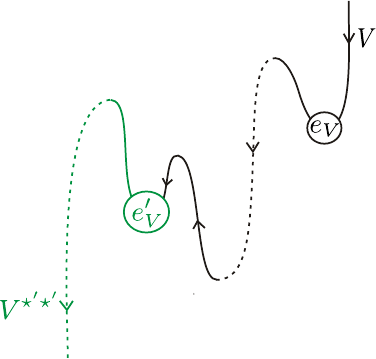} \ea \quad \text{(by construction)} \\
  & = \ba \ig{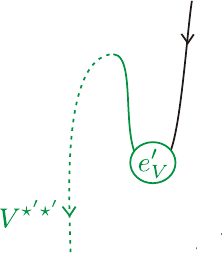} \ea  = \gamma_V'.
 \end{align*}
In other words, this bijection is natural with respect to changing the system of right duals.
\end{proof}

\subsection{The coherence theorem of Barrett and Westbury\label{cohbw}}
In \cite{ref:barrett_westbury}, Barrett and Westbury\footnote{We remark that here we are actually using the notation and slightly-different-but-equivalent reformulation given by M\"{u}ger \cite{ref:müger_from_subfactors_to_categories_and_topology_I}.} essentially defined a {\em strict pivotal category} as a strict monoidal category where every object $V$ has a chosen ambidextrous dual $(\overline{V}, \eta_V, \epsilon_V, \eta_{\overline{V}}, \epsilon_{\overline{V}})$ with $\overline{\overline{V}} = V$ in the obvious sense, which is compatible with the monoidal structure and which satisfies the even-handed equation --- that is, it does not matter if one computes the dagger of a morphism using the {\em right} unit and counit maps $(\eta_V, \epsilon_V)$ or the {\em left} unit and counit maps $(\eta_{\overline{V}}, \epsilon_{\overline{V}})$ . They proved a useful coherence theorem which had as one of its consequences that every pivotal category is equivalent to a strict one. Their motivation was to use spherical categories to define 3-manifold invariants --- although most natural examples of spherical categories are not strict, it is simpler to do calculations with a strict pivotal category than with a non-strict one. We also note that M\"{u}ger outlined how to define a strict pivotal 2-category in the same way \cite{ref:müger_from_subfactors_to_categories_and_topology_I}.

\subsection{Monoidal 2-categories with duals} In the context of their work on 2-tangles, Baez and Langford \cite{ref:baez_langford} defined a {\em monoidal 2-category with duals} as a semistrict monoidal 2-category equipped with a coherent notion of duals for objects, morphisms and 2-morphisms. To compare this notion with that of an even handed 2-category, we must first strip out from their definition those pieces of data which have to do with the monoidal structure on the 2-category and the notion of duals for objects, which are not needed for our purposes. We arrive at the notion of {\em strict 2-category with duals for morphisms and 2-morphisms}; the main feature is that every 2-morphism $\theta \colon F \Rightarrow G$ has a specified dual 2-morphism $\theta^* \colon G \Rightarrow F$ with $\theta^{**} = \theta$, and that every 1-morphism $F \colon A \rightarrow B$ has a specified right adjoint $(F^\star, \eta_F, \epsilon_F)$; by applying $*$ to $(\eta_F, \epsilon_F)$ we see that $F^\star$ is also a left adjoint of $F$. The main equations relating the duals for 2-morphisms and for 1-morphisms are that
 \[
  (\eta_{F^\star}, \epsilon_{F^\star}) = (\epsilon_F^*, \eta_F^*)
  \]
for every 1-morphism $F$ and
\[
 (\theta^\dagger)^* = (\theta^*)^\dagger, \;\; \text{or equivalently } \theta^\dagger = {}^\dagger \! \theta
 \]
for every 2-morphism $\theta$, where $\theta^\dagger, {}^\dagger \! \theta \colon G^\star \Rightarrow F^\star$ are calculated using the structure of $F^\star$ as a right and left adjoint of $F$ respectively. Clearly then a strict 2-category with duals for morphisms and 2-morphisms gives rise to an even-handed structure $\Psi$ on its 1-morphisms, defined first on the chosen duals by
 \[
 \Psi (\eta_F, \epsilon_F) = (\epsilon_F^*, \eta_F^*)
 \]
and then extended to {\em all} duals by naturality, as in the proof of Proposition \ref{pivev}. An important point to bear in mind though is that the resulting even-handed structure does not in {\em general} send
 \[
 (\eta, \epsilon) \mapsto (\epsilon^*, \eta^*)
 \]
because this formula is not natural (it is only natural with respect to {\em unitary} isomorphisms). In other words, in a 2-category with duals for morphisms and 2-morphisms it is {\em not} the case in general that `the natural way to turn an arbitrary right adjoint into a left adjoint is to take the stars of the unit and counit maps'; one must first apply a normalization factor as in the proof of Proposition \ref{pivev}.

\subsection{Comparison with even-handed structures}
In this subsection we compare the even-handed structures framework for dealing with duality in 2-categories to the frameworks of pivotal 2-categories and monoidal 2-categories with duals.

For our purposes, one of the main difficulties with these approaches is their starting point, namely that every 1-morphism $F$ comes equipped with a {\em fixed} chosen right dual $F^\star$. This is an unfortunate requirement, for a number of reasons. Firstly it generates awkward coherence diagrams which do not mesh comfortably with the string diagram notation, as we saw in Section \ref{pivssec}. We are not allowed to draw $(G \circ F)^\star$ as two separate lines $G^\star$ and $F^\star$ and this makes things awkward. One way round this is to use the coherence theorem of Barrett and Westbury \cite{ref:barrett_westbury} which justifies making these harmless identifications. However we could have avoided this entire issue from the start by not making the fixed choices of duals; moreover this would allow us to seamlessly make string diagram calculations in a {\em non-strict} context (such as those which actually occur in examples) without being led to make the unnatural step of first passing to a strictified situation and then appealing to coherence. Although this latter strategy is not actually required by the Barrett and Westbury coherence theorem, in practice this has often been the way in which it has been used (see \cite{ref:müger_from_subfactors_to_categories_and_topology_I, ref:barrett_westbury2} and also \cite{ref:hagge_hong}).

In the next chapter we will take advantage of this fact --- that the notion of an even-handed structure meshes well with the string diagram calculus --- in order to prove some new results about pivotal structures on fusion categories. Our proofs will be entirely elementary, yet since they appear to be new we are led to conclude that the awkwardness of pivotal structures is not just aesthetic but actually inhibits progress.

Another difficulty with having fixed chosen right adjoints is that this is invariably not what actually crops up in examples, at least in a 2-categorical situation where oftentimes adjoints are not given explicitly but only up to isomorphism. Our contention is that it is more natural to ask for a {\em function} which converts arbitrary right adjoints structure maps into left adjoint structure maps. Experience shows that this stripped-down requirement more naturally fits the examples we have in mind, such as the 2-category $\THilb$ of 2-Hilbert spaces. Here adjoints are not given canonically, but as we will see in Section \ref{EvenTHilb} there is at least a canonical {\em function} $\Psi$ which turns a right adjoint into a left adjoint, given at the level of adjunction isomorphisms $\{\phi \colon \Hom(Fx, y) \rightarrow \Hom(x, F^*y)\}$ by
 \[
  \Psi(\phi) = * \,\, \phi^* *
 \]
where $*$ is the $*$-structure on the morphisms in a 2-Hilbert space and $\phi^*$ is the Hilbert space adjoint of $\phi$. A similar analysis applies to the 2-category $\CYau$ of Calabi-Yau manifolds and integral kernels between them\footnote{A technical point which might be mentioned here is that in `large' 2-categories such as $\THilb$ and $\CYau$ where duals are not given canonically, there could be a set-theoretic issue if one wishes to apply a `chosen-dual paridigm' by initially arbitrarily choosing duals $F^*$ for each morphism and proceeding from there. Such a choice would range over a proper {\em class} of elements and is difficult to justify.}.

Even when adjoints {\em are} given canonically, this might not always be the most natural way to proceed, since these formulas are invariably not monoidal on-the-nose. For instance in the monoidal category of vector spaces, each vector space $V$ has a canonical left-and-right linear dual $V^\vee := \Hom(V, \mathbb{C})$. But unfortunately
 \[
  (V \otimes W)^\vee \neq W^\vee \otimes V^\vee
 \]
so if we insist that the dual of a vector space $V$ {\em must} be given by $\Hom(V, \mathbb{C})$, we will be in for an awkward time. A better option is to use a horses-for-courses approach and decide that $W^\vee \otimes V^\vee$ is {\em just as good} a dual for $(V \otimes W)$ as $(V \otimes W)^\vee$. This frees one up to use the most convenient dual for the calculation at hand, rather than being forced to use a dual chosen right at the outset.

\section{Even-handedness in terms of adjunction isomorphisms\label{catsection}}
Our goal in this section is to reformulate the notion of an even-handed structure on a 2-category in the context where the objects, morphisms and 2-morphisms of the 2-category are categories, functors and natural transformations respectively.

Suppose that $F\colon A \rightarrow B$ is a functor between categories $A$ and $B$, and that $F^*\colon B \rightarrow A$ is a right adjoint of $F$. Instead of expressing this via the unit and counit natural transformations $\eta$ and $\epsilon$, we will find it more convenient to express it in terms of the adjunction isomorphisms
 \[
  \{\phi \colon  \Hom(Fx, y) \rightarrow \Hom(x, F^*y)\}_{x \in A, y \in B}.
 \]
Recall the translation between these two pictures --- given a morphism $f \colon  Fx \rightarrow y$ we set $\phi(f) = F^*(f) \circ
\eta_x$ and similarly for $g \colon  x \rightarrow F^*y$ we have $\phi^{\mi}(g) = \epsilon_y \circ F(g)$. We will write $\Adj_{\Hom} (F \dashv F^*)$ for the set of all adjunction isomorphisms which exhibit $F^*$ as a right adjoint of $F$:
  \[
  \Adj_{\Hom} (F \dashv F^*) = \left\{ \begin{array}{c} \text{sets of natural isomorphisms}  \\ \{\phi \colon \Hom(F x, y) \rightarrow \Hom(x, F^*y) \}_{x \in A, y \in B} \end{array} \right\}
  \]
Recall that the group $\Aut(F)$ of automorphisms of the functor $F$ acts freely and transitively on $\Adj_{\Hom}(F \dashv F^*)$ by precomposition,
  \[
   \delta \cdot \phi = \phi \circ \pre(\delta).
  \]
That is, if $f \colon Fx \rightarrow y$ and $\delta \colon F \Rightarrow F$ is an automorphism then $(\delta \cdot \phi)(f) = \phi(f \circ \delta_x)$.

We now translate the definition of an even-handed structure from the language of unit and counit natural transformations into the language of adjunction isomorphisms. The only tricky point is the even-handed equation $\theta^\dagger = {}^\dagger \! \theta_\Psi$. When working with hom-set isomorphisms, to identify a morphism (such a component $\theta^\dagger_y$ for some $y\ \in B$) it is more natural to give a formula for {\em post- or pre-composition} by that morphism; the magical Yoneda lemma tells us that knowing how to post-compose or pre-compose with a morphism is the same thing as knowing the morphism itself. (To see this, recall that the Yoneda lemma states that for any pair of objects $x, x' \in A$ the map
 \begin{align*}
  \Hom(x,x') &\rightarrow \Nat (\underline{x}, \underline{x'}) \\
   f & \mapsto \post(f)
  \end{align*}
is a bijection, where $\underline{x}, \underline{x'} \colon A^\text{op} \rightarrow \Set$ are the presheaves represented by $x$ and $x'$ respectively.)

So, suppose that $A$ and $B$ are categories, and that $F, G \colon A \rightarrow B$ are functors between them with ambidextrous adjoints $F^*$ and $G^*$ respectively, and that $\theta \colon F \Rightarrow G$ is a natural transformation. Then we have the following characterization of the right and left daggers of $\theta$.
\begin{lem} \label{ppcomp} For objects $x \in A$ and $y \in B$, post- and pre-composition with the components of the right and left daggers $\theta^\dagger_y, {}^\dagger\!\theta_y \colon G^*y \rightarrow F^*y$ can be expressed as
 \begin{align*}
  \post(\theta_y^\dagger) &= \phi_F \circ \pre(\theta_x) \circ \phi_G^\mi \\
  \pre({}^\dagger\!\theta_y) &= \psi_G^\mi \circ \post(\theta_x) \circ \psi_F
  \end{align*}
 where $\phi_F$ and $\phi_G$ are the adjunction isomorphisms expressing $F^*$ and $G^*$ as {\em right} adjoints of $F$ and $G$ respectively, and $\psi_F$ and $\psi_G$ are the adjunction isomorphisms expressing $F^*$ and $G^*$ as {\em left} adjoints of $F$ and $G$ respectively.
\end{lem}
\begin{proof} We prove the formula for the right dagger $\theta^\dagger_y$; the proof for the left dagger  $^\dagger\!\theta_y$ is
similar. Given a morphism $f \colon x \rightarrow G^*y$, we can construct the following diagram:
 \[
  \xymatrix{   x \ar[rr]^{\eta^F_x} \ar[d]_{g} && F^*Fx \ar[d]^{F^*F(g)} \ar[rr]^{F^*(\theta_x)} && F^*Gx \ar[d]^{F^*G(g)} \\
  G^*y \ar[rr]_{\eta^F_{G^*y}} && F^*FG^*y \ar[rr]_{F^*(\theta_{G^*y})} && F^*GG^*y \ar[rr]_{F^*(\epsilon^G_y)} &&
  F^*(y)}
 \]
The counter-clockwise composite is the definition of $\post(\theta^\dagger_y)(g)$ from its string diagram formula \eqref{adjoints1}, while the clockwise composite is the formula we are trying to prove. The left-hand square commutes by the naturality of $\eta$, while the right-hand square commutes by the naturality of $\theta$.
\end{proof}
To make use of this lemma, we will use the notation whereby if $f$ is a morphism in a category, and $\Omega_f$ is some formula for how to post-compose morphisms with $f$, then $f = \post^\mi(\Omega_f)$ will be shorthand for the statement `$f$ is the unique morphism such that post-composing with it is the same as applying $\Omega_f$'. We will use a similar notation for the `inverse' of pre-composition. We then have the following.

\begin{prop} \label{altform} An even-handed structure on a full sub-2-category $\mathcal{C} \subset \Cat$ with ambidextrous adjoints is the same thing the data of, for every pair of adjoint functors $F, F^*$, a map
 \begin{align*}
  \Psi \colon \Adj_{\Hom} (F \dashv F^*) & \rightarrow \Adj_{\Hom} (F^* \dashv F) \\
   \phi & \mapsto \Psi(\phi)
  \end{align*}
satisfying:
 \begin{enumerate}
     \item $\Psi^2 = \id$.
    \item For composable adjoints  $ \xymatrix @1 @C=0.25in  { A \ar@<0.6ex>[r]^F & B \ar@<0.2ex>[l]^{F^*} \ar@<0.6ex>[r]^G  & C \ar@<0.2ex>[l]^{G^*} }$,
     \[
       \Psi (\phi_F \circ \phi_G) = \Psi (\phi_G) \circ \Psi (\phi_F)
     \]
     for all $\phi_F \in \Adj(F \dashv F^*)$ and $\phi_G \in \Adj(G \dashv G^*)$, and moreover $\Psi(\id) = \id$ for all objects $A$, where $\id \colon \Hom(x,y) \rightarrow \Hom(x,y)$ is the trivial adjunction isomorphism.
    \item If $F \dashv F^*$, $G \dashv G^*$ and $\theta \colon F \Rightarrow G$ is a natural transformation, then
     \[
      \post^{-1} \bigl(\phi_F \circ \pre(\theta) \circ \phi_G ^{-1}\bigr) = \pre^{-1}\bigl(\Psi (\phi_G)^{-1} \circ \post(\theta) \circ \Psi(\phi_F) \bigr)
    \]
  for all $\phi_F \in \Adj(F \dashv F^*)$ and $\phi_G \in \Adj(G \dashv G^*)$.
 \end{enumerate}
\end{prop}
\begin{proof} The first two properties are manifestly the translation into the hom-set adjunction isomorphisms language of the corresponding properties in the definition of an even-handed structure, while the third is precisely Lemma \ref{ppcomp}.
\end{proof}
We will find this reformulation of the notion of an even-handed structure very useful, because it uses precisely the sort of language which arises when one has a {\em trace} on a linear category.

\section{Even-handed structures from traces\label{tracelinsec}} Our goal in this section is to define the notion of a {\em trace} on a linear category, and to show how these give rise to even-handed structures. When the categories involved are semisimple, we will show that an even-handed structure is actually the {\em same thing}, up to a global scale factor, as the data of a trace on each category.

\subsection{Traces on linear categories}
A {\em $k$-linear category} is a category enriched in finite-dimensional $k$-vector spaces, where $k$ is some field. The following definition abstracts the algebraic geometry notion of a {\em Serre trace} coming from a `trivial Serre functor' (see for example \cite{ref:caldararu_willerton} and the references therein). Note that we write $V^\vee := \Hom(V, k)$ for the linear dual of a vector space.
\begin{defn} A {\em trace} on a $k$-linear category consists of a linear map $\Tr_x \colon  \End(x) \rightarrow k$ for each object $x$ which is
 \begin{enumerate}
  \item nondegenerate, in the sense that for all objects $x$ and $y$ the associated map
   \begin{align*}
    s\colon \Hom(x,y) &\rightarrow \Hom(y,x)^\vee \\
    f &\mapsto (g \mapsto \Tr_x(gf))
   \end{align*}
    is an isomorphism of vector spaces, and
  \item cyclic, that is $\Tr_x(gf) = \Tr_y(fg)$ for all $f \colon  x \rightarrow y$ and $g \colon  y \rightarrow x$.
   \end{enumerate}
\end{defn}
If the hom-sets in the linear category come equipped with a $\mathbb{Z}$-grading, then we shall require instead that the trace be {\em graded commutative} in the sense that $\Tr(gf) = (-1)^{\deg(f) \deg(g)} \Tr(fg)$ whenever this makes sense.

For our purposes, the main point about a trace is that it gives the composition operation in the category a certain symmetry, namely that one can express pre-composition in terms of the linear dual of post-composition, and vice-versa. We record this in a lemma.

 \begin{lem} \label{tracesym} Suppose that a $k$-linear category is equipped with a trace.
 \begin{enumerate} \item Given a morphism $f \colon x \rightarrow y$ and some other object $z$, the following diagrams commute:
  \[
  \ba \xymatrix{ \Hom(y,z) \ar[r]^-{\pre(f)} \ar[d]_s & \Hom(x,z) \ar[d]^s \\ \Hom(z,y) \ar[r]_-{\post(f)^\vee} & \Hom(z,x)^\vee} \ea \quad \quad
  \ba \xymatrix{ \Hom(z,x) \ar[r]^-{\pre(f)} \ar[d]_s & \Hom(z,y) \ar[d]^s \\ \Hom(x,z)^\vee \ar[r]_-{\post(f)^\vee} & \Hom(y,z)^\vee} \ea
  \]
  That is,
  \[
    \pre(f) = s^\mi \circ \post(f)^\vee \circ s \quad \text{and} \quad \post(f) = s^\mi \circ \pre(f)^\vee \circ s.
 \]
 \item Moreover, the following diagram commutes:
 \[
  \xymatrix{ \Hom(x,y) \ar[dr]_{\cong} \ar[r]^{s_{x,y}} & \Hom(y,x)^\vee \ar[d]^{{s_{y,x}^\mi}^\vee} \\ & \Hom(x,y)^{\vee \vee}}
 \]
  \end{enumerate}
 \end{lem}
 \begin{proof} (i) For the first diagram, suppose that $g \colon y \rightarrow z$ and $h \colon z \rightarrow x$. Then $s \circ \pre(f) (g)$ sends $h \mapsto \Tr_x (hgf)$ while $\post(f)^\vee \circ s (g)$ sends $h \mapsto \Tr_y (fhg) = \Tr_x (hgf)$ since the trace is cyclic. Note that if the trace is only graded cyclic, then this triangle fails to commute by the factor $(-1)^{\deg(f)\deg(g) \deg{h}}$. For the second diagram, $s \circ \post(f) (h)$ sends $g \mapsto \Tr_z(hfg)$, and so does $\pre(f)^\vee \circ s$, so the fact that the diagram commutes does not even require  the trace to be cyclic.

 (ii) Given a morphism $k \colon y \rightarrow x$, the clockwise composite sends $f \mapsto \left(s(k) \mapsto \Tr_x (k f)\right)$ while the canonical diagonal isomorphism sends $f \mapsto \left( s(g) \mapsto \Tr_y( f \circ k) \right)$, so the diagram commutes because the trace is cyclic. Note that if the trace is only graded cyclic, then this triangle fails to commute by the factor $(-1)^{\deg(f)\deg(k)}$.
\end{proof}

\subsection{Even-handed structures from traces\label{popq}}
In this subsection we show how to obtain even-handed structures from linear categories equipped with traces. We say that a functor between linear categories is {\em linear} if it is linear on the hom-sets, and we write $\LCat_k$ for the 2-category of $k$-linear categories, linear functors and natural transformations.

Suppose that $A$ and $B$ are linear categories equipped with a trace, that $F \colon  A \rightarrow B$ is a linear
functor between them and that $F^* \colon  B \rightarrow A$ is a right adjoint of $F$ with accompanying adjunction
isomorphisms
 \[
 \phi \colon  \Hom(Fx, y)  \rightarrow \Hom(x, F^*y) \quad x \in A, \, y \in B.
 \]
We can then equip $F^*$ as a {\em left} adjoint of $F$ by defining the hom-set isomorphisms
\[
 \phi^t  \colon  \Hom(F^*y, x) \rightarrow \Hom(y, Fx)
\]
by the diagram
 \[
  \ba \xymatrix{\Hom(F^*y, x) \ar[d]_{s_A} \ar[r]^{\phi^t} & \Hom(y, Fx) \ar[d]^{s_B} \\
  \Hom(x,F^*y)^\vee \ar[r]_{\phi^\vee} & \Hom(Fx, y)^\vee} \ea .
 \]
That is, we set
 \[
  \phi^t := s_B^\mi \circ \phi^\vee \circ s_A.
 \]
Note that $\phi^t$ can be regarded as the {\em adjoint} of $\phi$ with respect to the trace maps, because if $g\colon F^*y \rightarrow x$ then $\phi^t(g) \colon  y\rightarrow Fx$ is the unique morphism satisfying
 \[
 \Tr_{y}(f\phi^t(g)) = \Tr_{F^*y} (\phi(f) g)
 \]
for all $f \colon  Fx \rightarrow y$.

\begin{prop} \label{tracepn}Suppose that a full sub 2-category $\mathcal{C} \subseteq \LCat_k$ has ambidextrous
adjoints and that each category comes equipped with a trace. Then the formula
  \[
   \Psi(\phi) = \phi^t
  \]
equips $\mathcal{C}$ with an even-handed structure.
\end{prop}
\begin{proof}
We verify the conditions for an even-handed structure as presented in Proposition \ref{altform}. Firstly, $\Psi^2 = \id$ since for morphisms $f \colon Fx \rightarrow y$ we know that $\phi^{tt}(f) \colon Fx \rightarrow y$  is the unique morphism such that
 \[
  \Tr_x (g \phi^{tt}(f)) = \Tr_{Fx} (\phi^t(g) f)
 \]
for all $g \colon F^*y \rightarrow x$. But then $\phi^{tt}(f) = \phi(f)$, because
 \[
  \Tr_x (g \phi(f)) = \Tr_{F^*y} (\phi(f) g) = \Tr_y (f \phi^t(g)) = \Tr_{Fx}(\phi^t(g) f).
 \]
If the hom-sets are graded and the trace is graded cyclic, the left and right hand sides will pick up signs $(-1)^{\deg(g) \deg (\phi(f))}$ and $(-1)^{\deg (\phi^t(g)) \deg(f)}$ respectively, but these cancel since the adjunction isomorphisms $\phi$ (and hence $\phi^t$) are required to preserve the grading.\XXX{Ask simon about this}
Secondly, the transpose formula behaves correctly under composition $A \stackrel{F}{\rightarrow} B \stackrel{G}{\rightarrow} C$, since
 \begin{align*}
 \Psi(\phi_F \circ \phi_{G}) & = s_C^\mi \circ (\phi_{F} \circ \phi_{G})^\vee \circ s_A \\
  & = s_C^\mi \circ \phi_{G}^\vee \circ \phi_{F}^\vee \circ s_A    \\
  &= s_C^\mi \circ \phi_{G}^\vee \circ s_B \circ s_B^\mi \circ \phi_{F}^\vee \circ s_A  \\
  & = \Psi(\phi_{G}) \circ \Psi(\phi_{F}).
 \end{align*}
Moreover, it behaves correctly on identity adjunctions, since $\Psi(\id) = s_A^\mi \circ \id^\vee \circ s_A = \id$. Finally, the transpose formula ensures that the left and right daggers of a natural transformation $\theta \colon  F \Rightarrow G$ between two functors $F,G\colon A \rightarrow B$ agree, because for any $x \in A$ and $y \in B$ we calculate: \XXX{Decide if you want the Psi on the top or bottom when you do the left dagger}
 \begin{align*}
   \post\left(({}^\dagger\!\theta_\Psi)_y\right) &= s_A^\mi \circ \pre\left(({}^\dagger\!\theta^\Psi)_y\right)^\vee \circ s_A \\
    &= s_A^\mi \circ [\Psi(\phi_{G})^\mi \circ \post(\theta_x) \circ \Psi(\phi_{F})]^\vee \circ s_A  \quad \text{(by Lemma \ref{ppcomp})} \\
    &= s_A^\mi \circ \Psi(\phi_{F})^\vee \circ \post(\theta_x)^\vee \circ (\Psi (\phi_G)^\vee)^\mi \circ s_A \\
    &= s_A^\mi \circ s_A \circ \phi_{F} \circ s_B^\mi \circ \post(\theta_x)^\vee \circ s_B \circ \phi_{G}^{-1} \circ s_A^\mi \circ s_A \quad \text{(by Lemma \ref{tracesym}b)} \\
    &= \phi_{F} \circ \pre(\theta_x) \circ \phi_{G}^\mi \\
    &= \post(\theta_y^\dagger) \quad \text{(by Lemma \ref{ppcomp}).}
 \end{align*}
Therefore $^\dagger\!\theta_\Psi = \theta^\dagger$.
\end{proof}
\subsection{Even-handedness for semisimple categories\label{evsec}} To understand the notion of a trace and how it relates to even-handed structures, one should first consider the simplest class of examples. We use the formulation of M\"{u}ger   \cite{ref:müger_from_subfactors_to_categories_and_topology_I}
and say that a $k$-linear category is {\em semisimple} if it has direct sums and subobjects and if there exist objects $X_i$ labelled by a set $I$ such that $\Hom(X_i, X_j) \cong \delta_{ij}k$ (such objects are called {\em simple}) and such that for any two objects $x$ and $y$ the composition map
  \[
   \bigoplus_{i \in I} \Hom(x, X_i) \otimes \Hom(X_i, y) \rightarrow \Hom(x,y)
  \]
is an isomorphism. We write $\SLCat_k$ for the full sub-2-category of $\LCat_k$ whose objects are semisimple $k$-linear categories.

In this section we prove that an even-handed structure on a collection of semisimple categories is the same thing as equipping each category with a trace. On should think of this as the `categorification' of the fact that giving a $*$-structure on the morphisms between a collection of vector spaces is the same thing, up to a global scale factor, as providing each vector space with an inner product.

\begin{lem} A trace on a semisimple $k$-linear category is the same thing as an assignment of a nonzero number $k_i$ to each isomorphism class of simple object $e_i$.
\end{lem}
\begin{proof}
If the category has a trace, then we can define the number assigned to a simple object $e_i$ by taking the trace of the identity endomorphism,
 \[
  a_i := \Tr_{e_i}(\id).
 \]
Conversely, given such an assignment of numbers $e_i \mapsto a_i \in k^\times$, we obtain a trace on the category as follows. For each $x$, choose a basis $\{a_{i,p} \colon  e_i \rightarrow x\}_{p=1}^{\dim \Hom(e_i, x)}$ for $\Hom(e_i, x)$, and let $\{a_i^p \colon x \rightarrow e_i\}_{p=1}^{\dim \Hom(e_i, x)}$ be the corresponding dual basis, defined by the requirement that
 \[
  a_i^q \circ a_{i, p} = \delta_{qp} \id_{e_i}.
 \]
Note that the dual basis exists because category is semisimple. An endomorphism $f \colon  x \rightarrow x$ can be expressed in terms of this basis as
 \[
  f = \sum_{i,p,q} \langle a_i^p f a_{i,q} \rangle \, a_{i,p} a_i^q ,
 \]
where we have used the notation whereby scalars are stripped from scalar multiples of identity morphisms on the simple objects via
 \[
  \langle \lambda \id_{e_i} \rangle := \lambda.
 \]
We define the trace map $\Tr \colon  \End(x) \rightarrow k$ by summing the diagonal part of $f$, weighted by the numbers $k_i$:
 \[
  \Tr_x(f) := \sum_i k_i \sum_p \langle a_i^p f a_{i,p} \rangle.
 \]
It is easy to check that this definition is independent of the choice of basis we made for the hom-sets $\Hom(e_i, x)$, and moreover that it is nondegenerate and cyclic.
\end{proof}
To proceed further, we need to get a firm grip on the even-handed equation $\theta^\dagger = {}^\dagger \!\theta$ in the semisimple context. To this end, suppose that $F, G \colon A \rightarrow B$ are linear functors between semisimple categories with ambidextrous adjoints $F^*$ and $G^*$ respectively, with $\phi_F$ and $\psi_F$ witnessing the adjunctions $F \dashv F^*$ and $F^* \dashv F$ respectively and the same for $G$, and that $\theta \colon F \Rightarrow G$ is a natural transformation. We want to compute the relationship between the right dagger $\theta^\dagger \colon G^* \Rightarrow F^*$ and the left dagger ${}^\dagger \! \theta \colon G^* \Rightarrow F^*$.
\begin{lem}\label{ambcomplem} The matrices representing the components of ${}^\dagger \! \theta$ can be expressed in terms of the corresponding matrices for $\theta^\dagger$ and the matrices representing the adjunctions $\psi_F$ and $\psi_G$ as follows:
 \[
  [^\dagger\!\theta_{\mu, i}] = [\psi_G]^\mi [\theta^\dagger_{\mu, i}] [\psi_F].
 \]
\end{lem}
\begin{proof}
The precise meaning of this formula will become clear during the course of this proof. The strategy is to use Lemma \ref{ppcomp} which expresses the post- and pre-composition behaviour of $\theta^\dagger$ and $^\dagger\!\theta$ respectively. Let $e_i \in A$ and $e_\mu \in B$ be simple objects. We choose bases
 \[
  a_p \colon e_\mu \rightarrow F(e_i) \text{ and }  b_q \colon e_\mu \rightarrow G(e_i)
 \]
for the vector spaces $\Hom(e_\mu, Fe_i)$ and $\Hom(e_\mu, Ge_i)$ respectively. The corresponding dual bases are written as
\[
 a^p \colon F(e_i) \rightarrow e_\mu \text{ and } b^q \colon G(e_i) \rightarrow e_\mu.
\]
Let $\phi_F$ and $\phi_G$ be the adjunction isomorphisms exhibiting $F^*$ and $G^*$ as right adjoints of $F$ and $G$. We apply these isomorphisms to obtain bases for $\Hom(e_i, F^*e_\mu)$ and $\Hom(e_i, G^*e_\mu)$,
 \[
  \hat{a}_p := \phi_F (a^p) \colon e_i \rightarrow F^*(e_\mu) \,\, \text{ and } \,\, \hat{b}_q := \phi_G (b^q) \colon e_i \rightarrow G^*(e_\mu),
 \]
and we define
 \[
 \hat{a}^p \colon F^*(e_\mu) \rightarrow e_i \text{ and } \hat{b}^q \colon G^*(e_\mu) \rightarrow e_i
 \]
to be the corresponding dual bases.

Now we are ready to compute. We begin with the matrix elements of the right dagger $\theta^\dagger : G^* \Rightarrow F^*$. Using $\langle \lambda \id \rangle = \lambda$ to strip off scalars from multiples of identity morphisms as before, we have
\begin{align*}
  \langle \hat{a}^p \circ \theta_{e_\mu}^\dagger \circ \hat{b}_q \rangle &= \left\langle \hat{a}^p \circ \phi_F (\phi_G^\mi(\hat{b}_q)\circ \theta_{e_i}) \right\rangle & \text{by Lemma \ref{ppcomp} } \\
   &= \langle \hat{a}^p \phi_F (b^q \theta_{e_i}) \rangle.
  \end{align*}
We expand $b^q \circ \theta_{e_i}$ in the $a^{p'}$ basis,
 \[
  b^q \circ \theta_{e_i} = \sum_{p'} \langle b^q \circ \theta_{e_i} \circ a_p \rangle \, a^{p'}
 \]
and substitute back in to obtain
 \begin{align*}
    \langle \hat{a}^p \circ \theta_{e_\mu}^\dagger \circ \hat{b}_q \rangle &= \sum_{p'} \langle \contraction{}{a^p}{ \circ }{\hat{a}_{p'}}a^p \circ \hat{a}_{p'} \rangle \langle b^q \circ \theta_{e_i} \circ a_{p'} \rangle \\
     &=  \langle b^q \circ \theta_{e_i} \circ a_p \rangle \, .
 \end{align*}
Now we compute the matrix elements of the left dagger $\theta^\dagger$, again using Lemma~\ref{ppcomp}. We have
 \[
  \langle \hat{a}^p {} ^\dagger\!\theta_{e_i} \hat{b}_q \rangle = \langle \psi_G^\mi ( \theta_{e_i} \circ \psi_F(\hat{a}^p)) \circ \hat{b}_q \rangle \, .
 \]
We expand $\psi_F(\hat{a}^p)$ in the $a_{p'}$ basis,
 \[
  \psi_F(\hat{a}^p) = \sum_{p'} \langle \psi_F(\hat{a}^p) \circ a^{p'}\rangle \, a_{p'}
  \]
and substitute back in to obtain
 \[
  \langle \hat{a}^p {} ^\dagger\!\theta_{e} \hat{b}_q \rangle = \sum_{p'} \langle a^{p'} \circ \psi_F(\hat{a}^p) \rangle \, \langle \psi_G^\mi (\theta_{e_i} \circ a_{p'}) \circ \hat{b}_q \rangle \, .
 \]
Again, we expand $\theta_{e_i} \circ a_{p'}$ in the $b_{q'}$ basis, substitute back in, and rearrange terms to finally obtain
 \[
   \langle \hat{a}^p \circ  ^\dagger\!\theta_{e_i} \circ \hat{b}_q \rangle = \sum_{q', p'} \langle \psi_G^\mi (b_q) \circ \hat{b}_{q'} \rangle \langle b^{q'} \circ \theta_{e_i} \circ a_{p'} \rangle \langle a^{p'} \circ \psi_F(\hat{a}^p) \rangle.
 \]
We substitute in our result for the matrix elements of $\theta^\dagger$ above, and we write this as
 \[
  [^\dagger\!\theta_{\mu, i}]_{qp} = \sum_{q', p'} [\psi_G]^\mi_{qq'} [\theta_{\mu, i}^\dagger]_{q'p'} [\psi_F]_{p'p}
 \]
which is the statement of the lemma.
\end{proof}
We are now ready to prove that an even-handed structure on a collection of semisimple categories is the same thing, up to a global scale factor, as providing each category with a trace. We say that a {\em weighting} on a collection of semisimple categories is an assignment of a nonzero number $e_i \mapsto k_{e_i}$ to each isomorphism class of simple object in each category; in other words a weighting is the same thing as a trace on each category.

\XXX{Decide on the use of this full sub-2category language, here and in previous prop. Here you need all functors and all nat trans basically.}
\begin{prop} There is a canonical bijection
\[
 \left\{ \ba \begin{array}{c} \text{even-handed structures on a}\\ \text{full sub-2-category $\mathcal{S} \subset \SLCat_k$} \end{array} \ea \right\} \cong \left\{ \ba \begin{array}{c} \text{weightings on $\mathcal{S}$, up to} \\ \text{a global scale factor} \end{array} \ea \right\}.
 \]
\end{prop}
\begin{proof}
We begin with the most restrictive requirement on an even-handed structure, that the right and left daggers coincide. Suppose that $F \colon A \rightarrow B$ is a linear direct sum preserving functor between semisimple categories, with right adjoint $F^* \colon B \rightarrow A$ and associated adjunction isomorphisms $\phi_F$, and that $\theta\colon F \Rightarrow F$ is a natural transformation. If $e_i$ and $e_\mu$ index the simple objects in $A$ and $B$, then we will write
 \[
  \phi_{i,\mu} \colon \Hom(Fe_i, e_\mu) \rightarrow \Hom(e_i, F^*e_\mu)
 \]
for the behaviour of $\phi$ on the simple objects (which completely determines it). We showed in Lemma \ref{ambcomplem} that, after choosing a basis, the equation $\theta^\dagger = {}^\dagger\!\theta_\Psi$ amounts to the assertion that for each pair of simple objects $e_i \in A$ and $e_\mu \in B$ the following diagram of matrices commutes:
 \[
   \ba \xymatrix{  k^m \ar[d]_{ [\Psi(\phi)_{i, \mu}]} \ar[r]^-{[\theta^\dagger]} & k^n \ar[d]^{ [\Psi(\phi)_{i, \mu}]} \\
   k^m \ar[r]_-{ [\theta^\dagger]} & k^n } \ea \quad \ba \begin{array}{rcl} m\! \! \! \! &=& \! \! \! \! \dim \Hom(e_\mu, Fe_i) \\ n \! \! \! \! &=& \! \! \! \!\dim \Hom(e_\mu, Ge_i) \end{array} \ea
  \]
Since this must hold for {\em all} functors $F$ and {\em all} natural tranformations $\theta$, we see that $[\Psi(\phi)_{i, \mu}]$ is a natural transformation of the identity on $\Vect$, and hence we must have
 \[
  [\Psi(\phi)_{i, \mu}] = k_{\mu, i} \id
 \]
for some scalar $k_{\mu, i} \in k^\times$. The condition that $\Psi$ is compatible with composition then requires that for any triple of simple objects $e_i, e_\mu, e_\kappa$ inside $\mathcal{C}$ (they can come from different or the same categories), we must have
 \[
  k_{\kappa, \mu} k_{\mu, i} = k_{\kappa, i} \text{ and } k_{i,i} = 1.
 \]
The only solution is that
 \[
  k_{i, \mu} = \frac{k_\mu}{k_i}
 \]
for some assignment of scalars $e_i \mapsto k_i$ to each simple object; clearly this assignment is only determined up to a global scale factor. Note that the condition $\Psi^2 = \id$ is automatically fulfilled.
  \end{proof}
That the weighting is only defined up to a scalar factor ties in well with the geometric analogy we gave in Section \ref{geominterp} between an even-handed structure and a spin structure on a Riemannian manifold, because it has some formal similarities with the fact that the dimension of the space of harmonic spinors on a spin manifold is invariant under changing the scale of the metric \cite{ref:hitchin}.

\section{2-Hilbert spaces\label{EvenTHilb}} In this section we use Proposition \ref{tracepn} to show that the 2-category $\THilb$ of 2-Hilbert spaces comes has a canonical even-handed structure. This is because each 2-Hilbert space comes equipped with a canonical trace,
 \[
  \Tr(f) = (\id, f).
 \]
For instance recall that if the 2-Hilbert space is $\Rep(G)$ for some finite group $G$ then the resultant weightings on the simple objects $V_i$ --- the irreducible representations --- are precisely their dimensions, weighted by the reciprocal of the order of the group.

\begin{prop} \begin{enumerate} \item The 2-category $\THilb$ of 2-Hilbert spaces comes equipped with a canonical even-handed structure $\Psi$.
 \item We can express $\Psi$ in terms of the $*$-structure on the hom-sets and the inner products as
  \[
  \Psi(\phi) = *  \, \, \phi^* *.
 \]
where $\phi^*$ is the Hilbert space adjoint of $\phi$.
\end{enumerate}
\end{prop}
\begin{proof} The first claim follows directly from Proposition \ref{tracepn}. The second claim says that if $F \colon A \rightarrow B$ is a linear $*$-functor between 2-Hilbert spaces, and $F^* \colon B \rightarrow A$ is a right adjoint for $F$ with accompanying adjunction isomorphisms
 \[
  \phi \colon \Hom(Fx, y) \rightarrow \Hom(x, F^*y),
 \]
then the transposed map $\phi^t \colon \Hom(F^*y, x) \rightarrow \Hom(y, Fx)$ can be computed as the composite
 \[
 \xymatrix{\Hom(F^*y, x) \ar[r]^{*} & \Hom(y, Fx) \ar[r]^-{\phi^*} & \Hom(Fx, y) \ar[r]^{*} & \Hom(x, F^*y)}
 \]
where $\phi^* \colon \Hom(x, F^*y) \rightarrow \Hom(Fx, y)$ is the adjoint (in the ordinary sense of maps between Hilbert spaces) of $\phi$. Indeed, we know that for a morphism $g \colon F^*y \rightarrow x$, the transposed morphism $\phi^t(g) \colon x \rightarrow F^*(y)$ is the unique morphism satisfying $\Tr (f \phi^t(g)) = \Tr(\phi(f)g)$ for all $f \colon Fx \rightarrow y$.  But then we must have $\phi^t(g) = [\phi^*(g^*)]^*$, because
 \begin{align*}
  \Tr( f [\phi^*(g^*)]^*) &:= (\id, f [\phi^*(g^*)]^*) \\
  &= (\phi^*(g^*), f) \\
  &= (g^*, \phi(f)) \\
  &= (\id, \phi(f) g ) \\
  &=: \Tr( \phi(f) g).
 \end{align*}
\end{proof}
In the light of the geometric interpretation for even-handed structures in Section \ref{geominterp}, it is interesting to ponder the formal similarity between this formula $\Psi(\phi) = *\, \, \phi^* *$ and the formula for the adjoint of exterior derivative in Riemannian geometry, $d^* = * d *$.

\section{Derived categories of Calabi-Yau manifolds\label{yausec}}
A 2-Hilbert space can be thought of as the category of hermitian vector bundles over a discrete set --- the set of isomorphism classes of simple objects. Instead of vector bundles over a discrete set, one might attempt to consider coherent sheaves over a smooth space. This line of thought is indeed promising; in this section we show that the resulting `smooth' version of the 2-category of 2-Hilbert spaces also carries an even-handed structure.

Let $X$ be an $n$-dimensional compact complex manifold, and let $\mathcal{O}_X$ be its structure sheaf, that is the sheaf of holomorphic functions. Recall that a {\em coherent sheaf} on $X$ is a sheaf of $\mathcal{O}_X$-modules which is locally a quotient of a finite-rank locally-free sheaf. The coherent sheaves form an abelian category and we write the corresponding bounded derived category as $D(X)$. To avoid the nuisance of shift functors, we will rather work with the {\em graded} derived category $\mathbf{D}(X)$, which is defined to have the same objects as $D(X)$ but with hom-sets the graded vector spaces
 \[
  \Hom_{\mathbf{D}(X)} (\mathscr{E}, \mathscr{F}) := \bigoplus_i \Hom_{D(X)} (\mathscr{E}, \mathscr{F} [i]),
 \]
where $\mathscr{F} [i]$ is the complex $\mathscr{F}$ shifted by $i$. Note that the hom-sets in the graded derived category are finite because the complexes are bounded.

If $\omega_X$ is the canonical line bundle of $X$, then Serre duality gives for each coherent sheaf $\mathcal{E} \in \mathbf{D}(X)$ a natural trace map
 \[
  \Tr \colon \Hom_{\mathbf{D}(X)}^n (\mathscr{E}, \omega_X \otimes \mathscr{E}) \rightarrow \mathbb{C}
 \]
which is graded-commutative in the sense that
\[
 \Tr(fg) = (-1)^{\deg f} \Tr(gf)
\]
whenever this is defined, and is nondegenerate in the sense that it sets up natural bifunctorial isomorphisms
 \[
  \Hom^*_{\mathbf{D}(X)}(\mathscr{E}, \mathscr{F}) \stackrel{\cong}{\rightarrow} \Hom^{n- * }_{\mathbf{D}(X)} (\mathscr{F}, \omega_X \otimes \mathscr{E})^\vee.
  \]
In other words, if we restrict ourselves to {\em Calabi-Yau} manifolds (which are characterized as precisely those manifolds for which the canonical line bundle $\omega_X$ is trivial), then Serre duality gives the graded derived category $\mathbf{D}(X)$ a graded trace in the sense we defined in Section \ref{tracelinsec}. To translate this into an even-handed structure, we first need to define the 2-category we will be working in.

Firstly recall that for complex manifolds $X$ and $Y$ there is a way to transform an object $\mathscr{G} \in \mathbf{D}(X)$ into a functor $F \colon \mathbf{D}(X) \rightarrow \mathbf{D}(Y)$. Namely, one considers the diagram of projections
 \[
  \xymatrix{ & X \times Y \ar[dl]_{\pi_X} \ar[dr]^{\pi_Y} \\ X & & Y}
  \]
and then one defines $F$ by pulling $\mathscr{G}$ back to $\mathbf{D}(X \times Y)$, tensoring with $\mathscr{G}$ and then pushing down to $\mathbf{D}(Y)$:
 \[
  F(\mathscr{E}) = R\pi_{Y *}( \pi_{X}^* (\mathscr{E}) \, \underline{\otimes} \, \mathscr{G}).
 \]
A functor $F \colon \mathbf{D}(X) \rightarrow \mathbf{D}(Y)$ of this form is called an {\em integral kernel}.
\begin{defn} The 2-category $\CYau$ is defined as follows. An object is a compact Calabi-Yau manifold $X$. A morphism $F \colon X \rightarrow Y$ is a functor $F \colon \mathbf{D}(X) \rightarrow \mathbf{D}(Y)$ which is naturally isomorphic to an integral kernel. A 2-morphism is a natural transformation.
\end{defn}
By regarding each space $X$ as a placeholder for the derived category $\mathbf{D}(X)$, we can think of $\CYau$ as living inside the 2-category $\LCat$ of linear categories, linear functors and natural transformations. Thus we can apply Proposition \ref{tracepn} and conclude that it comes equipped with a canonical even-handed structure.

\begin{prop} The 2-category $\CYau$ comes equipped with a canonical even-handed structure, arising from the Serre duality on each derived category $\mathbf{D}(X)$.
\end{prop}
\begin{proof} Firstly, every morphism in $\CYau$ has an ambidextrous adjoint, because if $F : \mathbf{D}(X) \rightarrow \mathbf{D}(Y)$ is an integral kernel corresponding to a complex $\mathscr{G} \in \mathbf{D}(X \times Y)$, then we can consider $\mathscr{G}$ as an object in $\mathbf{D}(Y \times X)$ and hence we have the corresponding integral kernel
\[
 F^* \colon \mathbf{D}(Y) \rightarrow \mathbf{D}(X).
 \]
It is well-known (see e.g. \cite{ref:caldararu_willerton}) that if the Serre functors on $X$ and $Y$ are trivial then $F^*$ can be canonically equipped as a simultaneous right and left adjoint for $F$. However, observe that since the morphisms in $\CYau$ are not integral kernels but only functors which are {\em isomorphic} to integral kernels, there is no longer a canonical choice of ambdiextrous adjoint for the morphisms in $\CYau$ --- but at least we are assured that they exist.

Since each derived category $\mathbf{D}(X)$ comes equipped with a graded trace the proposition follows from Proposition \ref{tracepn}.
\end{proof}
In Section \ref{evsec} we saw that an even-handed structure on a collection of semisimple categories (considered as a full sub-2-category of $\SLCat_k$) is the same thing as an assignment of a nonzero number to each simple object in each category, up to a global scale factor. It would  be interesting to calculate the corresponding `moduli space' of even-handed structures on a collection of Calabi-Yau spaces (considered as a full sub-2-category of $\CYau$), and how this would relate to the geometry of the underlying Calabi-Yau spaces.

\chapter{Even-handed structures on fusion categories\label{FusionChap}}
In this chapter we consider the notion of even-handed structures in the realm of {\em fusion categories}, which are semisimple linear monoidal categories where every object has a dual. These categories have been much studied in the field of quantum algebra (see the lecture notes of M\"{u}ger \cite{ref:müger_lectures} for a recent overview). One of the important papers in this regard has been that of Etingof, Nikshych and Ostrik \cite{ref:eno}, where the foundations of the theory were laid out. Our main motivation has been the following conjecture made in \cite{ref:eno}:

\begin{conjecture}[{Etingof, Nikshych and Ostrik} \cite{ref:eno}] Every fusion category admits a pivotal structure.
\end{conjecture}

Recall from Chapter \ref{pivssec} that an even-handed structure on a monoidal category is essentially the same thing as a pivotal structure except that it fits more seamlessly with the string diagram calculus. In this chapter we will do two things:
 \begin{itemize}
  \item We will show how many of the important results of Etingof, Nikshych and Ostrik have simple and elegant proofs in terms of string diagrams. This is significant because these results were originally derived using algebraic methods from the theory of weak Hopf algebras.
  \item We will use the framework of even-handed structures to identify a pivotal structure on a fusion category as a {\em twisted monoidal natural transformation of the identity functor} on the category. The twisting is governed by a set of signs $\epsilon^i_{jk}$ associated to each triple of simple objects in the category, which arise from the nature of the duality operation in the category. We call these signs the {\em pivotal symbols} since they are analogous to the {\em 6j symbols} \cite{ref:carter_flath_saito} which govern the associator information in a fusion category.
 \end{itemize}
 \XXX{Don't forget to run through every argument that uses the old isomorphism naturality condition}
The starting point of our translation of the results of Etingof, Nikshych and Ostrik into a graphical framework arises from an observation of Hagge and Hong \cite{ref:hagge_hong}. Hagge and Hong showed that one of the fundamental results of \cite{ref:eno} --- that there is a more-or-less canonical monoidal natural isomorphism $\gamma \colon \id \Rightarrow \star \star \star \star$ between the identify functor and the fourth power of the dualizing functor `$\star$' on a fusion category --- can be given a purely string diagrammatic proof.

Hagge and Hong worked in the framework of strictified skeletal fusion categories and pivotal structures and appealed to coherence theorems, and our approach is different in two ways. Firstly we saw in Chapter \ref{stchap} that {\em string diagrams can be applied directly to the category at hand} without having to pass to a strictified situation. Secondly we work in the framework of even-handed structures, whereby the isomorphism $\gamma \colon \id \Rightarrow \star \star \star \star$ translates into an {\em involution} $T^i_{jk}$ on the basic hom-sets $\Hom(X_i, X_j \otimes X_k)$ between the simple objects. The eigenvalues of $T^i_{jk}$ are therefore $\pm 1$, and we call this collection of signs the {\em pivotal symbols} of the fusion category. We will prove that the pivotal symbols control the existence of pivotal structures on the fusion category.

The layout of this chapter is as follows. We begin in Section \ref{wex} with a warm up exercise: we compute the set of even-handed structures on a group-like monoidal category $(G, \omega)$ associated to a group $G$ and a 3-cocycle $\omega \in Z^3(G, U(1))$. Our computation of this toy example is not new (it appears for instance in the lecture notes of Boyarchenko \cite{ref:boyarchenko}), but it neatly illustrates the main ideas.

In Section \ref{stss} we review the basic facts about fusion categories, laying emphasis on the manifest algebraic information contained in a fusion category: the {\em fusion ring} and the more subtle information of the {\em paired dimensions} of the simple objects first explicated by M\"{u}ger \cite{ref:müger_from_subfactors_to_categories_and_topology_I}. In Section \ref{epiv} we use the graphical calculus and the idea of Hagge and Hong to introduce the involution operators $T^i_{jk}$, the eigenvalues of which are recorded in the {\em pivotal symbols}. In Section \ref{Factspivtensor} we give string diagram proofs of many of the main results of \cite{ref:eno}. Finally in Section \ref{existence} we prove our main result identifying a pivotal structure on a fusion category as a twisted monoidal natural transformation of the identity functor on the category. We show that unless each involution operator $T^i_{jk}$ equals plus or minus the identity map, then the fusion category {\em cannot} carry a pivotal structure; moreover we show that the pivotal structure can be made spherical if and only if these signs can be consistently removed.

\section{Warm-up exercise\label{wex}}
To prepare ourselves for calculating the even-handed structure on fusion categories, we start with a warm-up exercise. Recall that if $G$ is a group and $\omega \in Z^3(G, U(1))$ a group 3-cocycle, then we can form the monoidal category $(G, \omega)$ whose objects are the group elements $g \in G$, whose hom-sets $\Hom(g,h)$ are empty unless $g=h$ when they are equal to $U(1)$, and whose associator is given by the cocycle $\omega$. The following result is well-known (see eg. \cite{ref:boyarchenko}) in the context of pivotal structures.
\begin{lem}\label{glem} If $G$ is a group and $\omega \in Z^3(G, U(1))$ is a normalized 3-cocycle, then there is a canonical bijection
\[
\left\{ \ba \begin{array}{c} \text{even-handed structures} \\ \text{on $(G, \omega)$} \end{array} \ea \right\} \cong \left\{ \ba \begin{array}{c} \text{group homomorphisms} \\ \text{$f \colon G \rightarrow U(1)$} \end{array} \ea \right\}.
\]
Moreover the even-handed structure is spherical precisely when $f(g) = \pm 1$ for all $g \in G$.
\end{lem}
\begin{proof} Firstly we need to check that $(G, \omega)$ has ambidextrous duals; this is the only point where we need to perform a calculation directly involving $\omega$. Given a group element $g \in G$, let us write $\omega_g \equiv \omega(g, g^\mi, g)$; the 3-cocycle equation together with the fact that $\omega$ is normalized implies that $\omega_{g^\mi} = \frac{1}{\omega_g}$. If we trace through the snake diagrams we see that $(g^\mi, \eta_g, \epsilon_g)$ is a right dual of $g$ if and only if $\epsilon_g \eta_g = \omega_g$. In other words, the set of ways $\Adj(g \dashv g^\mi)$ in which $g^\mi$ can be expressed as a right dual of $g$ forms a $U(1)$-torsor, for we can choose $\epsilon_g \in U(1)$ freely after which we must set $\eta_g = \omega_g / \eta_g$.

Similarly, tracing through the snake diagrams one finds that $(g^\mi, \eta'_g, \epsilon'_g)$ is a {\em left} dual of $g$ if and only if $\epsilon'_g \eta'_g  = 1/\omega_g$, and the set of solutions $\Adj(g^\mi \dashv g)$ to this equation is again a $U(1)$-torsor.

Since $\Adj(g \dashv g^\mi)$ and $\Adj(g^\mi \dashv g)$ are nonempty for all $g \in G$, we conclude that $(G, \omega)$ indeed has ambidextrous duals.  Having established this fact, we no longer need to directly involve $\omega$ in our calculations at all. To see this, observe that an even-handed structure on $(G, \omega)$ consists of a bijective map
 \[
  \Psi \colon \Adj(g \dashv g^\mi) \rightarrow \Adj(g^\mi \dashv g)
 \]
for every $g \in G$, satisfying the axioms we listed in Chapter \ref{defesec}. Without knowing any of the `internal details' of these $U(1)$-torsors, it is clear that such a map $\Psi$ is uniquely determined by the `dimension function'
 \begin{align*}
  f \colon G & \rightarrow U(1)\\
  \intertext{which sends}
  g & \mapsto \; \ba \ig{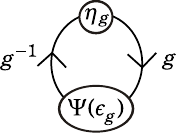} \ea
 \end{align*}
where $(\eta_g, \epsilon_g) \in \Adj(g \dashv g^\mi)$ is some choice of unit and counit maps which exhibit $g^\mi$ as a right dual of $g$; the resulting number $f(g)$ does not depend on this choice due to the naturality properties of $\Psi$ from Lemma \ref{narlem}.
Let us examine the requirements of the axioms of an even-handed structure on the function $f \colon G \rightarrow U(1)$. Although in this simple setting the monoidal axiom implies all the others, it is instructive to understand the requirements of each axiom in turn. The requirement that $\Psi^2 = \id$ implies that
 \[
  \ba \ig{y133.pdf} \ea \ba \ig{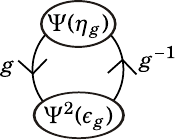} \ea = \ba \ig{y133.pdf} \ea \ba \ig{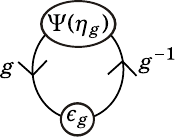} \ea .
 \]
We recognize the left hand side as $f(g) f(g^\mi)$. As for the right hand side, despite appearances it is actually a straightforward product of four numbers in $U(1)$, which we can parenthesize as $(\epsilon_g \eta_g) (\Psi(\epsilon_g) \Psi(\eta_g)) = \omega_g \, 1/\omega_g = 1$. (We remark that this calculation introduces the concept of {\em paired dimensions} which we will review in the next section). In other words,
 \[
  \Psi^2 = \id \quad \Leftrightarrow \quad f(g)f(g^\mi) = 1.
 \]
Similarly the requirement that $\Psi$ respects the monoidal structure means that
 \[
f(gh) = f(g) f(h),
 \]
and the requirement that $\Psi$ acts as the identity on the trivial adjunction means that $f(e) = 1$. Finally, the even-handed equation is automatically satisfied, since we can take scalar factors out of the diagrams and place them in front, so that we always have $\theta^\dagger = {}^\dagger\!\theta_\Psi$ for every morphism $\theta$. We conclude that an even-handed structure on $(G, \omega)$ is the same thing as a group homomorphism $f \colon G \rightarrow U(1)$.

Recall from Chapter \ref{defesec} that $\Psi$ is called {\em spherical} (as in the work of Barrett and Westbury \cite{ref:barrett_westbury}) if the left and right traces coincide, which in our case is the requirement that
 \[
 \ba \ig{y133.pdf} \ea = \ba \ig{y134.pdf} \ea.
\]
We recognize this as the equation $f(g) =f(g^\mi)$ for all $g \in G$, which means that $f(g)$ must take values in $\{1, -1\}$. This completes the proof.
\end{proof}
This example above neatly illustrates the kind of ideas involved in the notion of an `even-handed structure' on a monoidal category.

\section{Fusion categories\label{stss}}
In this section we begin our study of even-handed structures on fusion categories, by summarizing the basic facts about them which will be relevant for us. We will always work over the complex numbers $\mathbb{C}$. A {\em fusion category} $C$ is a semisimple linear rigid monoidal category with only finitely many isomorphism classes of simple objects indexed by a set $I$; moreover the unit object $1 \in C$ is required to be simple. Note that we will often assume that representative simple objects $\{X_i\}_{i \in I}$ have been chosen. The unit object takes the index zero, so that $X_0 = 1$.
Also, recall that rigidity is only the requirement that every object $V$ has a {\em right} dual $V^*$; however semisimplicity implies that $V^*$ is also a {\em left} dual of $V$, which we prove in Appendix \ref{AppendixFusion}. So a fusion category has ambidextrous duals; it therefore makes sense to talk about an even-handed structure on a fusion category.

We will write the collection of even-handed structures on a fusion category $C$ as $\EvenHanded(C)$, and we will write $\Aut_\otimes (\id)$ for the group of monoidal natural isomorphisms of the identity on $C$. The following is a very important fact to keep in mind.
\begin{lem} In the event that $\EvenHanded(C)$ is nonempty, it is acted on freely and transitively by $\Aut_\otimes(\id_C)$ via the formula
 \[
  (\delta \cdot \Psi) \left( \ba \ig{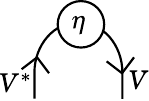} \ea, \ba \ig{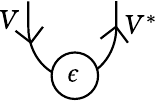} \ea \right) = \left( \ba \ig{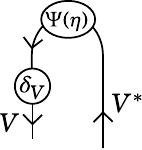} \ea, \ba \ig{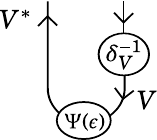} \ea \right).
  \]
\end{lem}
\begin{proof} Firstly, $\delta \cdot \Psi$ is indeed an even-handed structure on $C$ since the monoidal axiom and the even-handed axiom are clearly satisfied, while the following establishes that $(\delta \cdot \Psi)^2 = \id$:
 \[
 (\delta \cdot \Psi)^2 \left( \ba \ig{e340.pdf} \ea \right) =  \ba \ig{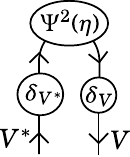} \ea = \ba \ig{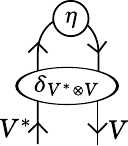} \ea = \ba \ig{e340.pdf} \ea.
 \]
In the last step above we used the naturality of $\delta$, and the fact that it is a {\em monoidal} natural isomorphism so that $\delta_1 = \id$.

Secondly, $\Aut_\otimes (\id_C)$ acts freely on $\EvenHanded(C)$, because if $\delta \cdot \Psi = \Psi$ then
 \[
  \ba \ig{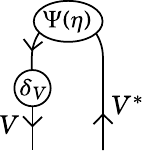} \ea = \ba \ig{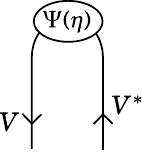} \ea \qquad \Rightarrow \qquad \ba \ig{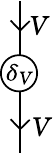} \ea = \ba \ig{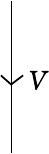} \ea.
\]
Finally, $\Aut_\otimes(\id_C)$ acts transitively, because if $\Psi'$ is another even-handed structure then we can define a monoidal natural transformation of the identity $\delta$ by setting
 \[
  \delta_V = \ba \ig{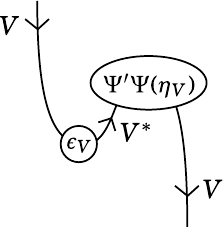} \ea
 \]
where $(V^*, \eta_V, \epsilon_V)$ is any choice of right dual for $V$; the formula is independent of this choice by the naturality properties of $\Psi$ expressed in Lemma \ref{narlem}. To check that $\delta$ is indeed a natural transformation, suppose that $\theta \colon V \rightarrow W$ is a morphism. Then we can use the even-handed equation twice, once for $\Psi$ and once for $\Psi'$, as follows:
 \begin{align*}
  \ba \ig{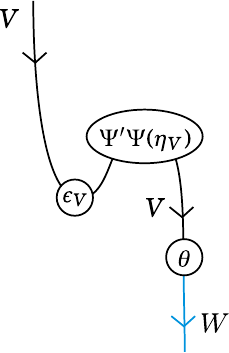} \ea =  \ba \ig{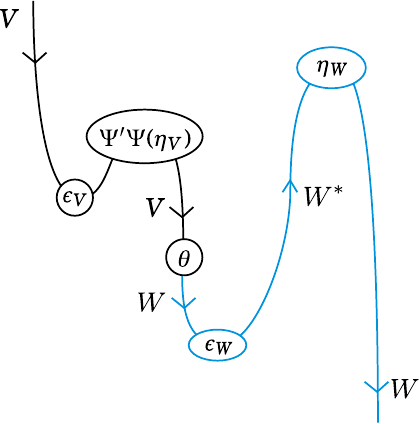} \ea  &= \ba \ig{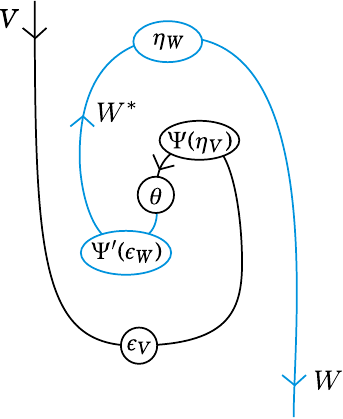} \ea \\  &= \ba \ig{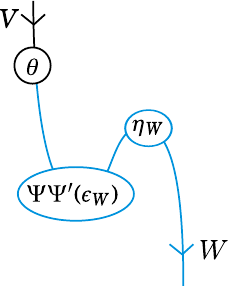} \ea.
\end{align*}
We recognize the right hand side as $\delta_W \circ \theta$ by setting $(\eta_W, \epsilon_W) = (\Psi' \Psi(\eta_W), \Psi' \Psi (\epsilon_W))$ and then using $\Psi'^2 = \id$. This establishes that $\delta$ is a natural transformation; it is clear that $\delta$ is a {\em monoidal} natural transformation because $\Psi$ respects the monoidal structure. This completes the proof.
\end{proof}
\noindent This means that every even-handed structure on a category can be obtained by acting via an automorphism of the identity on a given one. Combining this with Lemma \ref{natbugger}, we have the following corollary.
\begin{cor}\label{autcor} If $C$ is equipped with a braiding and a compatible twist (in particular, if $C$ is a symmetric monoidal category), then $\EvenHanded(C) \cong \Aut_\otimes(\id_C)$.
\end{cor}
There are two pieces of algebraic data one can manifestly extract from a fusion category $C$. The first is the {\em Grothendieck ring} which we will simply write as $[C]$ (instead of the more usual $K(C)$), defined as the ring of isomorphism classes of objects. This is not an arbitrary ring; firstly it has a distinguished basis $[X_i]$ given by the simple objects in which the expansion coefficients $[X_j] [X_k] = \sum_{i} N^i_{jk} [X_i]$ are nonnegative integers, and moreover  the Grothendieck ring has an anti-involution $* \colon [C] \rightarrow [C]$ satisfying $N^0_{jk} = \delta_{jk^*}$, where `$0$' refers to the index of the unit object. A ring with these properties is called a {\em fusion ring} \cite{ref:eno, ref:calaque_etingof}.

The other piece of manifest algebraic data one can extract from a fusion category is a bit more subtle and was first explicated by M\"{u}ger \cite{ref:müger_from_subfactors_to_categories_and_topology_I} --- to each set $\{X_i, X_i^*\}$ consisting of a simple object and its dual, we can assign the product of numbers
 \[
  d_{\{i, i^*\}} = \ba \ig{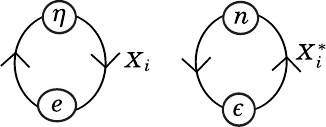} \ea,
 \]
where $(\eta, \epsilon)$ is some choice of unit and counit maps exhibiting $X_i^*$ as a right dual of $X_i$, and $(n, e)$ is some choice of unit and counit maps exhibiting $X_i^*$ as a {\em left} dual of $X_i$. We call $d_{\{i, i^*\}}$ the {\em paired norm} of $X_i$ and $X_i^*$. Observe that this product is independent of the choices we made because it is invariant under rescaling
 \[
 \eta \mapsto \lambda \eta, \epsilon \mapsto \frac{1}{\lambda} \epsilon, n \mapsto \mu n, e \mapsto \frac{1}{\mu} e.
\]
As an example, consider the group-like monoidal category $(G, \omega)$ from Section \ref{wex} as a fusion category, by allowing formal direct sums of group elements. Then we saw in the proof of Lemma \ref{glem} that the paired dimensions $d_{\{g, g^\mi\}}$ equal unity for all $g \in G$.

We will show in Section \ref{Factspivtensor} that the paired dimensions $d_{\{i, i^*\}}$ in a fusion category are positive and real; however for now observe that they are certainly nonzero since, for example, $\eta$ is zero if and only if $\epsilon$ is zero by semisimplicity, because $\Hom(1, X_i^* \otimes X_i)$ is in duality with $\Hom(X_i \otimes X_i^*, 1)$.

An even-handed structure $\Psi$ on a fusion category $C$ gives rise to an associated {\em dimension homomorphism} from the Grothendieck ring into $\mathbb{C}$,
 \begin{align*}
  \dim_{\Psi} & \colon K(C) \rightarrow \mathbb{C}, \\
  \intertext{defined by sending}
   [V] & \mapsto \ba \ig{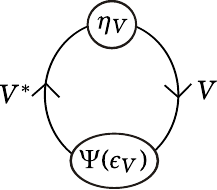} \ea
 \end{align*}
where $(V^*, \eta, \epsilon)$ is any choice of right dual for $V$. The fact that $\dim_\Psi$ is well-defined is a consequence of the naturality of $\Psi$.

\begin{lem} The dimension homomorphism $\dim_{\Psi} \colon K(C) \rightarrow \mathbb{C}$ associated to an even-handed structure $\Psi$ on a fusion category $C$ satisfies the following properties:
 \begin{itemize}
  \item It is a ring homomorphism,
  \item $\dim [X_i] \dim [X_i^*] = d_{\{i,i^*\}}$ for all simple objects $X_i$.
 \end{itemize}
\end{lem}
\begin{proof}
We have
 \begin{align*}
  \dim [V] \dim [W] &= \ba \ig{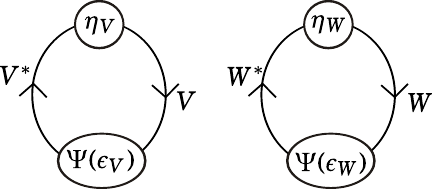} \ea \\
  &= \ba \ig{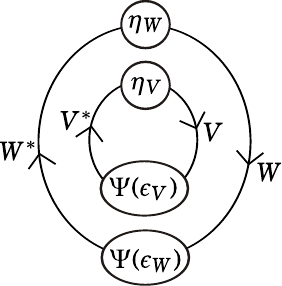} \ea = \ba \ig{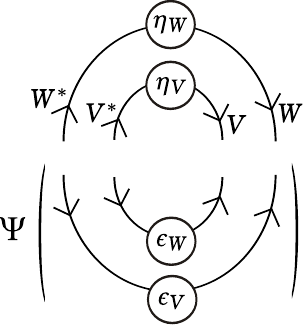} \ea  = \dim [V \otimes W].
  \end{align*}
(Note that the second equality above of `sliding a circle into another circle' does not use any properties of $\Psi$ but simply the fact that a morphism $1 \rightarrow 1$ is just a number.) To evaluate $\dim[X_i] \dim[X_i^*]$, make a choice $(\eta, \epsilon)$ of unit and counit maps exhibiting $X_i^*$ as a right dual of $X_i$, and then choose $(\Psi(\eta), \Psi(\epsilon))$ as the unit and counit maps exhibiting $X_i$ as a right dual of $X_i^*$. Then we have
 \[
 \dim[X_i] \dim[X_i^*] = \ba \ig{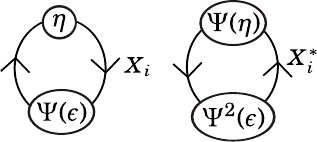} \ea = \ba \ig{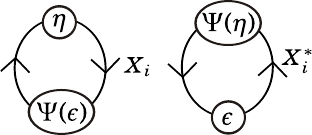} \ea = d_{\{i,i^*\}}.
 \]
\end{proof}
We call a function $f \colon K(C) \rightarrow \mathbb{C}$ from the Grothendieck ring of a fusion category to the ground field a {\em fusion homomorphism} if it satisfies the properties above. By analogy with the case of $(G, \omega)$, one might imagine that
 \[
  \operatorname{Even-Handed}(C) \cong \{\text{Fusion homomorphisms } f \colon K(C) \rightarrow \mathbb{C} \}.
 \]
Unfortunately we have not been able to obtain such a result. Instead, we will characterize the collection of even-handed structures on $C$ as the collection of {\em twisted} monoidal natural transformations of the identity, a statement we will explain below. However it is instructive to see that this characterization {\em does} hold for fusion categories of the form $\Rep(G)$ for a finite group $G$.
\begin{prop} The even-handed structures on $\Rep(G)$ are in 1-1 correspondence with the fusion homomorphisms from $[\Rep(G)]$ to $\mathbb{C}$.
\end{prop}
\begin{proof} We have
 \begin{align*}
  \operatorname{Even-Handed}(\Rep(G)) & \cong \Aut_\otimes(\id) \\
   & \cong Z(G) \\
   & = \{ g \in G \colon |\Tr_{V_i} (g) | = \dim V_i \text{ for all irreducibles $V_i$}\} \\
   & \cong \{\text{Fusion homomorphisms } f \colon [\Rep(G)] \rightarrow \mathbb{C} \}.
  \end{align*}
The first isomorphism uses the fact that $\Rep(G)$ is a symmetric monoidal category and Corollary \ref{autcor}. The second isomorphism is a result of M\"{u}ger \cite{ref:müger_center_group}. The equality in the third line is a basic result of representation theory (see Corollary 2.28 of \cite{ref:isaacs}). The final isomorphism uses two facts. Firstly, every ring homomorphism
 \[
  [\Rep(G)] \rightarrow \mathbb{C}
 \]
must take the form $V \mapsto \Tr_V(g)$ for some fixed $g \in G$ --- because a character like this certainly {\em is} a ring homomorphism, and there are as many such distinct characters as there are conjugacy classes in the group, which must exhaust all the ring homomorphisms since $[\Rep(G)]_\mathbb{C}$ is isomorphic to the space of functions on the conjugacy classes. Secondly, in $\Rep(G)$ the paired dimensions $d_{\{i,i^*\}}$ are just $\dim(V_i)^2$, so that a fusion homomorphism must satisfy $| f(V_i) | = \dim V_i$.
\end{proof}

\section{The pivotal symbols of a fusion category\label{epiv}}
In this section we use string diagrams to show that every fusion category comes equipped with involutions $T^i_{jk} \colon \Hom(X_i, X_j \otimes X_k) \rightarrow \Hom(X_i, X_j \otimes X_k)$ on the hom-sets between a single simple object and a tensor product of two other ones. Since $T^i_{jk}$ is an involution, its eigenvalues must equal $\pm 1$, and we call this collection of signs the {\em pivotal symbols} of the fusion category. They play a central role in the rest of this chapter.

We define a {\em root choice} on a fusion category as a symmetric choice $\{d_i\}$ of square roots of the paired dimensions, that is one which satisfies
 \[
  d_i^2 = d_{\{i,i^*\}} \text{ and } d_i = d_{i^*} \text{ for all $i$.}
 \]
Although we will show in Proposition \ref{realpositive} that the paired dimensions $d_{\{i,i^*\}}$ are positive and real, all that is necessary in the above definition is that they are nonzero.  The point about a root choice is that it allows one to establish the following convention.
\vskip 0.3cm
\framebox{\parbox[b]{12cm}{\textbf{Partnering Convention.} Whenever a unit and counit $\left(\ba \ig{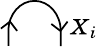} \ea, \ba \ig{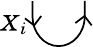} \ea \right)$ expressing $X_i^*$ as a {\em right} dual of $X_i$ appears together in some equation or statement with a unit and counit $\left(\ba \ig{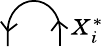} \ea, \ba \ig{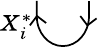} \ea \right)$ expressing $X_i^*$ as a {\em left} dual of $X_i$, it will always be understood that the former is arbitrary while the latter is determined uniquely by the requirement that
 \[
  \ba \ig{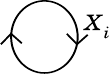} \ea = d_i \quad \text{(or equivalently } \ba \ig{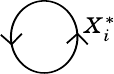} \ea = d_i \text{).}
 \]}}
\vskip 0.3cm
A root choice can thus be seen as a kind of `proto-even-handed structure' --- the game is to try and modify it to obtain a true even-handed structure.

Suppose that a root choice $\{d_i\}$ has been made. Using the Partner Convention, we can define the involution operators $T^i_{jk}$.
\begin{defn} The operators $T^i_{jk} \colon \Hom(X_i, X_j \otimes X_k) \rightarrow \Hom(X_i, X_j \otimes X_k)$ are defined by
 \[
  \ba \ig{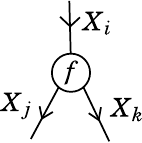} \ea \mapsto \ba \ig{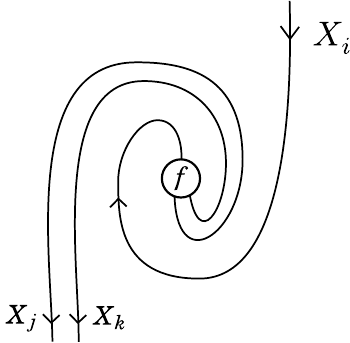} \ea.
 \]
\end{defn}
One might imagine that $T^i_{jk}$ is the identity operator, but that need not be the case --- observe for instance that changing the signs of the root choice $d_i \mapsto -d_i$ will have the effect of changing $T^i_{jk} \mapsto - T^i_{jk}$ (providing $j$ and $k$ are distinct from $i$ and $i^*$). If we choose a basis $\{a_p \colon X_i \rightarrow X_j \otimes X_k\}$ with corresponding dual basis $\{a^q \colon X_j \otimes X_k \rightarrow X_i\}$, we see that the matrix elements of $T^i_{jk}$ compute as
 \[
 a^q T^i_{jk}  a_p = \frac{1}{\ba \ig{e288.pdf} \ea} \ba \ig{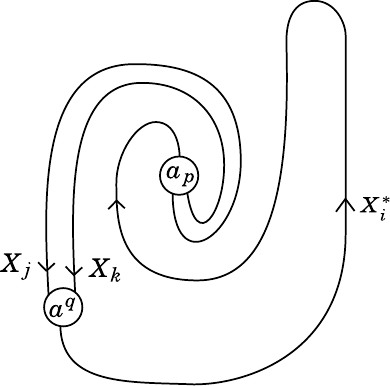} \ea = \frac{1}{\ba \ig{e288.pdf} \ea} \ba \ig{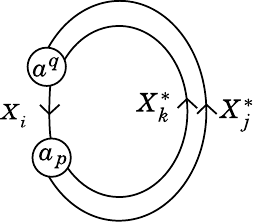} \ea
 \]
but we are unable to continue the calculation further. We record for later use here however one of the important properties of the Partner Convention --- it allows us to `flip over' closed loops, for instance:
\[
 \ba \ig{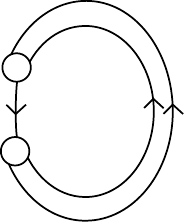} \ea = \ba \ig{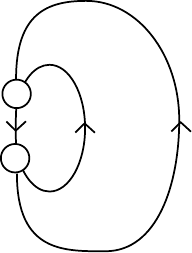} \ea = \ba \ig{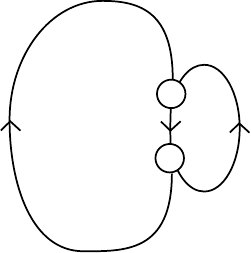} \ea
 \]
because
 \[
  \ba \ig{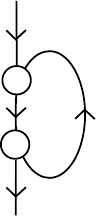} \ea = \lambda \ba \ig{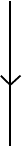} \ea \text{ for some $\lambda \in \mathbb{C}$, and } \ba \ig{e288.pdf} \ea = \ba \ig{e289.pdf} \ea.
 \]
To summarize: although we cannot show that $T^i_{jk}$ is the identity, one can show that its {\em square} is the identity. We credit the following graphical proof of this fact to Hagge and Hong \cite{ref:hagge_hong}. Our framework differs significantly from theirs: we do not work in a system of fixed chosen duals --- if we did, we could not conceive of $T^i_{jk}$ as an {\em involution}; this concept would translate into the quadruple dual $\ast \ast \ast \ast$ being monoidally naturally isomorphic to the identity functor. Moreover, we do not need to make assumptions about the fusion category being skeletal and strictified. Nevertheless the basic idea is the same.

\begin{lem}[{see Hagge and Hong \cite[Thm 3]{ref:hagge_hong}}] The operator $T^i_{jk}$ is an involution --- that is, $(T^i_{jk})^2 = \id$.
\end{lem}
\begin{proof} The operator $T^2$ sends
 \[
  \ba \ig{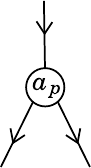} \ea \mapsto \ba \ig{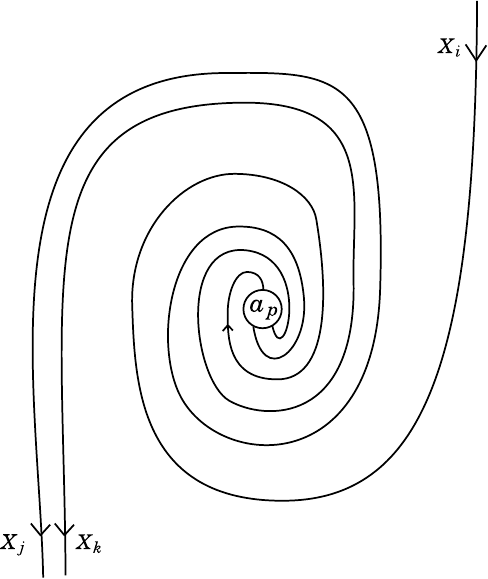} \ea .
 \]
Its matrix elements are thus:
\[
\ba \ig{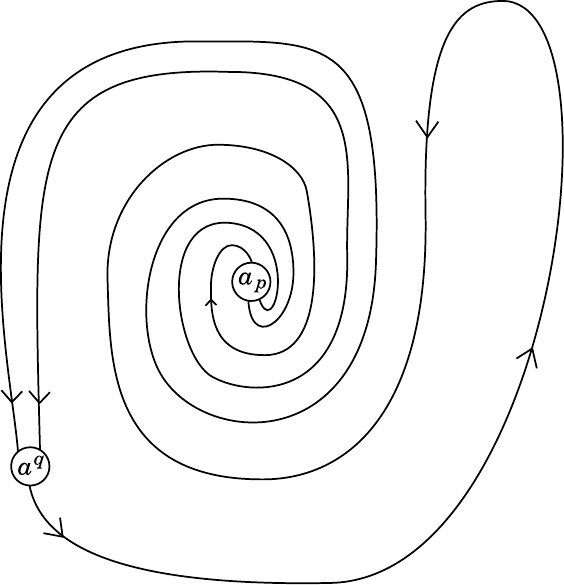} \ea = \ba \ig{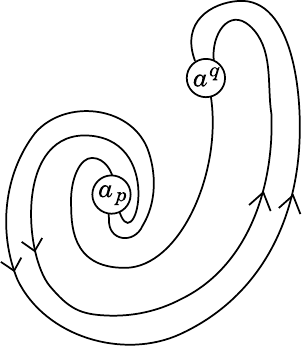} \ea
\]
\[
= \ba \ig{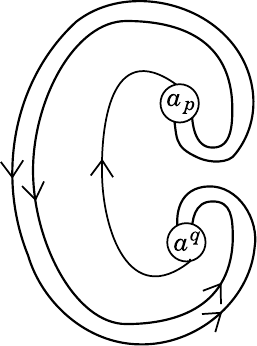} \ea = \ba \ig{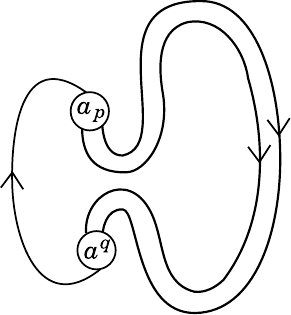} \ea = \ba \ig{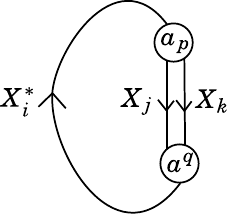} \ea = \delta^q_p.
\]
\end{proof}
This situation is somewhat reminiscent of the {\em Dirac belt trick} which proves that performing two full rotations is homotopic to the constant identity rotation. We acknowledge here that we have no good conceptual understanding of the involution operators $T^i_{jk}$. However, if the fusion category is the category of representations of a Hopf algebra, then the identity $(T^i_{jk})^2 = \id$ corresponds to {\em Radford's formula} for the fourth power of the antipode map (see \cite{ref:eno2} and the references therein). In fact, the string diagram argument above can be regarded as giving a graphical proof of Radford's formula.

In any event, since $T^i_{jk}$ is an involution, the vector space $\Hom(X_i, X_j \otimes X_k)$ has a basis of eigenvectors $a_p$ whose eigenvalues are $\pm 1$. We record these signs in the {\em pivotal symbols} $\epsilon^i_{jk, p}$:
 \[
  T^i_{jk} a_p = \epsilon^i_{jk, p} a_p.
 \]
We will see later that there is no hope of a pivotal structure on the category unless $T^i_{jk} = \pm \id$, so $\epsilon^i_{jk,p}$ will be required to be independent of $p$, whence we will write it as $\epsilon^i_{jk}$. We acknowledge that in fact this might {\em always} be the case, but we have been unable to show this. In any event, it is is this latter collection of signs that we {\em really} view as the pivotal symbols.

The pivotal symbols are analogous mathematical objects to the {\em 6j-symbols} \cite{ref:carter_flath_saito} which reflect the data of the associator. In the standard example of $\Rep(SL(2))$, the 6j-symbols consist of numbers
 \[
  \left\{ \begin{array}{ccc} i j k \\ p q r \end{array} \right\} \in \mathbb{C},
 \]
assigned to each sextuple of irreducible representations which reflect the way the tensor product of three irreducible representations can be decomposed in two natural ways. In a general fusion category, where the fusion rules are more complex than that of $\Rep(SL(2))$, the 6j-symbols become matrices instead of numbers, but the idea is the same.

\section{Facts about the pivotal tensor\label{Factspivtensor}}
In this section we use the graphical calculus to give alternative proofs of some of the main results of Etingof, Nikshych and Ostrik \cite{ref:eno}. We hope that this viewpoint offers a conceptually simpler and more streamlined way to understand them.

Our starting point is a fusion category $C$. We fix a root choice $\{d_i^2 = d_{i,i^*}\}$ of the paired dimensions, and we will use the resulting Partner Convention repeatedly. Moreover, whenever we write basis vectors $\{a_p \colon X_i \rightarrow X_j \otimes X_k\}$ for $\Hom(X_i, X_j \otimes X_k)$ it will always be understood that they an eigenbasis for the operators $T^i_{jk}$, with eigenvalues $\epsilon^i_{jk,p}$.

We start with the following lemma.
\begin{lem}\label{discomb} We have
 \[
  \ba \ig{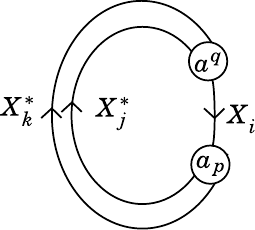} \ea = \epsilon^i_{jk, p}\ba \ig{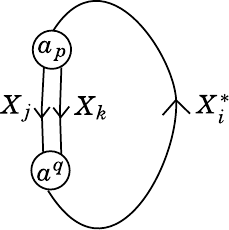} \ea \quad \left(= \delta^q_p \epsilon^i_{jk, p} d_i \right).
  \]
\end{lem}
\begin{proof}
\[
 \ba \ig{e209.pdf} \ea = \epsilon^i_{jk, p} \ba \ig{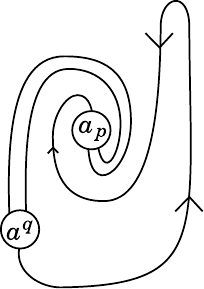} \ea = \epsilon^i_{jk, p} \ba \ig{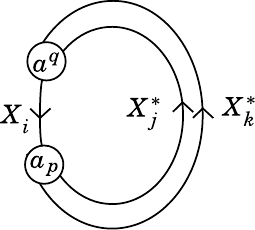} \ea = \epsilon^i_{jk, p} \ba \ig{e208.pdf} \ea.
\]
\end{proof}
The next proposition is very important: it establishes the {\em fusion rules} for the chosen roots $d_i$ in terms of the pivotal tensor. This allows the numbers $d_i$ to be compared with the {\em Frobenius-Perron} dimensions $d_+(X_i)$ of the simple objects, also introduced in \cite{ref:eno}, which can be defined as the unique homomorphism from the Grothendieck ring into the real numbers
 \[
 d_+ \colon [C] \rightarrow \mathbb{C}; \quad \text{i.e. } d_+(X_j) d_+(X_k) = \sum_i N^i_{jk} d_+(X_i)
\]
which takes positive values on the simple objects. The Frobenius-Perron dimension $d_+(X_i)$ of a simple object $X_i$ can also be characterized as the unique positive eigenvalue of the positive integer matrix representing left multiplication by $X_i$.

\begin{prop}\label{homprop} We have $d_j d_k = \displaystyle \sum_i \Tr (T^i_{jk}) \, d_i$.
\end{prop}
\begin{proof} \begin{align*}
 d_j d_k &= \ba \ig{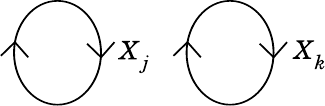} \ea \\
  &= \ba \ig{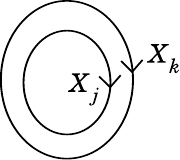} \ea \\
  &= \sum_{i,p} \ba \ig{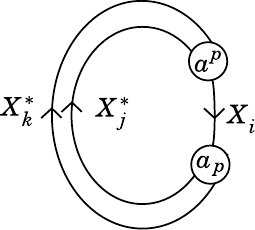} \ea \\
  &= \sum_{i,p} \epsilon^i_{jk,p} d_i \quad \text{(by Lemma \ref{discomb})} \\
  &= \sum_i \Tr(T^i_{jk}) \, d_i.
  \end{align*}
\end{proof}
In the previous section we emphasized the analogy between the pivotal symbols and the 6j-symbols in a fusion category. Now, it is known that the 6j-symbols have the following {\em symmetry properties} \cite[Lemma 2.7.8]{ref:carter_flath_saito}:
 \[
  \left\{ \begin{array}{ccc} i j k \\ p q r \end{array}\right\} = \left\{ \begin{array}{ccc} p q k \\ i j r \end{array}\right\}.
 \]
The next lemma establishes the corresponding symmetry properties for the pivotal symbols.
\begin{lem}\label{Tsymprop} The numbers $\Tr(T^i_{jk})$ have the following symmetry properties:
 \begin{enumerate}
  \item $\Tr(T^i_{jk}) = \Tr(T^{k^*}_{i^*j}) \quad$ \text{(Conjugate cyclic)} \\
  \item $\Tr(T^i_{jk}) = \Tr(T^{i^*}_{k^*j^*}) \quad$ \text{(Conjugate symmetric)}
 \end{enumerate}
 \end{lem}
 \begin{proof} The origin of these symmetries is that the presence of duals allows one to convert eigenvectors of the $T^i_{jk}$ operators amongst
 themselves, while preserving their eigenvalues.  To establish (i), suppose that $\{ a_p \colon X_i \rightarrow X_j \otimes X_k \}$ is an eigenbasis for the operator $T^i_{jk}$, so that $T^i_{jk} a_p = \epsilon^i_{jk, p} a_p$. Then the basis for $\Hom(X_k^*, X_i^*, X_j)$ given by
  \[
   b_p := \ba \ig{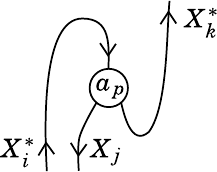} \ea
  \]
are eigenvectors of $T^{k^*}_{i^*, j}$ with the same eigenvalues as $a_p$, since
  \[
  T^{k^*}_{i^*, j} (b_p) = \ba \ig{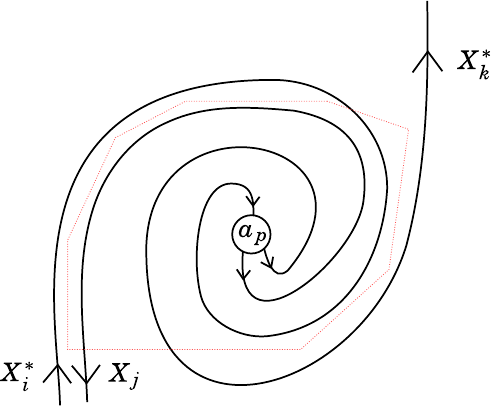} \ea = \epsilon^i_{jk, p} \ba \ig{e216.pdf} \ea = \epsilon^i_{jk,p} b_p.
  \]
The proof of (ii) is similar, except one uses the {\em dual} basis $\{a^p \colon X_j \otimes X_k \rightarrow X_i\}$ to define a basis $\{c^p\}$ for $\Hom(X_i^*, X_k^* \otimes X_j^*)$ via
 \[
  c^p := \ba \ig{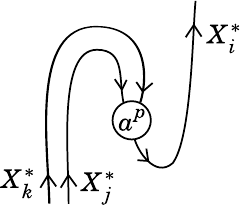} \ea.
 \]
A similar string diagram argument then establishes that $T^{i^*}_{k^*j^*} c^p = \epsilon^i_{j,k} c^p$.
 \end{proof}
We can use all these results to prove the following important fact.
\begin{prop}\label{realpositive} The paired dimensions $d_{\{i, i^*\}}$ are real and positive.
\end{prop}
\begin{proof} We organize the various roots $d_i \equiv \ba \ig{e288.pdf} \ea$ into a column vector:
 \[
   \mathbf{d} = \left( \ba \ig{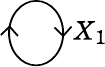} \ea, \ba \ig{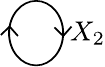} \ea, \ldots, \ba \ig{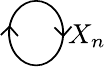} \ea \right)^T.
 \]
Define the matrices $A_j$ via $[A_j]_k^i = \Tr(T^i_{jk})$. Then Proposition \ref{homprop} says that $\mathbf{d}$ is a simultaneous eigenvector of each $A_j$ with eigenvalue $d_j$, i.e.
 \[
  A_j \mathbf{d} = d_j \mathbf{d}.
 \]
Thus we have $A_j A_{j^*} \mathbf{d} = d_j d_{j^*} \mathbf{d} = d_{\{j,j^*\}} \mathbf{d}$ and so $d_{\{j,j^*\}}$ is an eigenvalue of the matrix $A_j A_{j^*}$. But $A_{j^*} = A_j^T$, by the symmetry properties established in Lemma \ref{Tsymprop}:
 \begin{align*}
  [A_{j^*}]^i_k &= \Tr(T^i_{j^* k}) \\
   &= \Tr(T^{k^*}_{i^* j^*})  \\
   &= \Tr(T^k_{ji}) \\
   &= [A_j]^k_i.
  \end{align*}
Thus $d_{\{j,j^*\}}$ is an eigenvalue of the positive definite real matrix $A_j A_j^T$ and is therefore real and positive.
\end{proof}
We close this section by demonstrating how the presence of extra symmetry in the fusion category forces the pivotal tensor to be trivializable. We say that a fusion category has a {\em $*$-operation} is it is equipped with a contravariant monoidal endofunctor $*$ which squares to the identity\XXX{ understand why don't need positivity condition}. Many fusion categories arising from mathematical physics have this property, since they are often the categories of {\em unitary} representations of some or other structure, where the $*$-operation is `taking the adjoint' of a linear map.

\begin{prop} If a fusion category is equipped with a $*$-operation, then in the root choice where all the $d_i$ are positive each involution $T^i_{jk}$ must be the identity map.
\end{prop}
\begin{proof} The reason is that if the fusion category has a $*$-operation, then each involution $T^i_{jk}$ factors as $M^*M$ where $M$ is some other  linear operator; therefore $T^i_{jk}$ is positive definite and so it must be the identity map. Start by choosing the all the roots $d_i^2 = d_{\{i,i^\}}$ to be positive. If we have a $*$-operation, then we can {\em normalize} the units and counit maps at the level of the simple objects by the requirement that
 \[
  \ba \ig{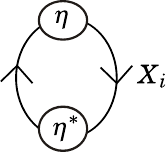} \ea = d_i \quad (\text{or equivalently } \ba \ig{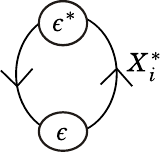} \ea = d_i ).
 \]
This normalization is determined up to a complex phase $e^{i \phi}$; we ask that these phases be fixed. This allows us to factor $T^i_{jk}$ as the composite $NM$,
 \[
  \ba \xymatrix{\Hom(X_i, X_j \otimes X_j) \ar@/^1.5pc/[r]^M & \Hom(X_k^* \otimes X_j^*, X_i^*) \ar@/^2pc/[l]^N} \ea,
 \]
where
 \[
  \ba \ig{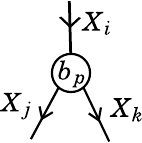} \ea \stackrel{M}{\mapsto} \ba \ig{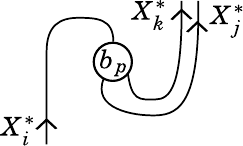} \ea \quad \text{and} \quad \ba \ig{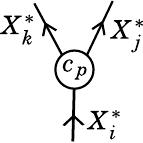} \ea \stackrel{N}{\mapsto} \ba \ig{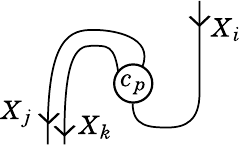} \ea.
 \]
Now, recall from Chapter \ref{evsec} that in a semisimple $*$-category we can always place an inner product on the hom-sets $\Hom(V,W)$ via the formula
 \[
  (f,g) = \Tr(f^* g)
 \]
where $\Tr \colon \End(V) \rightarrow \mathbb{C}$ is the trace operation on the endomorphisms in the category,
 \[
  \Tr(h) = \sum_{i,p} \langle b_{i,p}^\dagger h b_{i,p} \rangle
 \]
where $\{ b_{i,p} \colon X_i \rightarrow V\}$ is a `$*$-basis' which decomposes $V$ in terms of the simple objects $X_i$, that is $b_{i,p}^\dagger b_{i,q} = \delta_{pq}$ (see also Proposition 2.1 of \cite{ref:müger_galois}). By definition this inner product is compatible with the $*$-structure, that is
 \[
  (fg, h) = (g, f^*) \quad \text{and} \quad (fg, h) = (f, hg^*)
 \]
whenever this makes sense; in other words we have equipped our category with the structure of a {\em 2-Hilbert space} \cite{ref:baez_2_hilbert_spaces}.

We claim that with respect to these inner products, $N = M^*$ and hence $T= M^* M$ and hence $T$ must be the identity. Indeed, if we work with the basis for $\Hom(X_k^* \otimes X_j^*, X_i^*)$ given by
 \[
  c_p := \ba \ig{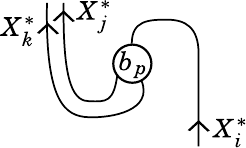} \ea
 \]
then
 \[
  (b_p, Nc_q) = \left( \ba \ig{e249.pdf} \ea, \ba \ig{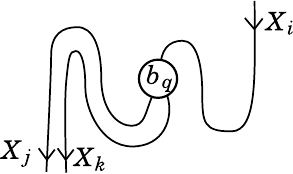} \ea \right) = \delta_{pq}
 \]
while
 \begin{align*}
  \left(Mb_p, c_q\right) = \left( \id_{X_i^*},  c_q (Mb_p)^\dagger\right) &= \left( \ba \ig{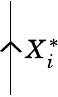} \ea, \;\ba \ig{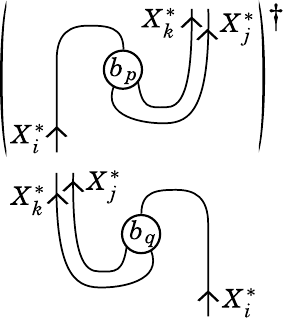} \ea\right)\\
   &= \left(\ba \ig{e257.pdf} \ea, \; \ba \ig{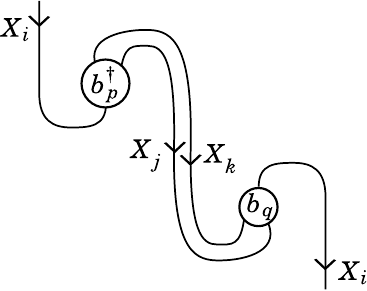} \ea\right) = \delta_{pq}.
 \end{align*}
Thus $N = M^*$ and the result has been proved.
\end{proof}

\section{Existence of pivotal and spherical structures\label{existence}}
Up to this point we have essentially been using Hagge and Hong's string diagram observation \cite[Thm 3]{ref:hagge_hong} to rederive the results of Calaque, Etingof and Nikshych \cite{ref:eno} about fusion categories purely in terms of string diagrams. In this section we go {\em beyond} these results in the following way. We show that an even-handed structure (equivalently, a pivotal structure) on a fusion category can only exist if all the involution operators $T^i_{jk}$ we introduced earlier equal plus or minus the identity map. This is precisely the requirement that the pivotal symbols $\epsilon^i_{jk, p}$ --- which are the eigenvalues of $T^i_{jk}$ --- are independent of $p$, so that we can write them simply as $\epsilon^i_{jk}$. In this situation ,we show that an even-handed structure is the same thing as a {\em $\epsilon$-twisted monoidal natural transformation of the identity} on the category. Moreover we show that the even-handed structure can be made spherical if and only if a root choice $\{d_i^2 = d_{\{i,i^*\}}\}$ of the paired dimensions exists in which all the involution operators $T^i_{jk}$ equal the identity map.

We apologize to the reader that it might well turn out that the involution operators $T^i_{jk}$ in a fusion category are {\em always} plus or minus the identity, and also that these signs might {\em always} be removable by making a different root choice $\{d_i^2 = d_{\{i,i^*\}} \}$. Nevertheless we hope our contribution at least clarifies the conditions necessary for the existence of a pivotal structure, and makes precise the equations which are involved. We hope that this result will help to decide the conjecture made in \cite{ref:eno} that a pivotal structure always exists.

Let us define what we mean by a `twisted monoidal natural transformation of the identity'. Firstly observe that a monoidal natural transformation of the identity on a semisimple monoidal linear category is the same thing as a collection of numbers $\{\theta_i\}$ assigned to the simple objects $X_i$ satisfying
 \[
  \theta_i = \theta_j \theta_k \quad \text{whenever $X_i$ appears in $X_j \otimes X_k$.}
  \]
This causes us to make the following definition.
\begin{defn} Let $C$ be a fusion category in which a root choice $\left\{d_i^2 = d_{\{i,i^*\}}\right\}$ of the paired dimensions has been made. Suppose that the resulting involutions $T^i_{jk} \colon \Hom(X_i, X_j \otimes X_k) \rightarrow \Hom(X_i, X_j \otimes X_j)$ are equal to plus or minus the identity map, that is $T^i_{jk} = \epsilon^i_{jk} \id$ for all $i,j,k$, where $\epsilon^i_{jk} = \pm 1$. Then an {\em $\epsilon$-twisted monoidal natural transformation of the identity} is a collection of numbers $\{t_i\}_{i \in I}$ satisfying
 \[
 t_j t_k = \epsilon^i_{jk} t_i \quad \text{whenever $X_i$ appears in $X_j \otimes X_k$.}
 \]
We write the collection of solutions to these equations as $\Aut_\otimes^\epsilon(\id_C)$.
\end{defn}
Observe that our claim relating even-handed structures to twisted natural transformations at least has the right symmetry group.
\begin{lem} The group $\Aut_\otimes(\id_C)$ acts freely and transitively on $\Aut_\otimes^\epsilon (\id_C)$.
\end{lem}
\begin{proof} A monoidal natural transformation of the identity $\theta$ acts on a twisted monoidal natural transformation by setting $t'_i = \theta_i t_i$. It is clear that this action is free and transitive.
\end{proof}
We will also need the following basic fact later.
 \begin{lem}\label{qqlem} If $t \in \Aut_\otimes^\epsilon(\id_C)$ is an $\epsilon$-twisted monoidal natural transformation of the identity, then we automatically have $t_1 = 1$ and $t_{i^*} = \frac{1}{t_i}$ for all $i$.
 \end{lem}
 \begin{proof} Since $1 \cong 1 \otimes 1$, we have $t_1 = \epsilon^1_{11} t_1^2$. But $\epsilon^1_{11} = 1$, by coherence in a monoidal category (every diagram constituted solely from the unit isomorphisms must commute), so that $t_1 = 1$. Similarly, since $1$ appears in $X_i^* \otimes X_i$, we have $t_i t_{i^*} = \epsilon^1_{i i^*} t_1$. But a string diagram argument shows that $\epsilon^1_{i i^*} = 1$; hence $t_{i^*} = \frac{1}{t_i}$.
 \end{proof}

\subsection*{Pivotal cohomology}
Our main result is that an even-handed structure on a fusion category is the same thing as an $\epsilon$-twisted monoidal natural transformation of the identity, but our second claim is that the even-handed structure can be made spherical structure if and only if the signs can be removed from the $\epsilon$-tensor by an appropriate choice of signs of the roots $d_i$. We can formalize this as follows. Let $A$ be a fusion ring; recall from Section \ref{stss} that this is a ring which looks like the Grothendieck ring of a fusion category. Suppose that the labels for the simple objects are parameterized by a set $I$. We define a {\em collection of pivotal symbols on $A$} as a choice of signs $\{\epsilon^i_{jk} = \pm 1\}_{i,j,k \in I}$ which satisfy the symmetry properties of Lemma \ref{Tsymprop}:
 \[
  \epsilon^i_{jk} = \epsilon^{k^*}_{i^* j}\, ,\;\; \epsilon^i_{jk} = \epsilon^{i^*}_{k^* j^*} \; \text{ and } \; \epsilon^1_{ii^*} = 1 \; \\; \text{for all $i,j,k \in I$.}
 \]
We say that two collections of pivotal symbols $\epsilon$ and $\epsilon'$ are {\em equivalent} if there is a function $f \colon I \rightarrow \{1, -1\}$ with $f_1 = 1$ and $f_i = f_{i^*}$ such that
 \[
  \epsilon'^i_{jk} = f_i f_j f_k \epsilon^i_{jk}.
 \]
We call the set of equivalence classes of collections of pivotal symbols
 \[
  H_\text{piv}(A, \mathbb{Z}/2) = \left\{ \text{candidate collections of pivotal symbols $\{\epsilon^i_{jk}\}$} \right\} / \sim
 \]
the {\em pivotal cohomology} of the fusion ring $A$. Observe that one can multiply collections of pivotal symbols together, so that $H_\text{piv}(A, \mathbb{Z}/2)$ is an abelian group. In summary, we see that a fusion category $C$ whose involution operators $T^i_{jk}$ are plus or minus the identity operator gives rise to a class $[\epsilon] \in H_\text{piv}([C], \mathbb{Z}/2)$ in the pivotal cohomology of $[C]$. Our claim will be that $C$ can carry spherical even-handed structure if and only if this class is trivial.

\subsubsection*{Examples for pivotal cohomology}
Let us calculate the pivotal cohomology of some of the examples of fusion rings listed in Section 1.3.2 of \cite{ref:calaque_etingof}.

\begin{enumerate}

\item The {\em Yang-Lee} fusion ring $A_1$ is generated by $1$ and $X$ with $X^2 = X$ and $X^* = X$. Thus the only nontrivial component of a pivotal tensor is $\epsilon^X_{XX}$, and its sign can be changed by setting $f_X = -1$. Thus $H_\text{piv} (A_1, \mathbb{Z}/2) = 1$.

\item The fusion rings $B_n$ are generated by $X_0, \ldots, X_{n-1}$ and $Y$, with $Y^2 = (n-1)Y + \sum_{i=0}^{n-1} X_i, XY = YX = Y, Y^* = Y, X_iX_j = X_{i+j}$ and $X_i^*=X_{-i}$ (indices are taken mod $n$). When $n=2$, the independent components can be taken to be $\epsilon^{X_1}_{YY}$ and $\epsilon^Y_{YY}$, and their signs can be changed independently by $f_{X_1}$ and $f_Y$, so $H(B_2, \mathbb{Z}/2) = 1$. When $n=3$, we find three independent components $\epsilon^{X_1}_{YY}, \epsilon^{X_2}_{X_1 X_1}$ and $\epsilon^Y_{YY}$, and only two degrees of freedom $f_Y$ and $f_{\{X_1, X_2\}}$ in adjusting the signs, so that $H_\text{piv}(B_3, \mathbb{Z}/2) = \mathbb{Z}/2$. Clearly as $n$ grows the number of independent components of $\epsilon^i_{jk}$ grows as $n^2$, while the degrees of freedom in adjusting the signs is only linear in $n$, so the pivotal cohomology will grow ever larger.

\item The fusion rings $R_G$ of the {\em Tambara-Yamagami categories} \cite{ref:tambara_yamagami} have the form $\mathbb{Z}[G] \oplus \mathbb{Z}[Y]$ for a finite abelian group $G$ with $Y^2 = \sum_{g \in G} g$, $gY = Yg = Y$, $gh = g\cdot h$, $g^* = g^{-1}$ and $Y^* = Y$. The same argument as we made for the fusion rings $B_n$ establishes that the size of $H_\text{piv}(R_G, \mathbb{Z}/2)$ grows ever larger with $|G|$.

\end{enumerate}
Nevertheless, although in these examples the pivotal cohomology of the fusion rings is generally nonzero, the {\em class} in that cohomology coming from the fusion category itself is always trivial, because all these categories are spherical (we remind the reader that {\em all} known examples of fusion categories can be given a spherical even-handed structure).

\subsection*{Main result}
We can now state the main result of this chapter.

\begin{thm} Let $C$ be a fusion category over $\mathbb{C}$ with representative simple objects $X_i$. Suppose that a choice of roots $d_i^2 = d_{\{i,i^*\}}$ of the paired dimensions has been made, with resulting involution operators $T^i_{jk} \colon \Hom(X_i, X_j \otimes X_k) \rightarrow \Hom(X_i, X_j \otimes X_k)$.  Then:
 \begin{enumerate}
  \item Unless $T^i_{jk}=\pm \id$ for all $i,j$ and $k$, the fusion category $C$ cannot carry an even-handed structure.
  \item Suppose that $T^i_{jk} = \epsilon^i_{jk} \id$ for all $i,j$ and $k$, where $\epsilon^i_{jk} = \pm 1$. Then an even-handed structure on $C$ is the same thing as an $\epsilon$-twisted monoidal natural transformation of the identity on $C$. That is, there is a canonical bijection of sets
       \[
        \EvenHanded(C) \cong \Aut_\otimes^\epsilon (\id_C).
       \]
   \item Furthermore, the even-handed structure can be made spherical if and only if $[\epsilon] = 0$ in $H_\piv([C], \mathbb{Z}/2)$.
  \end{enumerate}
\end{thm}
\begin{proof} (i) and (ii). Suppose that $\Psi$ is an even-handed structure on $C$. At the level of the simple objects, we can write (using the Partner Convention)
 \[
  \left(\ba \ig{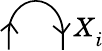} \ea , \ba \ig{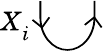} \ea \right) \stackrel{\Psi}{\mapsto} \left(\frac{1}{t_i} \,\, \ba \ig{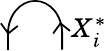} \ea, \; t_i \,\, \ba \ig{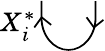} \ea \right)
 \]
for nonzero scalars $t_i$. Another way to write this formula is
\[
  \left(\ba \ig{e224.pdf} \ea , \ba \ig{e225.pdf} \ea \right) \stackrel{\Psi}{\mapsto} \sign(d_i) \left(\frac{1}{t_i} \ratone \, \ba \ig{e228.pdf} \ea, \; t_i \rattwo \,\, \ba \ig{e229.pdf} \ea \right)
 \]
where $\left(\ba \ig{e228.pdf} \ea, \ba \ig{e229.pdf} \ea \right)$ is an {\em arbitrary} choice of unit and counit maps expressing $X_i^*$ as a left dual of $X_i$; in this formula the signs of $d_i$ are made explicit. Observe also that since $\Psi^2 = \id$, we have $t_{i^*} = \frac{1}{t_i}$.

We now show that $\Psi$ is completely determined by these numbers $t_i$. Consider the behaviour of $\Psi$ at an arbitrary object $V \in C$,
 \[
  \left( \ba \ig{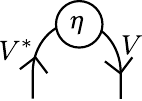} \ea, \ba \ig{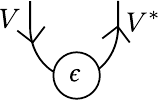} \ea \right) \stackrel{\Psi}{\mapsto} \left(\ba \ig{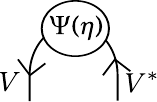} \ea , \ba \ig{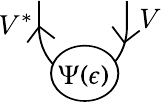} \ea\right).
 \]
By choosing bases
 \[
 \{v_{i,p} \colon X_i \rightarrow V\}_{p=1}^{\dim \Hom(X_i, V)}
 \]
for each hom-space $\Hom(X_i, V)$, with accompanying dual bases written as
 \[
  \{v_i^p \colon V \rightarrow X_i \},
\]
we find that
 \begin{align*}
  \ba \ig{e232.pdf} \ea = \sum_{i,p} \ba \ig{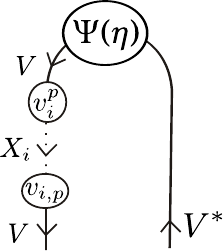} \ea &= \sum_{i,p} \ba \ig{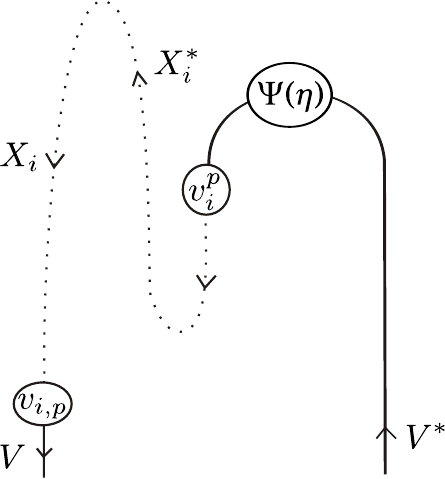} \ea
  \\ &= \sum_{i,p} \frac{1}{t_i} \ba \ig{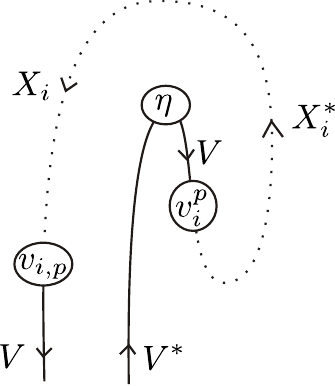} \ea.
 \end{align*}
where the last equality uses the even-handedness of $\Psi$. Similarly we have
 \begin{equation} \label{tjeq}
  \ba \ig{e233.pdf} \ea = \sum_{i,p} t_i \quad \ba \ig{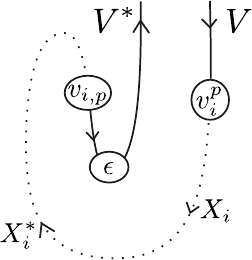} \ea.
    \end{equation}
Thus $\Psi$ is completely determined by the numbers $t_i$. Now, the fact that $\Psi$ respects tensor products means that
 \[
 \Psi\left( \ba \ig{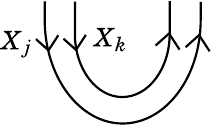} \ea \right) = t_j t_k \ba \ig{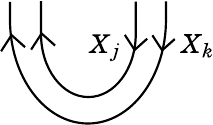} \ea
 \]
where we are again applying the Partner Convention. We can expand out this equation using \eqref{tjeq} to obtain
 \[
  \sum_{i,p} t_i \; \ba \ig{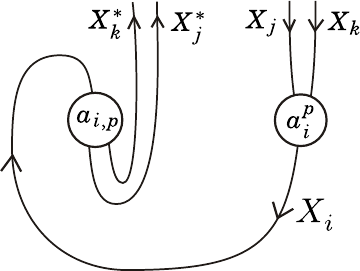} \ea  = t_j t_k \ba \ig{e239.pdf} \ea.
 \]
Thus we have
 \[
  \sum_{i,p} t_i \ba \ig{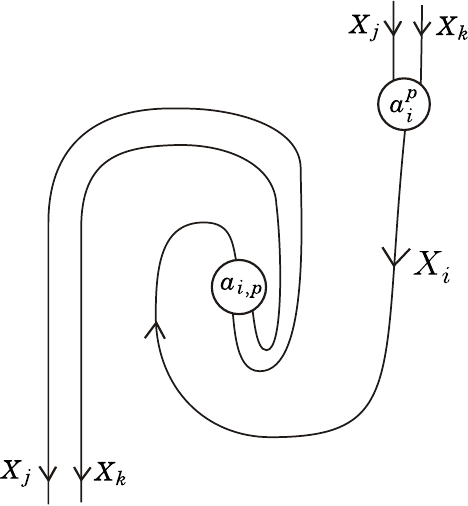} \ea = t_j t_k \ba \ig{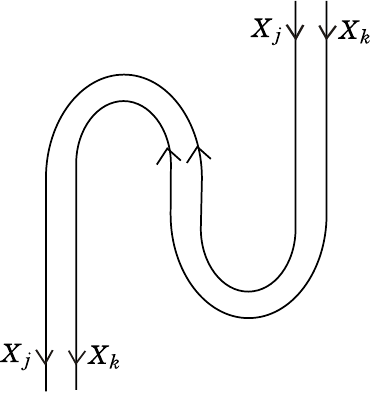} \ea
 \]
where we recognize the operator $T^i_{jk}$ acting on $a_{i,p}$, which is an eigenvector, so we have
 \[
  \sum_{i} t_i  \sum_p \epsilon^i_{jk, p}\ba \ig{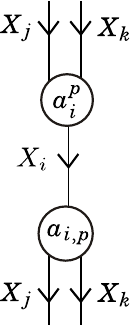} \ea = t_j t_k \;\; \ba \ig{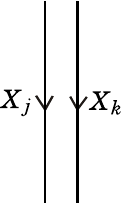} \ea.
 \]
In terms of matrices on the vector spaces $\Hom(X_i, X_j \otimes X_k)$, we therefore have
 \[
 t_i \left( \begin{array}{cccc} \epsilon^i_{jk, 1} & & & \\ & \epsilon^i_{jk, 2} & & \\ & & \ddots & \\ & & & \epsilon^i_{jk, n} \end{array} \right) = t_j t_k \id \quad \text{whenever $X_i$ appears in $X_j \otimes X_k$}.
\]
Thus for fixed $i,j,k$ all the pivotal symbols $\epsilon^i_{jk, p}$ must be equal --- whence we will write them as $\epsilon^i_{jk}$ --- and the numbers $\{t_i\}$ constitute an $\epsilon$-twisted monoidal natural transformation of the identity on $C$:
 \begin{equation} \label{etwist}
  t_i = \epsilon^i_{jk} t_j t_k \quad \text{whenever $X_i$ appears in $X_j \otimes X_k$.}
 \end{equation}
In the reverse direction, suppose that $\{t_i\}$ are some nonzero scalars satisfying \eqref{etwist}. Then we can define an even-handed structure $\Psi$ by the formula
 \[
  \left( \ba \ig{e230.pdf} \ea, \ba \ig{e231.pdf} \ea \right) \stackrel{\Psi}{\mapsto} \Bigg(\sum_{i,p} \frac{1}{t_i}  \ba \ig{e236.pdf} \ea , \; \sum_{i,p} t_i \ba \ig{e237.pdf} \ea\Bigg).
 \]
Let us verify that this formula satisfies the axioms. Firstly, we have $\Psi^2 = \id$ since a quick string diagram argument shows that this is equivalent to $t_{i^*} = \frac{1}{t_i}$, which we saw in Lemma \ref{qqlem} was a consequence of the equations \eqref{etwist} satisfied by the scalars $t_i$. We also saw that these equations required that $t_1 = 1$, so that $\Psi$ acts as the identity on the trivial adjunction on the unit object. Moreover $\Psi$ respects the monoidal structure, because we have seen that this is precisely equivalent to the equation $t_i = \epsilon^i_{jk} t_j t_k$ whenever $X_i$ appears in $X_j \otimes X_k$.

Finally, the even-handed equation is trivially satisfied, because it is trivial at the level of the simple objects. To see this, suppose $\theta \colon V \rightarrow W$ is a morphism. For each $i$ choose bases $\{v_{i,p} \colon X_i \rightarrow V\}$ and $\{w_{i,q} \colon X_j \rightarrow W\}$ for $\Hom(X_i, V)$ and $\Hom(X_i, W)$ respectively; the corresponding dual basis vectors are written as $v_i^p$ and $w_i^q$ as usual. Fix a right dual $(V^*, \ba \ig{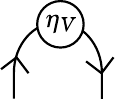} \ea, \ba \ig{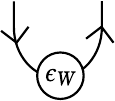} \ea)$ for $V$ and similarly for $W$. Now we also have the corresponding `$*$-dual' basis vectors
 \[
  {v_i^p}^* = \ba \ig{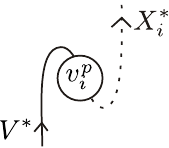} \ea, \quad v_{i,p}^* = \ba \ig{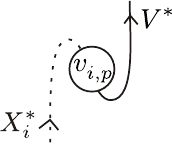} \ea, \quad v_{i,q}^* {v_i^p}^* = \delta^p_q
 \]
and the same for $W$. With this setup, the matrix elements of $\theta^*$ compute as the transpose of the matrix elements of $\theta$:
\[
 \langle v_{i,p}^* \theta^* {w_i^q}^* \rangle = \langle \ba \ig{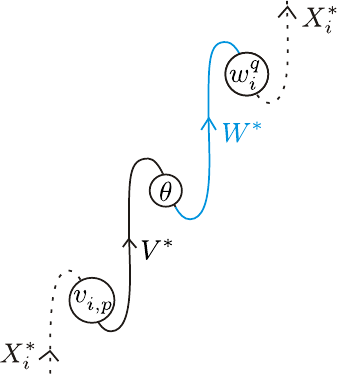} \ea \rangle = \langle \ba \ig{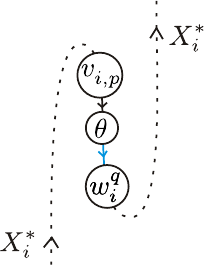} \ea \rangle = \langle w_i^q \theta v_{i,p} \rangle.
\]
And the same holds for $^* \! \theta_\psi$:
 \[
 \langle v_{i,p}^* {}^* \! \theta_\Psi {w_i^q}^* \rangle = \sum_{j,q'} \sum_{k,p'} \ba \ig{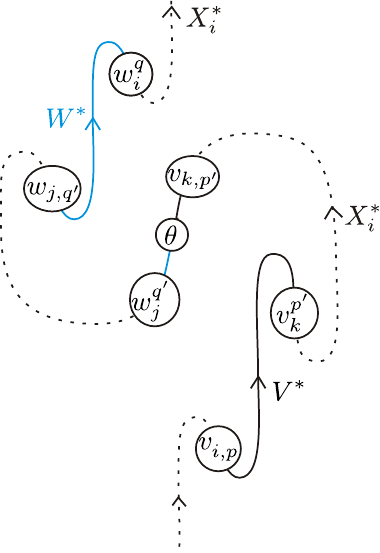} \ea = \ba \ig{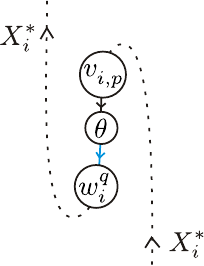} \ea = \langle w_i^q \theta v_{i,p} \rangle.
 \]

\noindent (iii) The pivotal structure specified by the numbers $t_i$ is spherical when $\dim [X_i] = \dim [X_i^*]$. This is the requirement that
 \[
  t_i \ba \ig{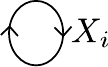} \ea = t_{i^*} \ba \ig{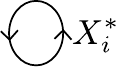} \ea,
 \]
which since $t_{i^*} = \frac{1}{t_i}$, is the requirement that $t_i^2 = 1$. Thus a spherical structure exists if and only if there is a solution of the pivotal equations \eqref{etwist} with each $t_i = \pm 1$. This means that
 \[
  \epsilon^i_{jk} = \sign(t_i) \sign(t_j) \sign(t_k).
 \]
But then we can remove all the signs from the $\epsilon$-tensor by choosing a new root system $d_i \mapsto \sign (t_i) d_i$, because under this transformation the the $\epsilon$ tensor will transform as
 \[
  \epsilon^i_{jk} \mapsto \sign(t_i) \sign(t_j) \sign (t_k) \epsilon^i_{jk} = 1.
 \]
Thus $[\epsilon] = 0$ in $H_\text{piv}(C, \frac{\mathbb{Z}}{2})$. Conversely, if $\epsilon$ can be made trivial then $t_i = 1$ is a spherical solution of the pivotal equations.
 \end{proof}

\chapter[2-representations and their 2-characters]{Unitary 2-representations and their 2-characters\label{2RepsChap}}
In this chapter we define and investigate unitary 2-representations of finite groups on 2-Hilbert spaces and their corresponding 2-characters.

The idea of a `categorified representation' of a group or 2-group as a group action on some sort of `categorified vector space' has already been studied by a number of authors, including Elgueta \cite{ref:elgueta}, Crane and Yetter \cite{ref:crane_yetter}, Barrett and Mackaay \cite{ref:barrett_mackaay}, Ganter and Kapranov \cite{ref:ganter_kapranov_rep_char_theory}, Ostrik \cite{ref:ostrik}, Freed \cite{ref:freed_remarks} and Baez, Baratin, Freidel and Wise \cite{ref:baez_baratin_freidel_wise}. These only represent the texts where group actions on {\em semisimple} linear categories were studied, because that is the area in which the present author is most familiar. There is also a body of work having to do with braid group actions on derived categories, starting with that of Deligne \cite{ref:deligne}, but also Seidel and Thomas \cite{ref:seidel_thomas}; a number of works having to do with the Khovanov homology research programme for categorifying tangle invariants \cite{ref:khovanov}; as well as works arising in geometric representation theory which are close in spirit to the Geometric Langlands programme \cite{ref:benzvi2}. 

The notion of the 2-character of a 2-representation which we will shortly discuss has been defined independently by Ganter and Kapranov \cite{ref:ganter_kapranov_rep_char_theory}, but there is also the important early work by Luztig on character sheaves \cite{ref:luztig} which is now being developed by Boyarchenko and Drinfeld amongst others in the setting of unipotent groups \cite{ref:boyarchenko_drinfeld}, and we mention also the work of Polesello and Waschkies \cite{ref:polesello_waschkies} in this regard.

Certainly there is a lot of work to be done in unifying the various pictures. As far as our own personal contribution here is concerned, we have already outlined in the introduction to this thesis what is novel in our approach, but for the benefit of the reader we briefly recall the following. Firstly, we study {\em unitary} 2-representations on the {\em 2-Hilbert spaces} of Baez \cite{ref:baez_2_hilbert_spaces}, which to our knowledge have not been explicitly defined or studied before. In this way we make contact with the Baez-Dolan research programme \cite{ref:baez_dolan_hda0} of extended topological quantum field theory, where {\em unitary} structures play the central role, because they provide the crucial ingredient of  duality. Secondly, we use {\em string diagrams} as our basic notation for working with 2-representations and their 2-characters, a strategy which considerably simplifies many constructions and computations. Thirdly, because a 2-Hilbert space can be regarded as a `categorified vector space equipped with an inner product', studying unitary 2-representations enables us to develop a {\em push-forward map} for the 2-character of a 2-representation. This enables one to see the 2-character as a {\em functor}
 \[
  \chi \colon [\TRep(G)] \rightarrow \Hilb_G (G)
 \]
from the homotopy category of unitary 2-representations to the category of vector bundles over the group equivariant under the conjugation action. One of our main results in this thesis, which we will prove in Chapter \ref{GerbesCharChap} after establishing the geometric correspondence between unitary 2-representations and {\em equivariant gerbes}, is that after one tensors the hom-sets in $[\TRep(G)]$ with $\mathbb{C}$ the resulting 2-character functor is {\em unitarily fully faithful}. This is the categorification of the fact that the ordinary character map is a unitary isomorphism
 \[
  \chi \colon [\Rep(G)] \rightarrow \Class(G)
 \]
from the set of isomorphism classes of unitary representations of a group to the Hilbert space of class functions becomes a unitary isomorphism after one tensors $[\Rep(G)]$ with $\mathbb{C}$.\XXX{Go over all your intros to chapters, changing notation to labelled notation for sections}

We should insert the disclaimer here that 2-Hilbert spaces are only a first approximation to `categorified Hilbert spaces' because they are {\em semisimple} and hence not rich enough to accommodate continuous geometry; this is the reason we must restrict ourselves to 2-representations of {\em finite} groups. We hope though that similar ideas will apply in the non-semisimple {\em derived} context that we mentioned above; the work of Costello \cite{ref:costello} appears to be a promising way to unify the two pictures. In addition, Baez, Baratin, Freidel and Wise have built on work of Crane and Yetter \cite{ref:yetter2, ref:crane_yetter} by accommodating 2-representations of Lie 2-groups using a notion of 2-Hilbert space where the objects are `measurable fields of Hilbert spaces' supported over a measurable space \cite{ref:baez_baratin_freidel_wise}.

We should also remark that not {\em all} facts about ordinary characters of representations carry over to 2-representations on 2-Hilbert spaces; for instance we will see in Chapter \ref{notdistinguish} that the 2-character does not always distinguish inequivalent 2-representations.

In Section \ref{utrep} we define unitary 2-representations and the 2-category $\TRep(G)$ which they constitute, firstly in a terse higher categorical way, and then by expanding out this definition in ordinary notation as well as in string diagram notation by introducing graphical elements to depict the various pieces of data involved. In Section \ref{exsec} we give a number of examples of unitary 2-representations, and the morphisms between them. In Section \ref{slem} we introduce new graphical elements into the string diagram notation, and prove some basic graphical identities. In Section \ref{evun} we show that unitary 2-representations are compatible with the canonical even-handed structure on $\THilb$, a fact which allows an important diagrammatic lemma regarding adjunctions to be established. Finally in Section \ref{charsecc} we define the notion of the 2-character of a 2-representation, and use the string diagram technology we have developed to give elegant proofs of the fact that the 2-character of a 2-representation produces an equivariant vector bundle over the group, and that the 2-character of a {\em morphism} of 2-representations produces a morphism between the corresponding equivariant vector bundles, in a manner which is well-defined and functorial with respect to composition of 1-morphisms in $\TRep(G)$.

\section{The 2-category of unitary 2-representations\label{utrep}}
In this section we define the 2-category $\TRep(G)$ of unitary 2-representations of a finite group $G$. We do this in two stages --- firstly we define it in a terse higher-categorical way, and then we expand out this definition explicitly in traditional notation as well as in string diagrams, introducing new graphical elements to depict the various pieces of data involved.

Since there are various conventions for terminology for 2-categories, we remind the reader we are essentially using those of Leinster \cite{ref:leinster_basic_bicategories}. Also, we use the notation that $\BG$ refers to the group $G$ thought of as a one object category, with the elements of $G$ as morphisms (we use bold face to distinguish $\BG$ from $BG$, the geometric realization of the nerve of $\BG$).

\begin{defn}[{compare \cite{ref:elgueta, ref:crane_yetter, ref:barrett_mackaay,  ref:ganter_kapranov_rep_char_theory, ref:ostrik}}] The 2-category $\TRep(G)$ of unitary 2-representations of a finite group $G$ is defined as follows. An object is a unitary weak 2-functor $\BG \rightarrow \THilb$, where $\BG$ is thought of as a 2-category with only identity 2-morphisms. A morphism is a
transformation whose coherence isomorphisms are unitary, and a 2-morphism is a modification.
\end{defn}

Recall from Chapter \ref{BigTCat} that a linear $*$=functor between 2-Hilbert spaces is called {\em unitary} if it is unitary at the level of hom-sets; by a `unitary' weak 2-functor $\alpha \colon G \rightarrow \THilb$ we mean one which sends each element of $G$ to a unitary endofunctor $\alpha_g \colon H \rightarrow H$ of a 2-Hilbert space.

We now expand this definition out.
\subsection{Unitary 2-representations\label{Untrep}}
A unitary 2-representation of $G$ consists of:
\begin{itemize}
   \item A finite-dimensional 2-Hilbert space $H$, drawn as
    \[
     \ba \ig{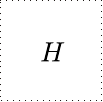}
   \ea \quad \text{or} \quad \ba \ig{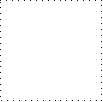} \ea \, \, \text{\; when $H$ is
   understood,}
   \]
   \item For each $g \in G$, a linear $*$-functor $H \la{\alpha_g} H$ which at the level of hom-sets is a unitary linear map, drawn as
   \[
     \ba \ig{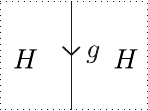} \ea \text{\; or simply \;} \ba \ig{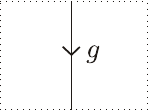} \ea,
    \]
   \item A unitary natural isomorphism $a(e) \colon \id_H \Rightarrow \alpha_e$ (where $e$ is the identity element of $G$), and for each $g_1, g_2 \in G$, a unitary natural isomorphism $a(g_2,g_1)
   \colon \alpha_{g_2} \circ \alpha_{g_1} \Rightarrow \alpha_{g_2g_1}$, drawn as
    \[
     \ba
     \ig{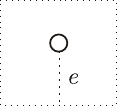} \ea \quad \ba \xymatrix @R=0.5cm {\id_A \ar@{=>}[d]^{a(e)}
     \\ \alpha_g} \ea \qquad , \qquad
    \ba \ig{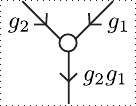} \ea \quad \ba \xymatrix @R=0.5cm {\alpha_{g_2} \circ
     \alpha_{g_1} \ar@{=>}[d]^{a(g_2,g_1)} \\ \alpha_{g_2 g_1}} \ea
     \]
\end{itemize}
such that
 \[
  \ba \xymatrix{ & \alpha_g \ar@{=>}[dl]_{a(e) \ast \id}
  \ar@{=>}[dr]^{\id \ast a(e)} \ar@{=}[dd] \\ \alpha_e \circ \alpha_g \ar@{=>}[dr]_{a(e,g)} & &
  \alpha_g \circ \alpha_e \ar@{=>}[dl]^{a(g,e)} \\ & \alpha_g} \ea \quad
  \text{and} \quad \ba
  \xymatrix @C = 0.1cm{ & \alpha_{g_3} \circ \alpha_{g_2} \circ \alpha_{g_1} \ar@{=>}[dr]^{\; a(g_3,g_2) \ast \id} \ar@{=>}[dl]_{\id \ast a(g_2,g_1)\;} \\ \alpha_{g_3 g_2} \circ \alpha_{g_1} \ar@{=>}[dr]_{a(g_3, g_2 g_1)\;} &
  & \alpha_{g_3} \circ \alpha_{g_2 g_1} \ar@{=>}[dl]^{\;a(g_3 g_2, g_1)} \\ & \alpha_{g_3 g_2 g_1}} \ea
 \]
commute, or in string diagrams,
  \be
  \label{br2} \ba \ig{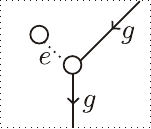} \ea = \ba \ig{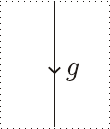} \ea = \ba \ig{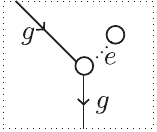} \ea \qquad \text{and} \qquad \ba \ig{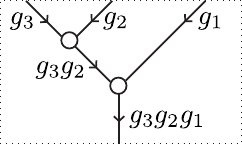} \ea = \ba \ig{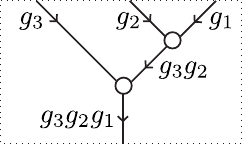} \ea.
  \ee
We will draw the inverse 2-isomorphisms $a(e)^\mi \colon \alpha_e \Rightarrow \id_H$
and $a(g_2,g_1)^\mi  \colon \alpha_{g_2 g_1} \Rightarrow \alpha_{g_2} \circ \alpha_{g_1}$ as
 \[
  \ba \xymatrix @R=0.5cm {\alpha_e \ar@{=>}[d]^{a(e)^*} \\ \id_A} \ea \quad \ba \ig{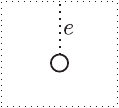} \ea \qquad , \qquad \ba \xymatrix @R=0.5cm {\alpha_{g_2 g_1}
  \ar@{=>}[d]^{a(g_2,g_1)^*} \\ \alpha_{g_2} \circ \alpha_{g_1}} \ea \quad \ba \ig{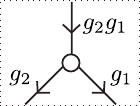} \ea .
 \]
Recall that these 2-isomorphisms are required to be {\em unitary}, so that $a(e)^\mi = a(e)^*$ and $a(g_2, g_1)^\mi = a(g_2, g_1)^*$. The fact that these satisfy $a(e)^\mi a(e) = \id$ and $a(e) a(e)^\mi = \id$, and similarly for the $a(g_2, g_1)$, is drawn as follows:
 \begin{eqnarray}
  \ba \ig{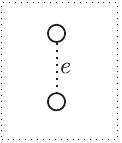} \ea = \ba \ig{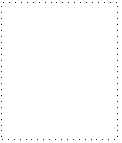} \ea \quad &,& \quad \ba
  \ig{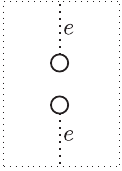} \ea = \ba \ig{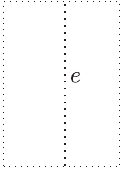} \ea \label{zr1}  \\
  \label{br1} \ba \ig{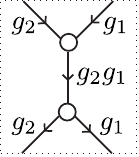} \ea = \ba \ig{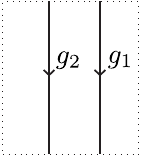} \ea \quad &,& \quad \ba
  \ig{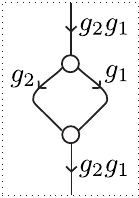} \ea = \ba \ig{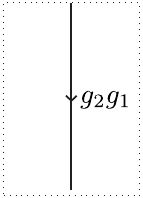} \ea. \label{zr2}
 \end{eqnarray}
We will abbreviate all of this data $(H, \{\alpha_g\}, a(e), \{a(g_2, g_1)\})$ simply as $\alpha$.

\subsection{Morphisms\label{Morphsec}} A morphism $\sigma \colon \alpha \rightarrow \beta$ of unitary 2-representations is a transformation from $\alpha$ to $\beta$ whose coherence isomorphisms are unitary. Thus, if $\alpha = (H_\alpha, \{\alpha_g\}, a(e), \{a(g_2, g_1)\})$ and $\beta = (H_\beta, \{\beta_g\}, b(e), \{b(g_2, g_1\})$, then it consists o:f
 \begin{itemize}
  \item A linear $*$-functor $\sigma \colon H_\alpha \rightarrow H_\beta$, drawn as
   \[
    \ba \ig{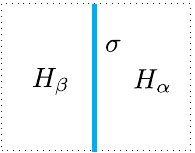} \ea \quad
    \text{or} \quad \ba \ig{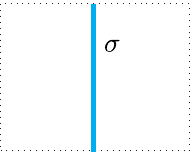} \ea \, \, \text{($H_\alpha$ and $H_\beta$
    understood)}
   \]
  The line for $\sigma$ is thick and coloured differently, so as to distinguish it from the lines for the functors $\alpha_g$ and $\beta_g$.
  \item For each $g \in G$ a unitary natural isomorphism $\sigma(g) \colon \beta_g
  \circ \sigma \Iso \sigma \circ \alpha_g$, drawn as
   \[
    \ba \ig{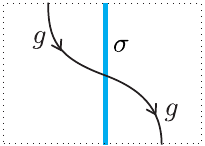} \ea \quad \ba \xymatrix{\beta_g \circ \sigma
    \ar@{=>}[d]^{\sigma(g)} \\ \sigma \circ \alpha_g} \ea
   \]
  \end{itemize}
such that
 \[
  \ba \xymatrix{ & \beta_{g_2} \circ \beta_{g_1} \circ \sigma \ar@{=>}[dl]_{b(g_2,g_1)
  \ast \id} \ar@{=>}[r]^{\id \ast \sigma(g_1)} & \beta_{g_2} \circ \sigma
  \circ \alpha_{g_1} \ar@{=>}[dd]^{\sigma(g_2) \ast \id} \\ \beta_{g_2g_1} \circ
  \sigma \ar@{=>}[dr]_{\sigma(g_2 g_1)} \\ & \sigma \circ \alpha_{g_2 g_1} & \sigma
  \circ \alpha_{g_2} \circ \alpha_{g_1} \ar@{=>}[l]^{\id \ast a(g_2,g_1)}} \ea \quad
  \ba \xymatrix{\sigma \ar@{=>}[d]_{b(e) \ast \id}
  \ar@{=>}[dr]^{\id \ast a(e)} \\ \beta_e \circ \sigma
  \ar@{=>}[r]_{\sigma(e)} & \sigma \circ \alpha_e} \ea
 \]
commute, or in string diagrams,
 \be \label{sig1}
  \ba \ig{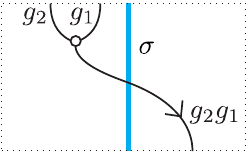} \ea = \ba \ig{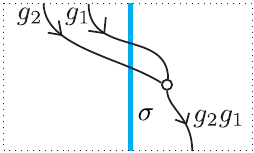} \ea \quad \text{and} \quad
  \ba \ig{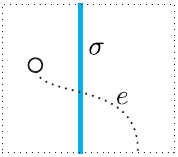} \ea = \ba \ig{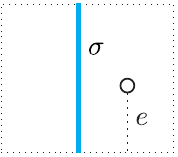} \ea.
 \ee
We will draw the inverse 2-isomorphisms $\sigma(g)^\mi \colon
\sigma \circ \alpha_g \Rightarrow \beta_g \circ \sigma$ as
 \[
  \ba \ig{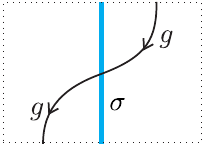} \ea \quad \ba \xymatrix{\sigma \circ \alpha_g
    \ar@{=>}[d]^{\sigma(g)^*} \\ \beta_g \circ \sigma} \ea \, .
 \]
These are required to be {\em unitary}, so that $\sigma(g)^* = \sigma(g)^\mi$. By definition we have $sigma(g)^\mi \sigma(g) = \id$ and $\sigma(g)
\sigma(g)^\mi = \id$, that is,
 \be \label{siginv}
 \ba \ig{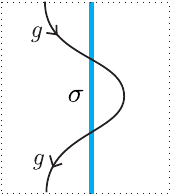} \ea = \ba \ig{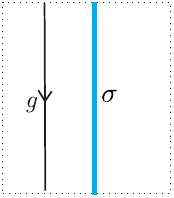} \ea \quad \text{and} \quad
 \ba \ig{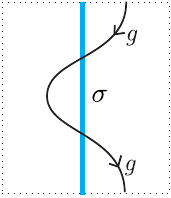} \ea = \ba \ig{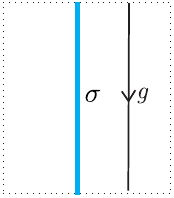} \ea \, .
 \ee
We will abbreviate all of this data $(\sigma, \{\sigma(g)\})$ simply as $\sigma$. Observe that a morphism of 2-actions of $G$, which might be called an {\em intertwiner}, really {\em does} have an `intertwining' aspect to it when expressed in terms of string diagrams.

\subsection{2-morphisms} Finally, if $\alpha$ and $\beta$ are unitary 2-representations of $G$, and $\sigma, \rho \colon \alpha \rightarrow \beta$ are morphisms between them, then a 2-morphism $\theta \colon \sigma \Rightarrow \rho$ is a modification from $\sigma$
 to $\rho$. Thus, $\theta$ is a natural transformation $\sigma$ to $\rho$, drawn as
 \[
  \ba \ig{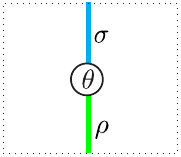} \ea \quad \ba \xymatrix{\sigma
  \ar@{=>}[d]^{\theta} \\ \rho} \ea \, ,
 \]
such that
 \[
 \xymatrix{ & \beta_g \circ \sigma \ar@{=>}[dl]_{\sigma(g)}
 \ar@{=>}[dr]^{\id \ast \theta} \\ \sigma \circ \alpha_g
 \ar@{=>}[dr]_{\theta \ast \id} & & \beta_g \circ \rho
 \ar@{=>}[dl]^{\rho(g)} \\ & \rho \circ \alpha_g}
 \]
commutes, or in string diagrams,
 \be \label{mod}
  \ba \ig{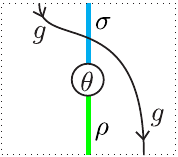} \ea = \ba \ig{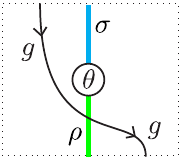} \ea \, .
 \ee
We trust that the simplicity of these diagrams has persuaded the reader the string diagrams are a useful notation for working with 2-representations. We will develop this notation further as we go along.

\section{Examples\label{exsec}}  We now give some examples to illustrate these ideas; we will say more about them in Section \ref{gexam} once we have established the geometric interpretation of unitary 2-representations in terms of equivariant gerbes. We encourage the reader to consult the paper of Ganter and Kapranov \cite{ref:ganter_kapranov_rep_char_theory} for examples of group actions on linear categories arising in a geometric context of a more advanced nature.

\subsection*{2-representations can be strictified}
Before we give the examples, let us first mention an important fact to keep in mind when thinking about 2-representations. A 2-representation $\alpha$ is called {\em strict} if all the
coherence isomorphisms are identities.
\begin{lem} Every 2-representation is equivalent inside $\TRep(G)$ to a strict \\ 2-representation.
\end{lem}
\begin{proof} The proof is essentially an application of the 2-Yoneda lemma, which can be found for instnace in \cite{ref:street_yoneda, ref:gurski} as well as \cite[page 60]{ref:fantechi}. Given a
2-representation of $G$ on a 2-Hilbert space $H$, we can define a corresponding strict 2-representation of $G$ on the 2-Hilbert space
 \[
  \Hom_{\TRep(G)} (\Hilb[G], \alpha),
 \]
where $\Hilb[G]$ is the category of $G$-graded Hilbert spaces on which $G$ acts by left multiplication.
\end{proof}
We warn the reader that this does not mean that there is no information in the coherence isomorphisms --- it just means that this information can always be shifted into the structure of a new and bigger category, if one wishes to do so. In other words, {\em a strict 2-representation on a `big' 2-Hilbert space} (such as `the category of {\em all} such and such') {\em is not always as trivial as it may seem} --- to decide this, one has to calculate the 2-cocycle of the corresponding equivariant gerbe, as we explain in Section \ref{extract}.

\subsection{Automorphisms of groups} Suppose $G \subseteq \Aut(K)$ is a subgroup of the automorphism group of a
finite group $K$. This gives rise to a unitary 2-representation of $G$ on the 2-Hilbert space $\Rep(K)$ by precomposition. That is, if
$V$ is a unitary representation of $K$, then $g \cdot V \equiv V^g$ has the same underlying vector space
except that the action of $k \in K$ on $V^g$ corresponds to the action of $g^\mi \cdot k$ on $V$. This is of course a {\em strict} 2-representation, but it is not necessarily frivolous, as we shall see in the next example. Also note that any 2-representation of this form will necessarily be {\em unitary}. That is because it can only permute irreducible representations of the same dimension amongst each other, and we saw in Chapter \ref{2HilbChap} that the scale factors on the irreducible representations in the 2-Hilbert space $\Rep(K)$ are precisely their dimensions divided by the order of the group, so that a functor which permutes the irreducibles is unitary if and only if it sends representations of the same dimension amongst themselves.

\subsection{The metaplectic representation} A good example of a nontrivial 2-representation of the above sort is the action of $SL_2(\mathbb{R})$ on $\Rep(\text{Heis})$, the category of representations of the Heisenberg group. Of course, $SL_2(\mathbb{R})$ is not a {\em finite} group, but all the definitions above still apply.

The Heisenberg group arises in quantum mechanics (see for example \cite{ref:carter_segal_macdonald}). It is the 3-dimensional Lie group with underlying manifold
$\mathbb{R}^2 \times U(1)$ --- with $\mathbb{R}^2$ thought of as phase space with elements being pairs $v = (v_z, v_p)$
--- and multiplication defined by
 \[
  (v, e^{i \theta}) \cdot (w, e^{i a}) = (v + w, e^{i\omega (v, w)} e^{i(\theta +
  a)}),
 \]
where $\omega(v, w) = \frac{1}{2}(v_z w_p - v_p w_z)$ is the canonical symplectic form on $\mathbb{R}^2$. Up to
isomorphism, there is only one irreducible representation of $\text{Heis}$ on a separable Hilbert space, with
the $U(1)$ factor acting centrally. Namely, the action on $L^2(\mathbb{R})$ given by
 \[
  (z \cdot f) (x) = e^{izx} f(x), \quad (p\cdot f)(x) = f(x-p).
 \]
Since there is only one irreducible representation, $\Rep(\text{Heis})$ is a one-dimensional 2-Hilbert space.

Now $SL_2 (\mathbb{R})$ is the group of symplectomorphisms of $\mathbb{R}^2$, hence it acts as automorphisms of
Heis, giving rise to a unitary 2-representation of $SL_2(\mathbb{R})$ on $\Rep(\text{Heis})$ via the standard
prescription $(g \cdot \rho)(v, e^{i \theta}) = \rho(g^\mi \cdot v, e^{i \theta})$.

This gives rise to a nontrivial projective representation of $SL_2 (\mathbb{R})$ on $L^2(\mathbb{R})$; the fact that the projective factor cannot be removed is known as the `metaplectic anomaly'. Indeed, the viewpoint of 2-representations elucidates somewhat the nature of this anomaly.  It might seem strange at
first that the action of $SL_2(\mathbb{R})$ --- the symmetry group of the classical phase space $\mathbb{R}^2$
--- does not survive quantization, becoming instead a projective representation. However $SL_2(\mathbb{R})$
{\em does} act on $\Rep(\text{Heis})$, the collection of {\em all} quantizations. From this we see that the
`anomaly' arose from an attempt to decategorify this action, by artificially choosing a fixed quantization
$\rho$.

\subsection{2-representations from exact sequences\label{exactseqsec}} We've seen how an action of $G$ on another group $K$ gives rise
to a unitary 2-representation of $G$ on $\Rep(K)$. The same can be said for a `weak' action of $G$ on $K$. Suppose
 \[
  1  \rightarrow K \stackrel{i}{\hookrightarrow} E \stackrel{\pi}{\twoheadrightarrow} G  \rightarrow 1
 \]
is an exact sequence of finite groups, which has been equipped with a set-theoretic section $s \colon G \rightarrow
E$ such that $s(e) = e$. We can think of this data as a homomorphism of 2-groups
 \[
  G \rightarrow AUT(K)
 \]
where $AUT(K)$ is the 2-group whose objects are the automorphisms of $K$ and whose morphisms are given by
conjugation (see \cite{ref:baez_lauda_2-groups, ref:baez_higher_schreier_theory}). Explicitly, one thinks of the group $K$ as being the morphisms of a one-object category (also
denoted $K$), and for each $g \in G$, $g \colon K \rightarrow K$ is the functor defined by conjugating in $E$,
 \[
  g \cdot k := s(g) k s(g)^\mi
 \]
where we have identified $K$ with its image in $E$. This determines a $K$-valued 2-cocycle $\varphi$
having the property that
 \[
 g_2 \cdot g_1 \cdot k = \varphi(g_2, g_1) [ (g_2 g_1)\cdot k ] \varphi(g_2, g_1)^\mi
 \]
for all $k \in K$.

This data gives rise to a unitary 2-representation $\alpha$ of $G$ on $\Rep(K)$, by precomposition. Explicitly,
if $\rho$ is a representation of $K$, and $g \in G$, then $\alpha_g(\rho)$ has the same underlying vector space
as $\rho$, with the action of $K$ given by
 \[
  \alpha_g (\rho) (k) = \rho(g^\mi \cdot k).
 \]
The coherence natural isomorphisms $a(g_2, g_1) \colon \alpha_{g_2} \circ \alpha_{g_1} \Rightarrow \alpha_{g_2
g_1}$ have components
 \[
  a(g_2, g_1)_\rho = \rho(\varphi(g_1^\mi, g_2^\mi))
 \]
while $a(e) \colon \id \Rightarrow \alpha_e$ is just the identity.

\subsection{Other examples of 2-representations} One expects to find similar examples of unitary 2-representations of groups arising from
automorphisms of other geometric or algebraic structures --- for instance, the automorphisms of a {\em rational vertex operator algebra} or of an {\em affine lie algebra} will act on their category of representations, which in good cases are 2-Hilbert spaces. We encourage the reader to view the slides of the talk of Mason in this regard \cite{ref:mason}.

\subsection{Morphisms of 2-representations from morphisms of exact sequences} We have seen how one obtains a
2-representation of $G$ from an exact sequence of groups (equipped with a set-theoretic section) with $G$ as the
final term, or equivalently from a weak action of $G$ on another group. A {\em morphism} of such a structure
gives rise to a morphism of 2-representations by induction. Indeed, suppose we have a map of exact sequences
 \[
  \ba \xymatrix{1 \ar[r] & K \ar[d]^{f_0} \ar[r] & E \ar[d]^{f_1} \ar[r] & G \ar[r] \ar[d]^\id & 1 \\
   1 \ar[r] & L \ar[r] & F \ar[r] & G \ar[r] & 1} \ea.
\]
In higher category language, this is essentially the same thing as a morphism inside the 2-category
 \[
  \Hom (\BG, \mathcal{G}\text{roups})
 \]
of weak 2-functors, transformations and modifications from $\BG$ (thought of as a 2-category with only identity 2-morphisms) to the 2-category of groups (objects are groups, morphisms are functors, 2-morphisms are
natural transformations). Then by inducing along $f_0$ we get a map
 \[
  \sigma \equiv \Ind(f_0) \colon \Rep(K) \rightarrow \Rep(L)
  \]
and also natural isomorphisms
 \[
  \sigma(g) \colon \beta_g \circ \sigma \Rightarrow \sigma \circ \alpha_g,
 \]
where $\alpha$ and $\beta$ are the associated 2-representations of $G$ on $\Rep(K)$ and $\Rep(L)$ respectively.
In other words, {\em a map of exact sequences gives rise to a morphism of 2-representations}.

\section{More graphical elements\label{slem}} The reason we have been drawing arrows on the strings representing the functors $\alpha_g$ involved in a 2-representation $\alpha$ is to conveniently distinguish group elements from their inverses: if a downward pointing section of a string is labeled `$g$' then it represents $\alpha_g$, and upward pointing sections of the same string represent $\alpha_{g^\mi}$. Using this convention we now construct some new graphical elements from the old ones. From now on we drop the bounding boxes on the diagrams.

Define $\eta_g \colon \id \Rightarrow \alpha_{g^\mi} \circ \alpha_g$ and $\epsilon_g \colon \alpha_g \circ \alpha_{g^\mi} \Rightarrow \id$ as:
 \[
 \eta_g = \ba \ig{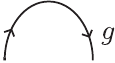} \ea := \ba \ig{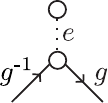} \ea  \ \ \ \equiv  \ba \xymatrix@R=0.5cm{ \id \ar@{=>}[d]^{\, a(e)} \\ \alpha_e
 \ar@{=>}[d]^{\, a(g^\mi, g)^\mi} \\ \alpha_{g^\mi} \circ \alpha_g} \ea
 \]
 \[
  \epsilon_g = \ba \ig{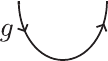} \ea := \ba \ig{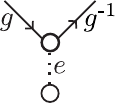}
 \ea  \ \ \ \equiv \ba \xymatrix@R=0.5cm{\alpha_g \circ \alpha_{g^\mi}
 \ar@{=>}[d]^{a(g, g^\mi)} \\ \alpha_e \ar@{=>}[d]^{a(e)^\mi} \\
 \id} \ea
 \]
These are indeed unitary natural transformations, since their inverses are clearly given by
 \[
  \eta_g^* = \ba \ig{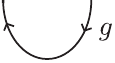} \ea := \ba \ig{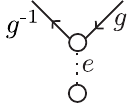} \ea \qquad
  \qquad \epsilon_g^* = \ba \ig{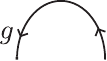} \ea := \ba \ig{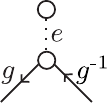}
  \ea \, .
 \]
In other words, we have the ``no
loops'' and ``merging'' rules
 \[
 \ba \ig{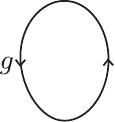} \ea = \ba \ig{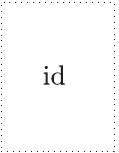} \ea \qquad \qquad \ba
 \ig{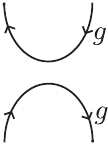} \ea = \ba \ig{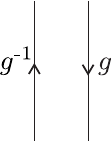} \ea
 \]
and similarly for the reverse orientations.

We now show that these new graphical elements behave as their string diagrams suggest. The first part of the following
lemma actually says in more orthodox terminology that `$\alpha_g$ is an ambidextrous adjoint equivalence from the underlying 2-Hilbert space to itself', or more precisely `for all $g \in G$, $\alpha_g \dashv \alpha_{g^\mi}$ via $(\eta_g,
\epsilon_g)$'. But it's the simple fact that these string diagrams can be manipulated in the obvious intuitive
fashion which is more important for us here.
\begin{lem} \label{ambilem} Suppose $\alpha$ is a 2-representation of $G$. The following graphical moves hold:
 \[
 \text{ \em (i)}
 \ba \ig{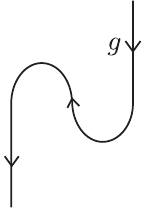} \ea = \ba \ig{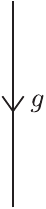} \ea = \ba \ig{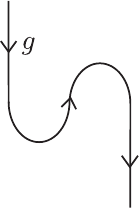} \ea
 \qquad \text{\em (ii)} \quad \ba \ig{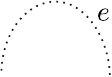} \ea = \ba \ig{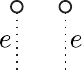} \ea
 \]
 \[
  \text{\em (iii)} \ba \ig{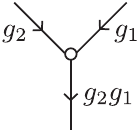} \ea = \ba \ig{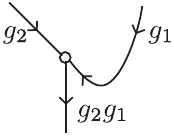} \ea = \ba \ig{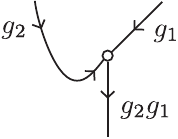} \ea  \qquad
  \text{\em (iv)} \ba \ig{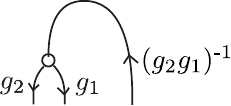} \ea = \ba \ig{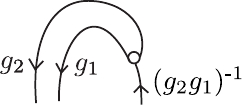} \ea
 \]
\end{lem}
\begin{proof} (i) The first equation as proved as follows,
 \[
 \ba \ig{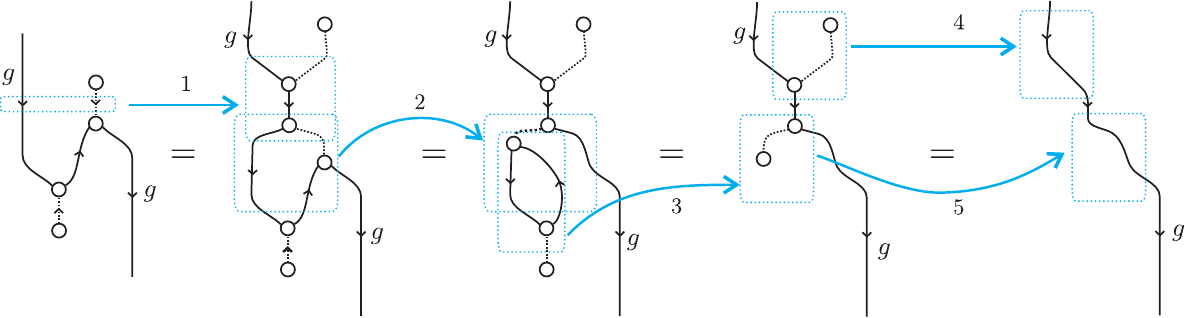} \ea .
 \]
In step 1 we zip together using the rule (\ref{zr1}a), in 2 we slide the button around using (\ref{br1}b), in 3
we unzip again using (\ref{zr1}a) and in 4 and 5 we contract the identity string using (\ref{br2}b). The other
equations are proved similarly.
 \end{proof}
Now we record for further use some allowable graphical manipulations for {\em morphisms} of 2-representations.
 \begin{lem} \label{onestlem} Suppose $\sigma \colon \alpha \rightarrow \beta$ is a morphism of 2-representations. The following graphical
 moves hold:
 \[
  \text{(i)} \ba \ig{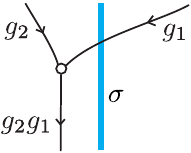} \ea = \ba \ig{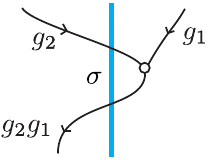} \ea \qquad \text{(ii)} \ba \ig{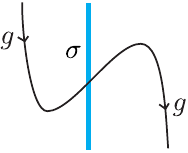} \ea = \ba
  \ig{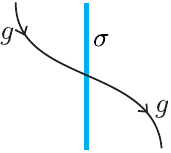} \ea
  \]
 \end{lem}
 \begin{proof} Equation (i) is proved as follows,
  \[
  \ba \ig{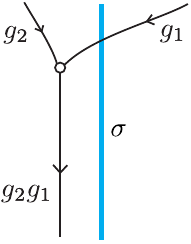} \ea \stackrel{\text{(a)}}{=} \ba \ig{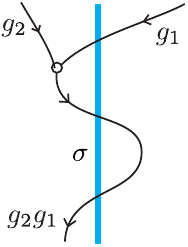} \ea
  \stackrel{\text{(b)}}{=} \ba \ig{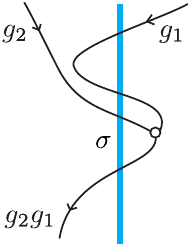} \ea
  \stackrel{\text{(c)}}{=} \ba \ig{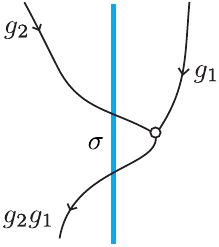} \ea \, ,
  \]
 where (a) uses the inverse rule (\ref{siginv}a), (b) uses the button-dragging rule
 (\ref{sig1}a), and (c) uses the reverse inverse rule
 (\ref{siginv}b). Equation (ii) is proved similarly, after one expands out the cups and caps $\eta_g$ and $\epsilon_g$ using their definitions we gave above.
 \end{proof}

\section{Even-handedness and unitary 2-representations\label{evun}} In this section we define what it means for a 2-representation to be compatible with a given even-handed structure, and we show that unitary 2-representations are compatible with the even-handed structure on $\THilb$. This enables certain important string diagram manipulations to be made.

Consider the category $\Rep(G)$ of unitary representations of a group $G$. If $\sigma \colon \rho_1 \rightarrow \rho_2$ is an intertwining map between two unitary representations, then the adjoint $\sigma^* \colon \rho_2 \rightarrow \rho_1$ is also an intertwining map, because
 \begin{align*}
  \sigma^* \circ \rho_2 (g) &= \sigma^* \circ \rho_2^*(g^\mi) \\
   &= (\rho_2 (g^\mi ) \circ \sigma)^* \\
   &= (\sigma \circ \rho_1 (g^\mi ))^* \\
   &= \rho_1(g) \circ \sigma^*.
  \end{align*}
The corresponding result for unitary 2-representations is a bit more subtle. Suppose that $\sigma \colon \alpha \rightarrow \beta$ is a morphism of unitary 2-representations of $G$, and that
 \[
  \sigma^* \colon H_\beta \rightarrow H_\alpha
 \]
is a functor between the underlying 2-Hilbert spaces which is adjoint to the underlying functor $\sigma \colon H_\alpha \rightarrow H_\beta$. We first need to equip $\sigma^*$ with the structure of a morphism of 2-representations
 \[
 \sigma^* \colon \beta \rightarrow \alpha.
 \]
This is easily done: choose unit and counit natural transformations $(\eta, \epsilon)$ expressing $\sigma^*$ as a right adjoint of $\sigma$, and then define the coherence isomorphisms $\sigma^*(g)$ as
 \[
 \sigma^*(g) := \sigma(g^\mi)^\dagger \colon \alpha_g \circ \sigma^* \Rightarrow \sigma^* \circ \beta_g.
 \]
In string diagrams,
 \[
 \sigma^*(g) \stackrel{\text{drawn as}}{\equiv} \ba \ig{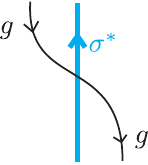} \ea := \ba \ig{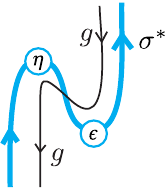} \ea = \ba \ig{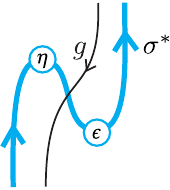} \ea
 \]
where the last simplification step uses Lemma \ref{onestlem} (ii). It is clear that this definition of $\sigma^*(g)$ satisfies the coherence equations \eqref{sig1}:
\begin{align*}
 \ba \ig{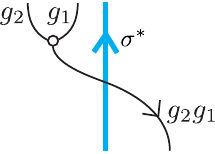} \ea := \ba \ig{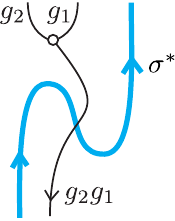} \ea &= \ba \ig{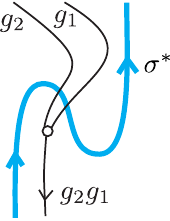} \ea \\
  &= \ba \ig{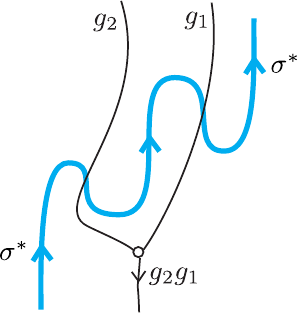} \ea =: \ba \ig{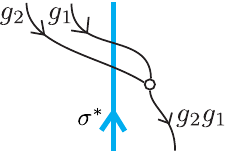} \ea.
 \end{align*}
Observe also that the inverse $\sigma^*(g)^\mi$ is given by
 \[
  \ba \ig{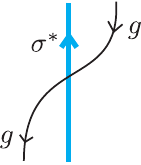} \ea := \ba \ig{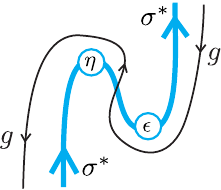} \ea.
 \]
But what if we had chosen to define $\sigma^*(g)$ via the structure of $\sigma^*$ as a {\em left} adjoint of $\sigma$ instead, using the unit and counit maps $(\Psi(\eta), \Psi(\epsilon)$ coming from the even-handed structure? We need the following definition.

\begin{defn} A 2-representation $\alpha$ of a group $G$ on an object inside an even-handed 2-category is said to be {\em compatible with the even-handed structure} if
 \[
  \Psi \left( \ba \ig{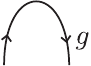} \ea, \ba \ig{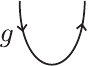} \ea \right) = \left( \ba \ig{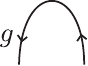} \ea, \ba \ig{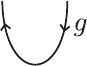} \ea \right)
 \]
for each $g \in G$, where $\left(\ba \ig{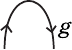} \ea, \, \ba \ig{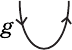} \ea\right)$ and $\left(\ba \ig{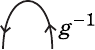} \ea, \, \ba \ig{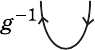} \ea\right)$ are the unit and counit maps introduced in Section \ref{slem}.
\end{defn}
Let us check that this property holds for unitary 2-representations inside $\THilb$.

\begin{lem} Unitary 2-representations are compatible with the standard even-handed structure on $\THilb$.
\end{lem}
\begin{proof} To simplify our notation, write $F \equiv \alpha_g$, $F^* \equiv \alpha_{g^\mi}$, $\eta \equiv a(g^\mi, g)^* \circ a(e)$ and $\epsilon \equiv a(e)^* \circ a(g, g^\mi)$. The standard even-handed structure $\phi \mapsto * \phi^* *$ on $\THilb$ is expressed in terms of the hom-set isomorphisms $\phi \colon \Hom(Fx, y) \rightarrow \Hom(x, F^* y)$ and not directly in terms of the unit and counit maps $\eta$ and $\epsilon$, so we first need to perform this translation. For $h \colon F^*y \rightarrow x$, we need to show that
 \[
  * \,\, \phi^* * (h) = F(h) \circ \epsilon_y^*, \quad \text{or equivalently that } \phi^*(h^*) = \epsilon_y \circ F(h^*).
 \]
For $f \colon Fx \rightarrow y$, we compute:
 \begin{align*}
  (\phi(f), h^*) &= (F^*(f) \circ \eta_x, h^*) \quad \text{(definition of $\phi$ in terms of $\eta$)} \\
   &= (FF^*(d) \circ F(\eta_x), F(h^*)) \quad \text{(since $\alpha_g$ is unitary on the hom-sets)} \\
   &= (\epsilon_y \circ FF^*(f) \circ F(\eta_x), \epsilon_y \circ F(h^*)) \quad \text{(since $\epsilon$ is unitary)} \\
   &= (f \circ  \epsilon_{Fx} \circ F(\eta_x), \epsilon_y \circ F(h^*)) \quad \text{(naturality of $\epsilon$)} \\
   &= (f, \epsilon_y \circ F(h^*)) \quad \text{(snake diagram, Lemma \ref{ambilem} (i))},
  \end{align*}
which is what we needed to show.
\end{proof}
The fact that unitary 2-representations are compatible with the even-handed structure on $\THilb$ means that we have
 \[
 \sigma^*(g) = \ba \ig{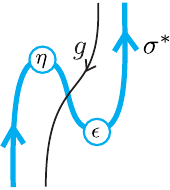} \ea = \ba \ig{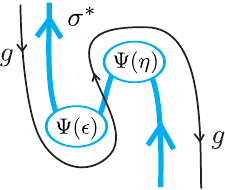} \ea
 \]
because the former is $\sigma(g^\mi)^\dagger$ and the latter is $^\dagger \! \sigma(g^\mi)_\Psi$. This allows us to establish the following important diagrammatic lemma, which we will need in Section \ref{funcchar}.
\begin{lem} \label{invproof} With this definition of $\sigma^*(g)$, the following equations hold:
  \begin{align*}
  &\text{(i)}\;\; \ba \ig{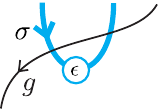} \ea = \ba \ig{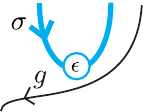} \ea &  &\text{(ii)} \;\; \ba \ig{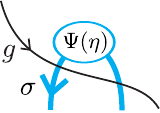} \ea = \ba \ig{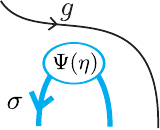} \ea \\
  \end{align*}
\end{lem}
\begin{proof} Equation (i) is almost a tautology:
 \[
  \ba \ig{e367.pdf} \ea = \ba \ig{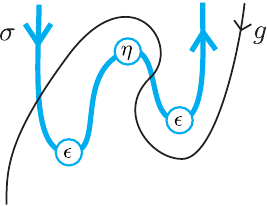} \ea = \ba \ig{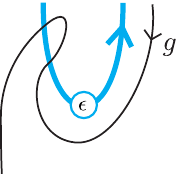} \ea = \ba \ig{e368.pdf} \ea.
  \]
Equation (ii) uses the fact that unitary 2-representations are compatible with the even-handed structure:
 \[
  \ba \ig{e369.pdf} \ea = \ba \ig{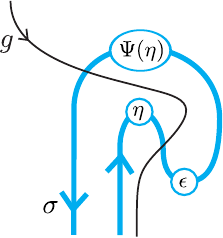} \ea = \ba \ig{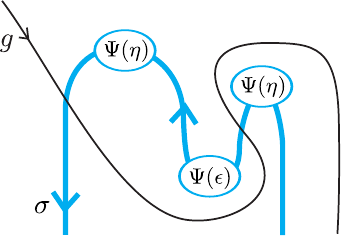} \ea = \ba \ig{e370.pdf} \ea.
  \]
\end{proof}

\section{2-characters of 2-representations\label{charsecc}}
\XXX{rewrite!} In this section we define and study the notion of the 2-character of a unitary 2-representation. The idea of the 2-character of a weak group action on a linear category has also been defined independently by Ganter and Kapranov \cite{ref:ganter_kapranov_rep_char_theory}.  Our treatment is novel in two important aspects. Firstly we show how 2-characters look especially simple when expressed in terms of string diagrams. Secondly, we use the canonical even-handed structure on $\THilb$ to show how the 2-character can be made {\em functorial} with respect to morphisms of 2-representations, as we explained in the introduction. This prepares the way for us to show, in Chapter \ref{GerbesCharChap}, how the 2-character of a unitary 2-representation corresponds to the `geometric character' of its associated equivariant gerbe, from which it will follow that the complexified 2-character is unitarily fully faithful.

\subsection{2-traces}
The basic idea of 2-characters, as we explained in the introduction, is that they categorify the notion of the character of an ordinary representation of a group. Ordinary characters are defined by taking traces, so we first need to define {\em 2-traces} (Ganter and Kapranov called this the {\em categorical trace}).

\begin{defn} The {\em 2-trace} of a linear endofunctor $F \colon H \rightarrow H$ on a 2-Hilbert space $H$ is the Hilbert space
 \[
 \ttr{F} = \Nat(\id_H, F) = \left\{ \ba \ig{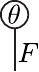} \ea \right\}.
 \]
 \end{defn}
If one thinks of $F$ via its associated matrix of Hilbert spaces $\Hom(e_j, Fe_i)$ where $e_i$ runs over a choice of simple objects for $H$, then the 2-trace corresponds to the direct sum of the Hilbert spaces along the diagonal, because a natural transformation $\id \Rightarrow F$ is freely and uniquely determined by its behaviour on the simple objects. Also recall from Chapter \ref{BigTCat} that the inner product on $\ttr{F}$ is defined as
 \[
   \langle \theta, \theta'\rangle = \sum_{i} k_i (\theta_{e_i}, \theta'_{e_i})
 \]
where $e_i$ runs over a choice of simple objects for $H$, and the weightings $k_i = (\id_{e_i}, \id_{e_i})$ are their scale factors.

\subsection{The loop groupoid\label{loopgdsec}}
\begin{figure}[t]
\centering
\ig{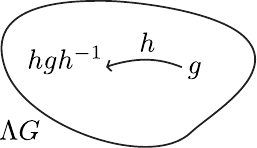}
\caption{\label{lgfig}The loop groupoid of a finite group}.
\end{figure}

In general the {\em loop groupoid} $\Lambda \mathcal{G}$ of a finite groupoid $\mathcal{G}$ is the category $\Fun(\mathbb{Z}, \mathcal{G})$ of functors and natural transformations from the group $\mathbb{Z}$ of integers, thought of as a one object category, to $\mathcal{G}$ (see \cite{ref:simon}). The reason for the name `loop groupoid' is that one thinks of $\mathbb{Z}$ as the homotopy group of the circle, so that the objects of the loop groupoid can be thought of as `loops' in $\mathcal{G}$. A special case is the loop groupoid $\Lambda \BG$ of a finite group $G$ (recall that $\BG$ refers to the group $G$ thought of as a one object category) which depicts the action of the group on itself by conjugation, since the objects can be identified with the elements $g \in G$, and the morphisms can be written as $hgh^\mi \sla{h} g$ (see Figure \ref{lgfig}). 

\subsection{2-characters\label{sec2char}}The loop groupoid is particularly convenient when discussing characters. The fact that the ordinary character $\chi_\rho$ of a representation $\rho$ of $G$ is conjugation invariant can be expressed by saying it is a function
 \[
  \chi_\rho \colon [\Lambda \BG] \rightarrow \mathbb{C}
 \]
from the set of connected components of the loop groupoid --- the conjugacy classes of $G$ --- to the complex numbers.
Similarly the 2-character $\chi_\alpha$ of a unitary 2-representation $\alpha$ will produce a {\em unitary equivariant vector bundle} over the group, that is, a {\em unitary representation} of the loop groupoid. By a unitary representation of the loop groupoid, we mean a functor
 \[
  \chi_\alpha \colon \Lambda \BG \rightarrow \Hilb
 \]
which sends morphisms in $\Lambda \BG$ to unitary maps in $\Hilb$. That is, we are defining a `unitary equivariant vector bundle over the group $G$' to be a unitary representation of the loop groupoid $\Lambda \BG$. We will write the category of unitary equivariant vector bundles over $G$ as $\Hilb_G(G)$.

\begin{defn} The {\em 2-character} $\chi_\alpha$ of a unitary 2-representation $\alpha$ of $G$ is the unitary equivariant vector bundle over $G$ given by
 \[
  \begin{split}
   \chi_F(g)  = \ttr{\alpha_g} &= \left\{ \ba\ig{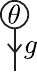} \ea \right\} \\
    \chi_\alpha(hgh^\mi \sla{h} g)\left( \ba \ig{c3.pdf}\ea \right) &=
    \ba\ig{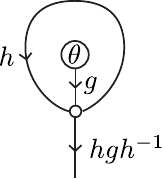} \ea.
  \end{split}
  \]
  \end{defn}
Let us verify that this definition makes sense.
\begin{prop} The 2-character $\chi_\alpha$ is indeed a unitary equivariant vector bundle over the group.
\end{prop}
\begin{proof} Using our graphical rules from Lemma \ref{ambilem}, we have
\begin{align*}
\chi_\alpha ( \sla{h_2} h_1 g h_1^\mi) \chi_\alpha (\sla{h_1} h) \left ( \ba
\ig{c3.pdf} \ea \right) &= \ba \ig{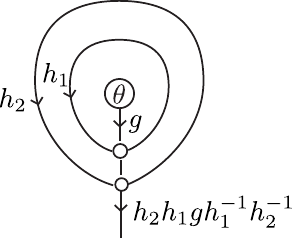} \ea =  \ba \ig{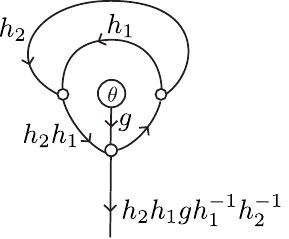} \ea \\ &= \ba\ig{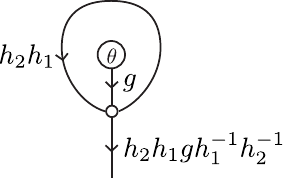} \ea
= \chi_\alpha( \la{h_2h_1} g)\left( \ba \ig{c3.pdf} \ea \right)
\end{align*}
and also
\[
 \chi_\alpha ( \sla{e} g) \left( \ba \ig{c3.pdf} \ea \right) = \ba \ig{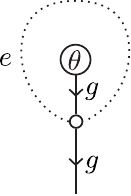} \ea = \ba \ig{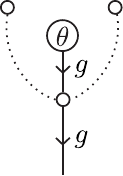} \ea =
 \ba \ig{c3.pdf} \ea.
\]
It is a {\em unitary} vector bundle because all the maps involved in its definition are unitary (see the proof of Lemma \ref{transglem} for an explicit formula).
\end{proof}

\subsection{Functoriality of the 2-character \label{funcchar}}
In this subsection we combine all the string diagram technology we have developed so far to define how to take the 2-character of a {\em morphism} of unitary 2-representations so as to obtain a morphism of the corresponding equivariant vector bundles over $G$. Moreover we show that this construction is functorial with respect to composition of morphisms in $\TRep(G)$.
\begin{defn} \label{2charmorphisms} If $\sigma \colon \alpha \rightarrow \beta$ is a morphism of unitary 2-representations, we define $\chi(\sigma) \colon \chi_\alpha \rightarrow \chi_\beta$ as the map of equivariant vector bundles over $G$ whose component at $g \in G$ is given by
 \be \label{stdiagmor}
 \begin{aligned}
  \chi(\sigma)_g \colon \chi_\alpha(g) & \rightarrow \chi_\beta (g) \\
  \ba \ig{c3.pdf} \ea & \mapsto \ba \ig{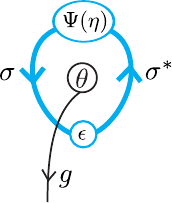} \ea
  \end{aligned}
 \ee
where $(\sigma^*, \eta, \epsilon)$ is some choice of right adjoint for $\sigma$, and $\Psi$ is the even-handed structure on $\THilb$ from Chapter \ref{EvenTHilb}.
 \end{defn}
We now verify that this definition makes sense.
We write $[\TRep(G)]$ for the {\em Grothendieck category} of $\TRep(G)$, defined as the category which has the same objects as $\TRep(G)$ but whose morphisms are {\em isomorphism classes} of 1-morphisms in $\TRep(G)$.
\begin{thm}\label{Thecharthm} In the situation above, the map $\chi(\sigma) \colon \chi_\alpha \rightarrow \chi_\beta$
 \begin{enumerate}
  \item does not depend on the choice of right adjoint $(\sigma^*, \eta, \epsilon)$ for $\sigma$,
  \item is indeed a morphism of equivariant vector bundles over $G$,
  \item does not depend on the isomorphism class of $\sigma$,
  \item is functorial with respect to composition of 1-morphisms in $\TRep(G)$,
  \end{enumerate}
 and hence $\chi$ descends to a functor
  \[
   \chi \colon [\TRep(G)] \rightarrow \Hilb_G (G).
  \]
 \end{thm}
 \begin{proof} (i) This follows from the naturality of $\Psi$. If $((\sigma^*)', \eta', \epsilon')$ is another right adjoint for $\sigma$, then there is a canonical natural isomorphism $\gamma \colon \sigma^* \Rightarrow (\sigma^*)'$ which transforms $(\eta, \epsilon)$ into $(\eta', \epsilon')$. Thus we have:
   \[
    \ba \ig{y28.pdf} \ea = \ba \ig{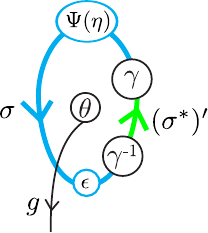} \ea \; = \ba \ig{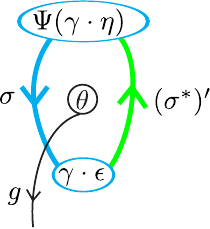} \ea.
   \]

   (ii) Using the graphical rules in Lemmas \ref{ambilem} and \ref{invproof}, we calculate:
   \[ \begin{aligned}
    \chi_\beta (\sla{h} g) \chi(\sigma)_g \left( \ba \ig{c3.pdf} \ea \right) = \ba \ig{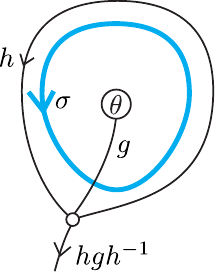} \ea = \ba
    \ig{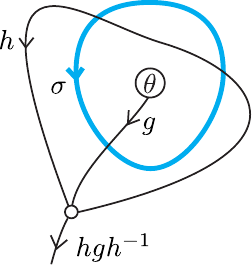} \ea = \ba \ig{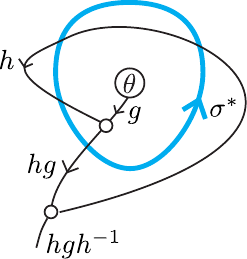} \ea \\ = \ba \ig{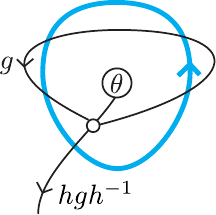} \ea = \ba \ig{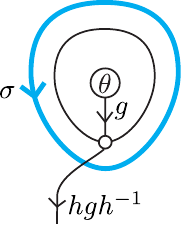} \ea = \chi(\sigma)_x
    \chi_\alpha(\sla{h} g) \left( \ba \ig{c3.pdf} \ea\right).
    \end{aligned}
   \]
   (iii) Suppose $\gamma \colon \sigma \Rightarrow \rho$ is an invertible 2-morphism in $\TRep(G)$, and that $(\rho^*, \eta', \epsilon')$ is a right adjoint for $\rho$. Then
 \begin{multline*}
    \ba \ig{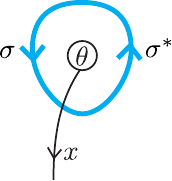} \ea = \ba \ig{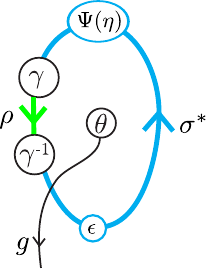} \ea = \ba \ig{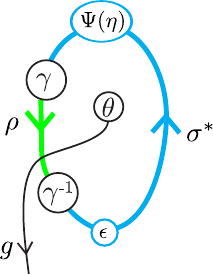} \ea \\ = \ba \ig{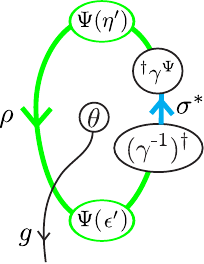} \ea = \ba \ig{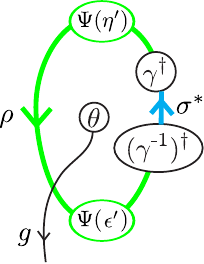}  \ea
         = \ba \ig{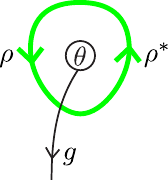} \ea
    \end{multline*}
 where the second step uses the fact that $\gamma^\mi$ is a 2-morphism in $\TRep(G)$ and the fourth step uses the even-handed equation.

  (iv)  This follows from the fact that an even-handed structure respects composition.
   \end{proof}
We hope that these diagrammatic proofs have convinced the reader of the utility of the string diagram notation. Explicit {\em expansions} of these diagrams in terms of concrete formulas can be found in Chapter \ref{2char2rep}. We remark here that the behaviour of the 2-character on morphisms really does use the canonical even-handed structure on $\THilb$ in an essential way. That is, if $\sigma \colon \alpha \rightarrow \beta$ is a morphism of 2-representations, then $\chi(\sigma) \colon \chi_\alpha \rightarrow \chi_\beta$ {\em depends} on the scale factors $k_i$ and $k_\mu$ on the simple objects in $H_\alpha$ and $H_\beta$. One way to see this is that our main theorem states that the 2-character of a unitary 2-representation corresponds to the geometric character of its associated equivariant gerbe --- and the behaviour of the geometric character on morphisms of equivariant gerbes really does depend on the metrics on the gerbes (see Chapter \ref{fgerbes}).

\chapter[Gerbes and their geometric characters]{Finite equivariant gerbes and their geometric characters\label{GerbesChap}}
In this chapter we introduce the main geometric actors of this thesis --- finite equivariant gerbes equipped with metrics, the 2-category which they constitute, and their geometric characters. As we discussed in the introduction, we encourage the reader to think of the notion of a finite equivariant gerbe equipped with a metric as the `categorification', at least in our finite discrete toy model setting, of the notion of an {\em equivariant hermitian holomorphic line bundle} over a compact hermitian manifold.

We hope that the notions we introduce in this chapter will have counterparts in the {\em smooth} setting, which would involve gerbes over smooth spaces equivariant under the action of a compact Lie group. In that setting one would probably want to use the language of {\em differentiable stacks}, and indeed the necessary technology to upgrade all the concepts we introduce in this chapter to the setting of stacks is most likely already available in \cite{ref:caldararu_willerton, ref:behrend_xu, ref:weinstein, ref:ganter_kapranov_rep_char_theory}. In addition, as mentioned in the Introduction, Sati, Schreiber, Sk\u{o}da and Stevenson have recently outlined a theory of nonabelian cohomology in a smooth $\omega$-groupoid setting \cite{ref:sati_et_al}, thereby generalizing the case $n=2$ (which is what we are considering here) to higher $n$. We leave this for future work.

In Section \ref{BigR} we define what we mean by a discrete equivariant gerbe equipped with a metric. In Section \ref{exeqgerbes} we give a useful example of such a gerbe, whose objects are the isomorphism classes of $U(1)$-bundles with connection over a connected smooth manifold $M$. In Section \ref{tired} we define the 2-category of finite equivariant gerbes, and we classify gerbes up to {\em isometric equivalence}. In Section \ref{uonebund} we introduce some more geometric language in the form of $U(1)$-bundles and line bundles on finite groupoids, and their associated spaces of sections.

We have encouraged the reader to think of an equivariant gerbe as the categorification of a hermitian line bundle; in Section \ref{tgess} we state the important theorem of Willerton \cite{ref:simon} which {\em links} these two pictures. It says that there is a unitary isomorphism --- the {\em twisted character map} --- between the space of isomorphism classes of equivariant vector bundles over an equivariant gerbe and the space of flat sections of a certain {\em transgressed line bundle}.  The last two sections define the geometric analogue of the 2-character of a 2-representation we introduced in Chapters \ref{sec2char} and \ref{funcchar}.  We define how to take the {\em geometric character} of a $G$-equivariant gerbe equipped with a metric in order to obtain a unitary equivariant vector bundle over $G$, and we show how to make this construction {\em functorial} with respect to morphisms of equivariant gerbes. We show that the geometric character descends to a functor from the Grothendieck category of $\Gerbes(G)$ to the category of equivariant vector bundles over $G$, and we use the theorem of Willerton to show that this functor is unitarily fully faithful after one tensors the hom-sets with $\mathbb{C}$. After we have established the correspondence between the 2-character of a unitary 2-representation and the geometric character of its associated equivariant gerbe in Chapter \ref{GerbesCharChap}, this will imply that the complexified 2-character is also unitarily fully faithful.

\section{Equivariant gerbes\label{BigR}} In this section we define discrete equivariant gerbes. Our definition agrees with the conventional notion of an `equivariant gerbe' as in the notes of Behrend and Xu \cite{ref:behrend_xu} after specializing their notion to this simplified setting --- though we also add in the idea of a {\em metric}.

\subsection*{$U(1)$-torsors and their tensor products} A $U(1)$-torsor is a set with a free and transitive left
action of $U(1)$. The tensor product $P \otimes Q$ of two $U(1)$-torsors is the torsor obtained from the
cartesian product $P \times Q$ by identifying $(e^{i \theta} p, q)$ with $(p, e^{i \theta} q)$ for any $p \in
P$, $q \in Q$ and $e^{i \theta} \in U(1)$; the equivalence class of $(p,q)$ is denoted $p \otimes q$.

\subsection*{Equivariant gerbes \label{EqGerbeSec}}
\begin{figure}
\center
\ig{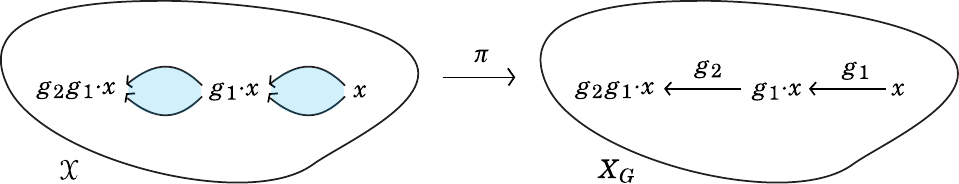}
\caption{\label{gerbefig}An equivariant gerbe $\X$. Each arrow in the action groupoid $X_G$ becomes a $U(1)$-torsor worth of arrows in $\X$, which are the shaded regions.}
\end{figure}

Let $X$ be a set acted on from the left by a group $G$. We think of $X$ via its associated {\em action groupoid} $X_G$, whose objects are the elements $x \in X$ and whose morphisms depict the action of $G$ on $X$, so that for each $x \in X$ and $g \in G$ we have a morphism $g \cdot x \sla{g} x$. Also let $\underline{U(1)}$ be the trivial `bundle of groups' on $X$ (see \cite{ref:moerdijk1}); as a groupoid its objects are the elements $x \in X$ and its hom-sets are given by $\Hom(x,y) = U(1)$ for $x = y$ and the empty set otherwise.
\newpage
\begin{defn} A {\em discrete $G$-equivariant gerbe} is a $U(1)$-central extension $\X$ (in the sense of Moerdijk \cite{ref:moerdijk1}) of the action groupoid $X_G$ associated to a $G$-set $X$:
 \[
   \underline{U(1)} \stackrel{i}{\hookrightarrow} \X \sra{\pi} X_G.
 \]
A {\em metric} on the gerbe is an assignment of a positive real number $k_x$ to each object $x \in \X$, invariant under the action of $G$.
\end{defn}
That is, when we say that $\X$ is `$U(1)$-central extension' of $X_G$ we mean that $\X$ has the same objects as $X_G$, $\pi$ is a full surjective functor and $i$ is an isomorphism onto the subgroupoid of arrows in $\X$ which project to identities in $X_G$ (see Figure \ref{gerbefig}). In other words, each arrow $\sla{g} x$ in the action groupoid $X_G$ becomes a {\em $U(1)$-torsor} worth of arrows in the equivariant gerbe $\X$; we write this torsor as
 \[
  \X_{\sla{g} x } := \pi^{-1} (\sla{g} x).
   \]
We shall use a non-calligraphic $X$ to refer to the underlying set of objects of an equivariant gerbe $\X$, and we say that the equivariant gerbe is {\em finite} if $X$ is a finite set.

\subsection*{Equivariant gerbes and cohomology}
Our notion of an `equivariant gerbe' is equivalent to the cocycle description of gerbes given by Willerton in \cite{ref:simon}. To see this, let us define a {\em section} of the gerbe as a set-theoretic map
\[
  s \colon \text{Arr}(X_G)  \rightarrow  \text{Arr}(\X) \\
\]
such that $\pi \circ s = \id$. Choosing a section gives rise to a  $U(1)$-valued 2-cocycle $c \in Z^2(X_G,
U(1))$ on the groupoid $X_G$, in the sense of Willerton \cite{ref:simon}. One defines $c$ as the correction term needed to balance the composition equation:
 \[
  s(\sla{g_2} g_1 \cdot x) \circ s(\sla{g_1} x) = c_x (g_2, g_1)
  s(\sla{g_2g_1} x).
 \]
Choosing a different section $s'$ will change $c$ by a coboundary, so that an equivariant gerbe $\X$ gives
rise to a cohomology class $[c_\X] \in H^2(X_G, U(1))$. In this thesis we have avoided working with cocycles, taking the viewpoint that the more fundamental geometric structure is the equivariant gerbe itself.

\subsection*{Tensor product of equivariant gerbes} Suppose $\Y$ and $\X$ are equivariant gerbes equipped with metrics. Note that we think of $\Y$ and $\X$ as structures in their own right; they have different underlying $G$-sets $Y$ and $X$ and we are not assuming them to be `gerbes over the same space'.  We define the {\em tensor product} $\Y \otimes \X$ as the equivariant gerbe with metric whose object set is the cartesian product $Y \times X$ on which $G$ acts diagonally, whose $g$-graded morphisms are the tensor product of those of $\Y$ and $\X$,
 \[
  (\Y \otimes \X)_{\sla{g} (y, x)} := \Y_{\sla{g} y} \otimes \X_{\sla{g} x},
  \]
and whose metric is the product metric. Also, if $\X$ is an equivariant gerbe, we write $\overline{\X}$ for the equivariant gerbe having the same underlying groupoid as $\X$ but with the conjugate action of $U(1)$ on its hom-sets.

\section{Examples of equivariant gerbes \label{exeqgerbes}} In this thesis, the main examples of discrete equivariant gerbes will be those arising from unitary 2-representations of $G$, but here is another example, which can be thought of as a reformulation of the ideas in the second chapter of Brylinski \cite{ref:brylinski}.

Suppose $M$ is a connected smooth manifold. Consider the groupoid $\mathcal{P}_M$ whose objects $(P, \nabla)$ are $U(1)$-bundles with
connection over $M$, and whose morphisms $f \colon (P, \nabla) \rightarrow (P', \nabla')$ are
diffeomorphisms $f \colon P \rightarrow P'$ which respect the action of $U(1)$ and which preserve parallel transport. If
there {\em is} an isomorphism from $(P, \nabla)$ to $(P', \nabla')$, then any other isomorphism must differ from
it by a constant factor in $U(1)$, since the maps must preserve parallel transport. Thus the nonempty hom-sets in
$\mathcal{P}_M$ are $U(1)$-torsors. In other words, $\mathcal{P}_M$ is an equivariant gerbe for the trivial group, i.e. an `ordinary gerbe'.

Now suppose a group $G$ acts on $M$ by diffeomorphisms. Let $\Pic(M)$ denote the set of isomorphism classes in
$\mathcal{P}_M$ (this set is isomorphic to a Deligne cohomology group), and suppose one chooses a distinguished
representative $(P, \nabla)_c$ for each isomorphism class $c \in \Pic^\nabla(M)$.  Now, the group $G$ acts on
$\mathcal{P}_M$ by push-forward. That is,
 \[
  g \ccdot (P, \nabla) = (g \ccdot P, g \ccdot \nabla)
 \]
where $g \ccdot P$ is the principal $U(1)$-bundle over $M$ whose fiber at $m$ is $P_{g^\mi \ccdot m}$, with $g \cdot \nabla$ the connection whose parallel transport map $(g \ccdot P)_m \rightarrow (g \ccdot P_{m'}$ along a path $\gamma \colon m \rightarrow m'$ is given by parallel transporting $P$ along $g^\mi \ccdot \gamma$.

Therefore $G$ also acts on $\Pic^\nabla(M)$. This gives rise to an associated equivariant gerbe $\X$ over the action groupoid $\Pic^\nabla(M)_G$ by a variant of the well-known Grothendieck construction for producing a fibration of categories from a weak 2-functor on a category taking values in $\Cat$  \cite{ref:grothendieck}. The objects of
$\X$ are the equivalence classes $c \in \Pic^\nabla(M)$ while the $g$-graded morphisms are given by
 \[
  \X_{\sla{g} c} := \Hom_{P_M} ( (L, \nabla)_{g \cdot c}, \,\, g_* (L, \nabla)_c).
 \]
Composition of $f_2 \in \X_{\sla{g_2} g_1 \cdot c}$ and $f_1 \in \X_{\sla{g_1} c}$ is defined by
 \[
 f_2 \diamond f_1 := \alpha_{g_2}(f_1) \, \circ \, f_2,
  \]
where we have suppressed the canonical isomorphisms ${g_2}_* {g_1}_* (L, \nabla) \cong (g_2 g_1)_* (L, \nabla)$. In
particular, it makes it clear that the stabilizer group of equivalence class of line bundle in the equivariant gerbe $\X$ is a central extension of the stabilizer group $H \subseteq G$ which fixes that class when it is considered as an element of $\mathcal{P}_M$. Many interesting central extensions
of groups arise in this way.

\section{The 2-category of equivariant gerbes\label{tired}} In this section we define the 2-category of equivariant gerbes, which is to be thought of as an equivariant version of the 2-category $\FinSpaces$ introduced in Chapter \ref{AGeomSec}.

\subsection*{Unitary vector bundles over equivariant gerbes} 
\begin{figure}[t]
\centering
\ig{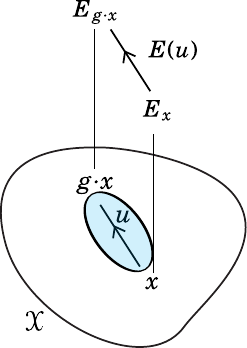}
\caption{\label{bungerbe}  A unitary equivariant vector bundle over an equivariant gerbe.}
\end{figure}

A {\em unitary equivariant vector bundle} $E$ over an
equivariant gerbe $\X$ is a functor $E \colon \X \rightarrow \Hilb$ which maps arrows in $\X$ to unitary maps in
$\Hilb$, and which preserves the $U(1)$-action on the hom-sets, so that $E(e^{i \phi} v) = e^{i \phi} E(v)$ for every arrow $v$ in the gerbe and $e^{i \phi} \in U(1)$; see Figure \ref{bungerbe}. A morphism $\theta \colon E \rightarrow E'$ of equivariant vector bundles over $\X$ is a natural
transformation; we write $\Hilb(\X)$ for the category of unitary equivariant vector bundles over $\X$.
If we choose a set-theoretic section $s$ of $\X$, giving rise to a groupoid 2-cocycle $c \in Z^2(X_G, U(1))$, then a
unitary equivariant vector bundle over the gerbe $\X$ can be regarded as a {\em $c$-twisted equivariant vector bundle}
$\hat{E}$ over the action groupoid $X_G$ in the sense of Willerton \cite{ref:simon}, using the prescription $\hat{E}(\sla{g} x) = E(s(\sla{g} x))$. Functoriality of $E$ means
that $\hat{E}$ is $c$-twisted functorial,
 \[
  \hat{E}(\sla{g_2} g_1 \cdot x) \hat{E}(\sla{g_1} x) = c_x(g_2, g_1) \hat{E}(\sla{g_2g_1} x).
 \]
This highlights one of the advantages of working with equivariant gerbes (as opposed to working with $G$-sets and cocycles) --- a {\em twisted} representation in the cocycle world is just an {\em ordinary} representation in the gerbes world. We also remark that, as shown by Willerton \cite{ref:simon}, a twisted equivariant vector bundle $\hat{E}$ over $X_G$ can also be thought of as a twisted representation of the {\em arrow algebra} $\mathbb{C}[X_G]$ of the groupoid. We follow Willerton and avoid using this algebraic description, for the simple reason that we prefer an inherently geometric formulation, but also because we want to avoid using cocycles.

\subsection*{The 2-category of equivariant gerbes}
\begin{figure}[t]
\centering
\ig{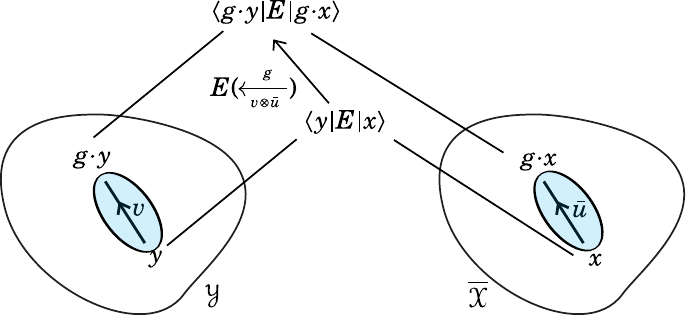}
\caption{\label{morfigg}A morphism of equivariant gerbes.}
\end{figure}

We now define the 2-category $\Gerbes(G)$ of finite $G$-equivariant gerbes equipped with metrics, which should be thought of as an equivariant version of the 2-category $\FinSpaces$ of finite spaces we introduced in Chapter \ref{AGeomSec}.
\begin{defn} For a finite group $G$, the 2-category $\Gerbes(G)$ of {\em finite $G$-equivariant gerbes} is defined as follows. An object is a
finite $G$-equivariant gerbe $\X$ equipped with a metric. A morphism $E \colon \X \rightarrow \Y$ of equivariant gerbes is a unitary equivariant vector bundle $E \colon \Y \otimes \overline{\X} \rightarrow \Hilb$, with composition given by convolution of vector bundles. A 2-morphism is a morphism of equivariant vector bundles.
\end{defn}
Let us explain this definition; see Figure \ref{morfigg}. A morphism of equivariant gerbes $E \colon \X \rightarrow \Y$ is to be thought of as the categorification of the idea of an {\em equivariant kernel} from Chapter \ref{GeometryOrdinaryRepChap} --- instead of the amplitude for going from $x \in \X$ to $y \in \Y$ being a mere {\em number} (more precisely, an element in a complex line), the amplitude is now a {\em vector space}, which we will write as
  \[
  \langle y | E | x \rangle.
 \]
Moreover, these amplitude vector spaces $\langle y | E | x \rangle$ are equivariant with respect to $x$ and $y$. That is, if $\bar{u} \in \overline{\X}_\mlag{x}{}{g}$ and $v \in \Y_\mlag{y}{}{g}$ are arrows in $\overline{\X}$ and $\Y$ which belong to the torsor of arrows emanating out of $x$ and $y$ in the direction of $g \in G$, then there is a prescribed unitary isomorphism
 \[
\langle y | E | x \rangle(\lag{g}{v \otimes \bar{u}}) \colon \langle y | E | x \rangle \rightarrow \langle g \ccdot y | E | g \ccdot x \rangle.
\]
In other words, a morphism of equivariant gerbes $E \colon \X \rightarrow \Y$ is an equivariant version of the notion of a morphism of finite spaces which we described in Chapter \ref{AGeomSec}. Just as in Chapter \ref{AGeomSec}, we compose morphisms $E \colon \X \rightarrow \Y$ and $F \colon \Y \rightarrow \Z$ of equivariant gerbes by setting $F \circ E$ to be the unitary equivariant vector bundle over $\Z \otimes \overline{\X}$ whose vector space of ways of going from $x \in \X$ to $z \in \Z$ is given by the weighted direct sum over all $y \in \Y$ of the tensor product amplitudes:
 \be \label{2kerncomp}
  \langle z | F \circ E | x \rangle := \hat{\displaystyle \bigoplus_{y \in \Y}} \, \frac{1}{k_y} \; \langle z | F | y \rangle \otimes \langle y | E | x \rangle.
 \ee
As an equivariant vector bundle over $\Z \otimes \overline{\X}$, the equivariant maps
\[
 \langle z | F \circ E | x \rangle(\lag{g}{w \otimes \bar{u}}) \colon \langle z | F \circ E | x \rangle \rightarrow \langle g \ccdot z | F \circ E | g \ccdot x \rangle, \quad w \in \Z_\mlag{z}{}{g} \,, \,  \bar{u} \in \overline{\X}_\mlag{x}{}{g}
\]
for $F \circ E$ are defined by simply taking the tensor product of the constituent equivariant maps:
\[
 \ba \xymatrix{ \langle z | F \circ E | x \rangle \ar[d]^{ (F \circ E)(\lag{g}{w \otimes \bar{u}}) } \\ \langle g \ccdot z | F \circ E | x \rangle} \ea
  :=
  \ba \xymatrix{ \hat{\displaystyle \bigoplus_{y \in \Y}} \langle z | F | y \rangle \otimes \langle y | E | x \rangle \ar[d]^{F(\lag{g}{w \otimes \bar{v}}) \otimes E(\lag{g}{v \otimes \bar{u}})} \\ \hat{\displaystyle \bigoplus_{y \in \Y}} \langle g \ccdot z | F | g \ccdot y \rangle \otimes \langle g \ccdot y | E | g \ccdot x \rangle} \ea
 \]
Note that this formula does not depend on the choices made for $v \in \Y_\mlag{y}{}{g}$ since it appears both linearly and antilinearly in the expression. Everything else works the same as in the 2-category $\FinSpaces$ we described in Chapter \ref{AGeomSec}; this completes our explanation of the 2-category $\Gerbes(G)$.

\subsection*{The inner product on the 2-morphisms} We can equip the 2-morphisms in $\Gerbes(G)$ with a natural inner product just as we did in equation \eqref{myip} for the 2-morphisms in $\FinSpaces$. That is, if $E, F \colon \X \rightarrow \Y$ are morphisms of finite equivariant gerbes equipped with metrics, and $f, g \colon E \Rightarrow F$ are 2-morphisms between them, then we set
 \be \label{myip2}
   \langle f, g \rangle := \sum_{x \in \X, y \in \Y} k_x k_y \Tr(f_{x,y}^* g).
 \ee

\subsection*{Classification of equivariant gerbes} We say two equivariant gerbes are {\em equivalent} if they are equivalent in the 2-category $\Gerbes(G)$; if they have metrics then we say they are {\em isometrically equivalent} if the 1-morphism $E\colon \X \rightarrow \Y$ furnishing the equivalence --- which is by definition an equivariant vector bundle over $\Y \otimes \overline{\X}$ --- is only supported on pairs of objects which share the same scale factor. The following result classifies equivariant gerbes up to isometric equivalence.
\begin{prop} \label{equivgerbes} Suppose $\X$ and $\Y$ are equivariant gerbes equipped with metrics. The following are equivalent:
 \begin{enumerate}
  \item $\X$ is isometrically equivalent to $\Y$.
  \item There exists an isomorphism of $G$-sets $f \colon X \rightarrow Y$, preserving the scale factors, such that $[c_{\X}] = [f^*(c_{\Y})]$ as cohomology classes in $H^2(X_G, U(1))$.
 \end{enumerate}
\end{prop}
\begin{proof} (i) $\Rightarrow$ (ii). Suppose we have an isometric equivalence
 \[
  \ba \xymatrix{\X \ar @/^1pc/[r]^E & \Y \ar
  @/^1pc/[l]^{F}} \ea.
 \]
Since $F \circ E \cong \id_\X$ and $E \circ F \cong \id_\Y$, the matrices
 \[
 [\dim E]_{y,x} = \dim \langle y | E | x \rangle, \quad [\dim F]_{x, y} = \dim \langle x | F | y \rangle
 \]
whose entries are the dimensions of the fibers of $E$ and $F$ must be inverses of each other:
 \[
  [\dim F] [\dim E] = 1 , \quad [\dim E] [\dim F] = 1.
 \]
Since the entries of these matrices are nonnegative integers, the only possibility is that they are permutation
matrices and hence define a map $f : X \rightarrow Y$ by setting $f(x)$ equal to the unique $y \in Y$
such that $\dim \langle y | E | x \rangle = 1$. Since $E$ is an equivariant vector bundle, we must have $f(g \ccdot
x) = g \ccdot f(x)$ so that $f$ is an isomorphism of $G$-sets; moreover since we are told that $\X$ is {\em isometrically} equivalent to $\Y$ we must have that $f$ preserves the scale factors on $X$ and $Y$.

Now choose sections $s \colon X_G \rightarrow \X$ and $s' \colon Y_G \rightarrow \Y$ of $\X$ and $\Y$, with resulting groupoid 2-cocycles $c \equiv c_{\X} \in Z^2(X_G, U(1))$ and $c' \equiv c_{\Y} \in Z^2(Y_G, U(1))$, and write $\hat{E}$ for the corresponding twisted equivariant vector bundle over $(Y \times X)_G$ associated to $E$. Thus for all $x \in \X$ and $y \in \Y$ we have
 \[
   \hat{E}(\sla{g_2} (g_1 \ccdot y, g_1 \ccdot x)) \circ \hat{E}(\sla{g_1} (y,x)) = \frac{c'_y(g_2, g_1)}{c_x(g_2, g_1)} \hat{E}(\sla{g_2 g_1} (y, x)).
 \]
The fibers of $\hat{E}$ are hermitian lines and have support given by the graph of $f$, that is points of the form $(f(x), x)$. If we choose a basis vector for each of these lines then we can identify the linear maps $\hat{E} (\sla{g} (f(x), x))$ with an element of $U(1)$. So we can rearrange the above equation into the form
 \[
  c_x (g_2, g_1) = c'_{f(x)} (g_2, g_1) \, \frac{\hat{E}(\sla{g_2 g_1} (f(x), x))}{\hat{E}(\sla{g_2} (g_1 \ccdot f(x), g_1 \ccdot x)) \circ \hat{E}(\sla{g_1} (f(x), x))}
  \]
which says precisely that $c$ differs from $c'$ by a coboundary $d \hat{E}$, in multiplicative notation.

(ii) $\Rightarrow$ (i). Write $c \equiv c_\X$ and $c' \equiv c_\Y$. Suppose we are given an isomorphism of $G$-sets $f : X \rightarrow X'$. To say that $[c] = [f^* c']$ as cohomology classes means that $c = (f^* c') d \gamma$ for some 1-cochain $\gamma$, ie. there must exist a map
 \[
  \gamma : \text{Arrows}(X_G) \rightarrow U(1)
 \]
such that
 \[
  c_x(g_2, g_1) = c'_{f(x)} (g_2, g_1) \frac{\gamma (\sla{g_2 g_1} x)}{\gamma(\sla{g_2} g_1 \ccdot x) \gamma (\sla{g_1} x)}.
 \]
We define vector bundles $\hat{E}$ and $\hat{F}$ over $Y \times X$ and $X \times Y$ respectively by setting
 \[
  \hat{E}_{y, x} = \delta_{y, f(x)} \mathbb{C}, \quad \hat{F}_{x,
 y} = \delta_{x, f^{\mi}(y)} \mathbb{C}.
 \]
We can use $\gamma$ to give $\hat{E}$ (respectively $\hat{F}$) the structure of a $c$- (respectively $c'$-) twisted equivariant vector bundle. We only need to do this on the graph of $f$, where we define
 \[
 \hat{E}(\sla{g} (f(x), x)) = \gamma (\sla{g} x), \quad \hat{F} (\sla{g} (x, f(x))) = \overline{ \gamma} (\sla{g} x).
 \]
It is easy to check that $\hat{E}$ and $\hat{F}$ are twisted equivariant vector bundles, which when considered as morphisms of equivariant gerbes furnish an equivalence between $\X$ and $\Y$ inside $\Gerbes(G)$. Moreover since the map $f$ preserved the scale factors on $X$ and $Y$, this equivalence is an isometric equivalence.
\end{proof}

\section{Line bundles and U(1)-bundles\label{uonebund}}
In this section we define what we mean by `$U(1)$-bundles' and `line bundles' on groupoids and their spaces of flat sections. All of this technology is essentially taken from the paper of Willerton \cite{ref:simon}; in particular we state the useful formula which computes the dimension of the space of flat sections of a line bundle as the integral of the transgression of the line bundle over the loop groupoid. We will use this theorem in Chapter \ref{eqref} to give elegant formulas for the dimensions of the hom-sets in $[\TRep(G)]$, the complexified Grothendieck category of unitary 2-representations of $G$.

\subsection*{Hermitian lines and $U(1)$-torsors} A {\em hermitian line} is a one-dimensional complex vector space with inner product, and a {\em $U(1)$-torsor} is simply a set equipped with a free and transitive action of $U(1)$. We write $\UOneTor$ for the category of $U(1)$-torsors and equivariant maps, and $\mathcal{L}$ for the category of hermitian lines and linear maps. To a $U(1)$-torsor $P$ we can associate a hermitian line $P_\mathbb{C}$ by taking the quotient of the cartesian product $P \times \mathbb{C}$ under the identifications $(e^{i \theta} p, \lambda) \sim (p, e^{i \theta} \lambda)$. We write the equivalence class of $(p, \lambda)$ as $p \otimes \lambda$, and the inner product on the line $P_\mathbb{C}$ is defined by $(p \otimes \lambda, p' \otimes \lambda') = \frac{p'}{p} \overline{\lambda}\lambda'$. Similarly, given a hermitian line $L$ we can associate a $U(1)$-torsor by taking the elements of unit norm.

\subsection*{$U(1)$-bundles and line bundles on groupoids}
We define a {\em $U(1)$-bundle with connection} over a finite groupoid $\mathcal{G}$ to be a functor $P \colon
\mathcal{G} \rightarrow \UOneTor$. Similarly a {\em hermitian line bundle with unitary connection} over
$\mathcal{G}$ is a functor $L \colon \mathcal{G} \rightarrow \mathcal{L}$, such that all the maps $L(\gamma)$ are
unitary, where $\gamma$ is an arrow in $\mathcal{G}$. We can use the conventions in the previous paragraph to
convert $U(1)$-bundles with connection into hermitian line bundles with unitary connection, and vice-versa.

\subsection*{$U(1)$-bundles and 1-cocycles}
A {\em trivialization} of a $U(1)$-bundle is a choice $\lambda_x \in P_x$ for each $x \in \mathcal{G}$. Choosing a
trivialization gives rise to a $U(1)$-valued groupoid 1-cocycle $\alpha \in Z^1(\mathcal{G}, U(1))$ (in the sense of \cite{ref:simon}) whose value on a morphism $\gamma$ in $\mathcal{G}$ is defined by the equation
 \[
  \alpha(\gamma) \lambda_{\text{target}(\gamma)} =
  P(\gamma)(\lambda_{\text{source}(\gamma)}).
 \]

\subsection*{Flat sections of line bundles}
A {\em flat section} of a line bundle $L \colon \mathcal{G} \rightarrow \mathcal{L}$ is a choice $s_x \in L_x$ for each $x
\in \mathcal{G}$, such that $s( \text{target}(\gamma)) = L(\gamma) s( \text{source}(\gamma))$ for all arrows
$\gamma \in \text{Arr} \, \mathcal{G}$. The space of flat sections of $L$ is denoted $\Gamma(L)$. If $s$ and
$s'$ are flat sections, then their fibrewise inner-product $(s, s')_x$ is a 0-form on $\mathcal{G}$, and hence can be integrated with respect to the natural measure on a groupoid (see \cite{ref:simon}, and also the more recent \cite{ref:weinstein}). So the space of sections $\Gamma(L)$ is endowed with an inner product via
 \[
  ( s, s') = \int_{x \in \mathcal{G}} (s, s')_x := \sum_{x \in \mathcal{G}} \frac{(s(x), s'(x))}{|x \rightarrow |},
 \]
where $|x \rightarrow| $ refers to the number of arrows emanating of $x$.
\subsection*{Transgression of a line bundle}
Recall the notion of the {\em loop groupoid} $\Lambda \mathcal{G} := \Fun(\mathbb{Z}, \mathcal{G})$ of a finite groupoid $\mathcal{G}$ from Chapter \ref{loopgdsec}. For instance, when $\mathcal{G}$ is the action groupoid $X_G$ corresponding to an action of $G$ on a finite set $X$, the objects of the loop groupoid are `loops' in $\mathcal{G}$ which we write as
\[
( \Fix{x}{g})
\]
and the morphisms are given by conjugation, which we write as
\[
(\Fixx{h \ccdot x}{hgh^\mi}{5.7}{3}) \sla{h} (\Fix{x}{g}).
\]
If $L$ is a line bundle over $\mathcal{G}$, then we define the {\em transgression} of $L$ to be the function
 \begin{align*}
  \tau(L) & \colon \Ob \Lambda \mathcal{G} \rightarrow \mathbb{C}\\
 \intertext{which assigns to each loop $\gamma$ in $\mathcal{G}$ the holonomy of $L$ around $\gamma$:}
  ( \Fix{x}{\gamma}) & \mapsto \Tr(L(\gamma) \colon L_x \rightarrow L_x).
 \end{align*}
Willerton has proved an elegant formula which computes the dimension of the space of flat sections of $L$ as the integral of the transgression of $L$ over the loop groupoid.
\begin{thm}[{Willerton \cite[Thm 6]{ref:simon}}]\label{FirstWill} If $L$ is a line bundle over a finite groupoid $\mathcal{G}$, then the dimension of the space of flat sections of $L$ computes as the integral of the transgression of $L$ over the loop groupoid:
 \[
  \dim \Gamma(L) = \int_{\Lambda \mathcal{G}} \tau(L).
 \]
\end{thm}
When we establish the correspondence between equivariant gerbes and unitary 2-representations in Chapter \ref{GerbesCharChap}, we will see that the {\em hom-sets} in the complexified Grothendieck category $[\TRep(G)]_\mathbb{C}$ can be thought of as the space of flat sections of a certain line bundle, so that the above formula then gives elegant compact formulas  (Corollary \ref{tre}) for the dimensions of these hom-sets.

\section{Transgression and twisted characters\label{tgess}} In this section we define the {\em transgressed line bundle} of an equivariant gerbe as a certain line bundle over the loop groupoid. Then we state the theorem of Willerton \cite{ref:simon} which identifies the space of isomorphism classes of equivariant vector bundles over the gerbe as the space of sections of this line bundle.

\subsection*{Line bundles from transgression of equivariant gerbes \label{transsec}}
\begin{figure}[t]
\centering
\ig{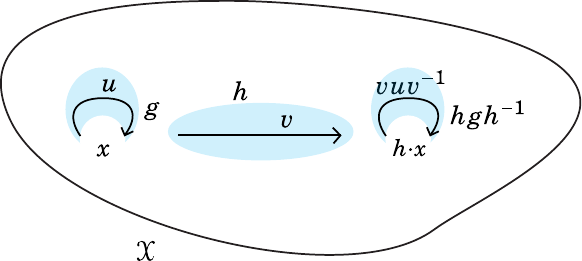}
\caption{\label{tbundle} The parallel transport operation in the transgressed $U(1)$-bundle of a gerbe $\X$ is obtained by conjugating `loops' in the gerbe by choosing arbitrary arrows $v \in \X_\mlag{x}{}{g}$; this operation does not depend on the choice of $v$.}
\end{figure}

Suppose $\X$ is an equivariant gerbe with underlying $G$-set $X$, thought of via its action groupoid $X_G$. We define the {\em transgressed $U(1)$-bundle} of $\X$ as the functor
 \[
 \tau(\X) \colon \Lambda X_G  \longrightarrow  \UOneTor
 \]
which sends an object of the loop groupoid to its corresponding $U(1)$-torsor in the gerbe,
\[
 ( \Fix{x}{g}) \quad  \mapsto  \quad \X\sFix{x}{g},
 \]
and an arrow in the loop groupoid to the conjugation operation on these torsors (see Figure \ref{tbundle}),
 \[
 (\Fixx{h\cdot x}{hgh^\mi}{5.7}{3}) \sla{h} (\Fix{x}{g})  \quad  \mapsto  \quad u \mapsto vuv^\mi.
\]
Here $v$ is an arbitrary choice of arrow in the torsor $\X_\mlag{x}{}{g}$; the formula is clearly independent of
this choice. The associated hermitian line bundle $\tau(\X)_\mathbb{C}$ is known as the {\em transgressed line bundle}. One should think of this construction as the line bundle arising from the {\em holonomy} of the gerbe as in Brylinski \cite{ref:brylinski}.

\subsection*{Transgression in terms of cocycles} In \cite{ref:simon}, Willerton defined transgression at the level of cocycles by using the `Parmesan map'; for 2-cocycles the resulting map
 \[
  \tau \colon Z^2 (X_G, U(1))  \rightarrow Z^1(\Lambda X_G, U(1))
 \]
works out as
 \[
  \tau(c)\left( \sla{h} (\Fix{x}{g})\right) = \frac{ c_x(hgh^\mi, h)}{c_x(h, g)}.
 \]
After passing to cocycles, the geometric definition of transgression we gave above reproduces this formula. To see this, choose a `normalized' section $s \colon X_G \rightarrow \X$ of the gerbe; that is one having the property that $s(\sla{e} x) = 1$ for all $x \in \X$ and $s(\sla{g^\mi} g \ccdot x) = s(\sla{g} x)^\mi$ for all $x \in \X, g \in G$. This section gives rise to a 2-cocycle $c \in Z^2(X_G, U(1))$ as we explained in Section \ref{BigR}. Now, the transgressed $U(1)$-bundle $\tau(\X)$ will act on the distinguished morphisms given by the section as
 \begin{align}
  \tau(\X)\left( \sla{h} (\Fix{x}{g}) \right) s(\Fix{x}{g}) &= s(\sla{h} x) \circ s(\Fix{x}{g}) \circ s( \sla{h^\mi} h \ccdot x) \nonumber \\
   &= c_x(h,g) c_{h \ccdot x}(hg, h^\mi) s(\Fixx{h\ccdot x}{hgh^\mi}{5.7}{3}) \nonumber \\
   &= \frac{c_x (h, g)}{ c_x(hgh^\mi, h)} s(\Fixx{h\ccdot x}{hgh^\mi}{5.7}{3}) \label{prect}
  \end{align}
where the last equality uses the cocycle equation $c_x (hgh^\mi, h) c_{h \ccdot x}(hg, h^\mi) = c_x(hg, e) c_x(h^\mi, h)$ and the fact that the section was normalized so the right hand side here is equal to unity. The formula \eqref{prect} is precisely the formula of Willerton for the transgression of $U(1)$-valued groupoid 2-cocycles \cite[Thm 3]{ref:simon}. In other words, if we think of equivariant gerbes as $U(1)$-valued 2-cocycles on an action groupoid $X_G$, and $U(1)$-bundles over groupoids as $U(1)$-valued 1-cocycles as we explained in Section \ref{uonebund}, then our direct geometric definition for transgression corresponds precisely to the cocycle formulation of transgression as defined by Willerton.

\subsection*{Twisted characters of equivariant vector bundles\label{chsec}}
Suppose $E \colon \X \rightarrow \Hilb$ is a unitary equivariant vector bundle over an equivariant gerbe $\X$. The
{\em twisted character} (or just {\em character} for short) of $E$ is a flat section of the transgressed line bundle,
 \[
  \chi_E \in \Gamma_{\Lambda X_G} (\tau(\X)_\mathbb{C}).
 \]
It is defined by setting
 \[
  \chi_E ( \Fix{x}{g} ) = u \otimes \Tr E(\lag{g}{u})^*
 \]
where $u$ is any morphism in $\X_\mlag{x}{}{g}$ (the choice of $u$ doesn't matter since the formula is
invariant under $u \mapsto e^{i \theta} u$) and the `$*$' refers to complex conjugation. We then have the following important result, which takes the formula from Theorem \ref{FirstWill} one step `higher'.
\begin{thm}[{Willerton \cite[Thm 11]{ref:simon}}] \label{Simthm} The twisted character map is a unitary isomorphism from the complexified Grothendieck group of isomorphism classes of unitary equivariant vector bundles over an equivariant gerbe $\X$ to the space of flat sections of the transgressed line bundle:
 \[
  \chi \colon [\Hilb_G (\X)]_\mathbb{C} \ra{\cong \, \, \, \,} \Gamma_{\Lambda X_G} (\tau(\X)_\mathbb{C}).
  \]
\end{thm}
This formula should be thought of as the second in a long line of formulas (Theorem \ref{FirstWill} is the first) which compute the space of flat sections of a `higher line bundle' using the {\em holonomy} of the line bundle. In Section \ref{fgerbes} we will use this formula to establish that the geometric character functor is unitarily fully faithful, which will in turn imply that the {\em 2-character} functor acting on unitary 2-representations is unitarily faithful (Corollary \ref{corc}), one of our main results in this thesis.

\section{The geometric character of an equivariant gerbe\label{gch}}
This section is the geometric analogue of Chapter \ref{sec2char}: we define how to take the {\em geometric character} of a $G$-equivariant gerbe equipped with a metric in order to obtain a unitary equivariant vector bundle over $G$.

\subsection{Pushing forward line bundles over loop groupoids}
\begin{figure}[t]
\centering
\ig{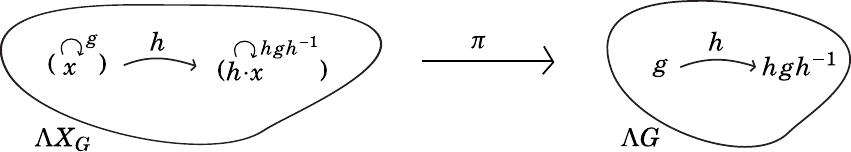}
\caption{\label{natfib}The fibration of loop groupoids associated to a $G$-set $X$.}
\end{figure}
For every $G$-set $X$ there is a natural fibration of loop groupoids
$\pi \colon \Lambda X_G \rightarrow \Lambda \BG$, where the fiber of $g \in \Lambda G$ is the set
 \[
  \left\{\Fix{x}{g}\right\}_{x \in \FFix (g)}
 \]
of fixed points of $g$ on $X$ (see Figure \ref{natfib}). If $X$ has a metric --- an equivariant assignment of a positive real number $k_x$ to each $x \in X$ --- then we can push-forward a unitary equivariant line bundle $L$ over $\Lambda X_G$ to a unitary equivariant vector bundle $\pi_* (L)$ over $G$ by taking the space of sections over the fixed points, as follows. The fiber of $\pi_*(L)$ at $g \in G$ is the Hilbert space
 \[
   \Sections (L|_{\FFix(g)}) := \hat{\displaystyle \bigoplus_{x \in \FFix(g)}} k_x \;  L \sFix{x}{g}
 \]
of sections of $L$ over the fixed points of $g$. That is, the fiber of the push-forward vector bundle $\pi_* (L)$ over $G$ at $g\in G$ is the space of sections of $L$ over the fixed points of $g$ on $X$, ie. a vector $\psi \in \pi_* (L)_g$ is an assignment
  \[
   ( \Fix{x}{g}) \mapsto \psi (\Fix{x}{g}) \in L \sFix{x}{g}
  \]
where $x$ ranges over the fixed points of $g$. The inner product on these sections is given by
 \[
  \langle \psi, \psi' \rangle = \sum_{x \in \FFix(g)} k_x \bigl( \psi(\Fix{x}{g}), \psi'(\Fix{x}{g}) \bigr).
 \]
As $g$ ranges over $G$, these vector spaces become a unitary equivariant vector bundle over $G$ via the natural $G$-action
 \[
 (h \cdot \psi) ( \Fixx{h\cdot x}{hgh^\mi}{5.7}{3}) = L(h) \psi(\Fix{x}{g}).
\]

\subsection{Definition of the geometric character\label{deftransg}}
In particular, if $\X$ is an equivariant gerbe equipped with a metric we write $\ch (\X) := \pi_* (\tau(\X)_\mathbb{C})$
for the push-forward of the transgressed line bundle of $\X$, and we call $\ch (\X)$ the {\em geometric character}
of $\X$. So, the geometric character of an equivariant gerbe is a certain equivariant vector bundle over the group $G$ constructed from the fixed point data of the gerbe. Its fiber at a group element $g \in G$ computes as
 \be \label{needlabel}
  \ch(\X)_g = \Sections (\tau(\X)_\mathbb{C}|_{\FFix (g)}).
 \ee
That is, a vector in the fiber of the geometric character bundle at $g \in G$ is something which assigns to each fixed point $x \in X$ of $g$ an element in the complex line $(\X \sFix{x}{g})_\mathbb{C}$. If $\X$ was just the `trivial' gerbe on a $G$-set $X$, then this fiber would just be the space of functions on the fixed points of $g$; in the general case it is simply the space of sections of the {\em transgressed line bundle} over the fixed points of $g$.

\subsection{Functoriality for the geometric character\label{fgerbes}} This section is the geometric analogue of Chapter \ref{funcchar}. We show how to make the geometric character construction {\em functorial} with respect to morphisms of equivariant gerbes, so that it descends to a functor from the Grothendieck category of $\Gerbes(G)$ to the category of equivariant vector bundles over $G$. Then we apply Theorem \ref{Simthm} to show that after one tensors the hom-sets in $[\Gerbes(G)]$ with $\mathbb{C}$, the resulting complexified geometric character functor is unitarily fully faithful. After we have established the equivalence between unitary 2-representations and equivariant gerbes in the next chapter, this result will imply that the complexified 2-character is also unitarily fully faithful.

\subsection*{The geometric character of a morphism of equivariant gerbes}
\begin{figure}[t]
\centering
\ig{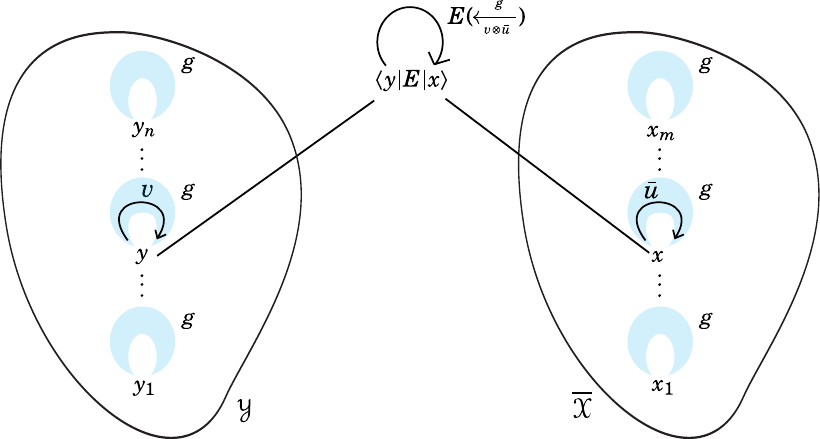}
\caption{\label{gcharmo}The fibrewise geometric character of a morphism $E \colon \X \rightarrow \Y$ is a linear map from the space of sections of the line bundle over the fixed points of $g$ on $X$ to space of sections of the line bundle over the fixed points of $g$ on $Y$. The matrix elements of this map are computed by taking the trace of the vector bundle $E$.}
\end{figure}

Suppose $E \colon \X \rightarrow \Y$ is a morphism of equivariant gerbes. Then we can define for each $g \in G$
a linear map
 \begin{align*}
  \ch(E)_g \colon \Sections(\tau(\X)_\mathbb{C} |_{\FFix_X(g)}) & \rightarrow \Sections(\tau(\Y)_\mathbb{C}
  |_{\FFix_{Y}(g)}) \\
   \psi & \mapsto \ch(E)_g (\psi)
  \end{align*}
by integrating the trace of $E$ over the simultaneous fixed points of $g$; see Figure \ref{gcharmo}.
It is easiest to give this map in terms of its matrix elements between orthonormal bases of sections $\{ \psi_x \in \ch(\X)_g\}_{x \in \FFix_X(g)}$ and $\{\phi_y \in \ch(\Y)_g\}_{y \in \FFix_Y(g)}$ consisting of sections which are localized over fixed points $x \in \FFix_X(g)$ and $y \in \FFix_{Y}(g)$ respectively:
\[
 \psi_x (\FixP{x'}{g}) = \delta_{xx'} \; u \otimes \frac{1}{\sqrt{k_x}} \; \in (\X \sFix{x}{g})_\mathbb{C}, \quad  \phi_y (\FixP{y'}{g}) = \delta_{yy'} \; v \otimes \frac{1}{\sqrt{k_y}} \in (\Y \! \sFix{y}{g})_\mathbb{C}.
\]
In terms of these basis vectors the matrix elements of $\ch(E)_g$ are defined as
 \be \label{matelements}
  \langle \phi_y , \, \ch(E)_g \psi_x \rangle = \Tr E (\lag{g}{v \otimes \bar{u}})^*.
  \ee
Note that once again this definition does not depend on the choices made for $v \in \Y \sFix{y}{g}$ and $u \in \X \sFix{x}{g}$.

\subsection*{Proof of functoriality}
Recall that we are writing $\Gerbes(G)$ for the 2-category of finite $G$-equivariant gerbes equipped with metrics. As in Chapter \ref{funcchar}, we will write $[\Gerbes(G)]$ for the {\em Grothendieck category} of $\Gerbes(G)$ --- the category obtained by identifying isomorphic 1-morphisms --- and $[\Gerbes]_\mathbb{C}$ for the category obtained from $[\Gerbes(G)]$ by tensoring the hom-sets with $\mathbb{C}$.
 \begin{thm}\label{geomchthm} The assignment
 \begin{align*}
  \ch \colon \Gerbes(G) & \quad \longrightarrow \quad \Hilb_G(G) \\
  \X & \quad \mapsto \quad \ch(\X) \\
  \X \sra{E} \Y & \quad \mapsto \quad \ch(\X) \ra{\ch(E)} \ch(\Y)
 \end{align*}
is functorial with respect to composition in $\Gerbes(G)$, and only depends on the isomorphism class of $E$, and
thus descends to a {\em geometric character functor}
  \[
   \ch \colon [\Gerbes(G)] \rightarrow \Hilb_G(G).
   \]
Moreover, after tensoring the hom-sets in $[\Gerbes(G)]$ with $\mathbb{C}$ the associated functor
  \[
   \ch_\mathbb{C} \colon [\Gerbes(G)]_\mathbb{C} \rightarrow \Hilb_G(G)
  \]
 is unitarily fully faithful.
\end{thm}
\begin{proof}
Functoriality follows from the fact that trace is multiplicative on tensor products and `localizes' on the fixed
points. In other words, if $\Z \sla{F} \Y \sla{E} \X$ are morphisms of equivariant gerbes, and if $u \in \X
\sFix{x}{g}$ and $w \in \Z \sFix{z}{g}$, then by the definition \eqref{2kerncomp} of composition in $\Gerbes(G)$,
\begin{align*}
 \Tr ( F \circ E) (\lag{g}{w \otimes \bar{u}}) &= \sum_{y \in \Y} \Tr F (\lag{g}{w \otimes \bar{v}}) \Tr E(\lag{g}{v
 \otimes \bar{u}}) \\
 &= \sum_{y \in \FFix_{Y}(g)} \Tr F(\lag{g}{w \otimes \bar{v}}) \Tr E(\lag{g}{v
 \otimes \bar{u}}).
\end{align*}
Also the trace of an equivariant vector bundle only depends on its isomorphism class; this gives the first
part of the proposition.

To prove the second part of the proposition, we show that the action of the geometric character on morphisms is
just a rearrangement of the twisted character map from Section \ref{chsec}, which is known to be unitary by Theorem \ref{Simthm}. Namely, we claim we have the following commutative diagram:
 \[
  \xymatrix{ \Hom_{[\Gerbes(G)]} (\X, \Y) := [\Hilb (\Y \otimes \overline{\X})] \ar[r]^-{\chi} \ar[dr]_{\ch} &
  \Gamma_{\Lambda (Y \times X)_G} (\tau(\Y \otimes \overline{\X})_\mathbb{C}) \ar[d]^\wedge \\
  & \Hom_{\Hilb_G(G)} (\ch(\X), \ch(\Y))}
 \]
The rearrangement map $\hat{}$ (the downwards arrow) works as follows. A section $\xi$ of the transgressed line bundle $\tau(\Y \otimes
\overline{\X})_\mathbb{C}$ is something which assigns to every simultaneous fixed point an element of the appropriate hermitian line:
 \[
  \Fixx{(y, x)}{g}{4}{3} \quad \stackrel{\xi}{\mapsto} \quad v \otimes \bar{u} \otimes \lambda \; \in \; \Y \sFix{y}{g} \otimes \overline{\X} \sFix{x}{g} \otimes \mathbb{C}.
 \]
Using the metrics on the gerbes, this can be regarded as a map of hermitian lines
\begin{align*}
 \hat{\xi} \colon (\X \sFix{x}{g})_\mathbb{C} \rightarrow (\X' \sFix{y}{g})_\mathbb{C}
 \intertext{via the formula}
 \frac{1}{\sqrt{k_x}} u \otimes 1 \mapsto \frac{1}{\sqrt{k_y}} v \otimes \lambda.
 \end{align*}
Using this correspondence, the section $\xi$ gives rise for each $g \in G$ to a linear map
 \[
  \hat{\xi}_g \colon \ch(\X)_g \rightarrow \ch(\Y)_g.
 \]
The fact that $\xi$ was a {\em flat} section translates into the statement that the collection of maps $\hat{\xi}_g$ is {\em equivariant} with respect to the action of $G$. Moreover one can check that the map $\xi \mapsto \hat{\xi}$ is unitary with respect to the natural inner products involved (note that the twisted character map $\chi$ does not use the metrics on the gerbes, but $\ch$ and $\;\hat{}\;$ do), and also that the above diagram indeed commutes. Now we apply Theorem \ref{Simthm}, which says that after tensoring the left hand side with $\mathbb{C}$ the
character map $\chi$ is a unitary isomorphism. This gives the second statement of the proposition.
\end{proof}
Let us summarize what we have done in this chapter. We have introduced the geometric language of finite equivariant gerbes equipped with metrics and their associated geometric characters. This framework allowed us to use the result of Willerton (Theorem \ref{Simthm}) to establish that the complexified geometric character is unitarily fully faithful functor. In the next chapter we will set up a correspondence between unitary 2-representations and equivariant gerbes, and show that under this correspondence the 2-character of a 2-representation corresponds to the geometric character of its associated equivariant gerbe. This will allow us to conclude that the complexified 2-character is a unitarily fully faithful functor, which is one of the main results of this thesis. 
\chapter{Geometric characters and 2-characters\label{GerbesCharChap}} 
In this chapter we prove our main theorem about 2-characters, namely that the complexified 2-character functor
 \[
  \chi \colon [\TRep(G)]_\mathbb{C} \rightarrow \Hilb_G(G)
 \]
is unitarily fully faithful. We do this by developing a correspondence between unitary 2-representations and finite equivariant gerbes equipped with metrics. We show that under this correspondence, the 2-character of a 2-representation corresponds naturally to the geometric character of its associated equivariant gerbe. The result then follows from Theorem \ref{geomchthm} which proved that the complexified geometric character functor is unitarily fully faithful.

We also show how this result has an interesting corollary: unlike the ordinary character of a representation, the {\em 2-character} does not always distinguish different 2-representations. This is likely an artifact of the generally accepted fact that semisimple categories are not quite the right notion of a `2-vector space'.

In Section \ref{extract} we show how to extract an equivariant gerbe from a marked unitary 2-representation, and we apply this procedure to some of the examples of unitary 2-representations we provided in Chapter \ref{exsec}. In Sections \ref{AMor} and \ref{BMor} we show how a morphism of unitary 2-representations gives rise to a morphism of equivariant gerbes, and the same for 2-morphisms. In Section \ref{eqref} we show that this correspondence is an equivalence of 2-categories. This allows us to give some elegant closed formulas for the dimensions of the hom-spaces in the complexified Grothendieck category $[\TRep(G)]_\mathbb{C}$ using the groupoid integration technology of \cite{ref:simon} which we reviewed in Chapter \ref{uonebund}. Finally in Section \ref{2char2rep} we prove our main theorem relating the 2-character to the geometric character, with the corollary that the complexified 2-character is unitarily fully faithful. We close the chapter in Section \ref{notdistinguish} by showing that the 2-character does not necessarily distinguish nonequivalent 2-representations.

\section{Equivariant gerbes from unitary 2-representations\label{extract}} In this section we show how to extract a $G$-equivariant gerbe equipped with a metric from a unitary 2-representation of $G$ on a marked 2-Hilbert space. Recall from Chapter \ref{2HilbChap} that a 2-Hilbert space is marked if each isomorphism class of simple objects has a distinguished representative $e_i$.

\begin{figure}[t]
\centering
\ig{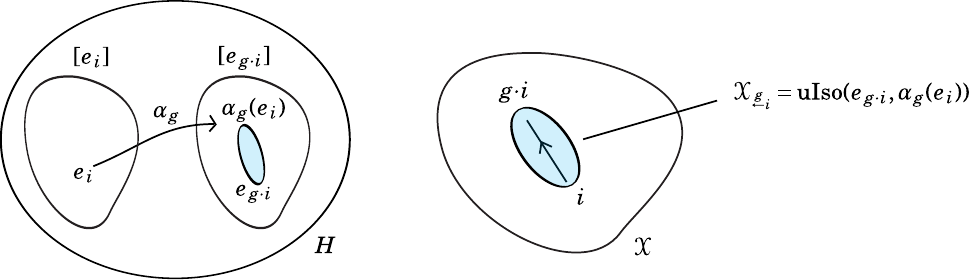}
\caption{\label{extractfig} Extracting an equivariant gerbe from a unitary 2-representation. The $U(1)$-torsor $\X_\mlag{i}{}{g}$ of arrows in the gerbe is the set of unitary arrows $u \colon e_{g \ccdot i} \rightarrow \alpha_g (e_i)$ in the 2-Hilbert space.}
\end{figure}

A unitary 2-representation $\alpha$ of $G$ on a marked 2-Hilbert space $H$ gives rise to an equivariant gerbe $\X$ by a variant of the {\em Grothendieck construction} (we learnt this approach from the paper of Cegarra et. al. \cite{ref:cegarra_et_al_graded_extensions}, but the original construction is due to Grothendieck \cite{ref:grothendieck}). Our variant of this construction works as follows. The base set $X$ of the gerbe is the set of isomorphism classes of simple objects in $H$; these are often referred to by their indices, so that $i \equiv [e_i]$ and so on. The set $X$ inherits a $G$-action via
 \[
  g \ccdot i := [\alpha_g (e_i)]
 \]
where $e_i$ is the distinguished simple object in the isomorphism class $[e_i]$. The equivariant gerbe $\X$ has the same objects as $X_G$, with the $G$-graded hom-sets given by the unitary isomorphisms in $H$ between the marked simple objects and the objects which have been acted on by $G$:
 \[
  \X_\mlag{i}{}{g} := \text{uIso}(e_{g \ccdot i}, \alpha_g (e_i)).
  \]
We refer the reader to Figure \ref{extractfig}.
Note that these unitary isomorphisms in $H$ do {\em not} lie on the unit circle in the hermitian line $\Hom(e_{g \ccdot i}, \alpha_g(e_i))$, but rather on the circle with radius $\sqrt{k_i}$.

The composite of $u \in \X_\mlag{i}{}{g}$ and $v \in \X_\mlag{g \ccdot i}{}{h}$ is defined as
 \be \label{gerbecomp}
  v \diamond u = a(h, g)_{e_i} \, \alpha_{h}(u) \, v,
 \ee
where we have used the `$\diamond$' symbol to stress that this is not {\em ordinary} composition of arrows in the 2-Hilbert space $H$. The identity morphisms $1_i \in \X_\mlag{i}{}{e}$ are given by the unit isomorphisms $a(e)_{e_i}$. It is well known (see eg. \cite{ref:cegarra_et_al_graded_extensions}) that this sort of construction produces a unital and associative composition law on the groupoid $\X$; we record the proof for completeness.
 \begin{lem} This construction equips $\X$ with a unital and associative composition law.
 \end{lem}
 \begin{proof} Suppose we have $u \in \X_\mlag{i}{}{g}$, $v \in \X_\mlag{g \ccdot i}{}{h}$ and $w \in \X_\mlag{hg \ccdot i}{}{k}$. Then
  \begin{align*}
   (w \diamond v) \diamond u &= a(kh, g)_{e_i} \alpha_{kh}(u) a(k,h)_{e_{g \ccdot i}} \alpha_k(v) w \\
    &= a(kh, g)_{e_i} a(k,h)_{\alpha_g(e_i)} \alpha_k \alpha_h (u) \alpha_k(v) w \\
    &= a(k, hg)_{e_i} \alpha_k(a(h,g)_{e_i}) \alpha_k \alpha_h (u) \alpha_k(v) w \\
    &= w \diamond (v \diamond u)
   \end{align*}
where the second equation uses naturality of $a(k,h)$ and the third equation uses the coherence equation \eqref{br2} from Chapter \ref{Untrep}. Similarly
 \[
  u \diamond 1_i = a(g,e)_{e_i} \alpha_g(a(e)_{e_i}) u = u \quad \text{and} \quad 1_{g \ccdot i} \diamond u = a(e,g)_{e_{g \ccdot i}} \alpha_e(u) a(e)_{e_{g \ccdot i}} = u.
 \]
\end{proof}
We we will write the inverse of $u \in \X_\mlag{i}{}{g}$ with respect to this composition law as $u^{\underline{\mi}} \in \X_\mlag{g \ccdot i}{}{\,g^{\text{-}1}}$. It is easy to check that $u^{\underline{\mi}}$ computes in terms of the {\em ordinary} inverse $u^\mi\equiv u^*$ as
 \be \label{inverselaw}
  u^{\underline{\mi}} = \alpha_{g^\mi}(u^\mi) a(g^\mi, g)^\mi_{e_i} a(e)_{e_i}.
 \ee
Finally, the metric on the equivariant gerbe $\X$ is defined as the scale factors on the simple objects $e_i$, that is $k_i = (\id_{e_i}, \id_{e_i})$.

\subsection*{Examples\label{gexam}} We now extract the associated equivariant gerbes from some of the examples of unitary 2-representations we gave in Section \ref{exsec}.
\subsubsection{Automorphisms of groups} The equivariant gerbe $\X$ arising from the 2-representation of $G \subseteq{\Aut(K)}$ on $\Rep(K)$ works as follows. Firstly one chooses distinguished irreducible representations $V_i$ of $K$. The underlying $G$-set of $\X$ is the set of isomorphism classes of irreducible representations (which identifies noncanonically with the conjugacy classes of $G$), and the graded hom-sets are
 \[
  \X_\mlag{i}{}{g} = \uIso_K(V_{g \ccdot i}, V_i^g).
 \]
In particular the $U(1)$-torsors above the fixed points are
 \[
  \X \sFix{i}{g} = \uIso_K(V_i, V_i^g) \subset U(V)
 \]
and there is no {\em apriori} preferred section of these torsors, unless $g$ is an inner automorphism. Let us summarize this discussion:
 \begin{quote}
 If $G$ acts as automorphisms of a group $K$, then each irreducible representation $\rho$ of $K$ will carry a projective representation of the subgroup $G_0 \subset G$ which fixes $\rho$.
 \end{quote}

\subsubsection{The metaplectic representation} In this way the equivariant gerbe arising from the 2-representation of $SL_2(\mathbb{R}) \subset \Aut(\text{Heis})$ on the category of representations of the Heisenberg group has a single object whose automorphism group is precisely the $\text{\em metaplectic}^c$ {\em group} $Mp^c(2)$ (the metaplectic group is the nontrivial double cover of $SL_2(\mathbb{R})$, see \cite[pg 2]{ref:robinson_rawnsley}). In other words, although the 2-representation is strict its corresponding equivariant gerbe and hence the associated projective representation are actually nontrivial.

\section{Morphisms of gerbes from morphisms of 2-representations\label{AMor}} Firstly, we define a morphism of marked 2-representations simply as an ordinary morphism of the underlying 2-representations which pays no attention to the distinguished simple objects, and similarly for the 2-morphisms, and we write the 2-category of marked unitary 2-representations of $G$ as $\TRep_m(G)$. Note that the obvious 2-functor $\TRep_\text{m} (G) \rightarrow \TRep(G)$ which forgets the marked simple objects is an equivalence of 2-categories.

A morphism $\sigma \colon \alpha \rightarrow \beta$ of marked 2-representations gives rise to a morphism of equivariant gerbes $\langle \sigma \rangle \colon \X_\alpha \rightarrow \X_\beta$ in the following way. The vector bundle $\langle \sigma \rangle$ over $\X_\beta \otimes \overline{\X}_\alpha$ is defined to
have fibers
 \[
  \langle \mu | \sigma | i \rangle  := \Hom(e_\mu, \sigma(e_i)).
 \]
For $v \in (\X_\beta)_\mlag{\mu}{}{g}$ and $\bar{u} \in (\overline{\X}_\alpha)_\mlag{i}{}{g}$, the unitary maps
 \[
 \langle \mu | \sigma | i \rangle (\lag{g}{v \otimes \bar{u}}) \colon  \langle \mu | \sigma | i\rangle
  \rightarrow \langle g \ccdot \mu| \sigma | g \ccdot i \rangle
  \]
send
 \be \label{eqdefnmor}
 \ba \xymatrix{\sigma(e_i) \\ e_\mu \ar[u]^\lambda} \ea \quad
 \mapsto \quad\ba \xymatrix{ & \sigma \alpha_g (e_i)
 \ar[dr]^{\sigma(u)^*} \\
 \beta_g \sigma(e_i)\ar[ur]^{\sigma_{e_i}(g)} & & \sigma(e_{g \ccdot i})
 \\
 \beta_g (e_\mu) \ar[u]^{\beta_g (\lambda)} & & e_{g \ccdot \mu}
 \ar[ll]^{v}} \ea
 \ee
We now show that this formula indeed defines an equivariant vector bundle $\langle \sigma \rangle \colon \X_\beta \otimes
\overline{\X}_\alpha \rightarrow \Hilb$; note that the conjugate of $\X_\alpha$ must be used because of the $\sigma(u)^*$ term occurring
above.
\begin{prop} This construction produces a unitary equivariant vector bundle $\langle \sigma \rangle \colon \X_\beta \otimes
\overline{\X}_\alpha \rightarrow \Hilb$.
\end{prop}
\begin{proof} Consider the diagram:
 \[
  \xymatrix @C=1.6cm { e_{g_2 g_1 \ccdot \mu} \ar[d]_{v_2} & & & \sigma(e_{g_2 g_1\ccdot i})  \\ \beta_{g_2} (e_{g_1 \ccdot \mu}) \ar[d]_{\beta_{g_2}(v_1)} & & \beta_{g_2} \sigma (e_{g_1 \ccdot i}) \ar[r]^{\sigma(g_2)_{e_{g_1 \ccdot i}}} & \sigma \alpha_{g_2} (e_{g_1 \ccdot i}) \ar[u]^{\sigma(u_2^*)} \\ \beta_{g_2} \beta_{g_1} (e_\mu) \ar[d]_{b(g_2, g_1)_{e_\mu}} \ar[r]^-{\beta_{g_2} \beta_{g_1}(\lambda)} & \beta_{g_2} \beta_{g_1} \sigma(e_i) \ar[d]^{b(g_1, g_1)_{\sigma(e_i)}} \ar[r]^-{\beta_{g_2}(\sigma(g_1)_{e_i})} & \beta_{g_2} \sigma \alpha_{g_1}(e_i) \ar[r]^-{\sigma(g_2)_{\alpha_{g_1}(e_i)}} \ar[u]^{\beta_{g_2} \sigma(u_1^*)} \ar[r] & \sigma \alpha_{g_2} \alpha_{g_1} (e_i) \ar[u]^{\sigma \alpha_{g_2}(u_1^*)} \\
  \beta_{g_2 g_1}(e_\mu) \ar[r]_-{\beta_{g_2 g_1}(\lambda)} & \beta_{g_2 g_1} \sigma(e_i) \ar[r]_-{\sigma(g_2 g_1)_{e_i}} & \sigma \alpha_{g_2 g_1} (e_i) \ar[ur]_{\sigma(a(g_2, g_1)^*_{e_i})}}
  \]
Taking the tour around the outside track is
\[
 \langle \mu | \sigma | i \rangle (\biglag{g_2 g_1}{(v_2 \diamond v_1) \otimes (\overline{u_2 \diamond u_1})})(\lambda),
 \]
while taking the tour around the inside track is
\[
\langle g \ccdot \mu | \sigma | g \ccdot i \rangle (\lag{g_2}{v_2 \otimes \bar{u_2}}) \circ \langle \mu | \sigma | i \rangle ( \lag{g_1}{v_1 \otimes \bar{u_1}})(\lambda).
 \]
The square at seven o' clock commutes by the naturality of $b(g_2, g_1)$, the irregular pentagon at five o' clock commutes because it is essentially the coherence equation \eqref{sig1} from Chapter \ref{Morphsec} for the natural isomorphisms $\sigma(g)$, and the square at four o' clock commutes by the naturality of $\sigma(g_2)$. Finally, the bundle is unitary precisely because all the morphisms involved in its construction are unitary (except the morphism $\lambda$ of course).
\end{proof}

\section{2-Morphisms\label{BMor}}
Similarly a 2-morphism $\theta \colon \sigma \rightarrow \rho$ between morphisms of marked 2-representations gives rise to a morphism of equivariant vector bundles $\langle \theta \rangle \colon \langle \sigma \rangle \rightarrow \langle \rho \rangle$ whose components are just given by postcomposition with the components of the natural transformation $\theta$:
 \begin{align*}
  \langle \theta \rangle_{\mu, i} \colon \langle \mu | \sigma | i \rangle & \rightarrow \langle \mu | \rho | i \rangle \\
   f & \mapsto \theta_{e_i} \circ f.
 \end{align*}

\section{Equivalence of 2-categories \label{eqref}} In this section we show that these constructions furnish an equivalence of 2-categories between  the 2-category of unitary 2-representations and the 2-category of equivariant gerbes equipped with metrics. This can be viewed as an equivariant version of the equivalence we established in Chapter \ref{2HilbChap} between the 2-category of 2-Hilbert spaces and the 2-category of finite spaces equipped with metrics.

We hasten to add that this kind of result --- obtaining an explicit understanding of unitary 2-representations and the morphisms and 2-morphisms between them --- is not really new, since related elements of it are present in the work of Elguera \cite[Cor 6.21]{ref:elgueta} and
some similar ideas also appear in the work of Barrett and Mackaay \cite[pg 17]{ref:barrett_mackaay}.  These references however do not use the language of equivariant gerbes; by using this language we hope our formulation expresses the {\em geometry} of the situation in a cleaner way:
  \begin{itemize}
   \item Firstly it shows how this result can be regarded as the `categorification' of the geometric correspondence between ordinary unitary representations and equivariant line bundles. 
   \item Secondly, the geometric language of equivariant gerbes is quite refined and enables us to understand the 2-category of 2-representations in a more succinct way. For instance, it allows us to identify the hom-sets in the complexified Grothendieck category $[\TRep(G)]_\mathbb{C}$ as spaces of flat sections of line bundles; moreover we can elegantly calculate the {\em dimension} of these hom-sets using the integration formula of Willerton (Theorem \ref{FirstWill}). 
   \item Finally, we essentially work directly with the 2-Hilbert spaces themselves and not some co-ordinatized skeleton of them, a strategy which is likely to be important in more advanced geometric situations.
   \end{itemize}
Let us state the result.
\begin{thm}\label{toodytheorem} The map
 \begin{align*}
  \TRep_\text{m}(G) & \quad \longrightarrow \quad \Gerbes(G) \\
  \alpha & \quad \mapsto \quad \X_\alpha \\
  \alpha \sra{\sigma} \beta & \quad \mapsto \quad \X_\alpha \sra{\langle \sigma \rangle} \X_\beta \\
  \sigma \sra{\theta} \rho & \quad \mapsto \quad \langle\sigma \rangle \sra{\langle \theta \rangle} \langle \rho
  \rangle
 \end{align*}
is functorial, and an equivalence of 2-categories. Moreover it is a unitary linear map at the level of 2-morphisms.
\end{thm}
\begin{proof} This map is just an equivariant version of the 2-functor
 \[
 \langle \ccdot \rangle \colon \THilb_\text{m} \rightarrow \FinSpaces
 \]
which we showed was an equivalence, and a unitary linear map at the level of 2-morphisms, in Proposition \ref{THilbprop}. The main thing to check is that if $\sigma \colon \alpha \rightarrow \beta$ and $\rho \colon \beta \rightarrow \gamma$ are morphisms of marked unitary 2-representations, then the compositor
 \[
    \langle \rho \rangle \circ \langle \sigma \rangle \rightarrow \langle \rho \circ \sigma \rangle
 \]
we defined in \eqref{mycompositor} is an equivariant map. That is, we need to check that it is compatible with the formula \eqref{eqdefnmor} for extracting equivariant vector bundles from morphisms of 2-representations. To this end, suppose that $e_i \in H_\alpha$, $e_\mu
\in H_\beta$ and $e_\kappa \in H_\gamma$ are the marked simple objects in the underlying 2-Hilbert spaces of $\alpha$, $\beta$ and $\gamma$ respectively, and that $u : e_{g \ccdot i} \rightarrow
\alpha_g(e_i)$, $v : e_{g \ccdot \mu} \rightarrow \beta_g(e_\mu)$ and $w : e_{g \ccdot \kappa} \rightarrow
\gamma_g(e_\kappa)$, are unitary morphisms. We must check that the following diagram commutes:
 \[
  \xymatrix{ \langle \kappa | \rho | \mu \rangle \otimes \langle \mu | \sigma | i \rangle \ar[d]_{ \circ}
  \ar[rrrr]^-{\displaystyle \langle \rho \rangle (\lag{g}{w \otimes \bar{v}}) \otimes \langle \sigma \rangle
   (\lag{g}{v \otimes \bar{u}})}
  &&&&
  \langle g \ccdot \kappa| \rho | g \ccdot \mu \rangle \otimes \langle g \ccdot \mu | \sigma | g \ccdot i \rangle
  \ar[d]^\circ \\
  \langle \kappa | \rho \sigma | i \rangle \ar[rrrr]_{\displaystyle \langle \rho \circ \sigma \rangle (
  \lag{g}{w \otimes \bar{u}})} &&&& \langle g  \ccdot\kappa | \rho \circ \sigma | g \ccdot i \rangle}
 \]
Indeed, for $\lambda_2 \in \langle \kappa | \rho | \mu \rangle$ and $\lambda_1 \in \langle \mu | \sigma | i
\rangle$, the counterclockwise and clockwise directions evaluate as
 \begin{align*}
 \ba \xymatrix @R=1pc @C=1pc {\ar[d]\\ \ar[r] & } \ea \left( \lambda_2 \otimes \lambda_1\right) &= \contraction{\rho\sigma(u)^* \circ \rho\bigl(\sigma_{e_i}(g)\bigr) \circ}{\rho_{\sigma(e_i)}(g)}{ \circ \gamma_g
  \bigl(}{\rho(\lambda_1}
{\rho\sigma(u)^* \circ \rho\bigl(\sigma_{e_i}(g)\bigr) \circ}{\rho_{\sigma(e_i)}(g)}{ \circ \gamma_g
  \bigl(}{\rho(\lambda_1}  \bigr) \circ \lambda_2)\circ u'' \\
  & = \rho \biggl(\sigma(u)^* \sigma_{e_i}(g) \beta_g(\lambda_1) u' \biggr) \circ \rho(u')^*
  \circ\rho_{e_\mu}(g) \circ \gamma_g(\lambda_2) \circ u'' \\
  & = \ba \xymatrix @R=1pc @C=1pc {\ar[r] & \ar[d] \\ {} & {}} \ea \left( \lambda_2 \otimes \lambda_1\right)
 \end{align*}
where the middle step uses the naturality of $\rho(g) : \gamma \circ \rho \Rightarrow \rho \circ \beta$, and the
$u'$ term drops out since $u' u'^* = \id$. This establishes that the map
 \[
  \TRep_\text{m}  \rightarrow \Gerbes(G)
 \]
is indeed a weak 2-functor. To see that it is an equivalence, the most direct method is to restrict attention to unitary 2-representations on skeletal 2-Hilbert spaces (all unitary 2-representations are equivalent to unitary 2-representations of this form, since every 2-Hilbert space is equivalent to a skeletal 2-Hilbert space). To give a skeletal 2-Hilbert space $\alpha$ the structure of a unitary 2-representation is the same thing as defining a $G$-action on the set $X_\alpha$ of simple objects, together with a groupoid 2-cocycle $c_\alpha \in Z^2(X_G, U(1))$. In this picture, a morphism $\sigma \colon \alpha \rightarrow \beta$ between skeletal 2-representations is the same thing as a $(c_\beta \times \overline{c_\alpha})$-twisted equivariant vector bundle over the groupoid $X_\beta \times X_\alpha$. Similarly a 2-morphism is just a morphism of equivariant vector bundles.  These data structures are precisely the same structures which define the objects, morphisms and 2-morphisms of $\Gerbes(G)$ once we pass to a cocycle description, as we explained in Chapter \ref{GerbesChap}. Therefore the 2-functor $\TRep_\text{m} \rightarrow \Gerbes(G)$ is indeed an equivalence of 2-categories.
\end{proof}
Combining this result with the classification of equivariant gerbes from Proposition \ref{equivgerbes} allows us to rederive some known results about 2-representations, but specialized to the unitary setting. We say that a 2-representation is {\em irreducible} if its associated equivariant gerbe has only a single orbit.
\begin{cor}[{Compare Ostrik \cite[Ex. 3.4]{ref:ostrik}, Elgueta \cite[Thm 7.5]{ref:elgueta}, Ganter and Kapranov \cite[Prop
7.3]{ref:ganter_kapranov_rep_char_theory}}] Irreducible unitary 2-representations of $G$ are classified up
to strong unitary equivalence by triples $(X, k, [\phi])$, where $X$ is a transitive $G$-set up to isomorphism, $k$ is a positive real number, and $[\phi]$ is an equivariant cohomology class $[\phi] \in H^2( X_G, U(1))$.
\end{cor}
\noindent We write $1$ for the trivial 2-representation of $G$ on $\Hilb$ (it is the unit object for the monoidal 2-category structure of $\TRep(G)$, but we will not discuss this here).
\begin{cor}[{Compare Elgueta \cite[Cor 6.21]{ref:elgueta}, Ganter and Kapranov \cite[Ex. 5.1]{ref:ganter_kapranov_rep_char_theory}}]  The endomorphism category of the unit object is monoidally equivalent to the category of representations of
$G$,
 \[
  \End_{\TRep(G)} (1) \simeq \Rep(G).
 \]
More generally, if $\alpha$ is any one-dimensional 2-representation of $G$, then there is a unitary equivalence
of 2-Hilbert spaces
 \[
  \Hom_{\TRep(G)} (1, \alpha) \simeq \Rep^\phi (G)
 \]
where $\phi \in Z^2(G, U(1))$ is the group 2-cocycle obtained from choosing a section of $\X_\alpha$.
\end{cor}
We would like to stress however that the geometric language of equivariant gerbes allows us to {\em go further} than these results, because the gerbal technology we developed in Chapter \ref{GerbesChap} allows us to give us a concrete geometric understanding of {\em all} the hom-sets in $[\TRep(G)]_\mathbb{C}$.
 \begin{cor} \label{tre}The space of morphisms between unitary 2-representations in $[\TRep(G)]_\mathbb{C}$ identifies as the space of flat sections of the transgressed line bundle over the loop groupoid of the product of their associated $G$-sets:
  \[
   \Hom(\alpha, \beta)_{[\TRep(G)]_\mathbb{C}} \cong \Gamma_{\Lambda (X_\beta \times X_\alpha)_G} (\tau (\X_\beta \otimes \overline{\X_\alpha})_\mathbb{C}).
  \]
Thus their dimensions can be calculated by integrating the transgression of the transgressed line bundle over the double loop groupoid:
  \[
   \dim \Hom (\alpha, \beta) = \int_{\Lambda^2 (X_\beta \times X_\alpha)_G} \tau^2(\X_\beta \otimes \overline{\X_\alpha}).
  \]
In particular, the dimension of the space of endomorphisms of an object computes as
 \[
  \dim \End (\alpha) = \frac{1}{|G|} |\{(i, j, g, h) : i,j  \in X_\alpha,  g,h \in G, i,j \in \FFix(g) \cap \FFix(h), gh=hg\}|.
 \]
\end{cor}
Having established the correspondence between unitary 2-representations and equivariant gerbes, we now turn to their {\em characters}.

\section{The 2-character and the geometric character\label{2char2rep}}
In this section we prove our main result in this chapter --- that the 2-character of a unitary 2-representation corresponds naturally to the geometric character of its associated equivariant gerbe, and hence the 2-character is a unitarily fully faithful functor at the level of the complexified Grothendieck category.
\begin{thm} \label{transthm} The 2-character of a marked unitary 2-representation is unitarily naturally
isomorphic to the geometric character (i.e. the push-forward of the transgression) of the associated equivariant
gerbe:
 \[
 \xymatrix @1 @C=0.1in{[\TRep_m(G)] \ar[dr]_{\chi} \ar[rr]^{\sim} && [\Gerbes(G)] \ar[dl]^{\ch} \\ & \Hilb_G(G)}
 \]
\begin{textblock}{0.3}(0.265,-0.0515)
     $\ba \ig{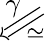} \ea$
\end{textblock}
\noindent That is, there are unitary isomorphisms $\gamma_\alpha \colon \chi_\alpha \stackrel{\cong}{\rightarrow} \ch(\X_\alpha)$,
natural in $\alpha$.
\end{thm}
Combining this with our knowledge of the geometric character functor from Theorem \ref{geomchthm} gives us one of the three main results in this thesis that we mentioned in the introduction.
\begin{thm}\label{corc} The complexified 2-character functor
 \[
  \chi_\mathbb{C} \colon [\TRep(G)]_\mathbb{C} \rightarrow \Hilb_G(G)
  \]
is a unitarily fully faithful functor from the complexified Grothendieck category of unitary representations of $G$ to the category of unitary equivariant vector bundles over $G$.
\end{thm}
Proving Theorem \ref{transthm} involves expanding out the abstract higher-categorical definitions for the 2-character and checking that they have the appropriate geometric behaviour. We do this in three steps --- firstly we define the isomorphisms $\gamma_\alpha$, then we show that they are indeed morphisms of equivariant vector bundles, and then we show that they are natural in $\alpha$.

\subsection{Defining the isomorphisms}
Given a marked 2-representation $\alpha$ the isomorphism of
equivariant vector bundles $\gamma_\alpha \colon \chi_\alpha \rightarrow \ch(\X_\alpha)$ is easy enough to write
down. For $g \in G$, the fibers of the 2-character compute, by definition, as
 \begin{align*}
  \chi_\alpha (g) &= \Nat (\id, \alpha_g) \\
   & \cong \{ (\theta_{e_i} \colon e_i \rightarrow \alpha_g(e_i))_{i \in \FFix(g)} \} ,
 \end{align*}
while the fibers of the geometric character are
 \begin{eqnarray*}
   \ch (\X_\alpha)(g) &=& \Sections \bigl(\tau(\X_{\alpha})_\mathbb{C} |_{\FFix(g)} \bigr) \\
    &=& \{ (\vartheta_i \in \text{uHom}(e_i, \alpha_g (e_i)) \otimes \mathbb{C})_{i \in \FFix(x)} \} .
  \end{eqnarray*}
So the fibrewise identification between these two complex lines is component-wise just the identification between a hermitian line and the line associated to its circle of radius $\sqrt{k_i}$. Recalling our conventions about $U(1)$-torsors and hermitian lines from Section \ref{uonebund}, the isomorphisms $\gamma_\alpha$ are given at the level of unitary morphisms $u \colon e_i \rightarrow \alpha_g (e_i)$ by
 \begin{align*}
  \gamma_\alpha \colon \chi_\alpha(g) & \rightarrow \ch (\X_\alpha)(g) \\
   u & \mapsto u \otimes \sqrt{k_i}.
  \end{align*}
Moreover, recalling our conventions about the inner products on $\Nat(\id, \alpha_g)$ and $\Sections \bigl(\tau(\X_\alpha)_\mathbb{C} |_{\FFix(g)} \bigr)$ from Chapters \ref{BigTCat} and \ref{deftransg} respectively, one sees that $\gamma_\alpha$ is indeed a unitary isomorphism.

\subsection{Verifying that the isomorphisms are equivariant}
The following lemma verifies that these fibrewise identifications are equivariant with respect to the action of
$G$, which for the 2-character is given by the string diagram formula from Chapter \ref{charsecc} and for the geometric
character by the transgression formula from Chapter \ref{transsec}. In other words, $\gamma_\alpha$ is indeed a morphism
in $\Hilb_G(G)$.

\begin{lem} \label{transglem} Let $\alpha$ be a marked 2-representation of $G$.
\begin{enumerate}
 \item The equivariant maps for the 2-character $\chi_\alpha$ compute as the
transgression, in the sense that
 \begin{align*}
  \chi_\alpha(hxh^\mi \sla{h} g) \colon \Nat(\id, \alpha_g) & \rightarrow \Nat(\id, \alpha_{hgh^\mi}) \\
  \intertext{evaluates as}
  \chi_\alpha (hgh^\mi \sla{h} g)(\theta)_{e_{h \ccdot i}} &= v \diamond \theta_{e_i} \diamond v^{\underline{\mi}}
  \end{align*}
where $v \colon e_{h \ccdot i} \rightarrow \alpha_h (e_i)$ is any unitary arrow in the underlying 2-Hilbert space
$H_\alpha$, and $\diamond$ is the twisted composition law from Chapter \ref{extract}.
 \item Therefore the following diagram commutes:
 \[
  \ba \xymatrix{\chi_\alpha (g) \ar[d]_{\gamma_\alpha(g)} \ar[rr]^{\chi_\alpha (\sla{h} g)} && \chi_\alpha (hgh^\mi) \ar[d]^{\gamma_\alpha
  (hgh^\mi)} \\
  \ch (\X_\alpha)(g) \ar[rr]_{\ch(\X_\alpha)(\sla{h} g)} && \ch (\X_\alpha)(hgh^\mi)} \ea\, .
 \]
\end{enumerate}
\end{lem}
\begin{proof} (i) We need to evaluate the string diagram formula for the map
$\chi_\alpha(\sla{h} g)$,
 \[
  \ba \ig{c3.pdf} \ea \mapsto \ba \ig{c4.pdf} \ea.
   \]
The right hand side computes as
\begin{eqnarray*}
 \bigl[\chi_\alpha (hgh^\mi \sla{h} g)(\theta)\bigr]_{e_{h \ccdot i}} &:=&  a(hg, h^\mi)_{e_{h \ccdot i}} a(h,g)_{\alpha_{h^\mi} (e_h \ccdot i)}
  \alpha_h (\theta_{\alpha_{h^\mi}(e_{h \ccdot i})}) \\ && \qquad \qquad \qquad \qquad \qquad a^*(h, h^\mi)_{e_{h \ccdot i}} a(e)_{e_{h \ccdot i}}  \\
   &\stackrel{\text{(a)}}{=}& a(hg, h^\mi)_{e_{h \ccdot i}} a(h,g)_{\alpha_{h^\mi} (e_h \ccdot i)}
  \alpha_h (\theta_{\alpha_{h^\mi}(e_{h \ccdot i})}) \alpha_h(v^{\underline{\mi}}) v  \\
  &\stackrel{\text{(b)}}{=}& a(hg, h^\mi)_{e_{h \ccdot i}} a(h,g)_{\alpha_{h^\mi} (e_h \ccdot i)}
  \alpha_h(\alpha_g(v^{\underline{\mi}})) \alpha_h(\theta_{e_i})v \\
 &\stackrel{\text{(c)}}{=}& a(hg, h^\mi)_{e_{h \ccdot i}} \alpha_{hg}(v^{\underline{\mi}}) a(h,g)_{e_i} \alpha_h(\theta_{e_i})  v \\
    &\stackrel{\text{(d)}}{=}& (v \diamond \theta_{e_i}) \diamond v^{\underline{\mi}}
 \end{eqnarray*}
where (a) uses the expression \eqref{inverselaw} for $v^{\underline{\mi}}$ with respect to the twisted composition law $\diamond$ from \eqref{gerbecomp}, (b) uses the naturality of $\theta$, (c)
uses the naturality of $a(h,g)$, and (d) again uses the  composition law $\diamond$ from
\eqref{gerbecomp}.

(ii) This is just (i), expressed more formally.
\end{proof}

\subsection{Verifying naturality}
It remains to show that the isomorphisms $\gamma_\alpha \colon \chi_\alpha \rightarrow \ch (\X_\alpha)$ are natural
with respect to morphisms of 2-representations $\sigma \colon \alpha \rightarrow \beta$. This amounts to computing
the string diagram formula \eqref{stdiagmor} for $\chi(\sigma)$ --- the action of the {\em 2-character} on
morphisms --- and observing that it corresponds to the formula \eqref{matelements} for the behaviour of the
geometric character on morphisms, which was defined in terms of the complex conjugate of the {\em
ordinary} character $\chi_{\langle \sigma \rangle}$ of the equivariant vector bundle $\langle \sigma \rangle$. In
a slogan, `the 2-character on morphisms is the ordinary character'.

\begin{lem} Let $\sigma \colon \alpha \rightarrow \beta$ be a morphism of marked 2-representations of $G$.
\begin{enumerate}
 \item The matrix elements of $\chi(\sigma)$ compute as the complex conjugate of the trace of the associated equivariant vector bundle $\langle \sigma \rangle$,
  \[
   \langle \theta_\mu, \chi(\sigma)_g \theta_i\rangle = \Tr \Bigl( \langle \mu | \sigma | i \rangle
   (\lag{g}{\theta_\mu \otimes \theta_i}) \Bigr)^*.
  \]
 \item Therefore the isomorphism $\gamma_\alpha \colon \chi_\alpha \sra{\cong} \ch(\X_\alpha)$ is natural in
 $\alpha$, that is for each $g \in G$ the following diagram commutes:
  \[
   \xymatrix{\chi_\alpha(g) \ar[r]^{\chi(\sigma)_g} \ar[d]_{\gamma_\alpha} & \chi_\beta(g) \ar[d]^{\gamma_\beta} \\
             \ch(\X_\alpha)(g) \ar[r]_{\ch(\langle \sigma \rangle)_g} & \ch(\X_\beta)(g)}
  \]
  \end{enumerate}
\end{lem}
\begin{proof}
(i) We need to evaluate the string diagram formula for the map $\chi(\sigma)_g$,
\[
\ba \ig{c3.pdf} \ea \mapsto \ba \ig{y28.pdf} \ea ,
\]
where $(\sigma^*, \eta, \epsilon)$ is some choice of right adjoint for $\sigma$ and $\Psi$ is the even-handed structure on $\THilb$. This gives
 \be \label{orrform}
 \chi(\sigma)_g(\theta)_{e_\mu} = \beta_g(\epsilon_{e_\mu}) \circ \sigma(g)^*_{\sigma^*(e_\mu)} \circ \sigma(\theta_{\sigma^*(e_\mu)}) \circ \Psi(\eta)_{e_\mu}.
 \ee
Our main task in evaluating this expression further is to compute the term $\Psi(\eta)_{e_\mu}$; to do this we need to translate the even-handed structure $\phi \stackrel{\Psi}{\mapsto} * \, \phi^* *$ on $\THilb$ from the language of adjunction isomorphisms into the language of unit and counit natural transformations.

Now, we know from Theorem \ref{Thecharthm} that the formula above does not depend on the right adjoint $(\sigma^*, \eta, \epsilon)$ that was chosen to evaluate it. If we fix the underlying functor $\sigma^*$, then the freedom in choosing the unit and counit maps $\eta$ and $\epsilon$ is given at the level of the adjunction isomorphisms between the hom-sets by the freedom in choosing arbitrary vector space isomorphisms
 \[
  \phi_{i, \mu} \colon \Hom(\sigma(e_i), e_\mu) \rightarrow \Hom(e_i, \sigma^*(e_\mu))
 \]
at the level of the simple objects. It is convenient to choose these isomorphisms $\phi_{i, \mu}$ to be `scaled unitary' in the sense that
 \[
  ( \phi_{i, \mu}(v), \phi_{i, \mu}(w)) = \frac{k_i}{k_\mu}  (v, w)
 \]
for $v,w \in \Hom(\sigma(e_i), e_\mu)$, or in other words
 \be \label{bt}
  \phi_{i, \mu}^* \phi_{i, \mu} = \frac{k_i}{k_\mu} \id.
 \ee
The advantage of this choice will presently become clear. Choose a {\em $*$-basis}
 \[
  \{ v^{(\mu, i)}_p \colon e_\mu \rightarrow \sigma(e_i)\}_{p=1}^{\dim \langle \mu | \sigma | i \rangle}
 \]
for each hom-space $\Hom(e_\mu, \sigma(e_i))$ --- that is, a basis satisfying
 \[
  v^{(\mu, i)*}_p \circ v^{(\mu, i)}_q = \delta_{pq} \id_{e_\mu} \quad \text{and hence} \quad \sum_{\mu, p} v^{(\mu, i)}_p \circ v^{(\mu, i)*}_p = \id_{\sigma(e_i)}.
 \]
In terms of such a basis, we claim that $\Psi(\eta)_{e_\mu}$ computes as
 \be \label{ieqn}
  \Psi(\eta)_{e_\mu} = \frac{1}{k_\mu} \sum_i k_i \sum_p \sigma(w^{(i, \mu)}_p) \circ v^{(\mu, i)}_p
 \ee
where $w^{(i, \mu)}_p \colon e_i \rightarrow \sigma^*(e_\mu)$ is defined as $\phi (v_p^{(\mu, i)*})$. To see this, firstly observe that the $w^{(i, \mu)}_p$ form a $*$-basis for $\Hom(e_i, \sigma^*(e_\mu))$ because
 \begin{align*}
  w^{(i, \mu)*}_p \circ w^{(i, \mu)}_q &= \frac{1}{k_i} (\id_{e_i},  w^{(\mu, i)*}_p \circ w^{(\mu, i)}_q) \; \id_{e_i} \\
   &= \frac{1}{k_i} ( w^{(i, \mu)}_p,  w^{(i, \mu)}_q) \; \id_{e_i} \\
   &= \frac{1}{k_i} ( \phi (v_p^{(\mu, i)*}),  \phi (v_q^{(\mu, i)*})) \; \id_{e_i} \\
   &= \frac{1}{k_i} \frac{k_i}{k_\mu} (v_p^{(\mu, i)*},  v_q^{(\mu, i)*}) \; \id_{e_i} \\
   &= \frac{1}{k_\mu} (\id_{e_\mu}, \; v_q^{(\mu, i)*}  v_p^{(\mu, i)}) \; \id_{e_i} \\
   &= \delta_{pq} \id_{e_i}.
 \end{align*}
(This is the utility of choosing the isomorphisms $\phi$ to be scaled unitary.) To prove \eqref{ieqn}, all we need to show is that when we translate this formula back into a hom-set isomorphism
 \[
  \psi \colon \Hom(\sigma^*(e_\mu), e_i) \rightarrow \Hom(e_\mu, \sigma(e_i))
 \]
via the standard translation between these two pictures then we must end up with $\psi = * \, \phi^* *$. Indeed, we compute:
 \begin{align*}
  \psi (w^{(i, \mu)*}_q) &:= \sigma (w^{(i, \mu)*}_q) \circ \Psi(\eta)_{e_\mu} \\
   &= \frac{1}{k_\mu} \sum_i k_i \sum_p \sigma(w^{(i, \mu)*}_q) \circ \sigma(w^{(i, \mu)}_p) \circ v^{(\mu, i)}_p \\
   &= \frac{k_i}{k_\mu} v^{(\mu, i)}_p \\
   &= * \, \phi^* * (w^{(i, \mu)*}_q) \quad \text{(by \eqref{bt})}.
  \end{align*}
A similar argument establishes that
 \be \label{step3}
  \epsilon_{e_\mu} = \sum_{i, p} v^{(\mu, i)*} \circ \sigma(w^{(i, \mu)*}).
 \ee
Moreover, by our formula from Lemma \ref{unique} we compute that
 \be \label{step4}
  \theta_{\sigma^*(e_\mu)} = \sum_{i,p} \alpha_g (w^{(i, \mu)}_p) \circ \theta_{e_i} \circ w^{(i, \mu)}_p.
 \ee
and similarly
 \be \label{step5}
  \sigma(g)^*_{\sigma^*(e_\mu)} = \sum_{i,p} \beta_g \sigma(w^{(i, \mu)}_p) \sigma(g)^*_{e_i} \sigma \alpha_g(w^{(i, \mu)*}_p).
 \ee
Substituting the expansions \eqref{ieqn}, \eqref{step3}, \eqref{step4} and \eqref{step5} into our original string diagram formula \eqref{orrform}, most of the indices pair off and one computes that for a natural transformation $\theta_i \in \chi_\alpha (g)$ supported on a single fixed point $x \in \FFix_X (g)$ (that is, the components $(\theta_i)_{e_j}$ of $\theta_i$ over the marked simple objects are zero unless $i=j$), we have:
 \[
    \chi(\sigma)_g(\theta_i)_{e_\mu} = \frac{k_i}{k_\mu} \sum_p \beta_x (v^{(\mu, i)*}_p) \circ \sigma(g)^*_{e_i} \circ \sigma(\theta_{e_i}) \circ v^{(\mu, i)}_p.
\]
Notice how the scale factors $k_i$ and $k_\mu$ have entered this description --- this underscores our point that this map {\em uses the even-handed structure in an intrinsic way}. We can identify this combination of terms as the complex conjugate of the trace of the associated equivariant vector bundle $\langle \sigma \rangle$ as follows. Fix an orthonormal basis $\{\theta_i \in \chi_\alpha (x)\}$ and $\{\theta_\mu \in \chi_\beta(x)\}$ of natural transformations supported exclusively over fixed points $i \in \FFix_X (g)$ and $\mu \in \FFix_{X'} (g)$. Recalling the relevant inner product from Chapter \ref{BigTCat}, the matrix elements in this basis are thus
 \begin{align*}
  \langle \theta_\mu, \, \chi(\sigma)_g \theta_i \rangle &= k_\mu (\theta_\mu, \, \chi(\sigma)_g \theta_i) \\
   &= \sum_p k_i (\theta_{e_\mu}, \,  \beta_g(v^{(\mu, i)*}_p)\sigma(g)_{e_i}^* \sigma(\theta_{e_i}) v^{(\mu, i)}_p) \\
   &= \frac{1}{k_\mu} \sum_p (v^{(\mu, i)}_p, \sigma(k_i \theta_{e_i}^*) \sigma(g)_{e_i} \beta_g(v^{(\mu, i)}_p) k_\mu \theta_{e_\mu})^* \\
   &= \Tr \Bigl( \langle \mu | \sigma | i \rangle
   (\lag{g}{\theta_\mu \otimes \bar{\theta_i}}) \Bigr)^*,
   \end{align*}
where the last step uses the definition of the equivariant vector bundle $\langle \sigma \rangle$ from
\eqref{eqdefnmor} (we needed to use $k_i \theta_{e_i}$ because $\theta_i$ was an orthonormal basis vector, so that $(\theta_{e_i}, \theta_{e_i}) = \frac{1}{k_i}$, and similarly for $\theta_{e_\mu}$). We also used the fact that the trace of a linear endomorphism $A$ of the Hilbert space $\langle \mu | \sigma | i \rangle$ can be expressed in terms of a $*$-basis $\{v_p\}$ as
 \[
  \Tr(A) = \frac{1}{k_\mu} \sum_p (v_p, Av_p)
 \]
since we must account for the fact that the basis vectors $v_p$ are orthogonal but not orthonormal:
 \[
  (v_p, v_q) = (\id_{e_\mu}, v_p^* c_q) = \delta_{pq} (\id_{e_\mu}, \id_{e_\mu}) = \delta_{pq} k_\mu.
 \]

(ii) This follows immediately from comparing the matrix elements of $\chi(\sigma)_g$ above to the matrix elements
of $\ch(\langle \sigma \rangle)_g$ given in equation \eqref{matelements}.
 \end{proof}

This completes the proof of Theorem \ref{transthm}.

\section{The 2-character does not distinguish 2-representations\label{notdistinguish}}
In this section we observe another corollary of our theorem that the 2-character of a 2-representation corresponds to the geometric character of its associated equivariant gerbe. Namely that, contrary to {\em ordinary} representation theory, the 2-character does not distinguish 2-representations.

\begin{cor}[{of Theorem \ref{transthm}}] The 2-character does not distinguish unitary 2-representations. That is, there exist nonequivalent unitary 2-representations $\alpha$ and $\beta$ having the property that $\chi_\alpha \cong \chi_\beta$ in $\Hilb_G (G)$.
\end{cor}
\begin{proof} We need only check this statement at the level of geometric characters of equivariant gerbes. Every finite $G$-set $X$ can be regarded as a `trivial' or `untwisted' equivariant gerbe $\X$ by replacing each arrow in the action groupoid $X_G$ with a copy of $U(1)$. The geometric character of such a gerbe is the equivariant vector bundle over $G$ whose fibers are the space of functions on the fixed points of $G$ on $X$:
 \[
  \ch(\X)_g = \{f \colon \FFix_X (g) \rightarrow \mathbb{C}\}.
 \]
The conjugation action on these fibers is given by conjugating the arguments of the functions; this completes the description of the geometric character $\ch(\X)$ as a representation of the loop groupoid $\Lambda \BG$. Now, a representation of a groupoid is determined up to isomorphism by its {\em character}, which is the function on the {\em loop groupoid} of the original groupoid defined by taking the trace of the representation over the endomorphisms in the groupoid (see Willerton \cite{ref:simon}). In our case, the character of the geometric character $\chi_{\ch(\X)}$ is a function on the {\em double loop groupoid} $\Lambda^2 \BG$, and unraveling the definition, we see that it assigns to every pair of commuting elements $g,h \in G$ the number of points $x \in X$ which are {\em simultaneously fixed by $g$ and $h$}:
 \begin{align*}
  \chi_{\ch(\X)} \colon \Lambda^2 \BG & \rightarrow \mathbb{C} \\
   (g,h) & \mapsto |\{ x \in X\colon g \ccdot x = h \ccdot x = x\}|.
  \end{align*}
We call this function the {\em double character} of the $G$-set $X$, and it can happen that two nonisomorphic $G$-sets $X$ and $Y$ have the same double character, in much the same way that two nonisomorphic $G$-sets can have the same {\em ordinary} character (that is, share the same number of fixed points for each group element). See \cite{ref:webster, ref:desmit_lenstra}. This is because two $G$-sets $X$ and $Y$ are isomorphic if and only if each {\em subgroup} of $G$ has the same number of fixed points in both $X$ and $Y$ --- this is basically the Yoneda lemma in the category of $G$-sets. This means it is not enough to know the number of fixed points for pairs of commuting elements, since these do not generate the most general subgroups that could occur.

Thus we have shown that nonequivalent 2-representations can give rise to isomorphic 2-characters.
\end{proof}

How are we to interpret this result? It seems best to understand it as a reaffirmation of the fact that {\em 2-Hilbert spaces, or semisimple categories in general, are not the final word on what a `2-vector space' should be}. They do not have continuous parameters and hence can only accommodate `permutation-type' 2-representations, and this is the phenomenon which is at play above. We can also interpret this result as stating that the functor $[\TRep(G)] \rightarrow [\TRep(G)]_\mathbb{C}$ does not reflect isomorphisms. 
\chapter{The higher categorical dimension of 2$\mathcal{R}\text{ep}(G)$.\label{Dim2RepGchap}}

In this chapter we compute $\Dim \TRep(G)$, the `higher-categorical dimension' of the 2-category $\TRep(G)$ of unitary 2-representations of $G$. By the `higher-categorical dimension' of a 2-category $\A$ we mean the category $\Hom(\id_\A, \id_\A)$ of transformations of the identity 2-functor on $\A$ and modifications between them; this is also known as the {\em generalized centre} of the 2-category $\A$ \cite{ref:street_conspectus, ref:baez_dolan_hda0}. By well-known higher-categorical reasoning which we will explain in detail shortly, the higher-categorical dimension of a 2-category forms a {\em braided monoidal} category, and our main result is as follows.
\begin{thmt} The higher-categorical dimension of the 2-category of unitary 2-representations of a finite group $G$ is equivalent, as a braided monoidal category, to the category of conjugation-equivariant hermitian vector bundles over $G$ equipped with the fusion tensor product. In symbols,
 \[
  \Dim \TRep(G) \simeq \Hilb_G^\fusion (G).
 \]
\end{thmt}
The `fusion tensor product' on $\Hilb^\text{fusion}_G (G)$ refers to the tensor product operation on the `category assigned to the circle' in three-tier topological quantum field theory, obtained by using the pair of pants to fuse $G$-bundles over $S^1$ together as in the lecture notes of Freed \cite{ref:freed2} (see also \cite{ref:freed_remarks, ref:freed_teleman_hopkins}). We have reviewed this braided monoidal structure on $\Hilb_G^\fusion (G)$ in Appendix \ref{AppFreed}.

Therefore, as we explained in the introduction, the relevance of this result is that it provides further evidence for the validity of the Extended TQFT Hypothesis of Baez and Dolan, at least in the case of the three-dimensional untwisted finite group model. This is because $\TRep(G)$ can be interpreted as the 2-category assigned to the point while $\Hilb_G^\fusion (G)$ can be interpreted as the category assigned to the circle. Hence this theorem confirms the `trace equation'
 \[
\Dim Z(M) \simeq Z(M \times S^1)
 \]
that one would expect to hold for all closed manifolds $M$ with $\dim M \leq 2$ if the Extended TQFT Hypothesis were to be true. Our theorem above deals with the case when $M$ is a point; the remaining cases are dealt with in Appendix \ref{AppendixTQFTFacts}.

Our strategy for proving this theorem is as follows. We do not calculate the higher-categorical dimension of $\TRep(G)$ directly; rather we use the fact that $\TRep(G)$ is equivalent to the 2-category $\Gerbes(G)$ of equivariant gerbes, and we instead calculate the higher-categorical dimension of $\Gerbes(G)$. Providing the assumption in the theorem about equivalences between 2-categories is indeed true, this will give the result. We will spell out our reasons for adopting this strategy in Section \ref{why} below, but the idea is essentially that the geometric framework of equivariant gerbes enables one to write down {\em precise formulas} for the various constructions involved {\em without making any choices}, a feature which is sadly not available if one works directly with the 2-category $\TRep(G)$.

In Section \ref{la} we recall the notion of the higher-categorical dimension of a 2-category, and we explain why it is a braided monoidal category, both by using the periodic table of Baez and Dolan as well by employing three-dimensional string diagrams. In Section \ref{lb} we illustrate the procedure we will use for the proof by carefully running through the same algorithm one categorical level down. We show step by step in a series of lemmas how the categorical dimension of the category $\Rep(G)$ of ordinary representations of a group computes as $Z(\mathbb{C}[G])$, the centre of the group algebra. This is not a new result, but the idea is to provide a service for the reader: the rest of the chapter `categorifies' these lemmas one by one, so that if confused the reader can return to this section to regain his or her bearings.

In Section \ref{why} we explain why we have been forced into the indirect route of calculating $\Dim \Gerbes(G)$ instead of calculating $\Dim \TRep(G)$ directly. In the remaining sections we proceed with our calculation of $\Dim \Gerbes(G)$, by establishing categorified analogues of each lemma in Section \ref{lb} in the world of equivariant gerbes. We begin in Section \ref{lc} by showing that the hom-category $\Hom(\EG, \X)$ of morphisms emanating out of the `regular equivariant gerbe' $\E G$ is generated by certain specific morphisms $U_x \colon \EG \rightarrow \X$ as $x$ ranges over $\X$. In Section \ref{ld} we show how to extract an equivariant vector bundle from a transformation of the identity on $\Gerbes(G)|_{\E G}$. In Section \ref{extendd} we show how, if one is given a transformation of the identity on $\Gerbes(G)|_{\E G}$, then one can {\em extend} that transformation to the entire 2-category. In Section \ref{le} we show that restricting and then extending transformations in this way recovers back the original transformation. In Section \ref{lf} we show that a {\em morphism} of transformations is uniquely and freely determined by its behaviour at the regular equivariant gerbe $\E G$. Finally in Section \ref{lg} we wrap up the proof.

\section{Higher categorical dimensions of 2-categories\label{la}}
From every 2-category $\mathcal{A}$ one can produce a braided monoidal category $\Dim \mathcal{A}$ which we call the {\em higher categorical dimension} (or just {\em dimension}) of $\mathcal{A}$. One defines it by setting
 \[
  \Dim \mathcal{A} := \Hom(\id_\mathcal{A}, \id_\mathcal{A}),
 \]
the category whose objects are transformations of the identity 2-functor on $\mathcal{A}$, and whose morphisms are modifications. This category is more commonly known as the {\em centre} or {\em generalized centre} of $\mathcal{A}$ (see eg. \cite{ref:street_conspectus, ref:baez_dolan_hda0}), because it can be thought of as a many-object generalization of notion of the {\em centre} $Z(C)$ of a monoidal category $C$, due to Joyal and Street\cite{ref:joyal_street_3}. Recall that the centre $Z(C)$ of a monoidal category is the braided monoidal category whose objects are pairs
 \[
  \underline{V} = (V, \{\theta^V_W \colon V \otimes W \rightarrow W \otimes V\}_{W \in C})
 \]
consisting of an object $V \in C$ together with a prescribed isomorphism $\theta^V_W \colon V \otimes W \rightarrow W \otimes V$ for every $W \in C$, satisfying various coherence diagrams. It turns out that if one thinks of the monoidal category $C$ as a one-object 2-category, then $\underline{V}$ corresponds precisely to a transformation of the identity 2-functor on $C$, at least if one reverses the normal convention for directions of transformation cells.

That is to say, so far in this thesis we have followed the widely-held convention (see eg. Leinster \cite{ref:leinster_basic_bicategories}) that the data of a transformation $\sigma \colon F \Rightarrow G$ between weak 2-functors $F, G \colon \mathcal{A} \rightarrow \mathcal{B}$ consists of, for each object $a \in \mathcal{A}$, a 1-morphism $\sigma_a \colon Fa \rightarrow Ga$ in $\mathcal{B}$, and for each 1-morphism $a \stackrel{f}{\rightarrow} b$ in $\mathcal{A}$, a 2-isomorphism $\sigma(f) \colon Gf \circ \sigma_a \Rightarrow \sigma_b \circ Ff$.
\vskip 0.3cm
\framebox{\parbox[b]{11cm}{ In order to align with the standard description of the centre of a monoidal category, in this chapter we will be using the opposite convention for transformations.\label{jimbo}}}
\vskip 0.3cm

That is, $\sigma(f)$ is taken to be a 2-isomorphism $\sigma(f) \colon \sigma_b \circ Ff \Rightarrow Gf \circ \sigma_a$. In the case of a monoidal category (with objects $V, W$ etc.) thought of as a 2-category, and where $\sigma \equiv \underline{V} \colon \id \Rightarrow \id$ is a transformation of the identity 2-functor, this means that the coherence cells $V(W)$ run in the direction
 \[
  \underline{V}(W) \colon V \otimes W \rightarrow W \otimes V
 \]
which is the ordinary convention for the centre of a monoidal category.

Before we explain in more detail the notion of the higher-categorical dimension of a 2-category $\mathcal{A}$, let us explain why we have called it the `dimension' of $\mathcal{A}$:
 \begin{itemize}
  \item The kind of 2-categories we are concerned with are the {\em higher-vector spaces} which arise in extended TQFT. We think of $\Dim \mathcal{A}$ as the higher-categorical trace of the identity 2-functor on $\mathcal{A}$, hence the word `dimension'.
  \item We would like a notation whereby the important formula
   \[
    Z(M \times S^1) \simeq \Dim Z(M)
    \]
  of extended TQFT, which we are trying to explicate in this thesis, continues to hold no matter which categorical level one is working at.
  \item We would like to avoid confusion with the more well-known concept of the `centre' of a {\em monoidal} higher-categorical structure, For instance, the 2-category $\TRep(G)$ is indeed a monoidal 2-category (a feature we have not discussed in this thesis), so that were we to use the word `centre' the uninitiated reader might imagine we were constructing a braided monoidal {\em 2-category} from $\TRep(G)$ as in \cite{ref:crans}, and not a braided monoidal {\em category}.
  \end{itemize}

Terminology issues aside, what {\em is} $\Dim \mathcal{A}$ and why does it produce a braided monoidal category?

\subsection*{Periodic table argument}
The most abstract answer to why the higher-categorical dimension of a 2-category produces a braided monoidal category is as follows: think of the 2-category $\mathcal{A}$ as an {\em object} in the 3-category $2\mathcal{C}\text{at}$ of 2-categories, weak 2-functors, transformations and modifications. Then $\Dim \mathcal{A} :=  \Hom(\id_\mathcal{A}, \id_\mathcal{A}) \subset 2\mathcal{C}\text{at}$ is a 3-category with one object and one 1-morphism and is hence a {\em braided monoidal category}. To see this, we need to recall the celebrated {\em periodic table} of Baez and Dolan \cite{ref:baez_dolan_hda0}:

  \vskip 0.5em
\begin{center}
{\small
\begin{tabular}{|c|c|c|c|}  \hline
      & $n = 0$   & $n = 1$    & $n = 2$          \\     \hline
$k = 0$  & sets      & categories & 2-categories     \\     \hline
$k = 1$  & monoids   & monoidal   & monoidal         \\
      &           & categories & 2-categories     \\     \hline
$k = 2$  &commutative& braided    & braided          \\
      &  monoids  & monoidal   & monoidal         \\
      &           & categories & 2-categories     \\     \hline
$k = 3$  &`       & symmetric  & weakly involutory \\
      &           & monoidal   & monoidal         \\
      &           & categories & 2-categories     \\     \hline
$k = 4$  &`'      & `'         &strongly involutory\\
      &           &            & monoidal         \\
      &           &            & 2-categories     \\     \hline
$k = 5$  &`'      &`'          & `'               \\
      &           &            &                  \\
      &           &            &                  \\      \hline
\end{tabular}} \vskip 1em
\end{center}
\vskip 0.5em
Recall that this table describes `$k$-tuply monoidal $n$-categories' --- that is, $(n+k)$-categories with only one
$j$-morphism for $j < k$.  The idea is that as one descends each column, the $n$-categories first acquire a `monoidal' or tensor
product structure, which then becomes increasingly `commutative' in character with increasing $k$, and Baez and Dolan conjectured that this process stabilizes at $k = n+2$.

The entry in this table which corresponds to a `3-category with one object and one 1-morphism' is found at $k=2$ and $n=1$ --- a braided monoidal category. For a formal {\em proof} that a tricategory with one object and one 1-morphism corresponds to a braided monoidal category, see \cite{ref:cheng_gurski}.

\subsection*{String diagrams argument} The most intuitive way to understand that the higher categorical dimension of a 2-category is a braided monoidal category is by using string diagrams to depict the relevant higher-dimensional algebra in the 3-category $2\mathcal{C}\text{at}$. Recall from Chapter \ref{stchap} that string diagrams for 2-categories are Poincar\'{e} dual to the ordinary globular notation. We follow the same method here. We depict an object $\mathcal{A} \in 2\mathcal{C}\text{at}$ --- that is, a 2-category --- as a solid cubical region:
 \[
  \ba \ig{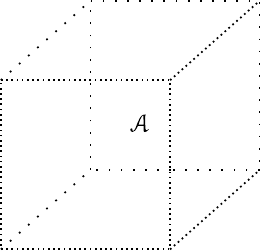} \ea
  \]
A weak 2-functor $F \colon \mathcal{A} \rightarrow \mathcal{B}$ is then drawn as a vertical plane, viewed from right to left as usual:
 \[
  \ba \ig{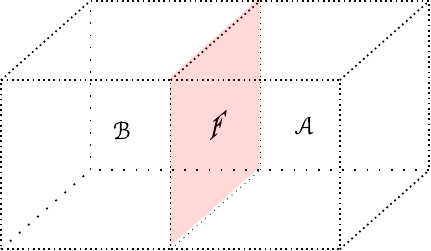} \ea
 \]
A transformation $\sigma \colon F \Rightarrow G$ is drawn as a vertical line:
 \[
  \ba \ig{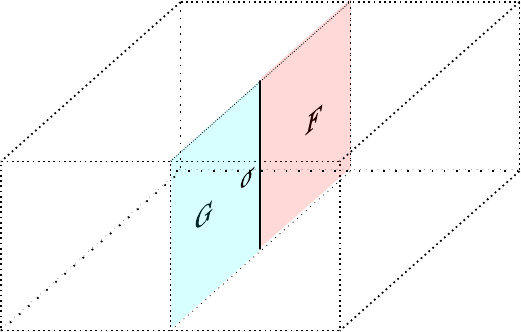} \ea
 \]
Note that we are reading these transformations from the back to the front. Finally, a modification $\theta \colon \sigma \threearrow \rho$ is drawn as a point (at least, a small circle):
 \[
 \ba \ig{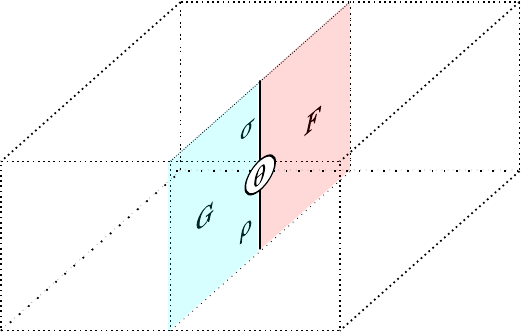} \ea
 \]
Now, we want to emphasize that these diagrams are not just a convenient mnemonic but are a {\em perfectly precise notation} for dealing with weak 2-functors, transformations and modifications between 2-categories. To see this, one should understand that each 3d diagram above can be viewed as a `movie' of transformations, obtained by taking sequential horizontal cross sections of the diagram from top to bottom. Every time a modification is crossed, a new frame is created. For instance, the movie corresponding to the diagram above is:
 \[
  \ba \ig{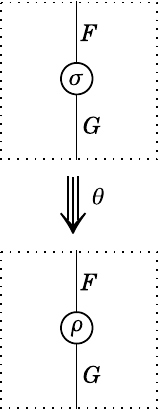} \ea
 \]
In all our movies (such as those in Appendix \ref{movieapp}), the source and target weak 2-functors for the transformations in each individual frame will be clearly understood from the context. We can thus appeal to the coherence theorem of Gurski \cite[Cor 10.2.5]{ref:gurski} to conclude that the resulting modification specified in a movie of transformations is independent of how each individual transformation has been parenthesized.

Let us apply this three-dimensional notation to $\Dim \mathcal{A} = \Hom(\id_{\mathcal{A}}, \id_\mathcal{A})$. A transformation $\T \colon \id_\mathcal{\A} \Rightarrow \id_\mathcal{\A}$ (`$\T$' stands for `transformation') is drawn simply as a vertical line, because the identity planes $\id_\mathcal{\A}$ are not labelled:
 \[
  \ba \ig{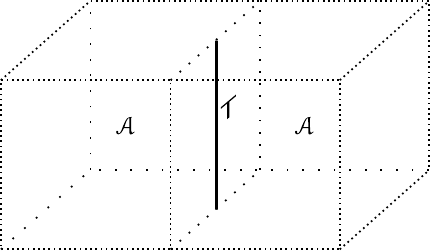} \ea
 \]
Observe that if $\Ss$ is another transformation of the identity, then there is a modification
 \be \label{braidingt}
  c \colon \T \circ \Ss \threearrow \Ss \circ \T
  \ee
whose component at $a \in \mathcal{A}$ is given by the coherence cell of $\T$ evaluated at $\Ss_a$:
\[
  \ba
\xy
(-11.5,7.5)*+{a}="1";
(-11.5,-7.5)*+{b}="2";
(11,7.5)*+{a}="3";
(11,-7.5)*+{b}="4";
{\ar_-{\T_a} "1";"2"};
{\ar^-{\Ss_a} "1";"3"};
{\ar^-{\T_a} "3";"4"};
{\ar_-{\Ss_a} "2";"4"};
{\ar@2{->}^-{\T(\Ss_a)} (0,2)*{};(-4,-2)*{}};
\endxy
 \ea
   \]
We choose to draw this modification in string diagrams as follows:
 \[
  \ba \ig{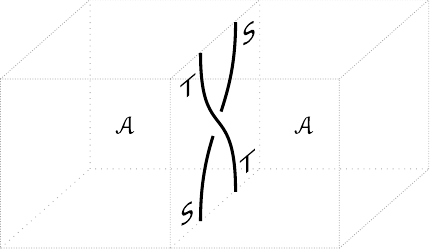} \ea
 \]
Then the coherence and naturality properties of the component 2-isomorphisms $\T(\Ss_a)$ tell us that the strings labelled $\T$, $\Ss$ etc. really do behave as braids in the corresponding string diagrams. Moreover, the coherence equation for the components of a modification $\theta \colon \T \threearrow \T'$ tell us that one can drag the circles along the strings:
\[
 \ba \ig{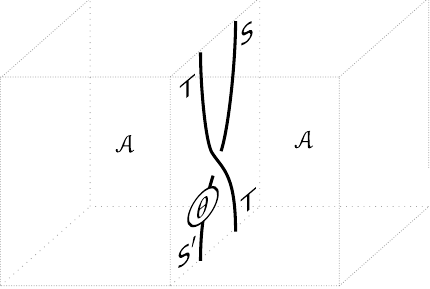} \ea = \ba \ig{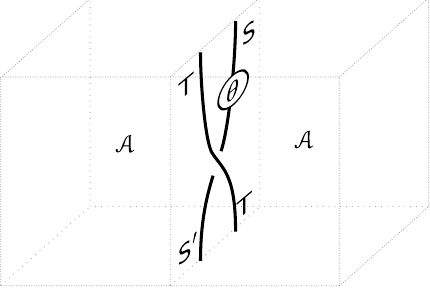} \ea.
\]
We hope this pictures make it clear why $\Dim \mathcal{\A}$ is a braided monoidal category. In fact, the same string diagrams give an elegant explanation for why a {\em general} tricategory with one object and one 1-morphism is a braided monoidal category, that is for Theorem 2.8 of Cheng and Gurski \cite{ref:cheng_gurski}. The braiding is obtained by creating `bubbles' using the unit coherence isomorphisms, which enables the strands to move around each other. However we will not deal with this general case here.

\section{The argument, one level down\label{lb}}
Having explained higher-categorical dimensions of 2-categories, we are now return to our task of calculating the higher-categorical dimension of $\TRep(G)$. In this section, which will be used as a blueprint for the later sections, we work out the corresponding problem one level down. We show that the higher-categorical dimension of the category of representations of a group computes as the center of the group algebra. This is not a new result (see for instance \cite{ref:simon}), but our purpose is to lay out the proof in a series of lemmas in such a way that the rest of the chapter is seen to be a direct categorification of each of these lemmas. Of course in conventional terms the proof is `obvious': a natural transformation of the identity on $\Rep(G)$ is entirely determined by what it does at the regular representation $\mathbb{C}[G]$. But things which appear obvious at this level become a bit more subtle when categorified; this is what forces us to go about the proof in this strange way.

Here then is the scene: we are trying to understand the nature of a natural transformation $\ttt$ (`$\ttt$' stands for `transformation') of the identity functor on $\Rep(G)$. Firstly, let us establish our notation. We write
 \[
 \mathbb{C}[G] = \left\{ \text{formal sums }\sum_{g \in G} \lambda_g \, g \right\}
\]
for the left regular representation of $G$, and we equip it with the inner product
 \[
  \langle g | h \rangle = \delta_{gh}.
 \]
We start with a lemma giving us a grip on the intertwiners in $\Rep(G)$ emanating out of the regular representation.

\begin{lem}\label{myfirst} An intertwiner of representations $u \colon \mathbb{C}[G] \rightarrow V$ is the same thing as a vector $v \in V$.
\end{lem}
\begin{proof}
Set $v = u(e)$. Then since $u$ is equivariant, we must have $u(g) = g \ccdot v$. Conversely, given any $v \in V$ we may clearly define an intertwiner by setting $u(g) = g \ccdot v$.
\end{proof}
Given a vector $v \in V$, we will write $u_v$ for the intertwiner $\mathbb{C}[G] \rightarrow V$ defined in the above lemma, that is $u_v(g) = g \ccdot v$.

Now, we will use the notation that $\Rep(G)|_{\mathbb{C}[G]}$ refers to the category of representations of $G$ `restricted' to $\mathbb{C}[G]$; that is it has one object $\mathbb{C}[G]$ and its morphisms are just intertwiners $\mathbb{C}[G] \rightarrow \mathbb{C}[G]$. Suppose $\ttt$ is a transformation of the identity on $\Rep(G)|_{\mathbb{C}[G]}$. By the above lemma, its solitary component $\ttt_{\mathbb{C}[G]} \colon \mathbb{C}[G] \rightarrow \mathbb{C}[G]$ is entirely determined by its value on the unit element $e \in \mathbb{C}[G]$. We encapsulate this information via its matrix elements:
 \[
  \ttt_g := \langle g | t_{\mathbb{C}[G]} | e \rangle.
  \]
In other words,
 \[
  \ttt_{\mathbb{C}[G]} (e) = \sum_g \ttt_g \; g.
 \]
Then we have the following lemma.
\begin{lem} \label{g6} \begin{enumerate}
 \item The numbers $\{\ttt_g\}_{g \in G}$ are invariant under conjugation, that is $\ttt_{hgh^\mi} = \ttt_g$ for all $g,h \in G$.
 \item All told, the map
 \begin{align*}
  \mathbf{v} \colon \Dim \left( \Rep(G)|_{\mathbb{C}[G]} \right) &\rightarrow Z(\mathbb{C}[G]) \\
  \ttt & \mapsto \sum_{g \in G} \ttt_{g^\mi} \; g
 \end{align*}
is an isomorphism of commutative algebras.
\end{enumerate}
\end{lem}
\begin{proof} (i) The map $\ttt_{\mathbb{C}[G]} \colon \mathbb{C}[G] \rightarrow \mathbb{C}[G]$ must commute with all intertwiners, and these are generated by the maps $u_h$ we defined in Lemma \ref{myfirst} as $h$ ranges over $G$. In terms of matrix elements, we must therefore have
 \begin{align*}
  \langle hg | \ttt_{\mathbb{C}[G]} \circ u_h | e \rangle & = \langle hg | u_h \circ  \ttt_{\mathbb{C}[G]} | e \rangle \\
  \Leftrightarrow \langle hg | \ttt_{\mathbb{C}[G]} | h \rangle & = \langle hgh^\mi |\ttt_{\mathbb{C}[G]} | e \rangle \\
  \Leftrightarrow \langle g | \ttt_{\mathbb{C}[G]} | e \rangle & = \langle hgh^\mi |\ttt_{\mathbb{C}[G]} | e \rangle.
 \end{align*}
That is, $\ttt_g = \ttt_{hgh^\mi}$.

(ii) We saw in (i) that a natural transformation $\ttt \in \Dim \left( \Rep(G)|_{\mathbb{C}[G]} \right)$ gives rise to an element $\sum_g \ttt_g \, g$ in the center of the group algebra. The only subtle point is that to make this into an {\em algebra} homomorphism, we have to use inverses as in the statement of the proposition. Indeed, we have
 \begin{align*}
 \mathbf{v}(\sss) \mathbf{v}(\ttt) &= \left( \sum_{g \in G} \langle g^\mi | \sss | e \rangle \right) \left( \sum_{h \in G} \langle h^\mi | \ttt | e \rangle |h \rangle \right) \\
  &= \sum_{g,h \in G} \langle g^\mi | \sss | e \rangle \langle h^\mi | \ttt | e \rangle |h \rangle |gh \rangle \\
  &= \sum_{g,h \in G} \langle h^\mi g^\mi | \sss | h^\mi \rangle \langle h^\mi | \ttt | e \rangle |gh \rangle \\
  &= \sum_{a, b \in G}  \langle a^\mi | \sss | b \rangle \langle b | \ttt | e \rangle |a \rangle \\
  &= \sum_{a \in G} \langle a^\mi | \sss \circ \ttt | e \rangle |a \rangle \\
  &= \mathbf{v}(\sss \circ \ttt).
 \end{align*}
\end{proof}
Now we establish that a natural transformation of the identity on $\Rep(G)$ is entirely determined by its behaviour at the regular representation. We do this in two steps. Firstly, we start with a transformation $\ttt$ of the identity on $\mathbb{C}[G]$, and then we define a {\em new} transformation $\hat{\ttt}$ of the identity on $\mathbb{C}[G]$, which is defined solely in terms of the data of $\ttt$ at the regular representation $\mathbb{C}[G]$. Then we show that $\hat{\ttt}$ is in fact equal to the original transformation $\ttt$. This is an admittedly convoluted way of looking at things, but it is the way our construction will run in the categorified setting.

Given $\ttt \in \Dim \Rep(G)$, we define the natural transformation $\hat{\ttt} \in \Dim \Rep(G)$ by defining its component at a representation $V$ by summing over the group action, weighted by the numbers $\ttt_g$ obtained by restricting $\ttt$ to the regular representation:
 \be \label{goodf}
  \ttt_V (v) = \sum_g \ttt_g \; \, g \ccdot v.
 \ee
Note that $\hat{\ttt}$ is defined solely in terms of the data of $\ttt$ at $\mathbb{C}[G]$. We have the following lemma.
\begin{lem} \label{mythird} With this definition, $\hat{\ttt}$ is indeed a transformation of the identity on $\Rep(G)$. Moreover, $\hat{\ttt}_{\mathbb{C}[G]} = \ttt_{\mathbb{C}[G]}$.
\end{lem}
\begin{proof}
 We must show that for each intertwiner $f \colon V \rightarrow W$, the following diagram commutes:
 \[
 \ba \xymatrix{ V \ar[r]^-{f} \ar[d]_{\hat{\ttt}_V} & W \ar[d]^{\hat{\ttt}_W} \\ V \ar[r]_-{f} & W} \ea.
\]
This follows from the fact that $f$ is an intertwiner:
 \begin{align*}
  \hat{\ttt}_W \circ f (v) &= \sum_g \ttt_g (g \ccdot f(v)) \\
   &= \sum_g \ttt_g f(g \ccdot v) \\
   &= f \circ \hat{\ttt}_V (v).
  \end{align*}
Moreover, it is clear that $\hat{\ttt}$ equals $\ttt$ at $\mathbb{C}[G]$. To check this, we only need to evaluate them at the unit element $e \in \mathbb{C}[G]$:
 \[
 \hat{\ttt}_{\mathbb{C}[G]} (e) = \sum_g \ttt_g \, g = \sum_g \langle g | \ttt_{\mathbb{C}[G]} | e \rangle \, g = \ttt_{\mathbb{C}[G]} (e).
 \]
\end{proof}
Finally we establish that $\hat{\ttt}$ actually recovers back our original transformation $\ttt$.
\begin{lem}\label{myfourth} In fact, $\ttt = \hat{\ttt}$.
\end{lem}
\begin{proof}
Indeed, this follows from the following naturality square, where $V$ is a representation and $v \in V$:
 \[
  \ba \xymatrix{ \mathbb{C}[G] \ar[r]^-{u_v} \ar[d]_{\ttt_{\mathbb{C}[G]}} & V \ar[d]^{\ttt_V} \\ \mathbb{C}[G] \ar[r]_-{u_v} & V} \ea.
 \]
 \end{proof}
In summary, we have shown the following.
\begin{lem}\label{g8} The restriction map
 \[
  \Dim \left(\Rep(G)\right) \rightarrow \Dim \left(\Rep(G)|_{\mathbb{C}[G]}\right)
 \]
is an isomorphism of commutative algebras.
\end{lem}
\begin{proof} This map is clearly an algebra homomorphism. Lemma \ref{mythird} establishes that it is surjective, and Lemma \ref{myfourth} establishes that it is injective.
\end{proof}
Combining Lemmas \ref{g6} and \ref{g8} then proves our result.
 \begin{prop} The categorical dimension of the category of unitary representations of a finite group $G$ is isomorphic, as a commutative algebra, to the center of the group algebra. That is, $\Dim \Rep(G) \cong Z(\mathbb{C}[G])$.
 \end{prop}
In the next sections we `categorify' each step in the above procedure, translated into the world of equivariant gerbes.

\section{Why the indirect method of equivariant gerbes?\label{why}}
Our strategy in this chapter is to calculate $\Dim \TRep (G)$ indirectly; that is we actually calculate $\Dim \Gerbes(G)$ and then we appeal to Theorem \ref{toodytheorem} which showed that $\TRep(G)$ and $\Gerbes(G)$ are equivalent 2-categories. We acknowledge that it would be more desirable to calculate $\Dim \TRep(G)$ directly; in this section we explain why we have been obliged to adopt the strategy we have.

The trouble is purely a technical one: it is difficult to write down choice-free formulas when working directly with 2-representations, and this causes problems in the proof. The difference between 2-representations and equivariant gerbes is analogous to the difference between an exact sequence of groups and a specific {\em 2-cocycle} which realizes that exact sequence, as we saw in Chapter \ref{exsec}.

Moreover, the language of categories and functors does not lend itself to writing down canonical choice-free formulas. Consider for instance attempting to categorify the formula \eqref{goodf} above. That is, suppose one is given a conjugation equivariant vector bundle $V$ over $G$ (one regards $V$ as a categorification of an element in the centre of the group algebra)\footnote{At this point, we do not want to think of a conjugation equivariant vector bundle over $G$ as the categorification of a {\em class function}, because as we saw in Lemma \ref{g6} the relevant algebra structure is that of $Z(\mathbb{C}[G])$ --- the {\em convolution product} --- and not the product of class functions.}, with which one wishes to define a transformation of the identity $\T$ on $\TRep(G)$. The natural formula to use is to set the component of $\T$ at a 2-representation $\alpha \in \TRep(G)$ to be the morphism of 2-Hilbert spaces $\T_\alpha \colon H_\alpha \rightarrow H_\alpha$ which acts on an object $v \in H_\alpha$ via
 \be \label{attempt}
  \T_\alpha (v) = \text{``} \, \bigoplus_{g \in G} V_g \odot \alpha_g (v) \, \text{''} .
 \ee
Besides the obvious fact that we have used the `tensoring with a vector space' operation $\odot$ which is not actually available in a 2-Hilbert space, at least in the way we have defined them, the problem with this formula is that there is no canonical way to `add' two objects in a 2-Hilbert space. That may seem like a trivial matter, and indeed it {\em is} merely a technicality, but it interferes with the proof. One may attempt to remedy this by passing to {\em skeletal} 2-Hilbert spaces; but even there one cannot proceed in an entirely choice-free fashion because at a certain point one will need to define linear $*$-functors between such spaces, and as we saw in Lemma \ref{charstlem} the value of such a functor on the hom-sets can be specified arbitrarily.

Another possible strategy one could adopt is to try and use an explicit {\em cocycle} description of 2-representations, similar to the treatment of Elgueta \cite{ref:elgueta}. The problem with cocycles is they too implicitly represent {\em choices}, and it proves very awkward to carry through with the calculation of $\Dim \TRep(G)$ in this manner, where one needs to accommodate the {\em twisted} structures in a canonical way.

This then is the advantage of equivariant gerbes: they constitute a {\em choice-free} geometric formalism for dealing with unitary 2-representations in which one can actually write down {\em canonical formulas}. For instance, we will see in Section \ref{extendd} that in the world of equivariant gerbes, \eqref{attempt} becomes the canonical formula
 \[
 \langle x' | \hat{\T}_\X | x \rangle = \bigoplus_\blag{x}{x'}{g} (\X_\mlag{x}{}{g})_\mathbb{C} \otimes  \T_g.
 \]
That is, one adds the vector spaces $V_g$ together after {\em twisting} them by the gerbe. We stress that this is a choice-free formula; even the direct sum $\bigoplus$ is not a {\em formal} direct sum of vector spaces (which would require an ordering) but is really defined as the `sections of the vector bundle over the relevant hom-sets in the gerbe twisted by the vector spaces $T_g$', just as in Chapter \ref{AGeomSec}.

The price we must pay is that we are no longer dealing directly with 2-representations. In Chapter \eqref{eqref} we established that there is a canonical 2-functor
 \[
  \TRep_\text{m} (G) \rightarrow \Gerbes(G)
 \]
from the 2-category of marked 2-representations to the 2-category of equivariant gerbes, which sets up an equivalence between these 2-categories. Notwithstanding the relatively harmless problem that there is no choice-free way to turn a {\em un}marked 2-representation into a marked one, the main problem is that there is no choice-free way to turn an equivariant gerbe into a 2-representation. To do that, one would first have to choose a {\em section} of the gerbe, which sets up a cocycle, which one then uses to define the coherence isomorphisms for the 2-representation. Thus, all we can truthfully claim is that the 2-category of 2-representations and the 2-category of equivariant gerbes are `merely' equivalent. To prove that this implies that their higher-categorical dimensions are equivalent, as we do in Appendix \ref{movieapp}, involves the assumption that this equivalence can be upgraded into a `coherent adjoint equivalence'. Nevertheless it is expected that this is indeed the case \cite{ref:gurski_private}.

\section{Understanding morphisms emanating out of $\E G$\label{lc}}
Having explained the necessary background and motivation, we can now commence with the task of computing $\Dim \Gerbes(G)$.  We begin in this section by finding the appropriate categorified analogue of Lemma \ref{myfirst}.

The replacement for the notion of the `regular 2-representation' in the world of equivariant gerbes is called $\E G$, the trivial equivariant gerbe over the free $G$-set $EG$. That is, the objects of $\E G$ are the elements $g \in G$, and its $G$-graded hom-sets are copies of $U(1)$,
 \[
  \E G_\mlag{g}{}{h} = U(1).
 \]
Suppose $\X \in \Gerbes(G)$ is an equivariant gerbe. For each $x \in \X$, we want to define a morphism
 \[
  U_x \colon \E G \rightarrow \X
 \]
which is the vector bundle analogue of the linear map $u_v \colon \mathbb{C}[G] \rightarrow \rho$ from Lemma \ref{myfirst}. Morally speaking, $U_x$ is the `map which sends $e \in G$ to $x \in \X$', but we need to express this in our vector bundle formalism. We define $U_x$ to be the equivariant vector bundle over $\X \otimes \overline{\E G}$ which is supported exclusively on the orbit
 \[
  \mathcal{O}_x = G \cdot (x, e) \subset \text{Ob} \, \X \otimes \overline{\E G}
 \]
 \begin{figure}[t]
\centering
\[
\begin{array}{ccc} \ba \ig{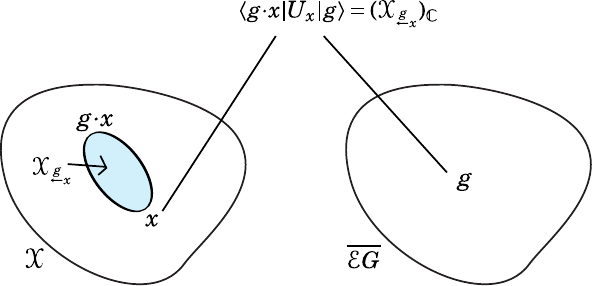} \ea &  \qquad \qquad & \ba \ig{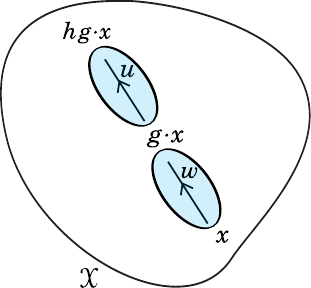} \ea \\ \text{(a)} & & \text{(b)}
\end{array}
\]
\caption{\label{morphx} (a) The equivariant vector bundle $U_x$ over $\X \otimes \overline{\E G}$ associated to a point $x \in \X$. (b) The equivariant maps for the bundle $U_x$ are simply given by composing arrows in the gerbe.}
\end{figure}
where its fibers are given by
 \[
  \langle g \ccdot x | U_x | g \rangle = (\X_\mlag{x}{}{g})_\mathbb{C},
 \]
as in Figure \ref{morphx}a, equipped with the equivariant maps
 \[
  U_x( \lag{h}{u \otimes \overline{1}} ) \colon \langle g \ccdot x | U_x | g \rangle \rightarrow \langle hg \ccdot x | U_x | hg \rangle
 \]
which are the linear extensions of the maps defined by `postcomposing with $u$':
 \begin{align*}
  \X_\mlag{x}{}{g} & \rightarrow \X_\mlag{x}{}{hg} \\
  w & \mapsto \post(u)(w) = u \circ w.
 \end{align*}
See Figure \ref{morphx}b. Then we have the following lemma.
\begin{lem} \label{aprelim}\begin{enumerate}  \item The morphisms $U_x \colon \E G \rightarrow \X$ defined above generate the hom-category $\Hom(\E G, \X)$, in the sense that every morphism is isomorphic to a direct sum of such morphisms.
 \item In particular, if $V \colon \X \rightarrow \Y$ is a morphism of equivariant gerbes, then there is a canonical isomorphism of equivariant vector bundles
  \[
  V \circ U_x \; \stackrel{\xi}{\longleftarrow} \; \bigoplus_{y' \in \Y} \langle y' | V | x \rangle \cdot U_{y'}.
 \]
 \end{enumerate}
\end{lem}
\begin{proof} (i) Since $G$ acts {\em freely} on $X \times EG$, every groupoid 2-cocycle on $(X \times EG)_G$ must be a coboundary, and thus no matter how twisted the gerbe $\X$ is, an equivariant vector bundle over $\X \otimes \overline{\E G}$ must be isomorphic to a direct sum of trivial line bundles over each orbit, which we can take to be the bundles $U_x$ we defined above.

(ii) We should first explain precisely what we mean by the vector bundle $\mathop{\oplus}_{y'} \langle y' | V | x \rangle \cdot U_{y'}$ --- it is the equivariant vector bundle over $\Y \otimes \overline{\E G}$ with fibers
  \begin{align*}
  \left \langle y \left| \mathop{\oplus}_{y'}  \langle y' | V | x \rangle \cdot U_{y'} \right| g \right\rangle &:= \bigoplus_{y'} \langle y' | V | x \rangle \otimes \langle y | U_{y'} | g \rangle \\
  &= \langle g^\mi \ccdot y | V | x \rangle \otimes (\Y_\mlag{}{y}{g})_\mathbb{C}
 \end{align*}
and with maps
 \[
  \left[\mathop{\oplus}_{y'}  \langle y' | V | x \rangle \cdot U_{y'} \right] (\lag{h}{w \otimes \bar{1}}) \colon \left \langle y \left| \mathop{\oplus}_{y'}  \langle y' | V | x \rangle \cdot U_{y'} \right| g \right\rangle \rightarrow \left \langle h \ccdot y \left| \mathop{\oplus}_{y'}  \langle y' | V | x \rangle \cdot U_{y'} \right| hg \right\rangle
 \]
defined by postcomposing with $w$:
 \begin{align*}
  \langle g^\mi \ccdot y | V | x \rangle \otimes (\Y_\mlag{}{y}{g})_\mathbb{C} & \rightarrow \langle g^\mi \ccdot y | V | x \rangle \otimes (\Y_\mmlag{}{h \ccdot y}{hg})_\mathbb{C} \\
  v \mapsto w \circ v.
 \end{align*}
We define the isomorphism $\xi$ fibrewise as follows:
 \[
   \xymatrix{  \displaystyle \left \langle y \left| \mathop{\oplus}_{y'}  \langle y' | V | x \rangle \cdot U_{y'} \right| g \right\rangle \ar@{=}[d]  \ar[rrr]^{\xi_{y, g}} &&& \langle y | V \circ U_x | g \rangle \ar@{=}[d]\\
   \langle g^\mi \ccdot y | V | x \rangle \otimes (\Y_\mlag{}{y}{g})_\mathbb{C} \ar[rrr]_-{W(\lag{g}{v \otimes \bar{u}}) \otimes (v \mapsto u)} &&& \langle y | V | g \ccdot x \rangle \otimes (\X_\mlag{x}{}{g})_\mathbb{C}   }
 \]
which as usual is independent of the choices of $v \in \Y_\mlag{}{y}{g}$ and $u \in \X_\mlag{x}{}{g}$. It is a simple matter of unravelling the definitions to see that $\xi$ is an isomorphism of equivariant vector bundles.
 \end{proof}

\section{Extracting equivariant vector bundles\label{ld}}
Our aim in this section is to establish the categorified analogue of Lemma \ref{g6}. We show that a transformation $\T$ of the identity on the 2-category $\Gerbes(G)_{\E G}$ is the same thing as a {\em vector bundle over $G$ equipped with isomorphisms lifting the action of conjugation} (another way to say this is that $\T$ is the same thing as a representation of the loop groupoid $\Lambda \BG$). Under this correspondence, composition of transformations gets sent to the {\em fusion product} of equivariant vector bundles over $G$ (not the {\em fiberwise} tensor product), and the same goes for the braiding. When we think of the category of conjugation equivariant vector bundles over a group as equipped with this braided monoidal structure, we write it as $\Hilb_G ^\text{fusion}(G)$, and we have recalled the definition in Appendix \ref{AppFreed}.

Consider that, by definition, a transformation of the identity $\T$ on $\Gerbes(G)_{\E G}$ consists of a morphism $\T_{\E G} \colon \E G \rightarrow \E G$, together with for each morphism $U \colon \E G \rightarrow \E G$, an isomorphism of equivariant vector bundles
\[
\ba
\xy
(-7.5,7.5)*+{\E G}="1";
(-7.5,-7.5)*+{\E G}="2";
(8,7.5)*+{\E G}="3";
(8,-7.5)*+{\E G}="4";
{\ar_-{\T_{\E G}} "1";"2"};
{\ar^-{U} "1";"3"};
{\ar^-{\T_{\E G}} "3";"4"};
{\ar_-{U} "2";"4"};
{\ar@2{->}^-{\T(U)} (1,2)*{};(-3,-2)*{}};
\endxy
\ea .
\]
Let us agree to write
 \[
  \T_g := \langle g | \T_{\E G} | e \rangle
 \]
for the fibers of the bundle $\T_{\E G}$ which `emanate out of the identity'. Now by Lemma \ref{aprelim}, $U$ is a direct sum of the generating line bundles $U_h$ as $h$ ranges over $G$, and since these isomorphisms $\T(U)$ are natural in $U$ it is sufficient to understand the coherence isomorphisms $\T(U_h)$ on these generators. At $(hg, e) \in \E G \otimes \overline{\E G}$, the fibers of $\T_{\E G} \circ U_h$ compute as
\begin{align*}
  \langle hg | \T_{\E G} \circ U_h | e \rangle &= \bigoplus_{a \in G} \langle hg | \T | a \rangle \otimes \langle a | U_h | e \rangle \\
   &= \langle hg | \T_{\E G} | h \rangle \otimes (\E G_\mlag{h}{h}{e})_\mathbb{C} \\
   &= \langle hg | \T_{\E G} | h \rangle.
 \end{align*}
In this last step, we have applied the following harmless convention, which we will use for the rest of this chapter:
\[
\text{\framebox{\parbox[b]{11cm}{ If a factor in a tensor product of vector spaces identifies canonically with $\mathbb{C}$, then it will not explicitly be written down.}}}\]
The corresponding fibers of $U_h \circ \T_{\E G}$ compute as
 \begin{align*}
 \langle hg | U_h \circ \T_{\E G} | e \rangle &= \bigoplus_{a \in G} \langle hg | U_h | a \rangle \otimes \langle a | \T_{\E G} | e \rangle \\
  &= \T_{hgh^\mi}.
 \end{align*}
We can therefore conveniently package the data of the coherence isomorphisms $\T(U_h)$ by defining the vector space isomorphisms
 \[
  \T(\sla{h} g) \colon \T_g \rightarrow \T_{hgh^\mmi}
 \]
as the composite
 \[
 \T_g \xrightarrow{=} \langle g | \T_{\E G} | e \rangle \xrightarrow{\T(\lag{h}{1 \otimes \bar{1}})} \langle hg | \T_{\E G} | h \rangle \xrightarrow{\T(U_h)_{hg, e}} \T_{hgh^\mmi}.
 \]
We then have the following result.

\begin{prop}\label{fyn} \begin{enumerate}\item The vector spaces $\{\T_g\}_{g \in G}$ together with the isomorphisms $\{ \T(\sla{\,h}  \! g) \colon \T_g \rightarrow \T_{hgh^\mmi}\}_{g,h \in G}$ furnish an equivariant vector bundle over $G$.
 \item Moreover, a modification $\theta \colon \T \rightarrow \Ss$ of transformations gives rise to a map of equivariant vector bundles
  \[
   \{ \theta_g \colon \T_g \rightarrow \Ss_g \}
  \]
 over $G$.
\item All told, the functor
 \begin{align*}
  \mathbf{V} \colon \Dim \left(\Gerbes(G)|_{\E G}\right) & \rightarrow \Hilb_G^\fusion(G) \\
   \mathbf{V}(\T)_g := \T_{g^\mmi}
 \end{align*}
is an equivalence of braided monoidal categories.
\end{enumerate}
\end{prop}
\begin{proof} (i) The reader will agree that, as long as we keep using our convention of ignoring factors in tensor products which identify canonically with $\mathbb{C}$, the composite vector bundle $U_{h_1} \circ U_{h_2}$ is equal to $U_{h_2 h_1}$. This allows us to use the accompanying coherence equation for the isomorphisms $\T(U_h)$:
 \[
 \ba
\xy
(-11.5,7.5)*+{\E G}="1";
(-11.5,-7.5)*+{\E G}="2";
(11,7.5)*+{\E G}="3";
(11,-7.5)*+{\E G}="4";
(33.5,7.5)*+{\E G}="5";
(33.5,-7.5)*+{\E G}="6";
{\ar_-{\T_{\E G}} "1";"2"};
{\ar^-{U_{h_2}} "1";"3"};
{\ar^-{\T_{\E G}} "3";"4"};
{\ar_-{U_{h_2}} "2";"4"};
{\ar^-{U_{h_1}} "3";"5"};
{\ar_-{U_{h_1}} "4";"6"};
{\ar^-{\T_{\E G}} "5";"6"};
{\ar@2{->}^-{\T(U_{h_2})} (-3,2)*{};(-7,-2)*{}};
{\ar@2{->}^-{\T(U_{h_1})} (23.5,2)*{};(19.5,-2)*{}};
\endxy
 \ea
 =
 \ba
\xy
(-7.5,7.5)*+{\E G}="1";
(-7.5,-7.5)*+{\E G}="2";
(12,7.5)*+{\E G}="3";
(12,-7.5)*+{\E G}="4";
{\ar_-{\T_{\E G}} "1";"2"};
{\ar^-{U_{h_2 h_{1}}} "1";"3"};
{\ar^-{\T_{\E G}} "3";"4"};
{\ar_-{U_{h_2 h_{1}}} "2";"4"};
{\ar@2{->}^-{\T(U_{h_{2}h_1})} (1,2)*{};(-3,-2)*{}};
\endxy
\ea
   \]
It is not hard to see that by expanding out this coherence equation fibrewise at $(h_2h_1g, e) \in \E G \otimes \overline{\E G}$ we obtain the equations
 \[
   \T(\sla{h_2} h_1 g h_1^\mi) \circ \T(\sla{h_1} g ) = \T( \sla{h_2 h_1} g). \\
 \]
Moreover since $U_e$ can be identified with the identify morphism $\E G \rightarrow \E G$, we have
 \[
  \T( \sla{e} g) = \id.
 \]
Thus the vector spaces $\{\T_g\}_{g \in G}$ together with the isomorphisms $\{\T(\sla{h} g)\}_{g,h \in G}$ furnish an equivariant vector bundle over $G$.

(ii) If $\theta \colon \T \rightarrow \Ss$ is a modification between transformations, we can define the fibrewise maps
  \[
   \theta_g := [\theta_{\E G}]_{g,e} \colon \T_g \rightarrow \Ss_g.
  \]
The condition that $\theta$ is a modification is the requirement that for all $h \in G$ the following pasting diagram commutes:
 \[
\xy
(-7.5,10)*+{\E G}="1";
(-11, 7)="1f"; (-5, 9)="1b";
(-7.5,-10)*+{\E G}="2";
(-11, -10)="2f"; (-5, -7)="2b";
(35,10)*+{\E G}="3";
(33.5, 7)="3f"; (38.5, 9)="3b";
(35,-10)*+{\E G}="4";
(33.5, -10)="4f"; (38.5, -7)="4b";
{\ar@/_0.8pc/_-{\T_{\E G}} "1f";"2f"};
{\ar@/^0.5pc/^-{\Ss_{\E G}} "1b";"2b"};
{\ar@/_0.8pc/_-{\T_{\E G}} "3f";"4f"};
{\ar@/^0.5pc/^-{\Ss_{\E G}} "3b";"4b"};
{\ar^-{U_h} "1";"3"};
{\ar_-{U_h} "2";"4"};
{\ar@2{->}@/{^0.5pc}/^>>>{{\Ss(U_h)}} (14.5,6)*+{};(12.5,0)*+{}};
{\ar@2{->}@/{_0.5pc}/_<<<<<{\T(U_h)} (11,2)*+{};(11,-4)*+{}};
{\ar@2{->}_-{\theta_{\E G}} (-13,-3)*+{};(-4,2)*+{}};
{\ar@2{->}_-{\theta_{\E G}} (31,-3)*+{};(38.5,2)*+{}};
\endxy
\]
The fibrewise components of this pasting diagram at $(hg, e)$ is precisely the requirement that $\{\theta_g \colon \T_g \rightarrow \Ss_g\}$ is a map of equivariant vector bundles:
   \[
   \xymatrix{\T_g \ar[r]^{\theta_g} \ar[d]_{\T(\sla{h}g)} & \Ss_g \ar[d]^{\Ss(\sla{h}g)} \\ \T_{hgh^\mmi} \ar[r]_{\theta_{hgh^\mmi}} & \Ss_{hgh^\mmi}}
  \]

(iii) Since an equivariant vector bundle over $\E G \otimes \overline{\E G}$ must be isomorphic to a constant vector bundle, we see that the only real information in a transformation $\T$ resides in the vector spaces $\T_g := \langle g | T | e \rangle$ together with the isomorphisms $\T(\sla{h}g)$ defined above, and the same goes for modifications between them. This establishes that the functor
 \[
  \mathbf{V} \colon \Dim \Gerbes(G)_{\E G} \rightarrow \Hilb_G^\fusion(G)
 \]
is essentially surjective and fully faithful, and hence an equivalence of categories.

The subtle point here is that we must use {\em inverses}, just as in Lemma \ref{g6}, in order for this to be a {\em monoidal} functor. For then we can define the fibrewise monoidal coherence isomorphisms as:
 \[
 \ba \xymatrix{ [\mathbf{V}(\Ss) \otimes \mathbf{V}(\T)]_g \ar[d] \\ [\mathbf{V}(\Ss \circ \T)]_g} \ea = \ba \xymatrix{ \bigoplus_{ab = g} \langle a^\mmi | \Ss | e \rangle \otimes \langle b^\mi | \T | e \rangle \ar[d]^{\Ss(g^\mmi a) \otimes \id} \\
  \bigoplus_{ab = g} \langle g^\mmi | \Ss | g^\mmi a \rangle \otimes \langle b^\mmi | \T | e \rangle \ar@{=}[d] \\ \bigoplus_{a} \langle g^\mmi | \Ss | a \rangle \otimes \langle a | \T | e \rangle} \ea
 \]
It is trivial to check that this data satisfies the coherence equations for a monoidal functor.

To see that it is a {\em braided} monoidal functor, recall from Section \ref{la} that the braiding $\Ss \circ \T \rightarrow \T \circ \Ss$ in $\Dim \left( \Gerbes(G)|_{\E G}\right)$ is given by the coherence isomorphism $\Ss(\T_{\E G})$,
\[
\ba
\xy
(-7.5,7.5)*+{\E G}="1";
(-7.5,-7.5)*+{\E G}="2";
(12,7.5)*+{\E G}="3";
(12,-7.5)*+{\E G}="4";
{\ar_-{\Ss_{\E G}} "1";"2"};
{\ar^-{\T_{\E G}} "1";"3"};
{\ar^-{\Ss_{\E G}} "3";"4"};
{\ar_-{\T_{\E G}} "2";"4"};
{\ar@2{->}^-{\Ss(\T_{\E G})} (1,2)*{};(-3,-2)*{}};
\endxy
\ea ,
\]
and that the braiding $\mathbf{V}(\Ss) \otimes \mathbf{V}(\T) \rightarrow \mathbf{V}(\T) \otimes \mathbf{V}(\Ss)$ in $\Hilb_G^\fusion(G)$ is given by {\em conjugation on the $\Ss$-factor} followed by rearrangement (see Appendix \ref{AppFreed}):
   \[
  \ba \xymatrix{ [\mathbf{V}(\Ss) \otimes \mathbf{V}(\T)]_g \ar[d] \\ [\mathbf{V}(\T) \otimes \mathbf{V}(\Ss)]_g
   }\ea = \ba \xymatrix{ \displaystyle \bigoplus_{ab = g} \Ss_{a^\mmi} \otimes \T_{b^\mmi} \ar[d]^{\Ss(\sla{b^\mmi}a^\mmi) \otimes \id} \\ \displaystyle \bigoplus_{ab = g} \Ss_{b^\mmi a^\mmi b} \otimes \T_{b^\mmi} \ar[d]^{\text{swap}} \\ \displaystyle \bigoplus_{ab = g} \T_{b^\mmi} \otimes \Ss_{b^\mmi a^\mmi b}} \ea
 \]
Now, $\T_{\E G}$ can be expressed as a direct sum of the generating line bundles \XXX{Big problem here! b?}
 \[
  \T_{\E G} \cong \bigoplus_{b \in g} \T_{b} \cdot U_b
 \]
just as in Lemma \ref{aprelim}. Thus by naturality, the coherence isomorphism $\Ss(\T_{\E G})$ can be expressed in terms of the basic isomorphisms $\Ss(U_b)$, whose fiberwise components are encoded by the conjugation isomorphisms $\Ss(\sla{b}a)$. In other words, the braiding isomorphism
 \[
  \Ss(\T_{\E G}) \colon \Ss \circ \T \rightarrow \T \circ \Ss
 \]
is {\em also} given by conjugation on the $\Ss$-factor and then rearrangment.  In this way we see that our functor recovers precisely the braiding on $\Hilb_G^\fusion(G)$.
 \end{proof}\XXX{Write down somewhere the fusion and braiding operations on HilbG(G). And change Repc(G) to Hilb.}

\section{Extending transformations\label{extendd}}
Our task in the next two sections is to show that a transformation $\T$ of the identity on $\Gerbes(G)$ is determined by its behaviour at $\E G$, the substitute for the `regular 2-representation' in the world of equivariant gerbes. In this section, the analogue of Lemma \ref{mythird}, we construct a transformation of the identity $\hat{\T}$ which {\em is} defined solely in terms of the data of $\T$ at $\E G$, and then in the next section we show that there is a canonical isomorphism
 \[
  \Omega^\T \colon \T \rightarrow \hat{\T}
 \]
inside $\Dim \Gerbes(G)$ --- which is the analogue of Lemma \ref{myfourth}.

\begin{figure}[t]
\centering
\ig{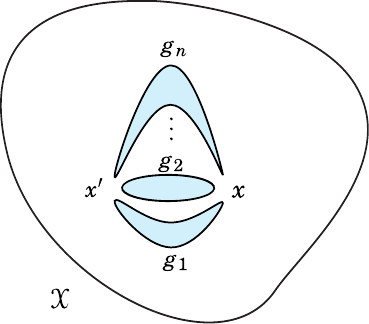}
\caption{\label{morphy} The fiber $\langle x' | \hat{\T}_\X | x \rangle$ of the equivariant vector bundle $\langle \hat{\T}_\X \rangle$ adds up all the vector spaces $\T_g$ where $g \ccdot x=x'$ together, after tensoring them with the corresponding hom-sets in the gerbe $\X$.}
\end{figure}

Given $\T$, we define $\hat{\T}$ as follows. Its component $\hat{\T}_\X \colon \X \rightarrow \X$ at an equivariant gerbe $\X$ is the equivariant vector bundle over $\X \otimes \overline{\X}$ whose fibers are
 \[
  \langle x' | \hat{\T}_\X | x \rangle = \bigoplus_\blag{x}{x'}{g} (\X_\mlag{x}{}{g})_\mathbb{C} \otimes  \T_g.
 \]
In other words, it collects each vector space $\T_g$ having the property that $g \ccdot x = x'$ together, after {\em twisting} them by the equivariant gerbe $\X$ (see Figure \ref{morphy}).
Being able to write down this formula so canonically is a great advantage of using the language of equivariant gerbes; had we opted for a cocycle description instead, or some other such scheme, we would not have been able to write down such a choice-free formula.

\begin{figure}[t]
\centering
\ig{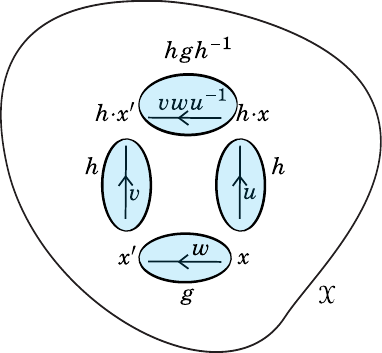}
\caption{\label{morphz} The adjoint operation $\ad_{v,u} \colon \X_\mlag{x}{x'}{g} \rightarrow \X_\mmlag{h \ccdot x}{h \ccdot x'}{hgh^\mmi}$ sends $w \mapsto vwu^\mi$.}
\end{figure}

As an equivariant vector bundle, the maps
 \[
\langle x' | \hat{\T}_\X | x \rangle (\lag{h}{v \otimes \bar{u}}) \colon  \langle x' | \hat{\T}_\X | x \rangle \rightarrow  \langle h \ccdot x' | \hat{\T}_\X | h \ccdot x \rangle
 \]
are given by
 \[
  \ba \xymatrix{\displaystyle \bigoplus_\blag{x}{x'}{g} (\X_\mlag{x}{}{g})_\mathbb{C} \otimes  \T_g \ar[d]^<<<<{\ad_{v,u} \otimes \T(\sla{h} g)} \\ \displaystyle \bigoplus_{\blag{x}{x'}{g}} (\X_\mmlag{h \ccdot x}{h \ccdot x'}{hgh^\mmi})_\mathbb{C} \otimes  \T_{hgh^\mi}} \ea
 \]
where
 \[
 \ad_{v, u} \colon  \X_\mlag{x}{x'}{g} \rightarrow \X_\mmlag{h \ccdot x}{h \ccdot x'}{hgh^\mmi}
 \]
is the `adjoint' operation which sends $w \mapsto vwu^\mi$ (see Figure \ref{morphz}), and the maps $\T(\sla{h} g) \colon \T_g \rightarrow \T_{hgh^\mmi}$ are the conjugation isomorphisms that we defined in Section \ref{ld}. Since both these operations compose appropriately, we see that $\hat{\T}_\X$ is indeed an equivariant vector bundle. The important point here though is that $\hat{\T}_\X$ is constructed {\em solely from the equivariant gerbe $\X$ and the data of $\T$ at $\E G$}.

\begin{figure}
\centering
\ig{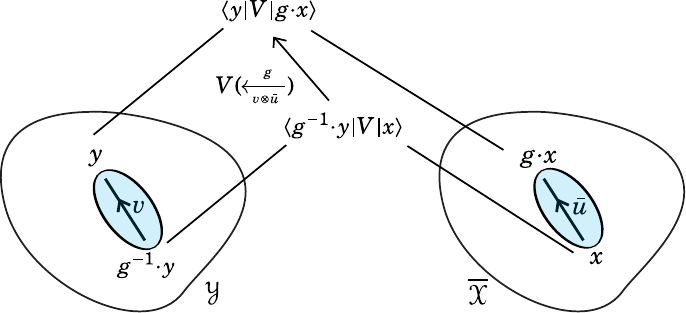}
\caption{\label{morphzz}The coherence isomorphisms $\hat{\T} (V)$ are defined component-wise by transporting $V$ along some chosen arrows arrows  $v \in \Y_\mlag{}{g}{g^\mi \ccdot y}$ and $\bar{u} \in \overline{\X}_\mlag{g}{}{x}$ while simultaneously sending $v \mapsto u$.}
\end{figure}

So far we have only defined the components $\hat{\T}_\X$ of the transformation $\hat{\T}$. Given a morphism $V \colon \X \rightarrow \Y$ of equivariant gerbes, we must also define the coherence isomorphisms
  \[
 \ba \xy
(-7.5,7.5)*+{\X}="1";
(-7.5,-7.5)*+{\X}="2";
(12,7.5)*+{\Y}="3";
(12,-7.5)*+{\Y}="4";
{\ar_-{\hat{\T}_\X} "1";"2"};
{\ar^-{V} "1";"3"};
{\ar^-{\hat{\T}_\Y} "3";"4"};
{\ar_-{V} "2";"4"};
{\ar@2{->}^-{\hat{\T}(V)} (1,2)*{};(-3,-2)*{}};
\endxy \ea.
\]
These are defined in components as follows:
 \[
 \ba \xymatrix{ \langle y | \hat{\T}_\Y \circ V | x \rangle \ar[d]^{\hat{\T}(V)_{y, x}} \\ \langle y | V \circ \hat{\T}_\X | x \rangle} \ea
  =
  \ba \xymatrix{ \displaystyle \bigoplus_{g \in G} (\Y_\mlag{}{y}{g})_\mathbb{C} \otimes \T_g \otimes \langle g^\mi \ccdot y | V | x \rangle \ar[d]^{\text{swap}} \\ \displaystyle \bigoplus_{g \in G}  \langle g^\mi \ccdot y | V | x \rangle \otimes (\Y_\mlag{y}{}{g})_\mathbb{C} \otimes \T_g \ar[d]^{V(\lag{g}{v \otimes \bar{u}}) \otimes (v \mapsto u) \otimes \id} \\ \displaystyle \bigoplus_{g \in G}  \langle y | V | g \ccdot x \rangle \otimes (\X_\mlag{x}{}{g})_\mathbb{C} \otimes \T_g} \ea
\]
That is to say, after initially rearranging things one chooses arrows $v \in \Y_\mlag{}{g}{g^\mi \ccdot y}$ and $u \in \X_\mlag{g}{}{x}$, which one uses to  parallel transport the vector bundle $V$ while simultaneously sending $v \mapsto u$ (see Figure \ref{morphzz}). This construction is canonical and does not depend on the choices of $v$ and $u$, because these two operations transform oppositely with respect to rescaling of $v$ and $u$. Observe also that the isomorphism $\hat{\T}(V)$ only uses the gerbe $\X$ and the vector bundle $V$; it does not use any information from the transformation $\T$ at all.

In the remainder of this section we will do two things. Firstly we will confirm that with this definition of the coherence isomorphisms $\hat{\T}(V)$, $\hat{\T}$ is indeed a transformation of the identity on $\Gerbes(G)$. Then we will use the coherence isomorphisms of $\T$ to establish the isomorphism $\Omega^T \colon \T \rightarrow \hat{\T}$.

\begin{lem}\label{Impy} With these definitions, $\hat{\T}$ is indeed a transformation of the identity on $\Gerbes(G)$. Moreover, its restriction $\hat{\T}|_{\E G}$ is isomorphic to $\T|_{\E G}$ as transformations in $\Dim \left(\Gerbes(G)|_{\E G}\right)$.
\end{lem}
\begin{proof} We first need to establish for a morphism $V \colon \X \rightarrow \Y$ of equivariant gerbes, the isomorphism of vector bundles
 \[
  \hat{\T}(V) \colon \hat{\T}_\Y \circ V \rightarrow V \circ \hat{\T}_\X
 \]
defined above is indeed a morphism of {\em equivariant} vector bundles --- that is, that the following diagram commutes:
 \[
  \xymatrix{ \langle y | \hat{\T}_\Y \circ V | x \rangle \ar[d]_{\hat{\T}(V)_{y, x}}  \ar[rrr]^{(\hat{\T}_\Y \circ V)(\lag{h}{t \otimes \bar{s}})} &&& \langle h \ccdot y | \hat{\T}_\Y \circ V | h \ccdot x \rangle \ar[d]^{\hat{\T}(V)_{h \ccdot y, h \ccdot x}} \\
  \langle y | V \circ \hat{\T}_\X | x \rangle \ar[rrr]_{(V\circ \hat{\T}_\X)(\lag{h}{t \otimes \bar{s}})} &&& \langle h \ccdot y | V \circ \hat{\T}_\X | h \ccdot x \rangle}
 \]
Using the definitions we gave above, we can expand out this diagram to obtain the following diagram,
\[
 \xymatrix{ \displaystyle \bigoplus_{g \in G} (\Y_\mlag{}{y}{g})_\mathbb{C} \otimes \T_g \otimes \langle g^\mmi \ccdot y | V | x \rangle \ar[d]_\zeta
  \ar[rrr]^-{\alpha} &&& \displaystyle \bigoplus_{g \in G} (\Y_\mmlag{}{h \ccdot y}{hgh^\mmi})_\mathbb{C} \otimes \T_{hgh^\mmi} \otimes \langle hg^\mmi \ccdot y | V | h \ccdot x \rangle \ar[d]^\xi \\
 \displaystyle \bigoplus_{g \in G} \langle y | V | g \ccdot x \rangle \otimes (\X_\mlag{x}{}{g})_\mathbb{C} \otimes \T_g
 \ar[rrr]_-{\beta} &&&  \displaystyle \bigoplus_{g \in G} \langle h \ccdot y | V | hg \ccdot x \rangle \otimes (\X_\mmlag{h \ccdot x}{}{hgh^\mmi})_\mathbb{C} \otimes \T_{hgh^\mmi}}
\]
where
\begin{align*}
 \xi &= \text{swap} \circ (v \mapsto u) \otimes \id \otimes V(\lag{g}{v \otimes \bar{u}}) \\
 \zeta &= \text{swap} \circ (tvt'^\mmi \mapsto s'us^\mmi) \otimes \id \otimes  V(\lag{hgh^\mmi}{tvt'^\mmi \otimes \overline{s'us^\mmi}}) \\
 \alpha &= \ad_{t, t'} \otimes \T(\sla{h}{g}) \otimes V(\lag{g}{t' \otimes \bar{s}}) \\
 \beta &= V(\lag{h}{t \otimes \bar{s'}}) \otimes \ad_{s', s} \otimes \T(\sla{h} g).
\end{align*}
The diagram above consists of three simultaneous `tracks' (each step is a tensor product of three operations). The track involving $\T$ clearly commutes, the track involving the gerbes $\X$ and $\Y$ is
 \[
 \ba \xymatrix{ \Y_\mlag{}{y}{x} \ar[d]_{v \mapsto u} \ar[r]^-{\ad_{t,t'}} & \Y_\mmlag{}{h \ccdot y}{hgh^\mmi} \ar[d]^{tvt'^\mmi \mapsto s'us^\mmi} \\
 \X\mlag{x}{}{h} \ar[r]_-{\ad_{s',s}} & \X_\mmlag{h \ccdot x}{}{hgh^\mmi}} \ea ,
\]
while the track involving the equivariant vector bundle $V$ is:
 \[
 \ba \xymatrix{ \langle x^\mmi \ccdot y | V | x \rangle \ar[rr]^-{V(\lag{h}{t'\otimes \bar{s}})} \ar[d]_{V(\lag{x}{v \otimes \bar{u}})} && \langle gx^\mmi \ccdot y | V | h \ccdot x \rangle \ar[d]^-{V(\xleftarrow[tvt'^\mmi \otimes \overline{s'us^\mmi}]{hgh^\mmi})} \\ \langle y | V  | x \ccdot x \rangle \ar[rr]_-{V(\lag{h}{t \otimes \bar{s'}})} && \langle h \ccdot y | V | gx \ccdot x \rangle} \ea.
 \]
Here we have been using the following arrows in the gerbes:
 \[
  \ba \xymatrix{h \ccdot y & gx^\mmi \ccdot y \ar[l]_-{\;\; tvt'^\mmi} \\ y \ar[u]^t & x^\mmi \ccdot y \ar[l]^v \ar[u]_{t'}} \ea \quad \quad
 \ba \xymatrix{gx \ccdot x & h \ccdot x \ar[l]_-{\;\;s'us^\mmi} \\ x \ccdot x \ar[u]^{s'} & x \ar[l]^u \ar[u]_{s}} \ea
 \]
In any event, both diagrams above commute, which is what we needed to show.

It remains to establish that the isomorphisms $\hat{\T}(V)$ obey the coherence equation for a transformation --- that is, if $W \colon \Y \rightarrow \Z$ is another morphism of equivariant gerbes then we must have the following equality:
 \[
 \ba
\xy
(-7.5,7.5)*+{\X}="1";
(-7.5,-7.5)*+{\X}="2";
(11,7.5)*+{\Y}="3";
(11,-7.5)*+{\Y}="4";
(29.5,7.5)*+{\Z}="5";
(29.5,-7.5)*+{\Z}="6";
{\ar_-{\hat{\T}_\X} "1";"2"};
{\ar^-{V} "1";"3"};
{\ar^-{\! \hat{\T}_\Y} "3";"4"};
{\ar_-{V} "2";"4"};
{\ar^-{W} "3";"5"};
{\ar_-{W} "4";"6"};
{\ar^-{\hat{\T}_\Z} "5";"6"};
{\ar@2{->}^-{\hat{\T}(V)} (1,2)*{};(-3,-2)*{}};
{\ar@2{->}^-{\hat{\T}(W)} (20.5,2)*{};(16.5,-2)*{}};
\endxy
 \ea
 =
 \ba
\xy
(-7.5,7.5)*+{\X}="1";
(-7.5,-7.5)*+{\X}="2";
(14,7.5)*+{\Z}="3";
(14,-7.5)*+{\Z}="4";
{\ar_-{\hat{\T}_\X} "1";"2"};
{\ar^-{W \circ V} "1";"3"};
{\ar^-{\hat{\T}_\Z} "3";"4"};
{\ar_-{W \circ V} "2";"4"};
{\ar@2{->}^-{\hat{\T}(W \circ V)} (1,2)*{};(-3,-2)*{}};
\endxy
\ea
\]
This equation indeed holds for the simple reason that $\hat{\T}(W \circ V)$ is {\em defined} in terms of the vector bundle $W \circ V$, which is in turn defined via the tensor product of $W$ and $V$, which is the requirement of the equation above. In other words, this equation is a routine unraveling of definitions, as the reader will confirm. This completes the proof that $\hat{\T}$ is indeed a transformation of the identity on $\Gerbes(G)$.

Finally, we need to show that the restriction of $\hat{\T}$ to $\E G$ is isomorphic in $\Dim\left( \Gerbes(G)|_{\E G}\right)$ to our original transformation $\T$. We only need to show that their corresponding equivariant vector bundles $\{\hat{\T}_g\}_{g \in G}$ and $\{\T\}_{g \in G}$ are isomorphic. Indeed, unraveling the definitions we see that
 \begin{align*}
  \hat{\T}_g &= \langle g | \hat{\T}_{\E G} | e \rangle  = (E G_\mlag{e}{}{g})_\mathbb{C} \otimes \langle g | \T_{\E G} | e \rangle = \langle g | \T_{\E G} | e \rangle \\
   &= \T_g
  \end{align*}
where in the second last step we once more used our convention of ignoring factors in tensor products which identify canonically with $\mathbb{C}$ (in this case, the factor {\em is actually} $\mathbb{C}$ by definition). So, the fibers of $\hat{\T}_g$ identify with those of $\T_g$. Moreover, the equivariant maps $\hat{\T}(\sla{h} g)$ identify with $\T(\sla{h} g)$, because if we unravel the definitions, we have:
 \begin{align*}
  \hat{\T}(\sla{h} g) &= \langle g | \T_{\E G} | e \rangle \stackrel{\hat{\T}(\lag{h}{1 \otimes \bar{1}})}{\longrightarrow} \langle hg | \hat{\T}_{\E G} | h \rangle \stackrel{\hat{\T}(U_h)_{hg, e}}{\longrightarrow} \langle hgh^\mi | \hat{\T}_{\E G} | e \rangle \\
  &= \T_g \stackrel{\T(\sla{h}g)}{\longrightarrow}  \T_{hgh^\mi} \stackrel{\id}{\longrightarrow} \T_{hgh^\mi} \\
  &= \T(\sla{h} g).
 \end{align*}
This completes the proof.
  \end{proof}

\section{A transformation is determined by its data at $\E G$\label{le}}
In this section, the analogue of Lemma \ref{myfourth}, we show that a transformation $\T$ of the identity on $\Gerbes(G)$ is indeed completely determined by its behaviour at the `regular equivariant gerbe' $\E G$. That is, we construct the components of an isomorphism $\Omega \colon \T \rightarrow \hat{\T}$, where $\hat{\T}$ is the transformation of the identity on $\Gerbes(G)$  we defined in the last section which was solely constructed out of the data of $\T$ at $\E G$. Then we confirm that $\Omega$ is indeed an isomorphism inside $\Dim \Gerbes(G)$.

To specify the isomorphism $\Omega$, we need to give for each equivariant gerbe $\X$ an isomorphism of equivariant vector bundles over $\X \otimes \overline{\X}$
  \[
   \Omega_\X \colon \T_\X \rightarrow \hat{\T}_\X.
  \]
We define these by observing that $\T$ is a transformation of the identity, so for each $x \in \X$ we have the coherence isomorphism
\[
\ba \xy
(-7.5,7.5)*+{\E G}="1";
(-7.5,-7.5)*+{\E G}="2";
(10.5,7.5)*+{\X}="3";
(10.5,-7.5)*+{\X}="4";
{\ar_-{\T_{\E G}} "1";"2"};
{\ar^-{U_x} "1";"3"};
{\ar^-{\T_\X} "3";"4"};
{\ar_-{U_x} "2";"4"};
{\ar@2{->}^-{\T(U_x)} (1,2)*{};(-3,-2)*{}};
\endxy \ea.
\]
At $(x', e) \in \X \otimes \overline{\E G}$, the fibres of the clockwise and counterclockwise composite vector bundles above compute as
 \begin{align*}
  \langle x' | \T_\X \circ U_x | e \rangle &= \bigoplus_{x'' \in \X} \langle x' | \T_\X | x'' \rangle \otimes \langle x'' | U_x | e \rangle \\
   &= \langle x' | \T_\X | x \rangle
 \end{align*}
and
 \begin{align*}
  \langle x' | U_x \circ \T_{\E G} | e \rangle &= \bigoplus_{g \in G} \langle x' | U_x | g \rangle \otimes \langle g | \T_{\E G} | e \rangle \\
   &= \bigoplus_\blag{x}{x'}{g} (\X_\mlag{x}{}{g})_\mathbb{C} \otimes \T_g \\
   &= \langle x' | \hat{\T}_\X | x \rangle
  \end{align*}
respectively, so we obtain fibrewise isomorphisms
 \[
  [\Omega_\X]_{x',x} \equiv \T(U_x)_{x', e} \colon \langle x' | \T_\X | x \rangle \rightarrow  \bigoplus_\blag{x}{x'}{g} (\X_\mlag{x}{}{g})_\mathbb{C} \otimes  \T_g.
 \]
We then have the main result of this section.
\begin{prop}\label{aTgim} Suppose that $\T$ is a transformation of the identity on $\Gerbes(G)$. Then:
 \begin{enumerate} \item For each equivariant gerbe $\X$, the isomorphism $\Omega_\X \colon \T_\X \rightarrow \hat{\T}_\X$ defined fibrewise above is indeed an isomorphism of equivariant vector bundles.
 \item Taken together, the isomorphisms
  \[
   \{ \Omega_\X \colon \T_\X \rightarrow \hat{\T}_\X \}
  \]
 satisfy the modification condition and hence $\Omega\colon \T \rightarrow \hat{\T}$ is indeed an isomorphism between these two transformations in $\Dim \Gerbes(G)$.
 \end{enumerate}
 \end{prop}
\begin{proof} (i) For this first part of the proof we regard the equivariant gerbe $\X$ as fixed, so we take the notational liberty of writing $\Omega \equiv \Omega_\X$. We want to show that the following diagram commutes:
  \be \label{ag}
  \xymatrix{ \langle x' | \T_\X | x \rangle \ar[rr]^-{\Omega_{x',x}} \ar[d]_-{\T_\X (\lag{h}{v \otimes \bar{u}})} && \displaystyle \bigoplus_\blag{x}{x'}{g} (\X_\mlag{x}{}{g})_\mathbb{C} \otimes  \T_g \ar[d]^-{\ad_{v,u} \otimes \T(\sla{h}g)} \\
  \langle h \ccdot x' | \T_\X | h \ccdot x \rangle \ar[rr]_-{\Omega_{h \ccdot x', h \ccdot x}} &&  \quad \quad \quad \displaystyle \bigoplus_\blag{x}{x'}{g} (\X_\mmlag{h \ccdot x}{}{hgh^\mmi})_\mathbb{C} \otimes  \T_{hgh^\mi}}
 \ee
We do this by combining two equations, one arising from the fact that
 \[
 \T(U_x) \colon \T_\X \circ U_x \rightarrow U_x \circ \T_{\E G}
 \]
is an isomorphism of {\em equivariant} vector bundles, and the other arising from the coherence equation on the isomorphisms $\T(U)$ when one composes $U$ and $V$. \XXX{rewrite this}.

To obtain the first equation, observe that the isomorphism $\T(U_x)$ has fibrewise components
 \[
  \T(U_x)_{h \ccdot x', h} \colon \xymatrix @R=0.3cm{\langle h \ccdot x' | \T_\X \circ U_x | h \rangle \ar@{=}[d] \ar[r] & \langle h \ccdot x' | U_x \circ \T_{\E G} | h \rangle \ar@{=}[d] \\ \langle h \ccdot x' | T_\X | h \ccdot x \rangle \otimes (\X_\mlag{x}{}{h})_\mathbb{C} & \smash{{\displaystyle \bigoplus_\blag{x}{x'}{g}}} (\X_\mlag{x}{}{hg})_\mathbb{C} \otimes \langle hg | T_{\E G} | h \rangle}.
 \]

\vspace{0.4cm}

\noindent This motivates us to define the isomorphisms
 \[
  \Omega^h_{x',x} :=  [\id \otimes \T_{\E G}(\lag{h^\mmi}{1 \otimes \bar{1}})] \circ \T(U_x)_{h \ccdot x', h},
 \]
in other words,
\[
 \Omega^h_{x',x} \colon \langle h \ccdot x' | \T_\X | h \ccdot x \rangle \otimes (\X_\mlag{x}{}{h})_\mathbb{C} \rightarrow \bigoplus_\blag{x}{x'}{g} (\X_\mlag{x}{}{hg})_\mathbb{C} \otimes \T_g.
 \]
Now, the fact that $\T(U_x)$ is an isomorphism of equivariant vector bundles translates into the following commutative diagram:
 \[
 \xymatrix{\langle x' | \T_\X \circ U_x | e \rangle \ar[rr]^-{\Omega_{x',x}} \ar[d]_-{(\T_\X \circ U_x) (\lag{h}{v \otimes \bar{1}})} && \langle x' | U_x \circ \T_{\E G} | e \rangle \ar[d]^-{ (U_x \circ \T) (\lag{h}{v \otimes \bar{1}})} \\
  \langle h \ccdot x' | \T_\X \circ U_x | h \rangle \ar[rr]_-{\T(U_x)_{h \ccdot x', x}}
  &&  \langle h \ccdot x' | U_x \circ \T_{\E G} | h \rangle}
\]
By expanding out this diagram we obtain our first equation, which relates $\Omega^h_{x',x}$ with $\Omega_{h \ccdot x', h \ccdot x}$:
 \be \label{ac}
 \ba \xymatrix{\langle x' | \T_\X | x \rangle \ar[rr]^-{\Omega_{x',x}} \ar[d]_-{\T_\X(\lag{h}{v \otimes \bar{u}}) \otimes u} &&
  \displaystyle \bigoplus_\blag{x}{x'}{g} (\X_\mlag{x}{}{g})_\mathbb{C} \otimes  \T_g \ar[d]^-{\post(v) \otimes \id}  \\
  \langle h \ccdot x' | \T_\X | h \ccdot x \rangle \otimes (\X_\mlag{x}{}{h})_\mathbb{C} \ar[rr]_-{\Omega^h_{x',x}} && \displaystyle \bigoplus_\blag{x}{x'}{h} (\X_\mlag{x}{}{hg})_\mathbb{C} \otimes  \T_g} \ea
 \ee
For the other equation, consider the coherence equation corresponding to the composite of the 1-morphisms $U_x$ and $U_h$:
 \be \label{acoh}
 \ba
\xy
(-7.5,7.5)*+{\E G}="1";
(-7.5,-7.5)*+{\E G}="2";
(11,7.5)*+{\E G}="3";
(11,-7.5)*+{\E G}="4";
(33.5,7.5)*+{\X}="5";
(33.5,-7.5)*+{\X}="6";
{\ar_-{\T_{\E G}} "1";"2"};
{\ar^-{U_{h}} "1";"3"};
{\ar^-{\T_{\E G}} "3";"4"};
{\ar_-{U_{h}} "2";"4"};
{\ar^-{U_x} "3";"5"};
{\ar_-{U_x} "4";"6"};
{\ar^-{\T_\X} "5";"6"};
{\ar@2{->}^-{\T(U_h)} (1,2)*{};(-3,-2)*{}};
{\ar@2{->}^-{\T(U_x)} (23.5,2)*{};(19.5,-2)*{}};
\endxy
 \ea
 =
 \ba
\xy
(-7.5,7.5)*+{\E G}="1";
(-7.5,-7.5)*+{\E G}="2";
(14,7.5)*+{\X}="3";
(14,-7.5)*+{\X}="4";
{\ar_-{\T_{\E G}} "1";"2"};
{\ar^-{U_x \circ U_h} "1";"3"};
{\ar^-{\T_\X} "3";"4"};
{\ar_-{U_x \circ U_h} "2";"4"};
{\ar@2{->}^-{\T(U_x \circ U_h)} (1,2)*{};(-3,-2)*{}};
\endxy
\ea
   \ee
At $(h \ccdot x', e) \in \X \otimes \overline{\E G}$,  \XXX{above there is a 2-morphism which strikes through something in the drawing}the left hand side of this coherence equation computes as:
 \be \label{ab}
  \ba \xymatrix{ \langle h \ccdot x' | \T_\X \circ U_x \circ U_h | e \rangle \ar[d]^{\T(U_x) \otimes \id} \\ \langle h \ccdot x' | U_x \circ \T_{\E G} \circ U_h | e \rangle \ar[d]^{\id \otimes \T(U_h)} \\ \langle h \ccdot x' | U_x \circ U_h \circ \T_{\E G} | e \rangle} \ea \quad = \quad \ba \xymatrix{ \langle h \ccdot x' | \T_\X | h \ccdot x \rangle \otimes (\X_\mlag{x}{}{h})_\mathbb{C} \ar[d]^{\Omega^h_{x',x}} \\ \displaystyle \bigoplus_\blag{x}{x'}{g} (\X_\mlag{x}{}{hg})_\mathbb{C} \otimes  \T_g \ar[d]^{\id \otimes \T(\sla{h} g)} \\ \displaystyle \bigoplus_\blag{x}{x'}{g} (\X_\mlag{x}{}{hg})_\mathbb{C} \otimes  \T_{hgh^\mi}} \ea .
 \ee
To evaluate the {\em right} hand of the coherence equation \eqref{acoh}, we need to relate $U_x \circ U_h$ to $U_{h \ccdot x}$. Indeed, for every $u \in \X_\mlag{x}{}{h}$ there is a canonical isomorphism
\[
\xy
(-15, 0)*+{\E G}="1";
(15,0)*+{\X}="2";
{\ar^-{U_{h \ccdot x}} "1";"2"};
{\ar@/_2pc/_{U_x \circ U_h} "1";"2"};
{\ar@2{->}^-{\, \pre(u)} (-2,-1.5)*{};(-2,-6.5)*{}};
\endxy
\]
which translates the orbit emanating out of $h \ccdot x$ into the orbit emanating out of $x$ by by precomposing with $u$ (see Figure XXX mention also that $g \cdot x = x'$ in figure):
 \[
  \ba \xymatrix{\langle h \ccdot x' | U_{h \ccdot x} | hgh^\mmi \rangle \ar[d]^{\pre(u)} \\ \langle h \ccdot x' | U_x \circ U_h | hgh^\mmi \rangle } \ea = \ba \xymatrix{\X_\mmlag{h \ccdot x'}{}{hgh^\mmi} \ar[d] \\ \X_\mlag{x}{}{hg}} \ea \quad \ba \xymatrix{t \ar@{|->}[d] \\ t \circ u} \ea.
  \]
By naturality of the transformation $\T$, this means the 2-cell on the right hand side of the coherence equation \eqref{acoh} computes as:
 \[
  \ba
\xy
(-7.5, 10)*+{\E G}="1";
(-7.5,-10)*+{\E G}="2";
(14,10)*+{\X}="3";
(14,-10)*+{\X}="4";
{\ar_-{\T_{\E G}} "1";"2"};
{\ar^-{U_x \circ U_h} "1";"3"};
{\ar^-{\T_\X} "3";"4"};
{\ar_-{U_x \circ U_h} "2";"4"};
{\ar@2{->}^-{\T(U_x \circ U_h)} (1,3)*{};(-3,-1)*{}};
\endxy
\ea
=
  \ba
\xy
(-7.5,10)*+{\E G}="1";
(-7.5,-10)*+{\E G}="2";
(20,10)*+{\X}="3";
(20,-10)*+{\X}="4";
{\ar_-{\T_{\E G}} "1";"2"};
{\ar_-{U_{h \ccdot x}} "1";"3"};
{\ar^-{\T_\X} "3";"4"};
{\ar^-{U_{h \ccdot x}} "2";"4"};
{\ar@/^2pc/^{U_x \circ U_h} "1";"3"};
{\ar@/_2pc/_{U_x \circ U_h} "2";"4"};
{\ar@2{->}^-{\T(U_{h \ccdot x})} (6,3)*{};(2,-1)*{}};
{\ar@2{->}^-{\, \pre(u^\mmi)} (2.5,15)*{};(2.5,11)*{}};
{\ar@2{->}^-{\, \pre(u)} (3,-12)*{};(3,-16)*{}};
\endxy
\ea
\]
At $(h \ccdot x', e) \in \X \times \overline{\E G}$, this expands out as:
 \[
 \ba \xymatrix{\langle h \ccdot x' | \T_\X \circ U_x \circ U_h | e \rangle \ar[d]^{\T(U_x \circ U_h)_{h \ccdot x', e}} \\ \langle h \ccdot x' | U_x \circ U_h \circ \T | e \rangle} \ea = \ba \xymatrix{\langle h \ccdot x' | \T_\X | h \ccdot x \rangle \otimes (\X_\mlag{x}{}{g})_\mathbb{C} \ar[d]^{\id \otimes \pre(u^\mi)} \\
 \langle h \ccdot x' | \T_\X | h \ccdot x \rangle \ar[d]^{\Omega_{h \ccdot x', h \ccdot x}} \\ \displaystyle \bigoplus_{\blag{x}{x'}{g}} (\X_\mmlag{h \ccdot x}{}{hgh^\mmi})_\mathbb{C} \otimes \T_{hgh^\mmi} \ar[d]^{\pre(u) \otimes \id} \\ \displaystyle \bigoplus_{\blag{x}{x'}{g}} (\X_\mlag{x}{}{hg})_\mathbb{C} \otimes \T_{hgh^\mmi}} \ea
 \]
Combining this with \eqref{ab}, we see that the coherence properties of $\T$ expressed in \eqref{acoh} amount to the commutative diagram \XXX{big problem here with i and j!}
  \be \label{ad}
  \ba \xymatrix @C=0.55in { \langle h \ccdot x' | \T_\X | h \ccdot x \rangle \otimes (\X_\mlag{x}{}{h})_\mathbb{C} \ar[dd]_{\Omega^h_{x',x}} \ar[r]^-{\id \otimes \pre(u^\mi)} & \langle h \ccdot x' | \T_\X | h \ccdot x \rangle \ar[dr]^-{\Omega_{h \ccdot x', h \ccdot x}} \\
   & & **[l] \displaystyle \bigoplus_{\blag{x}{x'}{g}} (\X_\mmlag{h \ccdot x}{}{hgh^\mmi})_\mathbb{C} \otimes \T_{hgh^\mi} \ar[dl]^-{\pre(u) \otimes \id}\\
    \displaystyle \bigoplus_{\blag{i}{j}{x}} (\X_\mlag{i}{}{gx})_\mathbb{C} \otimes \T_g \ar[r]_-{\id \otimes \T(\sla{h}g)} & \displaystyle \bigoplus_{\blag{x}{x'}{g}} (\X_\mlag{x}{}{hg})_\mathbb{C} \otimes \T_{hgh^\mmi}} \ea
   \ee
By combining the two equations \eqref{ac} and \eqref{ad} we have obtained, we can finally show that the diagram \eqref{ag} indeed commutes:
 \begin{align*}
  \Omega_{h \ccdot x', h \ccdot x} \circ \T_\X (\lag{h}{v \otimes \bar{u}}) &\stackrel{\text{\eqref{ac}}}{=} [\pre(u^\mmi)\otimes \id] \circ [\id \otimes \T(\sla{h}g)] \circ \Omega^h_{x',x} \circ [\id \otimes \pre(u)] \circ \T_\X (\lag{h}{v \otimes \bar{u}}) \\
  &\stackrel{\text{\eqref{ad}}}{=} [\pre(u^\mi) \otimes \id] \circ [\id \otimes \T(\sla{h} g)] \circ [\post(v) \otimes \id] \circ \Omega_{x',x} \\
  &= [\ad_{v,u} \otimes \T(\sla{h}g) ] \circ \Omega_{x',x}.
 \end{align*}

(ii) We need to show that for every $V \colon \X \rightarrow \Y$, the following pasting diagram commutes:
\[
\xy
(-7.5,10)*+{\X}="1";
(-11, 7)="1f"; (-5, 9)="1b";
(-7.5,-10)*+{\X}="2";
(-11, -10)="2f"; (-5, -7)="2b";
(35,10)*+{\Y}="3";
(33.5, 7)="3f"; (38.5, 9)="3b";
(35,-10)*+{\;\; \Y}="4";
(33.5, -10)="4f"; (38.5, -7)="4b";
{\ar@/_0.8pc/_-{\T_\X} "1f";"2f"};
{\ar@/^0.5pc/^-{\hat{\T}_\X} "1b";"2b"};
{\ar@/_0.8pc/_-{\T_\Y} "3f";"4f"};
{\ar@/^0.5pc/^-{\hat{\T}_\Y} "3b";"4b"};
{\ar^-{V} "1";"3"};
{\ar_-{V} "2";"4"};
{\ar@2{->}@/{^0.5pc}/^>>>{{\footnotesize \T(V)}} (14.5,6)*+{};(12.5,0)*+{}};
{\ar@2{->}@/{_0.5pc}/_<<<<<{\hat{\T}(V)} (11,2)*+{};(11,-4)*+{}};
{\ar@2{->}_-{\Omega_\X} (-13,-3)*+{};(-4,2)*+{}};
{\ar@2{->}_-{\Omega_\Y} (31,-3)*+{};(38.5,2)*+{}};
\endxy
\]
In terms of fibrewise components, this is the following diagram:
 \[
  \xymatrix{\langle y | \T_\Y \circ V | x \rangle \ar[d]_{\T(V)_{y,x}} \ar[rr]^-{\Omega_\Y \otimes \id} && \langle y | \hat{\T}_\Y \circ V | x \rangle \ar[d]^{\hat{\T}(V)_{y,x}} \\
   \langle y | V \circ \T_\X | x \rangle \ar[rr]_{\id \otimes \Omega_\X} && \langle y | V \circ \hat{\T}_\X | x \rangle}
  \]
We can expand this diagram out to the following requirement:
 \be \label{af}
 \ba \xymatrix{
  \displaystyle \bigoplus_{y'} \langle y | \T_\Y | y' \rangle \otimes \langle y' | V | x \rangle \ar[rr]^-{[\Omega_\Y]_{y, y'} \otimes \id} \ar[d]_{\T(V)_{y, x}} && \displaystyle \bigoplus_{y'} \bigoplus_\blag{y'}{y}{g} (\Y_\mlag{}{y}{g})_{\mathbb{C}} \otimes \T_g \otimes \langle y' | V | x \rangle \ar[d]^-{\text{swap} \circ (v \mapsto u) \otimes \id \otimes V(\lag{g}{v \otimes \bar{u}})} \\ \displaystyle \bigoplus_{x'} \langle y | V | x'\rangle \otimes \langle x' | \T_\X | x \rangle \ar[rr]_-{\id \otimes [\Omega_\X]_{x',x}} && \displaystyle \bigoplus_{x'} \bigoplus_\blag{x}{x'}{g} \langle y | V | x' \rangle \otimes (\X_\mlag{x}{}{g})_\mathbb{C} \otimes \T_g} \ea
  \ee
We claim that this follows from the following coherence equation on $\T$:
 \be \label{acoh2}
 \ba
\xy
(-7.5,7.5)*+{\E G}="1";
(-7.5,-7.5)*+{\E G}="2";
(11,7.5)*+{\X}="3";
(11,-7.5)*+{\X}="4";
(29.5,7.5)*+{\Y}="5";
(29.5,-7.5)*+{\Y}="6";
{\ar_-{\T_{\E G}} "1";"2"};
{\ar^-{U_x} "1";"3"};
{\ar^-{\T_\X} "3";"4"};
{\ar_-{U_x} "2";"4"};
{\ar^-{V} "3";"5"};
{\ar_-{V} "4";"6"};
{\ar^-{\T_\Y} "5";"6"};
{\ar@2{->}^-{\T(U_x)} (1,2)*{};(-3,-2)*{}};
{\ar@2{->}^-{\T(V)} (21,2)*{};(17,-2)*{}};
\endxy
 \ea
 = \ba
 \xy
(-7.5,7.5)*+{\E G}="1";
(-7.5,-7.5)*+{\E G}="2";
(14,7.5)*+{\Y}="3";
(14,-7.5)*+{\Y}="4";
{\ar_-{\T_{\E G}} "1";"2"};
{\ar^-{V \circ U_x} "1";"3"};
{\ar^-{\T_\X} "3";"4"};
{\ar_-{V \circ U_x} "2";"4"};
{\ar@2{->}^-{\T(V \circ U_x)} (1,2)*{};(-3,-2)*{}};
 \endxy
 \ea
 \ee
Indeed, at $(y,e) \in \Y \otimes \overline{\E G}$ the fibrewise component of the isomorphism of vector bundles given by the LHS above computes as:
 \[
  \xymatrix @C=1.8cm { \langle y | \T_\Y \circ V | x \rangle \mathclap{\quad \quad \qquad \quad \qquad \qquad =\langle y | \T_{\Y} \circ V \circ U_x | e \rangle} \ar[d]_-{\T(V)_{y, x}}  \\ \displaystyle \bigoplus_{x'} \langle y | V | x' \rangle \otimes \langle x' | \T_\X | x \rangle \ar[r]_-{\id \otimes [\Omega^\X]_{x',x}} & \displaystyle \bigoplus_{x'} \bigoplus_\blag{x}{x'}{g} \langle y | V | x' \rangle \otimes (\X_\mlag{x}{}{g})_\mathbb{C} \otimes \T_g}
 \]
To compute the RHS of \eqref{acoh2}, we need to expand out the vector bundle $V \circ U_x$ as a direct sum of the vector bundles $U_y$ as $y$ ranges over $\Y$, using  Lemma \ref{aprelim} (ii). That is, if we write
\[
V_x = \bigoplus_{y'} \langle y' |V|x \rangle U_y'
\]
then Lemma \ref{aprelim} (ii) gives us a canonical isomorphism $\xi \colon V_x \cong V \circ U_x$. By naturality of $\T$, we therefore compute the RHS of \eqref{acoh2} as the following composite:
\[
 \ba
 \xy
(-7.5,7.5)*+{\E G}="1";
(-7.5,-7.5)*+{\E G}="2";
(14,7.5)*+{\Y}="3";
(14,-7.5)*+{\Y}="4";
{\ar_-{\T_{\E G}} "1";"2"};
{\ar^-{V \circ U_x} "1";"3"};
{\ar^-{\T_\X} "3";"4"};
{\ar_-{V \circ U_x} "2";"4"};
{\ar@2{->}^-{\T(V \circ U_x)} (1,2)*{};(-3,-2)*{}};
 \endxy
 \ea
=
\ba \xy
(-7.5,10)*+{\E G}="1";
(-7.5,-10)*+{\E G}="2";
(20,10)*+{\Y}="3";
(20,-10)*+{\Y}="4";
{\ar_-{\T_{\E G}} "1";"2"};
{\ar_-{V_x} "1";"3"};
{\ar^-{\T_\Y} "3";"4"};
{\ar^-{V_x}  "2";"4"};
{\ar@/^2pc/^{V \circ U_x} "1";"3"};
{\ar@/_2pc/_{V \circ U_x} "2";"4"};
{\ar@2{->}^-{\T(V_x)} (6,3)*{};(2,-1)*{}};
{\ar@2{->}^-{\, \xi^\mi} (5,16)*{};(5,12)*{}};
{\ar@2{->}^-{\, \xi} (5,-12)*{};(5,-16)*{}};
\endxy
\ea
\]
At $(y, e) \in \Y \otimes \overline{\E G}$, this gives us:
 \[
  \xymatrix{\displaystyle \bigoplus_{y'} \langle y | \T_\Y | y' \rangle \otimes \langle y' | V \circ U_x | e \rangle \ar[r]^-{\id \otimes \xi_{y', e}^\mi = \id} & \displaystyle \bigoplus_{y'} \langle y | \T_\Y | y' \rangle  \otimes \langle y' | V | x \rangle \ar[d]^{\text{swap} \circ \; [\Omega_\Y]_{y, y'} \otimes \id} \\ &  \displaystyle \bigoplus_{y'} \bigoplus_\blag{y'}{y}{g} \langle y' | V | x  \rangle \otimes (\Y_\mlag{}{y}{g})_\mathbb{C} \otimes \T_g \ar[d]^{V(\lag{g}{v \otimes \bar{u}}) \otimes (v \mapsto u) \otimes \id} \\ & \displaystyle \bigoplus_{y'} \bigoplus_\blag{x}{x'}{g} \langle y | V | x' \rangle \otimes (\X_\mlag{x}{}{g})_\mathbb{C} \otimes \T_g}
 \]
Combining this with our previous expression for the LHS gives us the equation we need \eqref{af}.
\end{proof}

\section{The restriction map is fully faithful\label{lf}}
So far we have established that the operation of restricting a transformation of the identity $\T \in \Dim \left( \Gerbes(G) \right)$ to $\Dim \left(\Gerbes(G)|_{\E G}\right)$ is essentially surjective (Lemma \ref{Impy}) and also `essentially injective' (Proposition \ref{aTgim}). We have thus completed the categorified analogues of the lemmas from Section \ref{lb}. But there is a new categorical phenomenon which we must account for: we must show that the restriction functor is {\em fully faithful}.
\begin{prop}\label{afullyfaithful}Suppose that $\T, \Ss \in \Dim \Gerbes(G)$. Then for every morphism $\vartheta \colon \T_{\E G} \rightarrow \Ss_{\E G}$ in $\Dim \left(\Gerbes(G)|_{\E G}\right)$ there exists a unique extension
 \[
  \{ \theta_\X \colon \T_\X \rightarrow \Ss_\X \}
 \]
of $\vartheta$ to a morphism in $\Dim \Gerbes(G)$ such that $\theta_{\E G} = \vartheta$.
\end{prop}
\begin{proof} To prove uniqueness, we must show that a morphism of transformations
 \[
  \theta = \{ \theta_\X \colon \T_\X \rightarrow \Ss_\X \}
 \]
is determined by its component $\theta_{\E G}$ at $\E G$. Now, given an equivariant gerbe $\X$ the condition that $\theta$ is a modification between the transformations $\T$ and $\Ss$ implies that for each $x \in \X$ the following pasting diagram must commute:
\[
\xy
(-7.5,10)*+{\E G}="1";
(-11, 7)="1f"; (-5, 9)="1b";
(-7.5,-10)*+{\E G}="2";
(-11, -10)="2f"; (-5, -7)="2b";
(35,10)*+{\X}="3";
(33.5, 7)="3f"; (38.5, 9)="3b";
(35,-10)*+{\X}="4";
(33.5, -10)="4f"; (38.5, -7)="4b";
{\ar@/_0.8pc/_-{\T_{\E G}} "1f";"2f"};
{\ar@/^0.5pc/^-{\Ss_{\E G}} "1b";"2b"};
{\ar@/_0.8pc/_-{\T_\X} "3f";"4f"};
{\ar@/^0.5pc/^-{\Ss_\X} "3b";"4b"};
{\ar^-{U_x} "1";"3"};
{\ar_-{U_x} "2";"4"};
{\ar@2{->}@/{^0.5pc}/^>>>{{\footnotesize \Ss(U_x)}} (14.5,6)*+{};(12.5,0)*+{}};
{\ar@2{->}@/{_0.5pc}/_<<<<<{\T(U_x)} (11,2)*+{};(11,-4)*+{}};
{\ar@2{->}_-{\theta_{\E G}} (-13,-3)*+{};(-4,2)*+{}};
{\ar@2{->}_-{\theta_{\X}} (31,-3)*+{};(38.5,2)*+{}};
\endxy
\]
At $(x',e) \in \X \otimes \overline{\E G}$, this expands out to the commutative diagram
 \[
  \xymatrix{ \langle x' | \T_\X | x \rangle \ar[d]_-{[\Omega^{\T}]_{x',x}}  \ar[r]^-{[\theta_\X]_{x',x}} & \langle x' | \Ss_\X | x \rangle \ar[d]^-{[\Omega_{\Ss}]_{x',x}} \\ **[l] \displaystyle \bigoplus_\blag{x}{x'}{g} (\X_\mlag{x}{}{g})_\mathbb{C} \otimes \T_g \ar[r]_-{\id \otimes \theta_g} & **[r] \displaystyle \bigoplus_\blag{x}{x'}{g} (\X_\mlag{x}{}{g})_\mathbb{C} \otimes W_g  }
 \]
where
 \[
  \theta_g := [\theta_{\E G}]_{g,e} \colon \T_g \rightarrow \Ss_g
 \]
refers to the morphism of equivariant vector bundles $\{ \T_g \} \rightarrow \{ \Ss_g \}$ obtained by regarding $\T_{\E G}$ and $\Ss_{\E G}$ as equivariant vector bundles over $G$ (see XXX). In other words, the components of $\theta_\X$ can be expressed in terms of these morphisms $\theta_g$  as
 \be \label{apresc}
  [\theta_\X]_{x',x} = [\Omega_{\Ss}]^\mi_{x',x} \circ \biggl[ \bigoplus_\blag{x}{x'}{g} \id \otimes \theta_g \biggr] \circ [\Omega_\T]_{x',x},
 \ee
which establishes uniqueness.

To establish existence, suppose that $\T$ and $\Ss$ are transformations of the identity on $\Gerbes(G)$, and that we are given a map
 \[
  \vartheta \colon \T_{\E G} \rightarrow \Ss_{\E G}
 \]
in $\Dim \left( \Gerbes(G)_{\E G} \right)$ between their restrictions to $\E G$. Recall from XXX that $\vartheta$ is the same thing as a map of equivariant vector bundles over $G$,
 \[
  \{ \vartheta_g \colon \T_g \rightarrow \Ss_g\}.
 \]
We need to show that we can extend $\vartheta$ to a map $\theta \colon \T \rightarrow \Ss$ by defining each map $\theta_\X \colon \T_\X \rightarrow \Ss_\X$ in terms of $\vartheta$ by the above prescription \eqref{apresc}.

The main thing to check is that this definition of $\theta_\X$ is indeed a map of {\em equivariant} vector bundles. This is the requirement that the front face of the cube below commutes:
\[
\xymatrix @C=0.17cm { & \displaystyle \bigoplus_\blag{x}{x'}{g} (\X_\mlag{x}{}{g})_\mathbb{C} \otimes \T_g \ar'[d][dd]^-{\ad_{v,u} \otimes \T(\sla{h}g)} \ar[rr]^{\id \otimes \theta_g} && \displaystyle \bigoplus_\blag{x}{x'}{g} (\X_\mlag{x}{}{g})_\mathbb{C} \otimes \Ss_g \ar[dd]^{\ad_{v, u} \otimes \Ss(\sla{h}g)}\\
\langle x' | \T_\X | x \rangle \ar[ur]^-{\Omega^\T_{x',x}}  \ar[rr]^{\quad \quad \quad \quad \quad[\theta_\X]_{x',x}} \ar[dd]_{\T_\X(\lag{h}{v \otimes \bar{u}})} && \langle x' | \Ss_\X | x \rangle \ar[dd]^<<<<<<<<<{\Ss_\X(\lag{h}{v \otimes \bar{u}})} \ar[ur]_-{\Omega^\Ss_{x',x}}
\\ & \displaystyle \bigoplus_\blag{x}{x'}{g} (\X_\mmlag{h \ccdot x}{}{hgh^\mmi})_\mathbb{C} \otimes \T_{hgh^\mmi} \ar'[r][rr]_-{\id \otimes \theta_{hgh^\mmi} \qquad \qquad \qquad \qquad } && \displaystyle \bigoplus_\blag{x}{x'}{g} (\X_\mmlag{h \ccdot x}{}{hgh^\mmi})_\mathbb{C} \otimes \Ss_{hgh^\mmi} \\
\langle h \ccdot x' | \T_\X | h \ccdot x \rangle \ar[ur]^-{\Omega^\T_{h \ccdot x', h \ccdot x}} \ar[rr]_-{[\theta_\X]_{h\ccdot x', h\ccdot x}} && \langle h \ccdot x' | \Ss_\X | h \ccdot x \rangle \ar[ur]_{\Omega^\Ss_{h \ccdot x', h \ccdot x}}  }
\]
Now, the left and right faces commute since we showed in Proposition \ref{aTgim} (i) that $\Omega^\T$ and $\Omega^\Ss$ are equivariant.  The top and bottom faces are exactly our definition of $\theta_\X$ in terms of $\vartheta$. Finally the back face commutes, because it is precisely the condition that $\{ \vartheta_g \colon \T_g \rightarrow \Ss_g \}$ is a map of equivariant vector bundles over $G$. Hence the front face commutes.

Finally we must check that the formula \eqref{apresc} indeed satisfies the condition of a modification --- that is, that for every morphism $V \colon \X \rightarrow \Y$ the following pasting diagram commutes:
 \[
\xy
(-7.5,10)*+{\X}="1";
(-11, 7)="1f"; (-5, 9)="1b";
(-7.5,-10)*+{\X}="2";
(-11, -10)="2f"; (-5, -7)="2b";
(35,10)*+{\Y}="3";
(33.5, 7)="3f"; (38.5, 9)="3b";
(35,-10)*+{\Y}="4";
(33.5, -10)="4f"; (38.5, -7)="4b";
{\ar@/_0.8pc/_-{\T_\X} "1f";"2f"};
{\ar@/^0.5pc/^-{\Ss_\X} "1b";"2b"};
{\ar@/_0.8pc/_-{\T_\Y} "3f";"4f"};
{\ar@/^0.5pc/^-{\Ss_\Y} "3b";"4b"};
{\ar^-{V} "1";"3"};
{\ar_-{V} "2";"4"};
{\ar@2{->}@/{^0.5pc}/^>>>{{\footnotesize \Ss(V)}} (14.5,6)*+{};(12.5,0)*+{}};
{\ar@2{->}@/{_0.5pc}/_<<<<<{\T(V)} (11,2)*+{};(11,-4)*+{}};
{\ar@2{->}_-{\theta_\X} (-13,-3)*+{};(-4,2)*+{}};
{\ar@2{->}_-{\theta_\Y} (31,-3)*+{};(38.5,2)*+{}};
\endxy
\]
In components, this diagram is the requirement that the following diagram commutes:
 \[
  \xymatrix{ \langle y | \T_\Y \circ V | x \rangle \ar[d]_{\T(V)_{y, x}} \ar[r]^{\theta_\Y \otimes \id} & \langle y | \Ss_\Y \circ V | x \rangle \ar[d]^{\Ss(V)_{y, x}} \\
  \langle y | V \circ \T_\X | x \rangle \ar[r]_{\id \otimes \theta_\X} & \langle y | V \circ \Ss_\X | x \rangle}
 \]
To show that this equation holds, we use our formula \eqref{af} which expresses $\T(V)$ and $\Ss(V)$ in terms of $V$, $\Omega^\T$ and $\Omega^\Ss$ respectively. Substituting these expressions in, together with the our expression \eqref{apresc} for $\theta_\X$ and $\theta_\Y$ in terms of the maps $\{\vartheta_g \}$, and finally precomposing the above diagram with $[\Omega^\T_\Y]^\mi_{y,x}$ and postcomposing it with $[\Omega^\Ss_\X]_{y,x}$, we find that the $\Omega$ parts cancel and we are left with checking the following commutative diagram
 \[
  \xymatrix{ \displaystyle \bigoplus_{y'} \bigoplus_\blag{y'}{y}{g} (\Y_\mlag{}{y}{g})_\mathbb{C} \otimes \T_g \otimes \langle y' | V | x \rangle \ar[d]^-{\text{swap} \circ (v \mapsto u) \otimes \id \otimes V(\lag{g}{v \otimes \bar{u}})} \ar[rr]^-{\id \otimes \vartheta_g \otimes \id} && \displaystyle \bigoplus_{y'} \bigoplus_\blag{y'}{y}{g} (\Y_\mlag{}{y}{g})_\mathbb{C} \otimes \Ss_x \otimes \langle y' | V | x \rangle \ar[d]^-{\text{swap} \circ (v \mapsto u) \otimes \id \otimes V(\lag{g}{v \otimes \bar{u}})} \\
  \displaystyle \bigoplus_{x'} \bigoplus_\blag{x}{x'}{g} \langle y | V | x' \rangle \otimes (\X_\mlag{x}{}{g})_\mathbb{C} \otimes \T_g \ar[rr]_-{\id \otimes \id \otimes \theta_g} &&  \displaystyle \bigoplus_{x'} \bigoplus_\blag{x}{x'}{g} \langle y | V | x' \rangle \otimes (\X_\mlag{x}{}{g})_\mathbb{C} \otimes \Ss_g}
 \]
which clearly commutes.
\end{proof}

\section{Wrapping up the proof\label{lg}}
In this section we wrap up the proof of our main result in this chapter. Firstly, we can now conclude the following.
\begin{prop}\label{rrfun} The restriction functor
 \[
  \Dim(\Gerbes(G)) \rightarrow \Dim(\Gerbes(G)|_{\E G})
 \]
is an equivalence of braided monoidal categories.
\end{prop}
\begin{proof}
The restriction functor is clearly braided monoidal. Proposition \ref{aTgim} establishes that it is essentially surjective, while Proposition \ref{afullyfaithful} establishes that it is fully faithful.
\end{proof}
This allows us to establish the following.
\XXX{Sort out numbering here.}
\begin{thm}\label{joesot} The higher-categorical dimension of the 2-category of unitary 2-representations of a finite group $G$ is equivalent, as a braided monoidal category, to the category of conjugation-equivariant hermitian vector bundles over $G$ equipped with the fusion tensor product. In symbols,
 \[
  \Dim \TRep(G) \simeq \Hilb_G^\fusion (G).
 \]
\end{thm}
\begin{proof} This follows from the following sequence of equivalences of braided monoidal categories:
\begin{align*}
  \Dim \TRep(G) & \simeq \Dim \Gerbes(G) \\
   & \simeq \Dim \left(\Gerbes(G)|_{\E G}\right) \\
   & \simeq \Hilb_G^\fusion (G).
\end{align*}
The first equivalence uses Theorem \ref{toodytheorem} --- that the 2-category $\TRep(G)$ of unitary 2-representations is equivalent to the 2-category $\Gerbes(G)$ of finite equivariant gerbes --- together with Proposition \ref{natdimt} from Appendix \ref{movieapp} --- that if two 2-categories $\mathcal{A}$ and $\mathcal{B}$ are equivalent, then their higher-categorical dimensions are equivalent as braided monoidal categories. The second equivalence is Proposition \ref{rrfun}, and the third is Proposition \ref{fyn}.
\end{proof}
As we stated in the introduction to this thesis, the significance of this result is that it verifies, in the case of the untwisted 3d finite group model, the `crossing with the circle' equation of extended TQFT which is expected to hold if the Baez-Dolan Extended TQFT Hypothesis is indeed correct. Moreover, it shows that the braided monoidal structure on the category assigned to the circle can indeed be computed solely via `abstract nonsense' from the 2-category assigned to the point: it recovers precisely the fusion tensor product and braiding which Freed obtained using topological operations on the pair of pants cobordism \cite{ref:freed}.

\appendix

\chapter{Verifying the `crossing with $S^1$' equation for low codimension\label{AppendixTQFTFacts}}
In this appendix we record some computations using the groupoid technology of Willerton \cite{ref:simon} which verify the `crossing with the circle' equation
 \[
  Z(M \times S^1) \cong \Dim Z(M)
 \]
for closed manifolds $M$ in the $n$-dimensional twisted finite group extended TQFT, as we promised in the Introduction.  

Recall that the initial geometric data in the finite group model is a finite group $G$ and a group $n$-cocycle $\omega \in Z^n(\BG, U(1))$, and the quantum invariants $Z(M)$ are defined as the space of sections of the higher line bundles over $\mathcal{P}_M$ defined by transgression of the cocycle $\omega$.  We will verify the above formula in those dimensions $n$ where $Z(M)$ is a vector space or a category; that is when $\text{codim} \, M := n - \dim M$ is equal to 1 or 2.

Firstly, recall that the space of fields $\mathcal{P}_M$ on $M$ is defined as the groupoid of $G$-bundles on $M$; since such a bundle is determined by its holonomy we may set
 \[
  \mathcal{P}_M = \Fun (\Pi_1 (M), \BG)
 \]
where $\Pi_1(M)$ is the fundamental groupoid of $M$. In what follows one can assume that a finite number of points have been chosen on $M$ as the endpoints of these paths; this makes $\Pi_1(M)$ a finite groupoid.
  
Next, observe that the space of fields on $M \times S^1$ computes as the loop groupoid of the space of fields on $M$:
 \begin{align*}
  \mathcal{P}_{M \times S^1} &= \Fun(\Pi_1 (M \times S^1), \BG) \cong \Fun\left(\Pi_1 (S^1), \Fun(\Pi_1(M), \BG)\right)  \\
   = \Lambda \mathcal{P}_M.
  \end{align*}
Let $\omega_M \in Z^{\codim  M} (\mathcal{P}_M, U(1))$ be the transgression of the $n$-cocycle $\omega \in Z^n(\BG, U(1))$ to the space of fields $\mathcal{P}_M$. Under the identification of $\mathcal{P}_{M \times S^1}$ with $\Lambda \mathcal{P}_M$, the transgressed cocycle
 \[
  \omega_{M \times S^1} \in Z^{\codim M -1} (\mathcal{P}_{M \times S^1}, U(1))
 \]
identifies as the ordinary loop transgression of the cocycle $\omega_M$ (see Willerton \cite{ref:simon}), 
 \[
  \omega_{M \times S^1} = \tau(\omega_M).
 \]
This loop transgression can be given a concrete formula. In the case $\codim M = 1$, so that $\tau(\omega_M) \in Z^0(\Lambda \mathcal{P}_M, U(1))$, we have
 \[
  \tau(\omega_M)  ( \Fix{P}{\gamma} ) = \omega_M (\gamma),
 \]
where $P \in \mathcal{P}_M$ and $\gamma \in \Stab(P)$. Similarly when $\codim M = 2$, so that $\tau(\omega_M) \in Z^1(\Lambda \mathcal{P}_M, U(1))$, then
 \[
  \tau(\omega_M)\left(  \stackrel{\delta}{\leftarrow} (\Fix{P}{\gamma}) \right) = \frac{ \omega_M(\delta, \gamma)}{\omega_M(\gamma, \delta)}.
 \]
We now use the groupoid cocycle technology of Willerton \cite{ref:simon} to verify the `crossing with the circle' equation for these invariants. 

\begin{prop} Let $Z$ be the $n$-dimensional twisted finite group extended TQFT associated to a group $G$ and an $n$-cocycle $\omega \in Z^n(\BG, U(1))$. Then:
 \begin{enumerate}
  \item If $\dim M = n-1$, then we have the equality
   \[
    \int_{\mathcal{P}_{M \times S^1}} \tau(\omega_M) = \Dim \Gamma_{\mathcal{P}_M} \left((\omega_M)_\mathbb{C}\right)
   \]
  where $(\omega_M)_\mathbb{C}$ is the line bundle over $\mathcal{P}_M$ associated to the transgressed 1-cocycle $\omega_M \in Z^1(\mathcal{P}_M, U(1))$.
 \item If $\dim M = n-2$, then we have an isomorphism of vector spaces
  \[
   \Gamma_{\mathcal{P}_M} \left( \tau(\omega_M)_\mathbb{C} \right) \cong \Dim \Rep^{\omega_M} \left(\mathcal{P}_M\right)
  \]
  where $\omega_M \in Z^2(\mathcal{P}_M, U(1))$ is the transgressed 2-cocycle and $\Rep^{\omega_M} \left(\mathcal{P}_M\right)$ is the category of $\omega_M$-twisted vector bundles over $\mathcal{P}_M$.
\end{enumerate}
\end{prop}
\begin{proof} Formula (i) is precisely a specific case of the general integration formula \cite[Thm 6]{ref:simon}, while formula (ii) is precisely \cite[Thm 16]{ref:simon}.
\end{proof}
From the definition of the invariants $Z(M)$, we recognize this proposition as essentially verifying the claim that $Z(M \times S^1) \cong \Dim Z(M)$ in the case were $\codim \, M = 1$ or $2$. We use the word `essentially' because in item (ii) above we should really be recovering the {\em full Frobenius algebra structure} on
 \[
  Z(M \times S^1) = \Gamma_{\mathcal{P}_M} \left( \tau(\omega_M)_\mathbb{C} \right)
 \]
and not just an `isomorphism of vector spaces'. There are however some subtleties here regarding the inner products on these spaces and we leave this for later work. 
\chapter{Fusion categories have ambidextrous duals\label{AppendixFusion}}
In this appendix we record the proof of a lemma which was promised in Chapter \ref{stss}.

\begin{lem} A fusion category has ambidextrous duals.
\end{lem}
\begin{proof}
We are told that every object $V$ has a {\em right} dual $V^\star$; we must show that $V^\star$ is also a {\em left} dual of $V$. Observe that for every pair of objects $V$ and $W$ the presence of right duals gives rise to a bijection of hom-sets:
 \begin{align*}
  \Hom(1, V^\star \otimes W) &\cong \Hom(V, W) \\
   \ba \ig{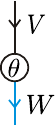} \ea & \mapsto \ba \ig{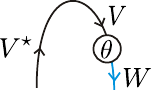} \ea \\
   \ba \ig{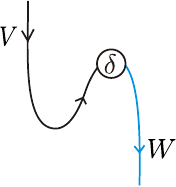} \ea & \mapsfrom \ba \ig{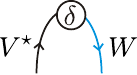} \ea
  \end{align*}
This means in particular (I learnt this argument from \cite[Prop 2.1]{ref:eno}) that the right dual $X^\star$ of a simple object $X$ must be the unique simple object up to isomorphism which has a nonzero map $1 \rightarrow X^\star \otimes X$, and moreover this map is unique up to a scalar. In exactly the same way, the {\em left} dual $^\star \!X$ of a simple object $X$ must be the unique simple object up to isomorphism which has a nonzero map $^\star \!X \otimes X \rightarrow 1$. But by semisimplicity, there must then also be a nonzero map in the other direction, $1 \rightarrow ^\star \!X \otimes X$, whence we must have $^\star \!X \cong X^\star$. Thus the simple objects have ambidextrous duals.

This means that {\em all} objects must have ambidextrous duals, by the following argument. Let $(X)_{i \in I}$ be some maximal choice of nonisomorphic simple objects. Given an object $V$, choose a basis
 \begin{align*}
  \{ v_{i,p} \colon X_i & \rightarrow V \}_{p=1}^{\dim \Hom(X_i, V)}
 \intertext{for each hom-set $\Hom(X_i, V)$. Since $C$ is semisimple, we have the corresponding dual bases}
  \{ v_i^p \colon V & \rightarrow X_i \}
 \end{align*}
for $\Hom(V, X_i)$ which satisfy
\[
 v_i^p \circ v_{i,q} = \delta^p_q \id_{X_i} \quad \text{and} \sum_{i,p} v_{i,p} \circ v_i^p = \id_V.
  \]
If $^\star V$ and $V^\star$ are left and right duals of $V$ respectively, then by making some choice of ambidextrous duals $X_i^*$ for the simple objects $X_i$ we can construct the following isomorphism $^\star V \rightarrow V^\star$:
 \[
  \sum_{i,p} \ba \ig{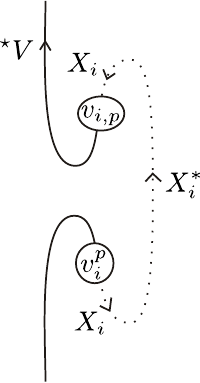} \ea
 \]
This is indeed an isomorphism, because its inverse $V^\star \rightarrow {} ^\star V$ is given by
 \[
 \sum_{i,p} \ba \ig{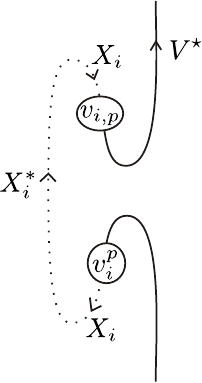} \ea,
 \]
as the reader will verify by an elementary calculation. Hence $^\star V \cong V^\star$ and all objects have ambidextrous duals.
\end{proof} 
\chapter[Naturality of the categorical dimension]{Naturality of the higher-categorical dimension\label{movieapp}}
In this appendix we show that if there is a {\em coherent adjoint equivalence} between two 2-categories $\A$ and $\B$  --- a notion we will shortly define --- then their higher-categorical dimensions are equivalent in a well-defined way, that is there is a braided monoidal functor $J \colon \Dim \A \rightarrow \Dim \B$ which is an equivalence of categories. We needed this result in our proof that $\Dim \TRep(G) \simeq \Hilb_G^\fusion (G)$ in Chapter \ref{lg}, since our strategy was to instead prove that $\Dim \Gerbes(G) \simeq \Hilb_G^\fusion (G)$, and then appeal to the equivalence of the 2-categories $\TRep(G)$ and $\Gerbes(G)$ which we established in Theorem \ref{toodytheorem}.

What we call a `coherent adjoint equivalence' is a special case of what Gurski calls a `biadjoint biequivalence' in \cite[Definition A.3.3]{ref:gurski}, with the following provisos:
\begin{itemize}
 \item We do not additionally require the modifications $b,b',s$ and $s'$ below to themselves satisfy the snake equations, simply because we do not need that property in the proof of Theorem \ref{joesot}, though we agree this should form part of the ultimate definition.
 \item We explicitly add in the extra triangulators $\Delta_1'$ and $\Delta_2'$ and the {\em horizontal cusp} coherence law, which is only implicitly present in \cite{ref:gurski}. This property is crucial in order for us to prove that the monoidal functor $J$ above is indeed a {\em braided} monoidal functor, so we decided to emphasize its importance here by listing it explicitly. This coherence law was explicated to us by J. Scott Carter \cite{ref:carter_private}.
 \end{itemize}
Gurski has shown that every equivalence between 2-categories can be upgraded into a biadjoint biequvalence (and hence a coherent adjoint equivalence) \cite{ref:gurski_new}. The definition of a coherent adjoint equivalence is as follows.
\begin{defn} A {\em coherent adjoint equivalence} between two 2-categories $\mathcal{A}$ and $\mathcal{B}$ consists of the following data:
 \begin{itemize}
  \item Weak 2-functors $F \colon \mathcal{A} \rightarrow \mathcal{B}$ and $G \colon \mathcal{B} \rightarrow \mathcal{A}$, drawn as:
   \[
    \ba \ig{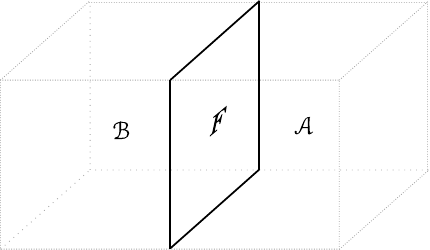} \ea \quad \text{and} \quad \ba \ig{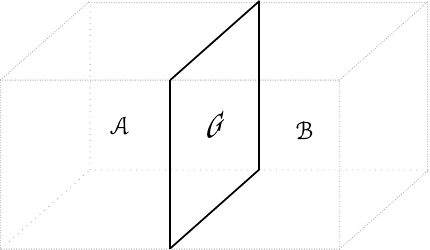} \ea
   \]
   \item Transformations $\eta \colon \id_\A \Rightarrow GF$, $\eta^* \colon GF \Rightarrow \id_\A$, $\epsilon \colon FG \Rightarrow \id_\B$ and $\epsilon^* \colon \id_\mathcal{B} \Rightarrow FG$, drawn as:
 \begin{align*}
 \eta = \ba \ig{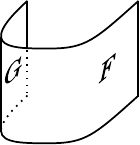} \ea, & \quad \eta^* = \ba \ig{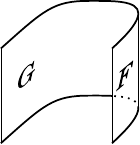} \ea \\
 \epsilon = \ba \ig{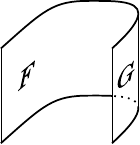} \ea, & \quad \epsilon^* = \ba \ig{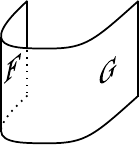} \ea
\end{align*}
  \item Invertible modifications
  \begin{align*}
  \Delta_1 = \ba \ig{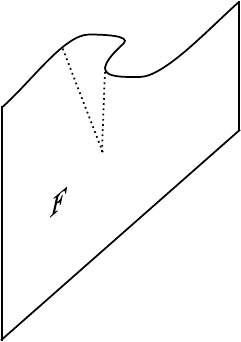} \ea, & \quad \quad \Delta_2 = \ba \ig{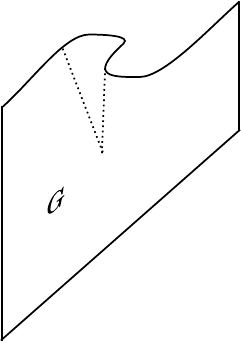} \ea \\
  \Delta_1' = \ba \ig{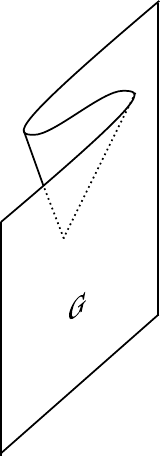} \ea, & \quad \quad \Delta_2' = \ba \ig{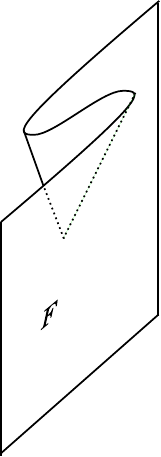} \ea
  \end{align*}
  called the {\em triangulators}, and
\item Invertible modifications $b \colon \id_{\id_\mathcal{A}} \threearrow \eta^* \eta$ (for `birth of a circle'), $s \colon \eta \eta^* \threearrow \id_{GF}$ (for `saddle), $b' \colon \id_{\id_\mathcal{B}} \threearrow \epsilon \epsilon^*$ and $s' \colon \epsilon^* \epsilon \threearrow \id_{FG}$, drawn as:
    \begin{align*}
    b = \ba \ig{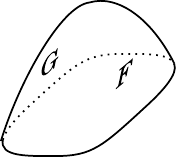} \ea, &\quad s = \ba \ig{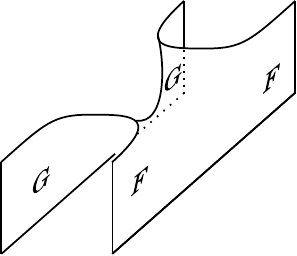} \ea \\
    b' = \ba \ig{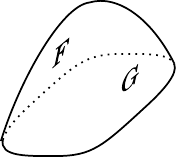} \ea, & \quad s' = \ba \ig{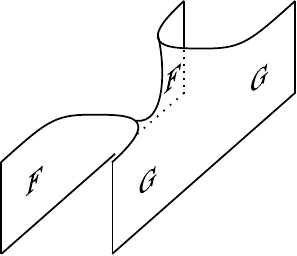} \ea \; \; .
    \end{align*}
  \end{itemize}
This data must satisfy the following equations (over and above the inverse laws for the modifications), drawn in movie moves as follows, looking `down from above':
 \begin{itemize}
  \item Swallowtail rules:
   \[
    \ba \ig{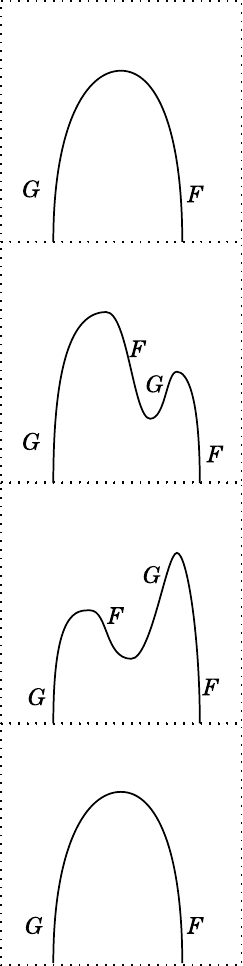} \ea \quad = \quad \ba \ig{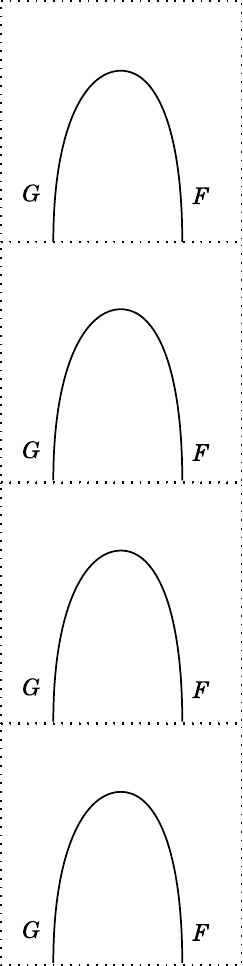} \ea
   \]
   together with the four other variants of this equation, and
  \item Horizontal cusp rules:
   \[
    \ba \ig{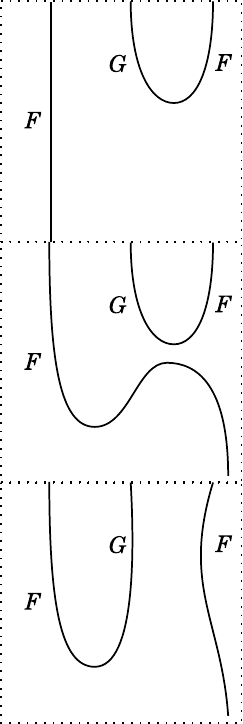} \ea \quad = \quad \ba \ig{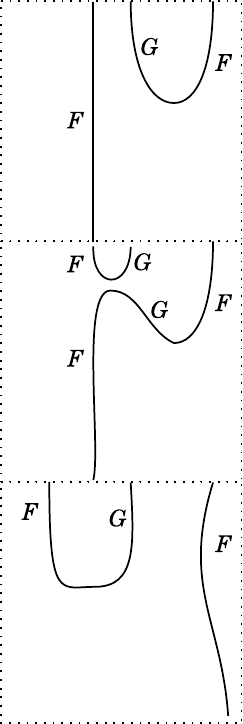} \ea
   \]
  together with the four other variants of this equation.
 \end{itemize}
\end{defn}

\begin{prop}\label{natdimt} Suppose that there is a coherent adjoint equivalence between two 2-categories $\A$ and $\B$. Then there is a well-defined equivalence of braided monoidal categories $K \colon \Dim \B \stackrel{\simeq}{\rightarrow} \Dim \A$.
\end{prop}
\begin{proof} Suppose that $\T$ is a transformation of the identity on $\B$. We define the transformation $K(\T) \in \Dim(\A)$ as
 \be \label{compt}
   K(\T) = \eta^* \circ (\id_G \ast \T \ast \id_F) \circ \eta.
  \ee
In string diagrams:
 \[
  K(T) = \ba \ig{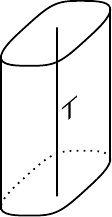} \ea.
 \]
Now, since $2\mathcal{C}\text{at}$ really {\em does} form a 3-category\footnote{Recall that we are using the weak terminology, so that for us the word `3-category' means `tricategory'.}\cite{ref:gordon_power_street, ref:gurski}, this composite indeed makes sense and is indeed a transformation $K(\T) \colon \id_\mathcal{A} \Rightarrow \id_\mathcal{A}$. The same formula works at the level of modifications as well; this shows that $K \colon \Dim \mathcal{A} \rightarrow \Dim \mathcal{B}$ is a well-defined functor between categories.

We need to establish that $K$ is monoidal. If $\T, \Ss \in \Dim \mathcal{B}$, we define the coherence cell
 \[
  \phi \colon K(\T) \circ K(\Ss) \threearrow K(\T \circ \Ss)
 \]
to be the modification $\id \ast e \ast \id$ where we are using the `saddle piece' $e$ we defined above:
 \[
 \phi = \ba \ig{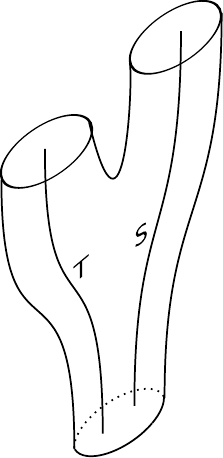} \ea
 \]
It is clear that $\phi$ satisfies the coherence equations for a monoidal functor, because of the interchange law for modifications in $2\mathcal{C}\text{at}$.

To check that $K$ is a {\em braided} monoidal functor, we must verify that the following diagram of modifications commutes:
 \[
  \begin{array}{rcl} K(\T)\circ K(\Ss) & \stackrel{c_\mathcal{A}}{\ba \ig{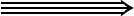} \ea} & K(\Ss) \circ K(\T) \\
    \phi \; \ba \ig{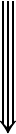} \ea & & \ba \ig{downthree.pdf}\ea \; \phi \\
     K(\T \circ \Ss) & \displaystyle \mathop{\ba \ig{acrossthree.pdf} \ea}_{K(c_\mathcal{A})} & K(\Ss \circ \T)
    \end{array}
  \]
In string diagrams, we must check the following equation:
 \[
  \ba \ig{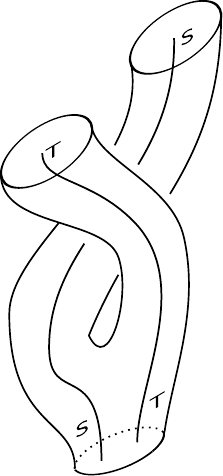} \ea = \ba \ig{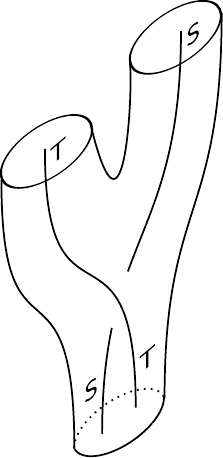} \ea.
 \]
We transform the left hand side to the right hand side in movie moves as follows. The author thanks J. Scott Carter for showing him this proof in private correspondence \cite{ref:carter_private}. The idea is to move the saddle point from the bottom to the top by using the horizontal cusp rule (H), as well as the swallowtail rule (S). Every other equality below simply uses the naturality of interchanging the order in which modifications are performed. \newpage
\[
\ba \ig{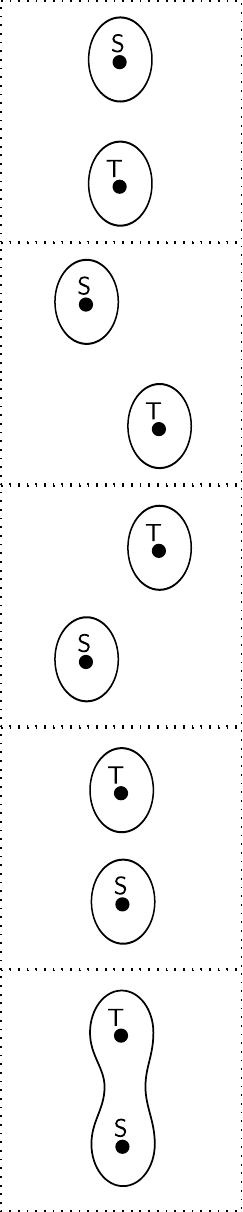} \ea \; \; \stackrel{(S)}{=} \; \; \ba \ig{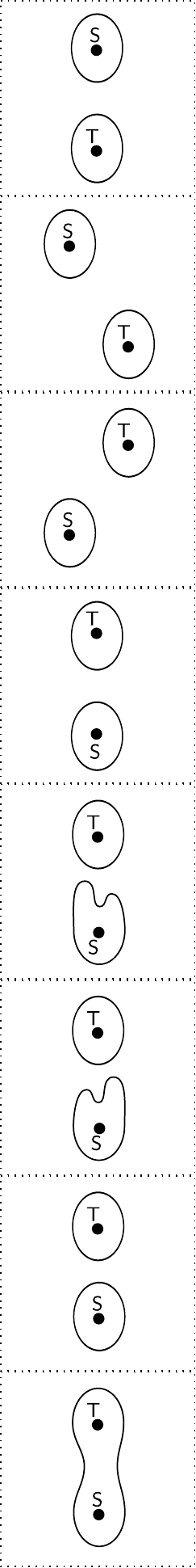} \ea  \; \; \stackrel{}{=} \; \; \ba \ig{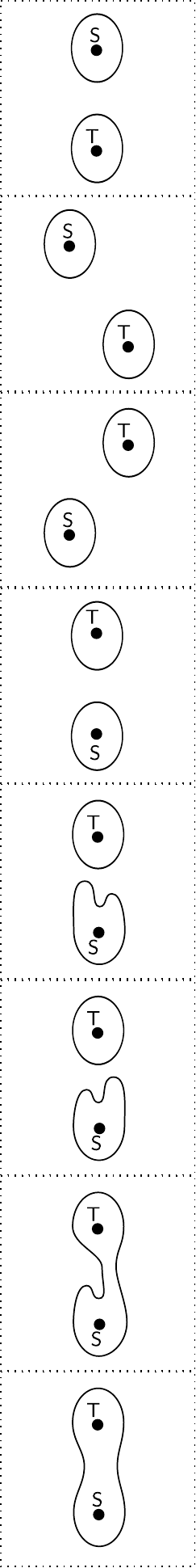} \ea \; \; \stackrel{}{=} \; \; \ba \ig{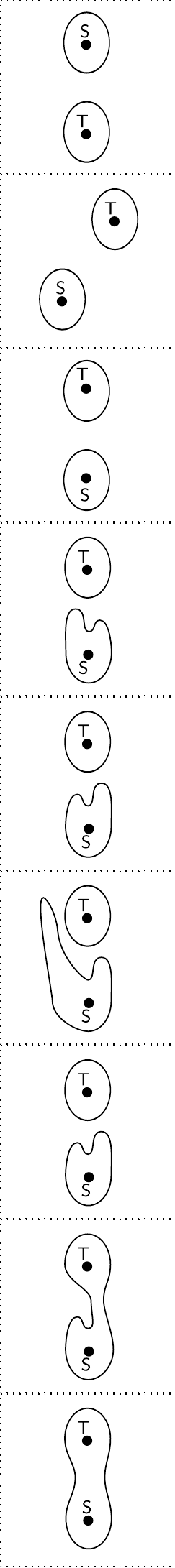} \ea
\]
\[
\stackrel{}{=} \; \; \ba \ig{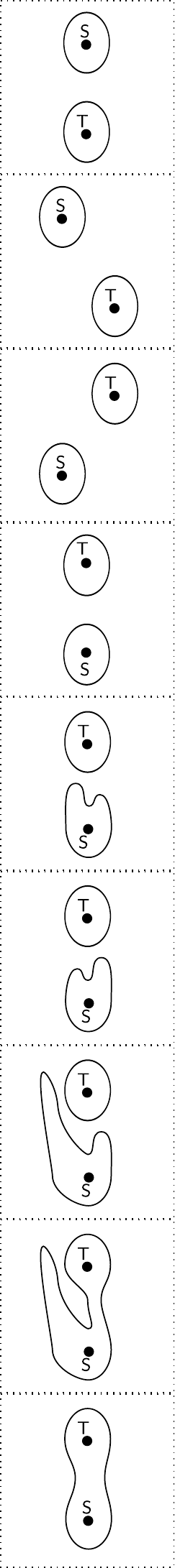} \ea \; \; \stackrel{}{=} \; \; \ba \ig{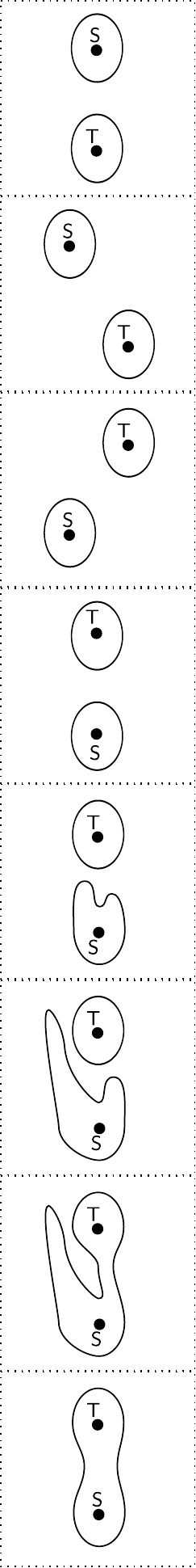} \ea  \; \; \stackrel{}{=} \; \; \ba \ig{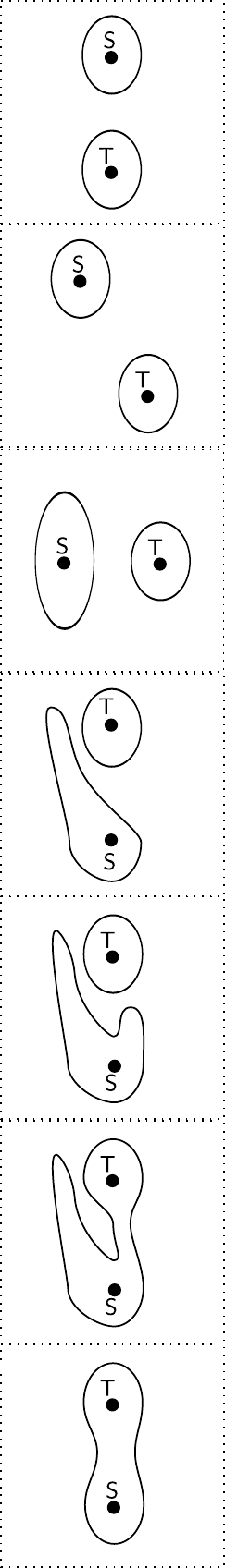} \ea \; \; \stackrel{(H)}{=} \; \; \ba \ig{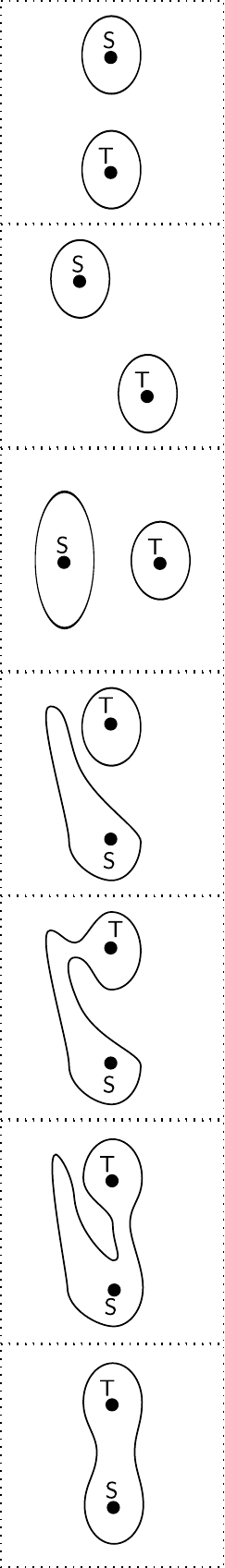} \ea
\]
\[
\stackrel{}{=}\; \;  \ba \ig{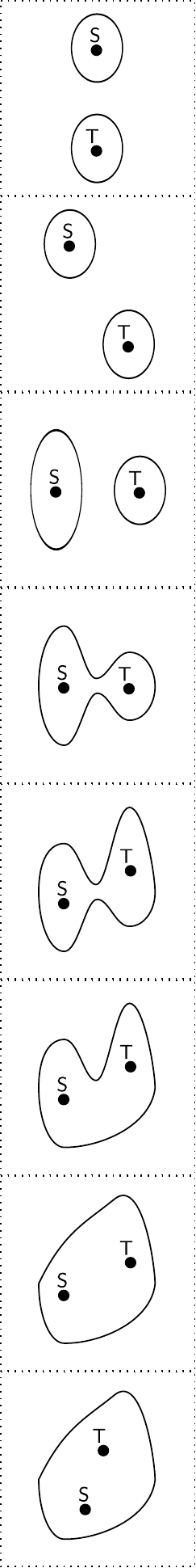} \ea \; \; \stackrel{}{=} \; \; \ba \ig{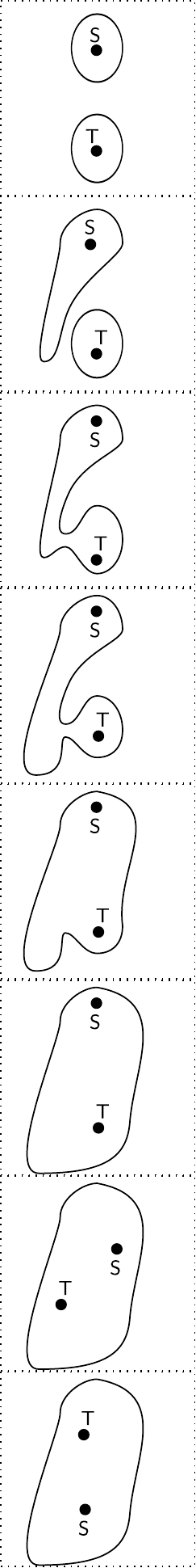} \ea  \; \; \stackrel{}{=} \; \; \ba \ig{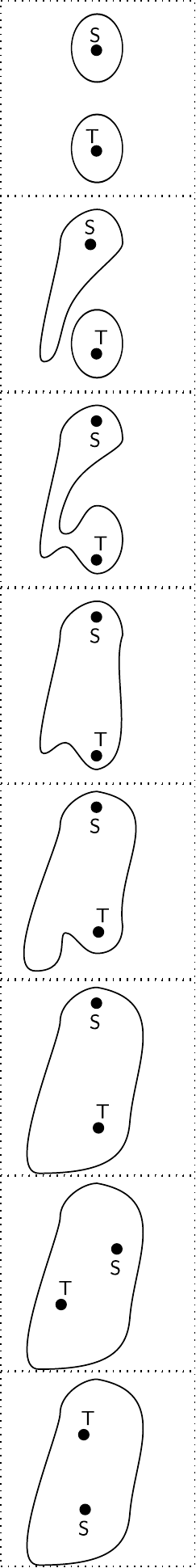} \ea \; \; \stackrel{(H)}{=} \; \; \ba \ig{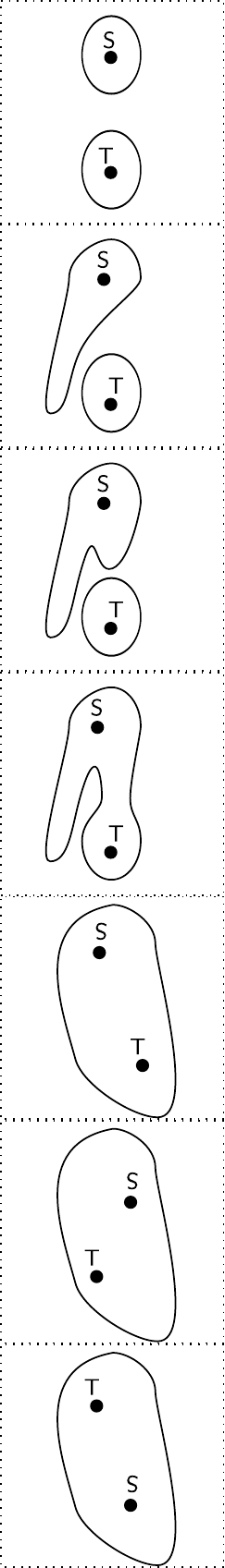} \ea
\]
\[
\stackrel{}{=}\; \;  \ba \ig{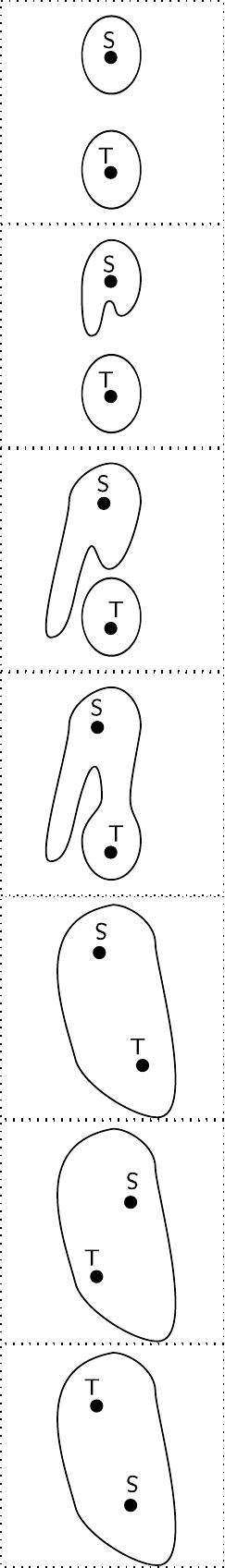} \ea \; \; \stackrel{(S)}{=} \; \; \ba \ig{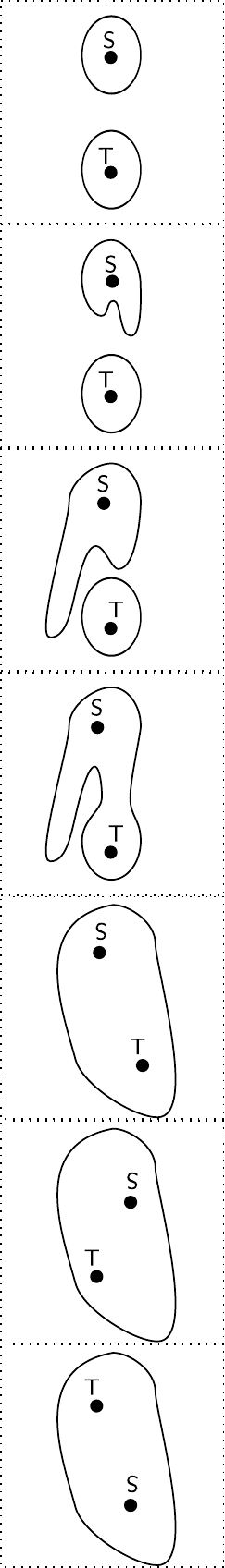} \ea  \; \; \stackrel{}{=} \; \; \ba \ig{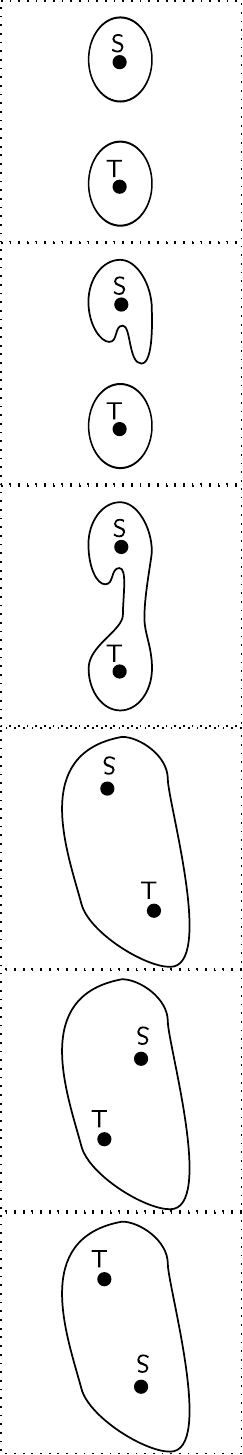} \ea \; \; \stackrel{}{=} \; \; \ba \ig{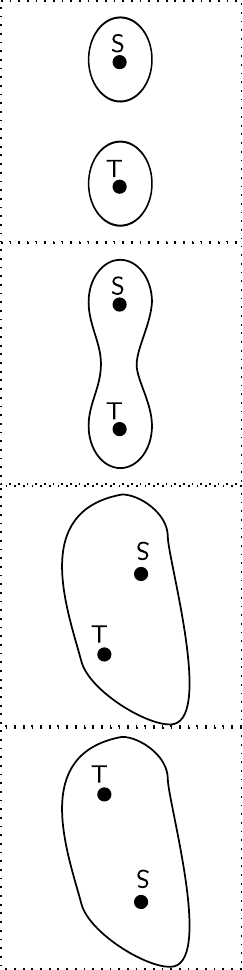} \ea \; .
\]
This completes the proof that $K \colon \Dim \mathcal{B} \rightarrow \Dim \mathcal{A}$ is a well-defined braided monoidal functor.

Finally, we must check that $K$ is an equivalence of categories. Indeed, the same constructions as above give us a well-defined braided moniodal functor $J \colon \Dim \A \rightarrow \Dim B$. Thus $K$ is essentially surjective, because for every $\T \in \Dim B$ we have an isomorphism $\T \rightarrow KJ(\T)$ given by the following movie:
\[
 \ba \ig{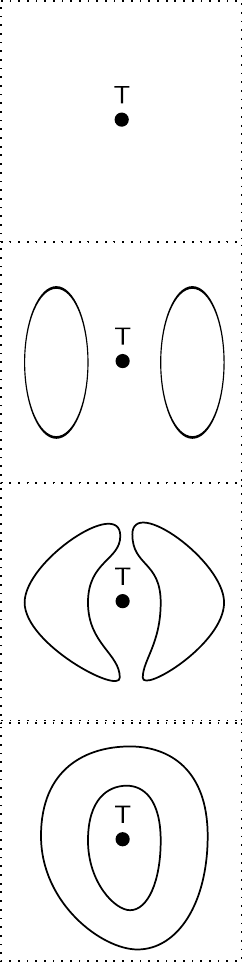} \ea
\]
The same construction shows that $K$ is fully faithful. This completes the proof that $\Dim \A$ is braided monoidally equivalent to $\Dim \B$.
\end{proof}

\chapter{Fusion tensor product for conjugation equivariant vector bundles\label{AppFreed}}
In this chapter we review the braided monoidal structure on the category $\Hilb^\text{fusion}_G (G)$ of $G$-equivariant unitary vector bundles over a finite group $G$, as in the lecture notes of Freed \cite{ref:freed2}. This tensor product is often known as the `fusion' product because it arises from considering how $G$-bundles on a circle fuse together via the pair-of-pants cobordism; it is the `category assigned to the circle' when one thinks of the untwisted finite group model as a three-tier extended TQFT.  The braiding arises from the diffeomorphism of the pair of pants which swaps the two input circles; see \cite{ref:freed2}.

\begin{defn} The braided monoidal category $\Hilb^\text{fusion}_G(G)$ is defined as follows. As a category, it is simply the category $\Hilb_G (G)$ of unitary equivariant vector bundles over the loop groupoid $\Lambda \BG$, as in Section \ref{loopgdsec}. The tensor product $V \otimes W$ of two equivariant vector bundles $V$ and $W$ is defined fibrewise by convolution. That is,
 \[
  (V \otimes W)_g = \oplus_{ab = g} V_a \otimes W_b.
 \]
This vector bundle $V \otimes W$ is equipped as an equivariant vector bundle by defining the equivariant maps fibrewise:
 \[
  (V \otimes W)(\lag{h} g) = \oplus_{ab = g} V(\lag{h} a) \otimes W(\lag{h} b).
 \]
The braiding
 \[
  c_{V, W} \colon V \otimes W \rightarrow W \otimes V
 \]
is defined fibrewise as by conjugating the first factor and then swapping the factors:
 \begin{align*}
  (c_{V, W})_g \colon \oplus_{ab = g} V_a \otimes W_b & \rightarrow \oplus_{ab = g} W_a \otimes V_b \\
    v_a \otimes w_b & \mapsto w_b \otimes V(\lag{b^\mi} a) (v_a).
  \end{align*}
\end{defn}
We remark that our convention for the braiding is opposite to that of Freed, who conjugates on the {\em second} factor and then swaps the factors. We have made this choice in order to recover the normal convention for the braiding on $\Dim \Gerbes(G)$, the category of transformations of the identity 2-functor on $\Gerbes(G)$, in Proposition \ref{fyn}(ii). On the other hand, this normal convention is itself {\em opposite} to the usual convention for transformations, as we pointed out in Chapter \ref{la}.

\newcommand{\arXiv}[1]{\href{http://arxiv.org/abs/#1}{{\tt arXiv:#1}}}

\newcommand{\hepth}[1]{\href{http://arxiv.org/abs/hep-th/#1}{{\tt arXiv:hep-th/#1}}}

\newcommand{\qalg}[1]{\href{http://arxiv.org/abs/q-alg/#1}{{\tt arXiv:q-alg/#1}}}

\newcommand{\quantph}[1]{\href{http://arxiv.org/abs/quant-ph/#1}{{\tt arXiv:quant-ph/#1}}}

\newcommand{\mathQA}[1]{\href{http://arxiv.org/abs/math.QA/#1}{{\tt arXiv:math.QA/#1}}}

\newcommand{\mathCT}[1]{\href{http://arxiv.org/abs/math.CT/#1}{{\tt arXiv:math.CT/#1}}}

\newcommand{\mathph}[1]{\href{http://arxiv.org/abs/math-ph/#1}{{\tt arXiv:math-ph/#1}}}

\newcommand{\mathDG}[1]{\href{http://arxiv.org/abs/math.DG/#1}{{\tt arXiv:math.DG/#1}}}

\newcommand{\Math}[1]{\href{http://arxiv.org/abs/math/#1}{{\tt arXiv:math/#1}}}

\newcommand{\mathAT}[1]{\href{http://arxiv.org/abs/math.AT/#1}{{\tt arXiv:math.AT/#1}}}

\newcommand{\mathKT}[1]{\href{http://arxiv.org/abs/math.KT/#1}{{\tt arXiv:math.KT/#1}}}

\newcommand{\arxiv}[1]{\href{http://arxiv.org/abs/#1}{{\tt arXiv:#1}}}

\XXX{Insert Tim Porter reference}


\begin{thebibliography}{77} 

\small              
\raggedright

\bibitem{ref:ambrose} W. Ambrose, Structure theorems for a special class of Banach algebras, Trans. Amer. Math. Soc. 57 (1945), 364-386.

\bibitem{ref:arnold} V. I. Arnold, {\em Catastrophe Theory}, 3rd edition, Berlin: Springer-Verlag (1992).

\bibitem{ref:atiyah} M. F. Atiyah, {\em The Geometry and Physics of Knots}, Lezioni Lincee, Cambridge University Press (1990).

\bibitem{ref:atiyah2} M. Atiyah, Topological quantum field theory, Publications Mathematiques de l'IHS 68 (1988) 176--186.

\bibitem{ref:atiyah_singer} M. F. Atiyah and I. M. Singer, The Index of Elliptic Operators: III, Annals of Mathematics, Vol 87 No. 3 (1968) 546--604.

\bibitem{ref:baez_2_hilbert_spaces} J. C. Baez, Higher-dimensional algebra II: 2-Hilbert
spaces, Adv. Math. 127 (1997), 125-189. Also available as \qalg{9609018}.

\bibitem{ref:baez_qq} J. C. Baez, Quantum Quandaries, available as \quantph{0404040}.

\bibitem{ref:baez_baratin_freidel_wise} J. C. Baez, A. Baratin, L. Freidel and D. Wise, Representations of 2-Groups on Higher Hilbert Spaces. Available as \arXiv{0812.4969}.

\bibitem{ref:baez_higher_schreier_theory} J. C. Baez, This Week's Finds in Mathematical Physics Week 223, \href{http://math.ucr.edu/home/baez/week223.html}{\tt http://math.ucr.edu/home/baez/week223.html}.

\bibitem{ref:baez_dolan_hda0} J. C. Baez and J. Dolan, Higher-dimensional Algebra and Topological Quantum Field Theory, J. Math. Phys. 36 (1995), 6073--6105. Also available as \qalg{9503002}.

\bibitem{ref:baez_dolan_1} J. C. Baez and J. Dolan, From finite sets to Feynman diagrams, {\em Mathematics unlimited --- 2001 and beyond}, 29--50, Springer Berlin (2001). Also available as \Math{0004133}.

\bibitem{ref:baez_langford} J.~C.~Baez, L.~Langford, Higher-dimensional algebra IV: 2-tangles, Adv. Math. 180 (2003), 705--764. Also available as \Math{9811139}.

\bibitem{ref:baez_lauda_2-groups} J. C. Baez and A. D. Lauda, Higher-dimensional algebra V: 2-Groups, Theory and Applications of Categories Vol 12 (2004) No. 14, 423--491. Also available as \mathQA{0307200}.

\bibitem{ref:baez_stay_rosetta} J. C. Baez and M. Stay, Physics, Topology, Logic and Computation: A Rosetta Stone. Available as \href{http://math.ucr.edu/home/baez/rosetta.pdf}{http://math.ucr.edu/home/baez/rosetta.pdf}.

\bibitem{ref:bakalov_kirillov} B. Bakalov and A. Kirillov, Jr., {\em Lectures on Tensor Categories and Modular Functors}, Amer. Math. Soc. University Lecture Series (2001).

\bibitem{ref:bar_pfaffle} C. B\"{a}r and F. Pf\"{a}ffle, Path integrals on manifolds by finite dimensional approximation. Available as \Math{0703272}.

\bibitem{ref:barrett_mackaay} J. W. Barrett and M. Mackaay,  Categorical
representations of categorical groups, Theory and Applications of Categories Vol 16 (2006) No. 20, 529--557. Also available as \mathCT{math.CT/0407463}.

\bibitem{ref:barrett_westbury} J. W. Barrett and B. W. Westbury, Spherical categories,
     Adv. Math. 143 (1999) 357. Also available as \hepth{9310164}.

\bibitem{ref:barrett_westbury2} J. W. Barrett and B. W. Westbury, Invariants of piecewise-linear 3-manifolds, Trans. Amer. Math. Soc. 348 503--510 (1995). Also available as \hepth{9311155}.

\bibitem{ref:bartlett} B. Bartlett, The geometry of unitary 2-representations of finite groups and their 2-characters. Available as \arXiv{0807.1329}.

\bibitem{ref:bartlett2} B. Bartlett, Categorical Aspects of Topological Quantum Field Theories, MSc thesis, Utrecht University (2005). Available as \mathQA{0512103}.

\bibitem{ref:behrend_xu} K. Behrend and P. Xu, Differentiable Stacks and Gerbes. Available as \Math{0605694}.

\bibitem{ref:benzvi} D. Ben-Zvi, Lecture series on Topological Field Theory and Geometric Langlands, KITP Santa Barbara, July 2008. Available as \href{http://www.math.utexas.edu/~benzvi/GRASP/lectures/benzvi/SBTFT.pdf}{\tt http://www.math.utexas.edu}{\tt /~benzvi/GRASP/lectures/benzvi/SBTFT.pdf}.

\bibitem{ref:benzvi2} D. Ben-Zvi, Lectures on The Geometric Langlands Correspondence, Mathematical Institute, Oxford, April 2007. Available at \href{http://people.maths.ox.ac.uk/~szendroi/langlands.html}{\tt http://people.maths.ox.ac.uk/~szendroi/langlands.html}.

\bibitem{ref:benzvi_francis_nadler} D. Ben-Zvi, J. Francis and D. Nadler, Integral Transforms and Drinfeld Centers in Derived Algebraic Geometry. Available as \arXiv{0805.0157}.

\bibitem{ref:berline_getzler_vergne} N. Berline, E. Getzler and M. Vergne, {\em Heat Kernels and Dirac Operators}, Springer-Verlag Berlin Heidelberg (2004).

\bibitem{ref:berman_et_al} R. Berman, B. Berndtsson and J. Sj\"{o}strand, A direct approach to Bergman kernel asymptotics for positive line bundles, Arkiv f\"{o}r Matematik 46 No. 2 (2008). Also available as \Math{math/0506367}.

\bibitem{ref:berman} R. Berman, Bergman kernels and local holomorphic Morse inequalities, Mathematische Zeitschrift, Volume 248, No. 2 (2004), pp. 325--344. Also available as \Math{0211235}.

\bibitem{ref:berndtsson} B. Berndtsson, Bergman kernels related to Hermitian line bundles over compact complex manifolds, in {\em Explorations in Complex and Riemannian Geometry}, Eds. J. Bland, K-T Kim and S. G. Krantz, Amer. Math. Soc.  Contemporary Mathematics Series 332, page 1--19.

\bibitem{ref:borceaux_dejean} F. Borceaux and D. Dejean, Cauchy completion in category theory, Cahiers topologie et g\'{e}om\'{e}trie diff\'{e}rentielle 27 (1986) No. 2, 133--146.

\bibitem{ref:boyarchenko} M. Boyarchenko, Introduction to modular categories, Lecture series at University of Chicago, Geometric Langlands Seminar. Available at \href{http://www.math.uchicago.edu/~mitya/langlands.html}{www.math.uchicago.edu/~mitya/langlands.html}.

\bibitem{ref:boyarchenko_drinfeld} M. Boyarchenko and V. Drinfeld, Character sheaves on unipotent groups in positive characteristic: foundations. Available as \arXiv{0810.0794}.

\bibitem{ref:brocker_dieck} T. Br\"{o}cker and T. tom Dieck, {\em Representations of Compact Lie Groups}, Springer Graduate Texts in Mathematics (1985).

\bibitem{ref:brylinski} J-L. Brylinski, {\em Loop Spaces, Characteristic Classes and Geometric Quantization},
Progress in Mathematics volume 107, Birkhauser-Boston (1993).

\bibitem{ref:brylinski_mclaughlin} J-L. Brylinski and D. McLaughlin, The geometry of degree-four characteristic classes and of line bundles on loop spaces I, Duke Math. J. 75 (1994) No. 3, 603--638.

\bibitem{ref:calaque_etingof} D. Calaque and P. Etingof, Lectures on tensor categories. Available as \mathQA{0401246v3}.

\bibitem{ref:caldararu_willerton} A. Caldararu and S. Willerton, The Mukai Pairing, I: A categorical
approach. Available as \Math{0707.2052}.



\bibitem{ref:carter_private} J. S. Carter, private correspondence.

\bibitem{ref:carter_flath_saito} J. S. Carter, D. E. Flath and M. Saito, {\em The Classical and Quantum 6j-symbols}, Princeton University Press (1995).

\bibitem{ref:carter_segal_macdonald} G. Segal, {\em Lie Groups}, taken from R. Carter, G. Segal and I. Macdonald, {\em Lectures on Lie Groups and Lie Algebras}, London Mathematical Society Student Texts 32 (1995).

\bibitem{ref:cegarra_et_al_graded_extensions} A.M. Cegarra, A.R. Garz\'{o}n and A.R.-Grandjean, Graded extensions of categories, J. Pure Appl. Algebra 154 (2000) 117-141. Also available as \href{http://www.ugr.es/~anillos/Preprints/total.ps}{\tt http://www.ugr.es/~anillos/Preprints/total.ps}.

\bibitem{ref:cheng_gurski} E. Cheng and N. Gurski, The periodic table of $n$-categories for low dimensions II: degenerate tricategories. Available as \arXiv{0706.2307}.

\bibitem{ref:cheng_gurski2} E. Cheng and N. Gurski, Towards an $n$-category of cobordisms, Theory and Applications of Categories, Vol. 18 (2007) No. 10 pp. 274-302.

\bibitem{ref:crane_sheppeard} L. Crane and M.D. Sheppeard, 2-categorical Poincare Representations and State Sum Applications. Available as \Math{0306440}.

\bibitem{ref:crane_yetter} L. Crane and D.N. Yetter, Measurable Categories and 2-Groups, Applied Categorical Structures 13 (2005) No. 5-6, 501--516. Also available as \mathQA{0305176}.

\bibitem{ref:crans} S. E. Crans, Generalized Centers of Braided and Sylleptic Monoidal 2-Categories, Adv. Math. 136 No. 2 (1998) 183--223. Also available as \href{crans.fol.nl/papers/gcbs.ps.gz}{\tt crans.fol.nl/papers/gcbs.ps.gz}.

\bibitem{ref:costello} K. Costello, Topological conformal field theories and Calabi-Yau categories, Adv. Math. 210 No. 1 (2007), 165--214. Also available as \Math{0412149}.

\bibitem{ref:deligne} P. Deligne, Action du groupe des tresses sur une cat\'{e}gorie, Invent. Math. 128 (1997), 159--175.

\bibitem{ref:desmit_lenstra} B. de Smit and H. W. Lenstra Jr., Linearly Equivalent Actions of Solvable Groups, J. Algebra Vol 228 (2000) 270--285.

\bibitem{ref:dijkgraaf_witten} R. Dijkgraaf and E. Witten, Topological gauge theories and group cohomology, Comm. Math. Phys. 129 (1990) 60-72.

\bibitem{ref:duistermaat} J. J. Duistermaat, {\em The Heat Kernel Lefschetz Fixed Point Formula for the Spin-c Dirac Operator}, Progress in Nonlinear Differential Equations and Their Applications Vol 18, Birkh\"{a}user.

\bibitem{ref:durhuus_jonsson} B. Durhuus and T. Jonsson, Classification and construction of unitary topological field theories in two dimensions, J. Math. Phys. 35 (1994) 5306-5313. Also available as \hepth{9308043}.

\bibitem{ref:elgueta} J. Elgeueta, Representation theory of 2-groups on Kapranov and Voevodsky's 2-category
2Vect, Adv. Math. 213 (2007) No. 1, 53--92. Also available as \Math{0408120v2}.

\bibitem{ref:elgueta2} J. Elgueta, A strict totally coordinatized version of Kapranov and Voevodsky's 2-category 2Vect. Available as \Math{0406475}.

\bibitem{ref:elgueta3} J. Elgueta, Generalized 2-vector spaces and general linear 2-groups, J. Pure and Applied Algebra Vol 212 Issue 9 (2008) 2069--2091. Also available as \Math{0606472}.

\bibitem{ref:eno} P. Etingof, D. Nikshych and V. Ostrik, On fusion categories, Annals of Mathematics, Volume 162, Number 2 (2005), 581--642. Also available as \Math{0203060v9}.

\bibitem{ref:eno2} P. Etingof, D. Nikshych and V. Ostrik, An analogue of Radford's $S^4$ formula for finite tensor categories, Int. Math. Res. Not. 54 (2004) 2915--2933. Also available as \Math{0404504}.

\bibitem{ref:fantechi} B.~Fantechi et al, {\em Fundamental Algebraic Geometry: Grothendieck's FGA Explained}, Amer. Math. Soc. Mathematical Surveys and Monographs Volume 123 (2006).

\bibitem{ref:freed} D. Freed, Higher algebraic structures and quantization, Comm. Math. Phys. 159 (1994) No. 2, 343--398.

\bibitem{ref:freed2} D. Freed, Quantum groups from path integrals. Available as \qalg{9501025}.

\bibitem{ref:freed_remarks} D. Freed, Remarks on Chern-Simons Theory. Available as \arXiv{0808.2507}.

\bibitem{ref:freed3} D. Freed, Locality and Integration in Topological Field Theory, in Group Theoretical Methods in Physics, Volume 2, Eds M.A. del Olmo, M. Santander and J.M. Guilarte, Ciemat, 1993, 35-54. Also available as \hepth{9209048}.

\bibitem{ref:freed_teleman_hopkins} D. Freed, M. Hopkins and C. Teleman, Loop Groups and Twisted K-Theory II. Available as \Math{0511232}.

\bibitem{ref:freed_quinn} D. Freed and F. Quinn, Chern-Simons theory with finite gauge group, Commun. Math. Phys. 156 (1993), 435-472.

\bibitem{ref:freyd_yetter} P.J. Freyd and D.N. Yetter, Coherence theorems via knot theory, J. Pure and App. Alg.
78 (1992) 49--76.

\bibitem{ref:freyd_yetter_2} P.J. Freyd and D.N. Yetter, Braided compact monoidal categories with applications to low dimensional topology, Adv. Math. 77 (1989), 156-182.

\bibitem{ref:ganter_kapranov_rep_char_theory} N. Ganter and M. Kapranov, Representation and character theory in 2-categories, Adv. Math. 217 (2008) No. 5, 2268--2300. Also available as \mathKT{0602510}.

\bibitem{ref:gordon_power_street} R. Gordon, A. J. Power and R. Street, Coherence for tricategories, Memoirs of the Amer. Math. Soc. 117 (1995) Number 558.

\bibitem{ref:grothendieck} A. Grothendieck, {\em Cat\'{e}gories fibr\'{e}es et descent (SGAI) expos\'{e}
VI}, Lecture Notes in Mathematics Vol. 224 (1971) 145--194, Springer Berlin.


\bibitem{ref:gurski} N.~Gurski, {\em An algebraic theory of tricategories}, Phd Thesis, University of Chicago (2006). Available as \href{http://gauss.math.yale.edu/~mg622/tricats.pdf}{\tt http://gauss.math.yale.edu/~mg622/tricats.pdf}.

\bibitem{ref:gurski_new} N.~Gurski, Biadjoint Biequivalences. Available as \href{
http://gauss.math.yale.edu/~mg622/biadjdraft.pdf}{\tt
http://gauss.math.yale.edu/~mg622/biadjdraft.pdf}.

\bibitem{ref:hagge_hong} T. J. Hagge and S-M Hong, Some non-braided fusion categories of rank 3. Available as \arXiv{0704.0208}.

\bibitem{ref:hitchin} N. Hichin, Harmonic spinors, Adv. Math. 14 (1974), 1--55.

\bibitem{ref:hollander} S. Hollander, A homotopy theory for stacks, Israel Journal of Math. 163 (2008), 93--124. Also available as \Math{0110247}.

\bibitem{ref:isaacs} I. M. Isaacs, {\em Character Theory of Finite Groups}, Dover Publications (1976).

\bibitem{ref:jones} V. Jones, A new knot polynomial and von Neumann algebras, Notices Amer. Math. Soc. 33 (1986), 219--225.

\bibitem{ref:joyal_street} A. Joyal and R. Street, The geometry of tensor calculus I, Adv. Math. 88 (1991), 55-112.

\bibitem{ref:joyal_street_btc} A. Joyal and R. Street, Braided tensor categories, Adv. Math. 102 (1993), 20--78.

\bibitem{ref:joyal_street_bmc} A. Joyal and R. Street, Braided monoidal categories, Macquarie Mathematics Reports No. 860081, Preprint, 1986. Available as \href{http://rutherglen.ics.mq.edu.au/~street/JS86.pdf}{http://rutherglen.ics.mq.edu.au/~street/JS86.pdf}.

\bibitem{ref:joyal_street_bmc2} A. Joyal and R. Street, Braided monoidal categories, Preprint (revised and expanded version of \cite{ref:joyal_street_bmc}).

\bibitem{ref:joyal_street_3} A. Joyal and R. Street, Tortile Yang-Baxter operators in tensor categories, J. Pure Appl. Algebra 71 (1991), 43--51.

\bibitem{ref:kapranov_voevodsky} M. M. Kapranov and V. A. Voevodsky, 2-categories
and the Zamolodchikov tetradedra equations, Proceedings of
Symposia in Pure Mathematics 56 (1994) 177--259.

\bibitem{ref:kassel} C.~Kassel, {\em Quantum Groups}, Springer-Verlag New York 1995.

\bibitem{ref:khovanov} M. Khovanov, A categorification of the Jones polynomial, Duke Math. J. 101 No. 3 (1999), 359--426. Also available as \mathQA{0207264}.

\bibitem{ref:kirwin} W. D. Kirwin, Coherent States in Geometric Quantization, J. Geom. Phys. 57 (2007) No. 2, 531-548. Also available as \Math{0502026}.

\bibitem{ref:kock} J. Kock, {\em Frobenius algebras and 2D topological quantum field theories}, London Mathematical Society Student Texts No. 59, Cambridge University Press (2003).

\bibitem{ref:kohno} T. Kohno, {\em Conformal Field Theory and Topology}, Translations of Mathematical Monographs, Iwanami Series in Modern Mathematics, Amer. Math. Soc. Vol 210 (1998).

\bibitem{ref:lauda} A.~Lauda, Frobenius algebras and planar open string topological field theories. Available as \Math{0508349}.

\bibitem{ref:lawrence} R. Lawrence, Triangulation, categories and extended field theories, in {\em Quantum Topology}, eds. R. Baadhio and L. Kaufmann, World Scientific, Singapore (1993), 191--208.

\bibitem{ref:lawvere} F. W. Lawvere, Metric spaces, generalized logic and closed categories, Rendiconti del Seminario Matematico e Fisico di Milano,  XLIII (1973), 135--166; republished in Reprints in Theory and Applications of Categories, No. 1 (2002),  1--37.

\bibitem{ref:leinster_basic_bicategories} T. Leinster, Basic Bicategories. Available as \mathCT{9810017}.

\bibitem{ref:luztig} G. Luztig, Character sheaves, I, Adv. in Math. Vol 56 (1985) 193--237.

\bibitem{ref:martins_porter} J. F. Martins and T. Porter, On Yetter's Invariant and an Extension of the Dijkgraaf-Witten Invariant to Categorical Groups, Theory and Applications of Categories, Vol. 18, (2007), No. 4, 118--150. Also available as \Math{0608484}.

\bibitem{ref:moerdijk1} I. Moerdijk, Introduction to the Language of Stacks and Gerbes. Available as \mathAT{0212266}.

\bibitem{ref:ma_marinescu} X. Ma and G. Marinescu, {\em Holomorphic Morse Inequalities and Bergman Kernels}, Birkhauser Verlag Ag Switzerland (2007).

\bibitem{ref:mason} G. Mason, Orbifold Conformal Field Theory and Cohomology of the Monster. Available as \href{www.newton.ac.uk/programmes/NST/Mason.pdf}{\tt www.newton.ac.uk/programmes/NST/Mason.pdf}.

\bibitem{ref:morton} J. Morton, Extended TQFT's and quantum gravity, PhD thesis University of California Riverside (2007). Available as \arXiv{0710.0032}.

\bibitem{ref:müger_from_subfactors_to_categories_and_topology_I} M. M\"{u}ger, From Subfactors to Categories
and Topology I. Frobenius algebras in and Morita equivalence of tensor categories. J. Pure Appl. Alg. 180,
(2003) 81--157. Also available as \Math{0111204}.

\bibitem{ref:müger_galois} M. M\"{u}ger, Galois Theory for Braided Tensor Categories and the Modular Closure, Adv. in Math., 150, Issue 2 (2000), 151--201. Also available as \mathCT{9812040}.

\bibitem{ref:mueger_abstract_duality} M. M\"{u}ger, Abstract duality for symmetric tensor
*-categories, appendix to H. Halversson, Algebraic Quantum Field Theory, in {\em Philosophy of Physics} (Handbook of the Philosophy of Science), Eds. J. Butterfield and J. Earman, Elsevier (2007). Also available as \mathph{0602036}.

\bibitem{ref:müger_lectures} M. M\"{u}ger, Tensor categories: A selective guided tour, \arXiv{0804.3587}.

\bibitem{ref:müger_center_group} M. M\"{u}ger, On the center of a compact group, Int. Math. Res. Not. (2004), 2751--2756.

\bibitem{ref:murray_singer} M. Murray and M. Singer, Gerbes,  Clifford Modules and the Index Theorem, Annals of Global Analysis and Geometry 26 (2004) No. 4, 355--367.  Also available as \mathDG{0302096}.

\bibitem{ref:ostrik} V. Ostrik, Module categories, weak Hopf algebras and modular invariants, Transformation Groups 8 (2003) No. 2, 177--206. Also available as \Math{0111139}.

\bibitem{ref:polesello_waschkies} P. Polesello and I. Waschkies, Homology Homotopy Appl. Vol 7 (1) (2005) 109--150. Also available as \Math{0407507}.

\bibitem{ref:poon} W. B. Poon, The Kodaira vanishing theorem and generalizations, Phd thesis, University of Hong Kong (2002). Available at \href{http://hub.hku.hk/handle/123456789/31640}{\tt http://hub.hku.hk/handle/123456789/31640}.

\bibitem{ref:quinn} F. Quinn, Lectures on axiomatic topological quantum field theory. Geometry and Quantum Field Theory, Volume 1, IAS/Park City Mathematical Series, Amer. Math. Soc., 325--453.

\bibitem{ref:resh} N. Y. Reshetikhin and V. G. Turaev, Invariants of 3-manifolds via link polynomials and quantum groups, Invent. Math. 103 (1991) 547--597.

\bibitem{ref:roberts_willerton} J. Roberts and S. Willerton, On the Rozansky-Witten weight systems. Available as \Math{0602653}.

\bibitem{ref:robinson_rawnsley} P.L. Robinson and J.H. Rawnsley, {\em The metaplectic representation, $Mp^\mathbb{C}$ structures and geometric quantization}, Memoirs of the American Mathematical Society Vol 81 No. 410.

\bibitem{ref:roe} J. Roe, {\em Elliptic operators, topology and asymptotic methods}, Chapman \& Hall/CRC, Second Edition (1998).

\bibitem{ref:sati_et_al} H. Sati, U. Schreiber, Z. Sk\u{o}da and D. Stevenson, Twisted differential nonabelian cohomology. Available as \href{http://www.math.uni-hamburg.de/home/schreiber/nactwist.pdf}{\tt http://www.math.uni-hamburg.de/home/schreiber/nactwist.pdf}.

\bibitem{ref:schreiber} U. Schreiber, Transgression of n-transport and n-connection. Personal notes, available at \href{http://www.math.uni-hamburg.de/home/schreiber/dgcaspaces.pdf}{\tt http://www.math.uni-hamburg.de/home/schreiber/dgcaspaces.pdf}.

\bibitem{ref:schreiber2} U. Schreiber, Transgression of n-transport and n-connection. Personal notes, available at \href{http://www.math.uni-hamburg.de/home/schreiber/dgcaspaces.pdf}{\tt http://www.math.uni-hamburg.de/home/schreiber/dgcaspaces.pdf}.

\bibitem{ref:schreiber_waldorf1} U. Schreiber and K. Waldorf, Connections on non-abelian Gerbes and their Holonomy. Available as \arXiv{0808.1923}.

\bibitem{ref:schreiber_waldorf2} U. Schreiber and K. Waldorf, Smooth Functors vs. Differential Forms. Available as \arXiv{0802.0663}.

\bibitem{ref:segal} G. Segal, The definition of conformal field theory, from {\em Topology, Geometry and Quantum Field Theory: Proceedings of the 2002 Oxford Symposium in Honour of the 60th birthday of Graeme Segal}, London Math. Soc. Lecture Note Series No. 308.

\bibitem{ref:segal2} G. Segal, `Stanford Notes' on Topological Field Theory. Available at \href{http://www.cgtp.duke.edu/ITP99/segal/}{\tt http://www.cgtp.duke.edu/ITP99/segal/}.

\bibitem{ref:seidel_thomas} P. Seidel and R. Thomas, Braid group actions on derived categories of coherent sheaves, Duke Math. J. Vol 108 No. 1 (2001) 37--108.

\bibitem{ref:simon} S. Willerton, The twisted Drinfeld double of a finite group via gerbes and
finite groupoids, Alg. and Geom. Top. Vol. 8 No. 3 (2008). Also available as \mathQA{0503266}.

\bibitem{ref:spera} M. Spera, On K\"{a}hlerian coherent states, Proceedings of the International Conference on
Geometry, Integrability and Quantization, Varna Bulgaria 1999 (eds. I. Mladenov and G. Naber), Coral Press (2000), 241--256.

\bibitem{ref:street_string_diagrams} R. Street, Categorical structures, in {\em Handbook
of Algebra}, Vol. 1, ed. M. Hazewinkel, Elsevier Amsterdam (1995), 529--577.

\bibitem{ref:street_yoneda} R. Street, Fibrations in bicategories, Cahiers topologie et g\'{e}om\'{e}trie diff\'{e}rentielle 21 (1980) 111--160.

\bibitem{ref:street_colimits} R. Street, Absolute colimits in enriched categories, Cahiers topologie et g\'{e}om\'{e}trie diff\'{e}rentielle 24 (1983) No. 4, 111--160.

\bibitem{ref:street_conspectus} R. Street, An Australian conspectus of higher categories, prepared for the Institute of Mathematics and its Applications Summer Program on {\em n-Categories: Foundations and Applications}, University of Minnesota (2004). Available as \href{www.maths.mq.edu.au/~street/Minneapolis.pdf}{\tt www.maths.mq.edu.au/~street/Minneapolis.pdf}.

\bibitem{ref:thom} R. Thom, {\em Structural Stability and Morphogenesis: An Outline of a General Theory of Models}, Westview Press (1989). Originally published as {\em Stablit\'{e} structurelle et morphog\'{e}n\`{e}se}, W. A. Benjamin Inc. (1972).

\bibitem{ref:tambara_yamagami} D. Tambara and S. Yamagami, Tensor categories with fusion rules of self-duality for finite abelian groups, J. Algebra 209 (1998), no. 2, 692-707.



\bibitem{ref:turaev2} V. Turaev, {\em Quantum Invariants of Knots and 3-Manifolds}, Walter de Gruyter, Berlin (1994).

\bibitem{ref:vicary} J. Vicary, Categorical formulation of quantum algebras. Available as \quantph{0805.0432}.

\bibitem{ref:webster} B. Webster, Small linearly equivalent $G$-sets and a construction of Beaulieu, Journal of Algebra, Vol 317 Issue 1 (2007) 306-323. Also available as \Math{1610205}.

\bibitem{ref:wells} R. O. Wells, {\em Differential Analysis on Complex Manifoldss}, Springer (1980).

\bibitem{ref:witten} E. Witten, Quantum field theory and the Jones polynomial, Comm. Math. Phys. Vol 121, No. 3 (1989), 351--399.

\bibitem{ref:witten2} E. Witten, Mirror manifolds and topological field theory, Essays on Mirror Manifolds ed. S.-T. Yau, International Press, Hong Kong, 1992, pp. 120–158. Also available as \hepth{9112056}.

\bibitem{ref:weinstein} A. Weinstein, The volume of a differentiable stack. Available as \arXiv{0809.2130}.

\bibitem{ref:woit} P. Woit, {\em Not Even Wrong: The Failure of String Theory and the Search for Unity in Physical Law}, Perseus Publishing (2007).

\bibitem{ref:woodhouse} N. M. J. Woodhouse, {\em Geometric Quantization}, Oxford Mathematical Monographs (1992).

\bibitem{ref:yetter} D. N. Yetter, Triangulations and TQFT's, taken from {\em Quantum Topology}, eds. R. Baadhio and L. H. Kaufmann, Series on Knots and Everything Vol. 3, World Scientific, (1993).

\bibitem{ref:yetter2} D. N. Yetter, Measurable Categories, Applied Categorical Structures 13 (2005) 469--500.

\end{thebibliography}
\end{document}